\input amssym.def
\input amssym.tex

Let $G$ be a connected algebraic reductive group over an algebraic closure of a prime field ${\Bbb F}_p$, defined over ${\Bbb F}_q$ thanks to a Frobenius $F$. Let $\ell$ be a prime different from $p$. Let $B$ be an $\ell$-block of the group of rational points $G^F$. Under mild restrictions on $\ell$, we  show the existence of an algebraic reductive group $H$ defined over ${\Bbb F}_q$ {\it via} a Frobenius $F$, and of a unipotent $\ell$-block $b$ of $H^F$ such that :

the respective defect groups of $b$ and $B$ are isomorphic, the
associated Brauer categories are isomorphic and there is a height preserving one-to-one map from the set of
irreducible representations  of $b$ onto the set of
irreducible representations  of $B$.

\def\AA{{\bf A}}
\def\BB{{\bf B}}
\def\DD{{\bf D}}
\def\CC{{\bf C}}
\def\EE{{\bf E}}
\def\FF{{\bf F}}
\def\GG{{\bf G}}
\def\XX{{\bf X}}

\def\FDD{{\frak D}} 


\def\tg{{\tilde g}}
\def\tG{{\tilde {G}}}
\def\tL{{\tilde L}}
\def\ts{{\tilde s}}


\def\lexp#1#2{\kern\scriptspace
\vphantom{#2}^{#1}\kern-\scriptspace#2}
\def\cent#1#2{{\rm C}_{#1}(#2)}	             
\def\cento#1#2{{\rm C}^\circ_{#1}(#2)}	      
\def\CF#1#2{{\cal CF}(#1,#2)}       
\def\nor#1#2{{\rm N}_{#1}(#2)}         
\def\zo#1{{\rm Z}^\circ(#1)}           
\def\z#1{{\rm Z}(#1)}                  


\def\Q{{\Bbb Q}}
\def\N{{\Bbb N}}
\def\Z{{\Bbb Z}}
\def\C{{\Bbb C}}
\def\R{{\Bbb R}}


\def\Ind#1#2{{\rm Ind}_{#1}^{#2}}


\def\dec{decomposition}
\def\irr{irreducible}
\def\al{\alpha}
\def\la{\lambda}
\def\ga{G_{\rm a}}


\def\a#1{{#1}_{\rm a}}
\def\b#1{{#1}_{\rm b}}

\def\F#1{{\Bbb F}_{#1}}
\def\gen#1{\langle #1\rangle}
\def\II#1{{\rm Irr}({#1})}

\def\Lu#1#2{{\rm R}_{#1}^{#2}\,}
\def\slu#1#2{{}^*\!{\rm R}_{#1}^{#2}}

\def\ser#1#2{{\cal E}({#1},#2 )}
\def\lser#1#2{{\cal E}_{\ell}\big({#1},#2 \big)}
\def\Res#1#2{{\rm Res}_{#2}^{#1}}
\def\Ind#1#2{{\rm Ind}_{#1}^{#2}}
\def\Ker{{\rm Ker}\,}

\def\scal#1#2#3{\langle{#1},{#2}\rangle_{#3}}

\def\bl#1#2{{\rm Bl}(#1;#2)}

\def\joinrel{\mathrel{\mkern-4mu}}

\def\surli#1{\vbox{\ialign{##\crcr\hrulefill\kern1pt\crcr\noalign
{\kern1pt\nointerlineskip}$\hfil\displaystyle{#1}\hfil$\crcr}}}
\def\surline#1{\vbox{{\kern.5pt\nointerlineskip}\ialign
{##\crcr\kern1pt\hrulefill\kern1pt\crcr\noalign
{\kern.5pt\nointerlineskip}$\hfil\scriptstyle{#1}\hfil$\crcr}}}
\def\mapright#1{\smash{\mathop{\relbar\joinrel\longrightarrow}\limits^{#1}}}
\def\mapleft#1{\smash{\mathop{\longleftarrow\joinrel\relbar}\limits^{#1}}}
\def\mapdown#1{\Big\downarrow\rlap{$\vcenter{\hbox{$\scriptstyle#1$}}$}}
\def\mapup#1{\Big\uparrow\rlap{$\vcenter{\hbox{$\scriptstyle#1$}}$}}
\def\diag#1#2#3#4#5#6#7#8#9{\def\normalbaselines
{\baselineskip20pt\lineskip3pt
\lineskiplimit3pt}\matrix {#1&\mapright{#2}&#3\cr \mapdown{#8}&#9&\mapdown{#4}\cr
#7&\mapright{#6}&#5\cr} }  
\font\ninerm=cmr9

\def\cF{{\cal F}}

\def\bull{\hfill\vrule height 1.6ex width .7ex depth -.1ex
\smallskip}

\def\piu{\pi_{\rm un}}

\def\preuve{\medskip\noindent {\it Proof. }}
\def\simex#1{{\sim}_{#1}}

\openup 1\jot

{

\bigskip\centerline{\bf TOWARDS A JORDAN DECOMPOSITION OF BLOCKS OF FINITE REDUCTIVE GROUPS}
}

\vskip1cm\bigskip\centerline{MICHEL E. ENGUEHARD}

\vskip1cm
{\narrower\narrower
\smallskip\noindent{\it Foreword. }{\ninerm In 2008 was published by Michel Enguehard in Journal of Algebra [0] "Vers une d\'ecomposition de Jordan des blocs des groupes r\'eductifs finis". So I really have to say some words on that [0]-paper and the present text, name it [$\infty$]-paper. Several years after the publication of [0] I received, as the presumed author, several quite judiciouses questions on some results and proofs contained in it. I discovered in  [0] inaccuracies, incomplete proofs, not to say more. To save my homonymous collegue reputation I decided to rewrite [0] in english, I hoped that bad english would be saved by better mathematics. In a sense this is a joint work from M. Enguehard and M. E. Enguehard.

{\frenchspacing Je continue en fran\c cais, pour \^etre mieux compris de l'auteur de [0) tout en esp\'erant ne pas \^etre lu par mes autres lecteurs, s'il y en a. Je l'imagine alternativement d\'esinvolte et besogneux et l'ai maudit \`a la fois de sa l\'eg\' eret\'e et de son acharnement. J'ai craint de sombrer dans les m\^eme travers, ou dans la schizophr\'enie math\'ematique. Mais, qui sait? il se peut que toutes ces Propositions soient exactes.}} 
\smallskip}

\bigskip\centerline{Dedicated to Edmond Lavergne who gave me my middle name}

\vskip1truecm
{\narrower\narrower{
\smallskip\noindent ABSTRACT  (that is English summary in [0]). Let $G$ be a reductive algebraic  group over an
algebraic closure of a prime field  $\F p$, defined over $\F q$, with Frobenius
endomorphism $F$. Let $G^F$ be the subgroup of rational points. The center of $G$ is not assumed to be connected. Let $\ell$ be a
prime number, different from $p$.  If ($G^*,F)$ is in duality with $(G,F)$, then,
by a theorem of M. Brou\'e and J. Michel  [9], for any $\ell$-block --- {\bf further ``block" means ``$\ell$-block" }---
$B$ of $G^F$ there exists  a unique
$(G^*)^F$-conjugacy class
$(s)$  of
$\ell$'-semi-simple
elements such that at least one  irreducible representation of $B$ belongs to the rational Lusztig's
series $\ser{G^F}s$ associated ([18], [22]) to the $G^F$-conjugacy class $ (s)$. If $s=1$, $B$ is said to be unipotent. If $G$ is not
connected, with identity component
$G^\circ$,  define the ``unipotent
$\ell $-blocks of
$G^F$" as the $\ell$-blocks that cover some unipotent $\ell$-block of $(G^\circ)^F$. 
From $(G,F)$ and $(s)$, we construct a reductive algebraic group $(G(s),F)$
defined over $\F q$ and, assuming $\ell$ good for $G$ and {\it  some other slight restrictions on $\ell$, see  our Theorem~1.4}, a one-to-one map from the set of unipotent $\ell$-blocks of
$G(s)^F$ onto  the set of $\ell$-blocks of
$G^F$ with associated class  $(s)$   such that

if $B$ corresponds to $b$, then
there is a height preserving one-to-one map from the set of
irreducible representations  $\II b$ onto the set $\II B$,
the respective defect groups of $b$ and $B$ are isomorphic, the
associated Brauer categories are isomorphic.}
 \bigskip}

\noindent{\bf Introduction }

Let $(G, \F q, F, G^*, \ell)$ as in the abstract above, and $s$ be a semi-simple element of $ G^{*\,F}$ (the notation $G^{*\,F}$ has to be understood as $(G^*)^F$). When the center of $G$ is connected, more generally when the centralizer $\cent{G^*}s$ is connected, there is a so called Jordan decomposition of \irr\ representations which is defined by a one-to-one map between Lusztig series $$\Psi_{G,s}\colon \ser{\cent{G^*}s^F}1\to \ser{G^F}s$$
with strong properties [26], see here section~1.3. A Jordan \dec\ of blocks would associate to any block $B$ of $G^F$ a unipotent block $b$ of a related ``finite reductive group", with strong similarities between $B$ and $b$.

Let $s$ be an $\ell'$-semi-simple element in $G^{*\,F}$ such that $B$ acts trivially on some element of $\ser{G^F}s$. By a theorem of Brou\'e and Michel [9] the $G^{*\,F}$-conjugacy class of $s$ is well defined, $B$ will be said {\it in series $(s)$}. 

Assume that $\cent{G^*}s$ is a Levi subgroup of $G^*$ and let $L$ be a Levi subgroup of $G$ in the dual $G^F$-conjugacy class of the $G^{*\,F}$-conjugacy class of $\cent{G^*}s$. Then $\Psi_{G,s}$ above may be defined from the Lusztig functor $\Lu LG$ and there is a ``perfect isometry" (see [8], [9],) between a unipotent block of $L^F$ and a block $B$ of $G^F$ with $\II B=\Psi_{G,s}(\II b)$ : that is our Jordan decomposition.

A deeper result of  Bonnaf\'e and Rouquier (see [5] or [16], chapters 10, 11) say that if $\cent{G^*}s$ is contained in an $F$-stable  Levi subgroup $L^*$ of $G^*$ with dual $L$ in $G$, then $\Lu LG$ induces a Morita equivalence between the sum of the blocks of $L$ in series  $(s)$ and the sum of the blocks of $G^F$ in series $(s)$. If furthermore $L^*=\cent{G^*}s$ there is a Morita  equivalence between unipotent blocks of $L^F$ and
blocks of $G^F$ associated to $(s)$. One may hope a similar result when the connected component $\cento{G^*}s$ of $\cent{G^*}s$ is a Levi subgroup of $G^*$.

In the more general case, $\cento{G^*}s$ being a Levi subgroup of $G^*$ or not, one may construct  a reductive group $G(s)$ defined over $\F q$, in duality with $\cent{G^*}s$ --- thus $G(s)$ is  connected only if $\cent{G^*}s$ is  connected. Then, as said in the abstract above,  there is a good one-to-one map between the set of unipotent blocks of $G(s)^F$ and the set of blocks of $G^F$ in series $(s)$. Some restrictions on $\ell$ and on the type of $G$ are required by our proof (Assumption~2.1.2 in Theorem~1.4). Thanks to the classification of blocks given in [15] by Cabanes and Enguehard the proof is of combinatorial type. Several of the properties we need for the final proof are proved inductively, reducing to minimal cases. 

In Part~1 we first precise our notations on finite and algebraic groups, and morphisms. In section 1.2 are proved elementary properties of the centralizers in $G$ of finite $\ell$-subgroups. In section 1.3 we collect properties of the Jordan decomposition of \irr\ representations of $G^F$, due to Lusztig [26], [27]. Our main Theorem is Theorem~1.4. In section 1.5 we present Jordan decomposition for $2$-blocks of groups of classical type in odd characteristic, a case excluded in 1.4 and sections 2 to 4.

In Part~2 we mix the main results of Cabanes and Enguehard [13---15], where blocks are classified by so-called ``cuspidal data", see section 2.1, with Generalized Harish-Chandra theory for unipotent representations from Brou\'e-Malle-Michel [10] and Jordan decomposition. We obtain a convenient description of the set $\II B$ and a form of 
Generalized Harish-Chandra theory for representations in $\ser{G^F}s$ when $\cent{G^*}s$ is connected (2.2.1, Propositions~2.2.4, 2.3.6). In section 2.4 we study Clifford theory between blocks in a so-called regular embedding $G\to H$ defined over $\F q$, and obtain combinatorial results on blocks in relation with the non-connexity quotient of $\cent{G^*}s$, Proposition~2.4.4.

In Part~3, sections 3.1, 3.2, 3.3, the group $G(s)$, in duality with $\cent{G^*}s$ is defined by a root datum. We prove the required properties to compute combinatorial parameters in  Clifford theory for unipotent blocks between $G(s)^F$ and $G(s)^{\circ\,F}$, where $G(s)^\circ$ is the connected component of $1$ in $G(s)$. These properties are ``Non-multiplicity conditions" (Proposition~3.1.1, (B), Proposition~3.1.2, (C)) and relations between defect groups and the quotient $G(s)^F/G(s)^{\circ\,F}$ (Proposition~3.1.2, (A)-(B)). Using results of section~2.4, we obtain  a one-to-one map sending a unipotent block $b$ of $G(s)^F$ to a block $B$ in series $(s)$ of $G^F$, with a height preserving one-to-one map $\II b\to \II B$, see Propositions~3.4.1, 3.4.2.

The part 4 is devoted to  Brauer's categories of the corresponding blocks $b$ and $B$. One may identify a defect group $D$ of $B$ with a defect group of $b$ (Proposition~4.1.2). Then for any subgroup $X$ of $D$, if $(X,b_X)$ and $(X,B_X)$ are Brauer subpairs of respectively $(1,b)$ in $G(s)^F$ and $(1,B)$ in $G^F$, the quotients $\nor{G(s)^F}{X,b_X}/\cent{G(s)^F}X$ and $\nor{G^F}{X,B_X}/\cent{G^F}X$ are isomorphic (Proposition~4.2.4).

In the Appendix, section 5.1, are collected results on Clifford theory for blocks. Section 5.2  is devoted to some useful remark on unipotent Generalized Harish-Chandra series, as computed in [10] (Brou\'e-Malle-Michel).  In section~5.3, assuming Mackey decomposition formula, we deduce for classical types with connected center a commutation formula between Jordan decomposition and Lusztig functor in any series from the similar one for unipotents. In section 5.4 we show that the description of $\II {G^F}$ for $G^F={\rm SL}(\pm q)$ by Bonnaf\'e [3] and Cabanes [12] imply the existence of Generalized $d$-Harish-Chandra series.

\bigskip\noindent{\bf 1. Backgrounds. A theorem}

\bigskip\noindent{\bf 1.1. Notations and terminology }

\medskip\noindent{\bf 1.1.1. Finite groups }

The cardinality of a finite set $X$, or the {\bf order} of a group $X$, is denoted by $|X|$. The unit for multiplication in groups is denoted by $1$ and if $X$ is a group such that $|X|=1$, we may write $X=1$. Group actions on sets or modules are on the left and conjugacy may be denoted exponentially : $\lexp xy=xyx^{-1}$. The commutator of elements $x,y$ of a group is  $[x,y]=\lexp xyy^{-1}$ and, for subsets $X,Y$ of a group , $\gen X$ is the subgroup generated by $X$ and $[X,Y]$ is the subgroup generated by the set of commutators $[x,y]$ for $(x,y)\in X\times Y$. 

If  $X$ is a subset of a group acting on a set $E$, $E^X$ is the subset of fixed points. If $y\in E$ or $y\subseteq E$,  and $X$ is a group acting on $E$, $X_y$ is the {\bf stabilizer} of $y$ in $X$ : $X_y=\{g\in X\mid g.y=y\}$, with traditional exceptions if action is induced by conjugacy inside a group $X$ :  if $Y, Y'\subseteq X$, $\cent XY=\cap_{y\in Y}X_y$ is the {\bf centralizer} of $Y$ in $X$, $\cent {Y'}Y=Y'\cap\cent XY$, $\z X=\cent XX$ is the {\bf center} of $X$, $\nor X Y=\{x\in X\mid \lexp xY=Y\}$ is the {\bf normalizer} of $Y$ in $X$. We may mix the two notations and write $\nor X {Y,\la}$ for $\nor X Y_\la$.

If $\pi$ is a set of primes, $\pi'$ is its complementary in the set of all primes, any integer $n$ is a product $n=n_\pi n_{\pi'}$ where any prime divisor of $n_\pi$ (resp. $n_{\pi'}$) belongs to $\pi$ (resp. $\pi'$). Let $X$ be a  finite group, $X$ is said a $\pi$-group if $|X|=|X|_\pi$, an element $g$ of a group is said a $\pi$-element if there is an integer $n$ such that $g^n=1$ and $n=n_\pi$; any $g\in X$ is a product $g=g_\pi g_{\pi'}=g_{\pi'}g_\pi$ where $g_\pi$ (resp. $g_{\pi'}$) is a $\pi$-element (resp. $\pi'$--element). We denote  $X_\pi$  the set of $\pi$-elements of $X$, hence $X$ is a $\pi$-group if and only if $X=X_\pi$.

Let $X$ be a finite group. If $\cal O$ is a commutative ring, ${\cal O} X$ is the group algebra of $X$ on $\cal O$. We denote by $\II X$ the set of \irr\ characters of $X$, i. e. the trace maps defined by simple $\C X$-modules. When $X$ is abelian the tensor product defines a product in $\II X$, the group so obtained is denoted $X^\wedge$. Any \irr\ character of $X$ is a central function on $X$ and $\II X$ is a basis of the space $\CF X \C$ of all central functions on $X$. If one consider a field $K$ of characteristic zero and containing a $|X|$th root of unity, one recover a set ${\rm Irr}_K (X)$ which is in bijection with $\II X$ (once a bijection between $|X|$th-roots of unity is given) and is a basis of $\CF X K$.  So we omit the subscript $K$ in that case. This applies to the field $K$ of any ``$\ell$-modular splitting system"  (see [28], \S 3.6) $({\cal O},K,k)$ of $X$. Such a triple allows to introduce $\ell$-blocks inside ${\cal O} X$ ---as well in $k X$---, called {\bf blocks of $X$},  and blocks idempotents, i.e. primitive idempotents in the center of ${\cal O} X$, see [28], \S 1.8. By extension from ${\cal O}$ to $K$ one obtains a partition
$$\II X=\cup_{B}\II B$$ where $B$  ranges over the set of blocks of $X$. Define, for $\xi\in\II X$, the block  $B_X(\xi)$ by $$\xi\in\II{B_X(\xi)}\,.$$

On the space $\CF X K$  is defined the usual scalar product 
$$\scal \phi\psi X=|X|^{-1}\sum_{g\in X}\phi(g)\psi(g^{-1})$$
and,  if $\phi,\psi\in\II X$, then $\scal\phi\psi X=\delta_{\phi,\psi}$. The associated {\bf norm} is denoted $||?||$, i.e.  $||\phi||^2=\scal \phi\phi X$ for $f\in\CF X K$.

Given a morphism $\sigma\colon Y\to X$ (or simply an inclusion of groups $Y\subseteq X$), the {\bf restriction} from  $X$ to  $Y$, applied to representations or to central functions,  is denoted $\Res{}\sigma$ (or $\Res XY$). {\bf Induction} from a subgroup $Y$ to a group $X$ is denoted $\Ind YX$. If $Y$ is a normal subgroup of $X$ and $\eta\in\II Y$ one defines 
$$\II{X\mid\eta}=\{\chi\in\II X\mid \scal{\Res XY\chi}\eta Y\neq 0\}\,.$$
Then we say that $\chi$ {\bf covers} $\eta$. Similarly a block $b$ of $X$ {\bf covers} a block $c$ of $y$ if there exist $\eta\in\II c$ and $\chi\in\II b$ such that $\chi$ covers $\eta$ [28], 5.1.

The notation ``tensor product" of representations or of central functions is used to produce representations or central functions on a central product, as well with one fixed group : if $\chi\in\II X$ and $\xi\in\II Y$, and $Z$ injects in the centers of $X$ and $Y$, if furthermore $\chi/\chi(1)$ and $\xi/\xi(1)$ are equal on $Z$, $\chi\otimes\xi$ may be considered as an element of  $\II{X\times_Z Y}$. But when $X=Y$ and $\chi/\chi(1)$ and $\xi/\xi(1)$ are conjugate complexes  on $Z$, $\chi\otimes \xi$ may be considered as an element of  $\II{X/Z }$ by restriction of the preceding one to the diagonal subgroup. We hope the good interpretation is given by 
the context.
 
\medskip\noindent{\bf 1.1.2. Algebraic groups }

 All along  $G$ is an algebraic group. The connected component of $1$ is denoted $G^\circ$, but we prefer 
 
 \medskip\noindent{\bf 1.1.2.0.  Notations. }{\sl $\cento G g$ is $\cent Gg^\circ$,  $\zo G$ is $\z G^\circ$, $G^{\circ\, F}$ is $(G^\circ)^F$, $G^{* \,F}$ is $(G^*)^F$.}
 
 \smallskip  Let $p$ be a prime number, different from $\ell$. We consider first connected reductive groups $G$ on an algebraic closure $\F{}$ of a prime field $\F p$, 
that are defined over a finite field $\F q$ ($q$ a power of $p$) thanks to an endomorphism $F\colon G\to G$. 

If $T$ is a maximal torus in $G$, then is defined a {\bf root datum} $$\FDD( G,T):=(X(T),Y(T),\Phi,\Phi^\vee)$$ (group of characters of $T$, group of one parameter subgroups of $T$, set of roots in $\R\otimes_\Z X(T)$, set of coroots) and that root datum defines $G$ up to  interior isomorphisms induced by $T$. To emphasize the choice of $T$ or/and to recall the algebraic group $G$, we may write $\Phi_G(T)$ instead of $\Phi$. The types of connected Dynkin diagrams of root system, or types of \irr\ algebraic groups are denoted as usually $\AA_n$, $\BB_n$,\dots, $\EE_8$. Assuming that $T$ is $F$-stable, $F$ acts on each of the 4 objects of the root datum, hence acts on the set of connected components of the Dynkin diagram of the root system $\Phi$. To every orbit of $F$ there corresponds an $F$-stable component of $[G,G]$, minimal as $F$-stable component, with a so-called {\bf rational type} $({\bf X},q^m)$, defined by one of the above types, twisted or not (so $\lexp 2\AA_n$, $\lexp 2\DD_n$, $\lexp 3\DD_4$,  $\lexp 2\EE_6$ appear) and an extension $\F{q^m}$ of $\F q$. $G$ is said {\bf rationally \irr}  if $[G,G]$ has only one such component. In reference  to ``Ennola's conjecture" we write $(\AA_n,-q^m)$ instead of $(\lexp 2\AA_n,q^m)$.

A {\bf dual root datum} of $\FDD
(G,T)$ is isomorphic to $(Y(T),X(T),\Phi^\vee,\Phi)$. The pair $(Y(T),X(T))$  defines a so-called {\bf dual torus} $T^*$, i.e. $X(T^*)=Y(T)$ and $Y(T^*)=X(T)$,  and the dual root datum defines, up to some isomorphisms, an algebraic reductive group usually denoted $G^*$, said to be in duality with $G$ :
so $T^*$ is a maximal torus in $G^*$ and $\FDD(G^*,T^*)\cong (Y(T),X(T),\Phi^\vee,\Phi)$. We may write then that {\bf $(T\subseteq G)$ and $(T^*\subseteq G^*)$ are in duality}. Any maximal torus of $G^*$ is in duality with some maximal tori of $G$, more generally any Levi-subgroup of $G^*$ is in duality with some Levi subgroup of $G$. The action of $F$ on $\FDD(G,T)$ gives an action on $\FDD(G^*,T^*)$, hence some endomorphism of $G^*$, we denote $F$ for simplicity. Then a $G^{*\,F}$-conjugacy class of $F$-stable Levi subgroups $L^*$ of $G^*$ corresponds to a $G^F$-conjugacy class of $F$-stable Levi subgroups $L$ of $G$ so that $L$ and $L^*$ may be defined by dual root data (the duality is not uniquely defined but is coherent with duality between $G$ and $G^*$, see formulas (1.3.2.1), (1.3.2.4), (1.3.2.5)). We say then that $(L,F)$ and $(L^*,F)$ (or $L$ and $L^*$) are in {\bf dual conjugacy classes}. The Weyl group $W(G,T):=\nor G T/T$ is a group of automorphisms of $\FDD(G,T)$, so is anti-isomorphic to $W(G^*,T^*)$.
Given once for all an isomorphism $\F{}^\times\to (\Q/\Z)_{p'}$ and an imbedding $\F{}^\times \hookrightarrow \bar \Q_\ell^\times$, one obtains isomorphisms $T^F\cong (T^{*\,F})^\wedge$, $(T^F)^\wedge\cong T^{*\,F}$.

Let $\pi\colon G_{\rm sc}\to [G,G]$ be a simply connected covering of the derived group of $G$. The kernel of $\pi$ is denoted $\cF (G)$. If $G$ is adjoint, $\cF (G)$ is the {\bf fundamental group} common to groups of same  type.

\medskip\noindent{\bf 1.1.3. Proposition. }{ \sl (a) Let $G$ and $G^*$ be in duality. The finite abelian groups $\cF(G)$ and $\z {G^*}/\zo {G^*}$ are in duality, hence are isomorphic.

(b) Let $L$ be a Levi subgroup of $G$. Then $\z L/\zo L$ is isomorphic to a quotient of $\z  G/\zo G$ and $\cF (L)$ to a subgroup of $\cF(G)$.} 

\preuve (a) Let $\FDD(G,T)=(X,Y,\Phi,\Phi^\vee)$ and $\FDD(G^*,T^*)$ be root data in duality, defining the duality between $G$ and $G^*$. Let $\Omega$ be the group of weights, dual over $\Z$ of $\gen \Phi$. Then $\cF (G)=\cF([G,G])$ is is isomorphic to the dual (as a finite abelian group ...) of the $p'$-torsion group of $\Omega/(X(T)\cap[G,G])$ hence isomorphic to the $p'$-torsion group of $Y(T\cap[G,G])/\Z\Phi^\vee$ (one has $\Z\Phi^\vee\subseteq Y(T\cap[G,G])\subseteq Y$). One  knows that the finite group $\z {G^*}/\zo {G^*}$ is isomorphic to the dual of the $p'$-torsion of $Y/\Z\Phi^\vee$ ([17] 4.5.8). As $Y/Y(T\cap [G,G])$ has no torsion, the torsion groups of $Y/\Z\Phi^\vee$ and  $Y(T\cap[G,G])/\Z\Phi^\vee$ are isomorphic. 

Hence our claim.

(b) We may assume that $T\subseteq L$ and $L$ is defined by a subroot datum of $\FDD(G,T)$, so that the set $\Phi_L$   of roots of $L$ is contained in $\Phi$. As $\Z\Phi/\Z\Phi_L$ has no torsion, the torsion group of $X/\Z\Phi_L$ injects in the torsion group of $X/\Z\Phi$, hence $\z L/\zo L$ is isomorphic to a quotient of $\z  G/\zo G$.
\bull

\medskip\noindent{\bf 1.1.4. On morphisms }

{\sl (a) A morphism $$\sigma\colon G\to H$$ between algebraic reductive groups is said to be {\bf isotypic} if its kernel is central in $G$ and $[H,H]=\sigma([G,G])$. If $ G, H$ and $\sigma$ are defined on $\F{q}$ by some endomorphisms $F$, we may write $\sigma\colon (G,F)\to (H,F)$. Between groups $G^*$, $H^*$ in duality with resp. $G$ and $H$, there exists ``dual morphisms" $$\sigma^*\colon H^*\to G^*\, ,$$  where $\sigma^*$ is isotypic and eventually defined over $\F q$.  

Then $\cF(G)$ (resp. $\cF(H^*)$) is isomorphic to a subgroup of $\cF(H)$ (resp. $\cF(G^*)$) and $\z{G^*}/\zo{G^*}$ (resp. $\z H/\zo H$) is isomorphic to a quotient of  $\z{H^*}/\zo{H^*}$ (resp.  $\z G/\zo G$).

(b) An isotypic morphism $\sigma\colon G\to H$ with kernel $\{1\}$ will be called an {\bf embedding}. Then $\cF (G)$ is isomorphic to $\cF(H)$. If $\sigma^*\colon H^*\to G^*$ is a dual morphism of an embedding $\sigma$, then $\sigma^*(H^*)=G^*$ and $\z{G^*}/\zo{G^*}$  is isomorphic to $\z{H^*}/\zo{H^*}$.

 If furthermore the center of $H$ is connected, $\sigma$ is called a {\bf regular embedding}. An isotypic morphism $\sigma\colon G\to H$ is a regular embedding if and only if one (and then every) dual morphism $\sigma^*\colon H^*\to G^*$ satisfies : $\sigma^*(H^*)=G^*$, the kernel of $\sigma^*$ is a torus  and $\cF(H^*)=\{1\}$.
 
 (c) For any $(G,F)$ there exists a regular embedding $\sigma\colon (G,F)\to (H,F)$. Given two regular embeddings $\sigma_j\colon (G,F)\to (H_j,F)$ ($j=1,2$) there exists a third one $(G,F)\to (H,F)$ that factorizes through $\sigma_1$ and $\sigma_2$.
 
 (d) For any isotypic morphism $\sigma\colon G\to H$ (resp . $(G,F)\to (H,F)$) there exist regular embeddings $G\to G_0$, $H\to H_0$ (resp.$(G,F)\to (G_0,F)$, $(H,F)\to (H_0,F)$) and an isotypic morphism $\sigma_0\colon G_0\to H_0$ (resp. $(G_0,F)\to (H_0,F)$) whose kernel is a torus and extending $\sigma$ (resp. $(\sigma,F)$).
 }  
 
\medskip\noindent{\it On proofs. } The existence of $\sigma^*$ in (b) follows from the existence of a morphism of dual  root data, so $\sigma^*$ is not unique. The construction of regular embeddings is based on products as follows :

 Let $\sigma\colon  (G,F)\to (H,F)$ be isotypic as in (d). 
Recall the center of $G$ is contained in any maximal torus. Let $T$ be an $F$-stable subtorus of $G$ that contains  the kernel $K$  of $\sigma$. Let 
$$G_1=G\times_KT=G\times T/\{(k,k^{-1})\mid k\in K\}$$ 
The composed map $G\to G\times T\to G_1$  is an embedding (if $K=\z G$, $\z {G_1}$ is a torus hence that embedding is regular, a way to prove (c)). Then $\sigma$ is the restriction of some $\sigma_1\colon G_1\to H$, where $\Ker \sigma_1$ is the central torus $T_1:=K\times_KT$, isomorphic to $T$. A regular embedding $G_1\to G_0$ defines a regular embedding   $H\to
H_0:=G_0/T_1$ such that the following diagram of isotypic morphisms is commutative
$$\def\normalbaselines
{\baselineskip20pt\lineskip3pt
\lineskiplimit3pt}\matrix
{K&&T_1&\mapright{1}&T_1\cr
\mapdown{}&&\mapdown{}&&\mapdown{}\cr
G&\hookrightarrow{}&G_1&\hookrightarrow{}&G_0\cr
\mapdown{\sigma}&&\mapdown{\sigma_1}&&\mapdown{\sigma_0}\cr
H&\mapright{1}&{H}&\hookrightarrow{} &H_0\cr}   $$
\bull

 \medskip\noindent{\bf 1.1.5. $\ell$,  $q$ and $d$. The decomposition $G=\a G\b G$}
 \smallskip We have assumed once for all
that $\ell$ is a prime number, $\ell\neq p$.
 
 \smallskip \noindent(1.1.5.1) $\quad\quad $ If $\ell$ is odd, let $d=d_{q,\ell}$ be the order of $q$ mod $\ell$ and
 $E=\{n\in \N\mid n_{\ell'}=d\}$\hfill
 
 $E$ is denoted $E_{q,\ell}$ in [16] Theorem~21.7. 
 
 For any positive integer $a$, let $\phi_a\in\Z [X]$ be  the $a$-th cyclotomic polynomial. If $\ell>2$ and $a\in\N^\times$, then  $\ell$ divides $\phi_a(q)$ if and only if $a\in E$. Any algebraic group $(G,F)$ we consider, if defined on $\F q$, has a {\bf  polynomial order} $P_{G,F}(X)\in \Z [X]$ such that $|G^F|=P_{G,F}(q)$ [16] Section~13.1. If $G$ is a torus its polynomial order is a product of cyclotomic polynomials. Let $A$ be a set of positive integers. A {\bf $\phi_A$-subgroup ($\phi_a$-subgroup if $A=\{a\}$) $S$ of $(G,F)$} is an $F$-stable torus in $G$ such that $P_{S,F}$ is a product of  $\phi_a$ where $a\in A$. One has ``Sylow's theorems" for $\phi_A$-subgroups in $G$. For any torus $T$, we denote $T_{\phi_A}$ its maximal $\phi_A$-subgroup. By definition of $E$ with respect to $\ell$ one has $T^F_\ell\subseteq T_{\phi_E}$. A Levi subgroup of $(G,F)$ is said to be {\bf $A$-split} ({\bf $a$-split} if $A=\{a\}$) if it is the centralizer in $G$ of a $\phi_A$-subgroup. 
 For more results on these notions see [10] or [16] Chapter 13.
 
  \medskip\noindent{\bf 1.1.5.2. Definition. }[16] 22.4, 22.5.  {\sl Let $(G,F)$ be defined on $\F q$ as above and assume $\ell>2$. Let $\a G$ be the product in $G$ of $\zo G$ and all rationally \irr\ components of $[G,G]$ of type $(\AA_n,r)$ ($r=\pm q^m$) where $\ell$ divides $|r-1|$. Let  $\b G$ be the  product in $G$ of all rationally \irr\ components of $[G,G]$ which are not included in $\a G$.
  
  One has $G=\a G.\b G$ (central product), $\z{\b G}^F$ and $G^F/\a G^F.\b G^F$ are commutative $\ell'$-groups.}
  
  The choice of components of  $\b G$ is made so that $\cF ({\b G})^F_\ell=1=(\z {\b G}/\zo{\b G})^F_\ell$. For any central $F$-stable subgroup $A$ of $G$ one has $\a{(G/A)}=\a G/A\cap\a G$ and $\b{(G/A)}=\b G/A\cap\b G$. For any $F$-stable Levi subgroup $L$ of $G$, one has $L\cap\a G\subseteq \a L$ and $\b L\subseteq\b G$. Currently $\b L\neq L\cap\b G$ nevertheless one has $\z{L\cap \b G}^F_\ell=\zo{L\cap\b G}^F_\ell$ by Proposition~1.1.3. In inductive proofs the following properties are frequently used  [16] Proposition~22.5, Theorem~22.2:
  
   \medskip\noindent{\bf 1.1.5.3.  }{\sl Assume $\ell>2$. If $Y$ is an $\ell$-subgroup of $G^F$ such that $\z{\cent {G^F}Y}_\ell\subseteq \z G\a G$, then $Y\subseteq \a G$.
   
   If $G=\b G$ then any proper $E$-split Levi subgroup of $G$ is contained in a proper $d$-split Levi subgroup of $G$.}

  The last assertion is an immediate consequence of an easy to verify fact when $G=\b G$ : if $L$ is an $F$-stable Levi subgroup of $G$ such that $\zo L_{\phi_d}\subseteq \zo G$, then $\zo L_{\phi_E}\subseteq \zo G$.
  
  Let $(G^*,F)$ be in duality with $(G,F)$. Then one has isotypic morphisms $$\a{(G^*)}\to (\a G)^* ,\quad \b {(G^*)}\to (\b G)^*\leqno{(1.1.5.4)}$$
 
 \bigskip\noindent {\bf 1.2. Centralizers and connexity }
 
 For any subgroup $X$ of $G$ or $x\in G$ we denote
 $${\rm A}_G(X)=\cent G X/\cento GX,\quad {\rm A}_G(x)={\rm A}_G(\gen x)\leqno{(1.2.0)}$$  Proposition~1.1.3 (a) gives a relation between ${\rm A}_G(G)$ and $\cF(G^*)$. Proposition~1.2.6 below relies ${\rm A}_G(X)$ to $\cF( G)$ and $\cF(G^*)$ when $X$ is a finite $\ell$-subgroup of $G$.
  
 We first note some elementary results for later use.
 
 \medskip\noindent{\bf 1.2.1. Lemma. }{\sl  Let $\rho\colon H\to K$ be a surjective morphism between
groups,  whose kernel $\Ker \rho$ is central in $H$. Let $X$ be a  subgroup of
$H$ with finite exponent. 

The exponent of
$\rho^{-1}(\cent K{\rho(X)}/\cent HX$ divides the exponent of $X$ and the exponent of $\Ker
\rho\cap [H,H]$. If $\pi$ is a set of primes,  $X$ is a finite $\pi$-group and
$\Ker
\rho\cap [H,H]$ is a finite $\pi'$-group, then $\rho(\cent HX)=\cent
K{\rho(X)}$.}

\preuve  The last assertion is a direct consequence of the first one.

 One defines a bi-morphism $(\rho^{-1}(\cent K{\rho(X)})/\cent HX)\times X\to \Ker
\rho\cap [H,H]$ by
 restriction of $(h,x)\mapsto [h,x]$ : indeed, when  $[h,x]\in\z H$ --- and this happens if  $\rho(h)\in\cent K{\rho(X)}$ and
$x\in X$ --- one has, for $h', y\in H$, $[h,x][h,y]=[h,xy]$ and
$[h,x][h',x]=[hh',x]$.
Let $k$ be the exponent of $X$, then $[h^k,x]=[h,x]^k=[h,x^k]=1$ so that $h^k\in\cent HX$. 
\bull
\vfill\eject 
 \noindent{\bf 1.2.2. Proposition. }{\sl Let $X$ be a subgroup of $G$.
 
 (a) If $X.\zo G=X'.\zo G$  one has ${\rm A}_G(X)={\rm A}_{G}(X')$. One has  exact sequences of morphisms
 $$1\to \zo
G/(\zo G\cap
\cento{[G,G]}X)\to \cent{G}X/\cento{[G,G]}X\to {\rm A}_G(X)\to 1$$
$$1\to \cent{[G,G]}X/\cento{[G,G]}X\to \cent{G}X/\cento{[G,G]}X\to \zo G/(\zo
G\cap [G,G] )\to 1
$$
hence $|\cent{[G,G]}X/\cento{[G,G]}X|=|{\rm A}_G(X)|.|\zo G\cap[G,G]/\zo G\cap\cento{[G,G]}X|$.

 (b) If $H$ is an algebraic subgroup of $G$ and $X\subseteq H$, there exists an exact sequence of morphisms $$1\to \cento GX\cap H/\cento HX\to {\rm A}_H(X)\to {\rm A}_G(X)\to \cent GX/\cento GX.\cent H X\to 1$$
 
 (c) If $X$ is a subgroup of $Y$, there exists an exact sequence of morphisms
 $$1\to \cento GX\cap\cent GY/\cento
GY\to {\rm A}_G(Y)\to {\rm A}_G(X)\to \cent GX/\cento GX.\cent GY\to 1$$ }
\preuve The first assertion in (a) follows from the equality $G=\zo G.[G,G]$.

The four exact sequences are given by isomorphisms theorems, knowing that $\cento HX\subseteq \cento GX$ and $\cent HX=H\cap \cent GX$ in (b) and that $\cento GY\subseteq \cento GX$ and $\cento GX.\cent GY\subseteq \cent GX$ in
(c).
\bull

\medskip\noindent{\bf 1.2.3. Proposition. }{\sl  Let $X$ be a finite $p'$-subgroup of $G$,

(a) Let $Z$ be a central subgroup of $G$ 
and $\rho\colon G\to G/Z$
 the quotient morphism. Then $$\rho(\cento
GX)=\cento{G/Z}{\rho(X)}$$ so that ${\rm A}_G(X)$ is isomorphic to a subgroup of ${\rm A}_{\rho(G)}(\rho(X))$.

(b) Let $\sigma\colon G\to H$ be an embedding. Then ${\rm A}_H(\sigma(X))$ is isomorphic to a quotient of ${\rm A}_{G}(X)$.

(c) An isotypic morphism $\sigma\colon G\to H$ defines a morphism ${\rm A}_G(X)\to {\rm A}_H(\sigma(X))$.}

\preuve (a) The property is well known in case 
$X$ is cyclic or more generally contained in a torus,  $\cento GX$ being described by a root datum [16] Proposition~13.13. Assume $X$ is generated by a finite set of semi-simple elements
$x_j$, one has 
$$\cento GX\subseteq \cap_j\cento G{x_j}\subseteq \cent
GX$$  
and indices are finite. Applying $\rho$ we obtain, knowing that $Z\subseteq\cent G{x_j}$,
$$\rho(\cento {G}X)\subseteq \rho(
\cap_j\cento G{x_j})\subseteq  \cap_j\rho(\cento G{x_j})\subseteq  \cap_j\rho(\cent G{x_j})=\rho(\cap_j\cent G{x_j})=\rho(\cent GX)$$ 
with finite indices, and, on $G/Z$-side
$$\cap_j\rho(\cento G{x_j})=\cap_j\cento{G/Z}{\rho(x_j)}\subseteq
\cap_j\cent{G/Z}{\rho(x_j)}=\cent{G/Z}{\rho(X)}$$ 
with finite indices. But
$\rho(\cento GX)$ is connected hence $\rho(\cento GX)=\cento{G/Z}{\rho(X)}$. 

(b) This is a special case of (a) in 1.2.2. Directly : one has $H=\zo H.\sigma(G)$, hence $\cent H{\sigma(X)}=\zo H.\cent{\sigma(G)}{\sigma(X)}$ and $\cento H{\sigma(X)}=\zo H.\sigma(\cento GX)$. By isomorphism theorems ${\rm A}_H(\sigma(X))$ is a quotient of ${\rm A}_{\sigma(G)}(\sigma(X))\cong {\rm A}_G(X)$.

(c) Any isotypic morphism may be obtained by composition of an embedding and a quotient morphism.
\bull

As we are interested by centralizers of $\ell$-subgroups and $\ell$ is different from the characteristic of the base field, the key-property is given by one of Steinberg's theorems in [31], essentially that ${\rm A}_G(s)$ is 1 if $s$ is semi-simple and $G$ simply connected.

\medskip\noindent{\bf 1.2.4. Proposition.
}{\sl  (a) Let $s$ be a semi-simple element of $G$. The group ${\rm A}_{G}(s)$ is
isomorphic to a subgroup of ${\cF}(G)$ and its exponent divides the exponent of $s$.

(b) If $L$ is a Levi subgroup of $G$ and $s\in L$, $A_L(s)$ is isomorphic to a subgroup of $A_G(s)$.

(c) Let $s$ and $t$ be commuting semi-simple elements of $G$ with coprime orders. One has $$\cento Gs\cap\cento Gt=\cento G{st}$$ hence
${\rm A}_{G}(st)$ is isomorphic to ${\rm A}_{\cento Gs}(t)\times {\rm A}_{\cento Gt}(s)$.

(d) If 
$\cent Gs$ is connected for any semi-simple $\ell$-element
$s$ in $ G$, then ${\cF}(G)$ has order prime to $\ell$.

(e) Assume that ${\cF}(G)_\ell=\{1\}$. If the order of  a semi-simple element $s$ of $G$  is prime to $\ell$, $|\cF(\cento G s)|$ is prime to $\ell$.}

\preuve (a) See [31] 9.5 and Proposition~1.2.3.

(b) is a special case of (b) in Proposition~1.2.2, with $(\gen s,L)$ instead of $(X,H)$ : one has $L=\cent G{\zo L}$ and $L\cap\cento Gs=\cent{\cento G s}{\zo L}$.

(c) If $\cF(G)=\{1\}$, then $\cento G{st}=\cent G{st}=\cent Gs\cap\cent Gt=\cento Gs\cap\cento Gt$, thanks to (a). 

When  $\cF(G)\neq\{1\}$ let $\rho\colon H\to G$ be a covering such that $\cF (H)=\{1\}$ (see Proposition~1.1.4.1). There exist semi-simple elements $s',t'$ in $H$ with coprime orders such that $\rho(s')=s$, $\rho(t')=t$ so that $\rho(\cento H{s'})=\cento Gs$, $\rho(\cento H{t'})=\cento Gt$, $\rho(\cento H{s't'}=\cento G{st}$. One obtains $\cento Gs\cap\cento Gt=\cento G{st}$. But $\cento Gs\cap \cent Gt=\cent{\cento Gs}t$ and $\cent Gs\cap \cento Gt=\cent{\cento Gt}s$. The isomorphism  we claim follows.

(d) Let $\rho\colon H\to [G,G]$ be a simply connected covering and
 $G_1=H/{\cF}(G)_{\ell}$. Then ${\cF}(G_1)$ is isomorphic to ${\cF}(G)_{\ell}$ and $\rho$  
factors through
$\rho_1\colon G_1
\to G$. By (a) if 
  $s_1\in G_1$ is semi-simple and of order prime to $\ell$,
then $\cent{G_1}{s_1}$ is connected.  If $s_1\in G_\ell$,
then $\rho_1(\cent{G_1}{s_1})\subseteq \cent G{\rho_1(s_1)}=\cento
G{\rho_1(s_1)}=\rho_1(\cento{G_1}{s_1})$, hence $\cent{G_1}{s_1}$ is connected. 
Any semi-simple element $t$ of $G_1$ is a product $t=t_\ell.t_{\ell'}$ and, using (b) one has 
$\cent {G_1}t=\cent {G_1}{t_\ell}\cap \cent {G_1}{t_{\ell'}}=\cento
{G_1}{t_\ell}\cap \cento {G_1}{t_{\ell'}}=\cento Gt$. By [31]
9.9 we have ${\cF}(G_1)=\{1\}$.

(e) Let $t\in\cento Gs_\ell$. By (a) $\cent Gt$ is connected. We have $\cent{\cento Gs}t=\cento Gs\cap\cento Gt=\cento G{st}$ by (b). By (c) $\cF (\cento Gs)_\ell=\{1\}$.
\bull

When $\ell$ is good for $G$  (see [16]
section 13.2) and $s\in G_\ell$, then $\cento Gs$ is a Levi subgroup of $G$. If furthermore $\z G^F_\ell=\zo G^F_\ell$, then $s\in\zo{\cento Gs}$ (see Proposition~1.1.3)  hence  $\cento Gz=\cent G{\zo{\cento Gs}_{\phi_E}}$. Thus $E$-split Levi subgroup of $G$ are good examples of (connected) centralizers of abelian $\ell$-subgroups of $G^F$. In the following Proposition, (b) applies when $G=\zo G\b G$ :

 \medskip\noindent{\bf 1.2.5. Proposition. }[15] Propositions~2.4, 3.2, [16] Proposition 13.19. {\sl Assume the prime $\ell$ is good for $G$. 
 
 (a) If an $F$-stable Levi subgroup $K$ of $G$ satisfies $K=\cento G{\z K^F_\ell}$, then it is $E$-split.
 
 (b) Assume that $\ell$ does not divides $(\z G/\zo G)^F$. Let $S$ be a $\phi_{E}$-subgroup of $G$. Let $K=\cent G S$ be the associated $E$-split Levi subgroup of $G$  (1.1.5). Then
 
 (b.1) $\zo K^F_\ell=\z K^F_\ell$.
 
 (b.2) $K=\cento G{S^F_\ell}=\cento G{\z K^F_\ell}$ and $K^F=\cento G{\z K^F_\ell}^F=\cent{G^F}{\z K^F_\ell}$.
 }
 
\medskip\noindent{\bf 1.2.5.1. On Proposition 1.2.5. }   There exist examples of $E$-split Levi subgroup $K$ that does not satisfy the equality $K=\cento G{\z K^F_\ell}$ as in assertion (b.2) of Proposition~1.2.5, even with $\ell$ good, but they are deduced from one single extremal case, we want to describe here. 

Indeed a good prime may divide the order of $\z G/\zo G$ only in rational types $\AA$, $\lexp 2\AA$. Assume  $K$ is an $E$-split Levi subgroup of $G$ such that $K\neq \cento G{\z K^F_\ell}$. Let $M=\cento G{\z K^F_\ell}$. Then $M$ is an $E$-split Levi subgroup of $G$ by assertion (a). Furthermore $K$ is an $E$-split Levi subgroup of $M$ and $\z K^F_\ell\subseteq \z M$. So we are reduced to check the $E$-split Levi subgroups $K$ of $G$ with $\z K^F_\ell\subseteq \z G$.

Let $K_1=K\cap [G,G]$. Then $K_1$ is an $E$-split Levi subgroup of $[G,G]$ and $\z{K}^F_\ell\subseteq \zo G.\z{K_1}^F_\ell$, hence $\cento G{\z K^F_\ell}=\zo G.\cento{[G,G]}{\z{K_1}^F_\ell}$ so that ${\z {K_1}^F_\ell}\subseteq \z{[G,G]}$. The same situation occurs in $[G,G]$, we may assume $G$ semi-simple. In  a simply connected covering  $\pi\colon \hat G\to G$, $\hat K=\pi^{-1}(K)$ is an $E$-split Levi subgroup  of $\hat G$ and $\pi(\z{\hat K}^F_\ell)\subseteq  \z K^F_\ell$ so that $G=\cento G{\z K^F_\ell}\subseteq \pi(\cento{\hat G}{\z {\hat K}^F_\ell})$, hence $\z {\hat K}^F_\ell\subseteq \z {\hat G}$. Assume $G$ simply connected. $G$ is a direct product of rationnally \irr\ components and the inclusion occurs  in each component. Assume $(G,F)$ of type $\AA_n(\epsilon q)$, $\epsilon \in\{-1,1\}$. The $E$-split Levi subgroups of such a group are easy to describe. Recall that the $E$-split Levi subgroups of $(G={\rm GL}_n,F)$ where $G^F={\rm GL}_n(\epsilon q)$ are direct products  of the form $S\times \AA_{N}(\epsilon q)$ where the polynomial degree of the torus $S$ may be written $\prod_j(X^{d\ell^{\al(j)}}-1))$ so that $|S^F|=\prod_j|(\epsilon q)^{d\ell^{\al(j)}}-1|$ and $n=N+d.\sum_j\ell^{\al(j)}$, with the special case $N=0$ and $|S^F|=\prod_j|(\epsilon q)^{d\ell^{\al(j)}}-1|/(q-\epsilon)$. One sees that, when $G^F={\rm SL}_n(\epsilon q)$, the condition $\z K^F_\ell\subseteq \z G$ is satisfied only when $\ell$ divides $(q-\epsilon)$, $n$ divides $(q-\epsilon)_\ell$ and $K$ is a so called ``Coxeter torus",  that is a maximal $F$-stable torus $T$ such that $|T^F|=|(\epsilon q)^n-1/\epsilon q-1|$ (Coxeter tori are $G^F$-conjugate).

\medskip

In the following Proposition on centralizers of $\ell$-subgroups, we assume $\ell$ good and use the \dec\ $G=\a G.\b G$ in an inductive proof. But in view of a similar results for $2$-groups in classical types, see Proposition~1.2.7, we notice that the hypothesis ``$\ell$ good" is used only in parts (B) and (D) of the proof.

\medskip\noindent{\bf 1.2.6. Proposition. }{\sl  Assume $\ell$ is good for $G$. Let
$Y$ be a finite
$\ell$-subgroup of
$G$.

(a)
$\cento GY$ is an algebraic reductive group. 

(b)  For any (eventually $F$-stable if $Y\subseteq G^F$) maximal torus $T_Y$ of $\cento GY$, there exists
a maximal  (eventually $F$-stable)  torus
$T$ of
$G$ such that $T_Y\subseteq T$ and $Y\subseteq\nor GT$. Thus
${\rm W}(\cento GY,T_Y)$ is isomorphic
to a subgroup of ${\rm W}_{G}(\cent G{T_Y})\cong\nor G{\cent G{T_Y}}/\cent G{T_Y}$. 

(c) Assume $Y\subseteq G^F$. Let $(T^*\subseteq G^*,F)$ in duality with $(T\subseteq G,F)$. There exists a finite 
$\ell$-subgroup $Y'$ of $\nor{G^*}{T^*}$ and a maximal $F$-stable torus $T^*_{Y'}$ of $\cento{G^*}{Y'}$  such that $T^*_{Y'}\subseteq T^*$ and the Levi subgroups $\cent G{T_Y}$ of $G$
and
$\cent{G^*}{T^*_{Y'}}$ of $G^*$ are in dual conjugacy classes. Moreover $(T^*_{Y'}\subseteq\cento {G^*}{Y'},F)$ and $(T_Y\subseteq\cento GY,F)$ are in duality. The groups
$\nor {G} {\cento GY,T_Y}$ and $\nor {G^*}
{\cento{G^*}{Y'},T^*_{Y'}}$ act on the root datum $\FDD(\cento GY,T_Y)$ by contragredient actions. 

(d) The group $\z{\cento GY}/\zo
{\cento GY}$ is isomorphic to a quotient of $\z G/\zo G$ and $\cF(\cento GY)$ is
isomorphic to a subgroup of
$\cF(G)$, with $F$-action when $Y\subseteq G^F$. 

(e) ${\rm A}_G(Y)$ is an $\ell$-group. If ${\cF}(G)_\ell=\{1\}$ and $(\z G/\zo G)_\ell=\{1\}$ then
${\rm A}_G(Y)=\{1\}$.  If $Y\subseteq G^F$,
${\cF}(G)^F_\ell=\{1\}$ and $(\z G/\zo G)^F_\ell=\{1\}$ --- and that is the case when $G=\zo G\b G$ ---  then
${\rm A}_G(Y)^F=\{1\}$. 

(f) If $L$ is a  Levi subgroup of $G$ and $Y\subseteq L$, ${\rm A}_L(Y)$ is
 isomorphic to a subgroup of ${\rm A}_G(Y)$. If furthermore $Y\subseteq G^F$, ${\rm A}_L(Y)^F$ is
 isomorphic to a subgroup of ${\rm A}_G(Y)^F$.
  }

\medskip\noindent{\it Comments. }  In case $Y$ is abelian, specially  contained in a torus, as $\ell$ is good $\cento G Y$ is a Levi subgroup of $G$ and all these properties are well known (see [16] Proposition~13.16 and [13] 2.1). There are similar results for some automorphisms  : if $\sigma$ is a quasi-simple automorphism and $(T_\sigma)$ is a maximal torus of $(G^\sigma)^\circ$, then $\cent G{T_\sigma}$ is a maximal and $\sigma$-stable torus of $G$ (see [31] or [21] Theorem~1.8).

When $\ell$ is not good (c) and (d) may fail.

\medskip\noindent{\it Proof of Proposition~1.2.6. }

Note that any assertion refering to some Frobenius $F$ in the Proposition applies to any finite $\ell$-subgroup $Y$ of $G$ thanks to the fact that there exists a power  $F_1$ of $F$ such that $Y\subseteq G^{F_1} $. Thus to prove (a), (b), (d) and (f) one may assume that $Y\subseteq G^F$. In (e) it is sufficient to prove that 

\centerline{$Y\subseteq G^F$,
${\cF}(G)^F_\ell=\{1\}=(\z G/\zo G)^F_\ell$ implies ${\rm A}_G(Y)^F=\{1\}$.}

In (c) it is said that the morphism $X(T)\to X(T_Y)$ given by the inclusion  $T_Y\subseteq T$ (see (b)) is direct, with  a section ${\rm Y}({T_Y}^*)\to {\rm Y}(T^*)$ that is a  $\nor{{\rm
W}(G,T)}{{\rm W}(\cent G{T_Y},T)}$-morphism. Note that if $Y\subseteq \nor GT$, then
$Y$ is a split extension of $Y\cap T$ by $Y/Y\cap T$, so $Y$ is defined by the root datum $\FDD(G,T)$. The torus $T_Y$ is equally defined by $\FDD(G,T)$ : one has $T_Y=(T^Y)^\circ=T\cap\cento GY$.

On the connexity of centralizers of finite nilpotent $p'$-subgroups (resp. $\ell$-subgroups) assertion (e) says that $G$ is a ``good guy" (resp. ``$\ell$-good guy") if $\z G$ is connected and $\cF(G)=1$ (resp. $ |\z G/\zo G|$ and $|\cF(G)|$ are prime to $\ell$). As $\ell$ is good, if $[G,G]$ has no component of type $\AA$, $G$ is an $\ell$-good guy. It appears that Proposition~1.2.6 is easily proved  using induction in the semi-simple rank for $G$ an $\ell$-good guy (see (D) in the proof below). Independantly  we verify all assertions when $G$ is the good guy ${\rm GL}_n$ and acceed to any
$G$ of type $\AA$ by a standard way. Thus some groups are considered twice, some partial results are proved twice.

(A) Some general implications : 

All the properties in Proposition~1.2.6 go from two groups to their direct product.

\smallskip (A.i) {\sl (a) implies (f)}.

Indeed $L\cap\cento GY=\cent{\cento GY}{\zo L}$ is a Levi subgroup of the algebraic reductive group $\cento GY$ so is connected: one has  $L\cap\cento GY\subseteq\cento LY$. As 
$\cent LY=L\cap \cent GY$ and $\cent GY/\cento GY$ is finite, one has  $L\cap\cento GY=\cento LY$ and obtains an injective map $\cent LY/\cento LY\to\cent GY/\cento
GY$. If $Y\subseteq G^F$, that map commute with $F$-action.

\smallskip (A.ii) {\sl The first assertion in (b) implies the second one. }

 One has 
  $T_Y=\nor {\cento GY}{T_Y}\cap \cent G{T_Y}$, hence by  isomorphisms theorems  
  
  \noindent $W(\cento G Y,T_Y)=\nor{\cento GY}{T_Y}/T_Y\cong \nor {\cento GY}{T_Y}.\cent G{T_Y}/\cent G{T_Y}$, 
   a subgroup of $\nor G{\cent G{T_Y}}/\cent G{T_Y}$.

 \smallskip (A.iii) {\sl The last
assertion of (c) is a direct consequence of the preceding one.}

By the anti-isomorphism $W(G,T)\to W(G^*,T^*)$, $\nor G{\cento
G{Y},T_Y}/\nor{\cento GY}{T_Y}$ maps onto 

\noindent $\nor {G^*}{\cento
{G^*}{Y'},T_{Y'}}/\nor{\cento {G^*}{Y'}}{T_{Y'}}$. 

\smallskip (A.iv) {\sl  (d) follows from 

(d') If $\z G=\zo G$ and $\cF (G)=\{1\}$, then $\z{\cento GY}=\zo{\cento G Y}$ and $\cF(\cento GY)=\{1\}$.}

Given $G$, there exists a commutative diagram of isotypic morphisms defined over $\F q$
$$\diag {(H,F)} {} {(H_0,F)}{\pi_0}{(G_0,F)}{} {(G,F)}{\pi} {}$$
where horizontal maps are regular embeddings, $\pi$ and $\pi_0$ are quotients by central torii and $\cF(H)=\cF(H_0)=\{1\}$ (see Proposition~1.1.4). 

Let $X$, $Y$ be $\ell$-subgroups of $G^F$, $H^F$ with $\pi(X)=Y$.  
(d') applies to $H_0$ : let $D_0:= \cento{H_0}X$, we have $\cF(D_0)=\{1\}$ and $\z D=\zo D$. 

Let $C_0=\cento{G_0}Y$. By Proposition~1.2.3, $\pi_0(D_0)=C_0$, so that $\pi_0(\z{D_0})=\z {C_0}$. As $\z{D_0}$ is connected, so is $\z{C_0}$. Let $C=\cento G Y$. Clearly $C=C_0\cap G$ and $C_0=\z{G_0}.C$, hence $\z {C_0}=\z {G_0} .\z C$, $\z C=\z{C_0}\cap G$. As $\z {C_0}$ and $\z{G_0}$ are connected and $\z {G_0}.\zo C$ is connected with finite index in $\z{C_0}$, $\z {C_0}=\z {G_0}.\zo C$. Finally $\z C=\z G.\zo C$. 
 By isomorphisms theorems $\z C/\zo C$ is isomorphic to $ \z G/\zo C\cap\z G$. As $\zo G\subseteq \zo C$, $ \z G/\zo C\cap\z G$ is a quotient of $\z G/\zo G$. All morphisms commute with $F$ and $(\z C/\zo C)^F
=\z C^F/\zo C^F$ is a quotient of $\z G^F/\zo G^F$.

 Let $D=\cento HX$. Clearly $D_0=\z{H_0}D$, so that $[D,D]=[D_0,D_0]$ and $\cF(D_0)=\{1\}$ (by (d')) implies $\cF(D)=\{1\}$. As $\pi(X)=Y$, one has $\pi(D)=C$, $\pi([D,D])=[C,C]$ and $[C,C]$ is isomorphic to $[D,D]/{\ker \pi}\cap [D,D]$. As $[D,D]$ is simply connected, $\cF (C)= {\ker \pi}\cap [D,D]$, a subgroup of $\cF (G)=\ker \pi\cap[H,H]$.

When (c) is satisfied, one sees that $\cF(\cento GY)$ is isomorphic to $\z{\cento{G^*}{Y'}}/\zo{\cento{G^*}{Y'}}$ by Proposition~1.1.3 , hence 
a section of $\z{G^*}/\zo{G^*}$, isomorphic to a section of $\cF(G)$.

(A.v) {\sl A short exact sequence :
$$1\to (\zo
G\cap [G,G] )/(\zo G\cap
\cento{[G,G]}Y\to \cent{[G,G]}Y/\cento{[G,G]}Y\to {\rm A}_G(Y)\to 1
\leqno{ (\cal K)}$$ If $Y$ is $F$-stable, the groups in  ($\cal K$) are $F$-stable and morphisms are $F$-morphisms.}

We have $G=\zo G.[G,G]$ hence $\cent GY=\zo G.\cent{[G,G]}Y$ and  $\cento GY=\zo G.\cento {[G,G]}{Y}$. The two equalities imply that ${\rm A}_G(Y)$ is a quotient of  $\cent{[G,G]}Y/\cento{[G,G]}Y$. By isomorphism theorems the kernel is isomorphic to $(\zo
G\cap \cent{[G,G]}Y )/(\zo G\cap
\cento{[G,G]}Y$ and $\zo
G\cap \cent{[G,G]}Y =\zo
G\cap [G,G] $. Note that there exists a finite $\ell$-subgroup $Y'$ of $[G,G]$ such that $\zo G.Y=\zo G.Y'$. 

Then $\cent{[G,G]}Y/\cento{[G,G]}Y={\rm A}_{[G,G]}(Y')$ and ${\rm A}_G(Y)={\rm A}_G(Y')$.

(B)  Assume first $Y\subseteq G^F$, $\cF(G)^F_\ell=\{1\}=\z G^F_\ell/\zo G^F_\ell$, a property satisfied when $G=\zo G.\b G$. The property goes from $(G,F)$  to $(G^*,F)$ by (a) in Proposition~1.1.3. 

 We may assume $ [G,G]\neq 1$ and $Y\not\subseteq\zo G$, if not there is nothing to prove.

Let $z\in\z{\zo G^F_\ell.Y}\cap [G, G]$,
 $z\neq 1$. As $\ell$ is good, $\cento Gz$ is  a Levi  subgroup of $G$,
hence 
${\cF}(\cento Gz)^F_\ell=1$ and $\z{\cento Gz}^F_\ell=\zo{\cento Gz}^F_\ell$ ((b) in Proposition
1.1.3). By Proposition~1.2.4 (a) 
$\cento Gz^F=\cent {G}z^F$, hence $Y\subseteq \cento Gz$. The  semi-simple rank of
${\cent Gz}$ is less that the one of $G$ and 
induction applies in $\cent Gz$. Let $T_Y$ be a maximal $F$-stable torus in $\cento{\cent Gz}Y$. As $z\in T_Y$, $\cent
G{T_Y}\subseteq \cent Gz$. But
$\cento{\cent Gz}Y=\cento GY$ and $\cent{\cent Gz}Y=\cent GY$.   A maximal torus $T$ in $\cent Gz$ is maximal in $G$. Thus one obtains (a) and (b) for $Y$ in $G$ by (a) and (b) for $Y$ in $\cent Gz$, and ${\rm A}_G(Y)^F=\{1\}$ by (e) in $\cent G z$. Assertion (d) for $Y$ in $G$ is deduced from  (d) for $Y$ in $\cent G z$ and Proposition~1.1.3.

To apply assertion (c) in $\cent Gz$, choice  $L^*:=\cent Gz^* $ as a Levi subgroup of $G^*$ : there exists an $\ell$-subgroup $Y'$ of $(L^*)^F$  and a maximal $F$-stable torus $T^*$ of $L^*$ such that $Y'\subseteq \nor{L^*}{T^*}$ and $T^*_{Y'}= T^*\cap \cent{L^*}{Y'}$ maximal in $\cento{L^*}{Y'}$ with  dualities between $(T\subseteq \cent {\cent Gz}{T_Y})$ and $(T^*\subseteq \cent{L^*}{T^*_{Y'}})$ and between $(T_Y\subseteq \cent{\cent Gz}Y)$ and $(T^*_{Y'}\subseteq \cent{L^*}{Y'})$. As $T^*$ is a maximal torus in $G^*$ one has (c) for $Y$ in $G$. 

There is a ``non rational version" of the preceding proof, with same steps, and one obtains:

{\sl Assume $Y\subseteq G_\ell$ and $\cF(G)_\ell=\{1\}=\z G_\ell/\zo G_\ell$. One has (a), (b), (c), $\cF(\cento GY)_\ell=\{1\}$, $\z{\cento GY}_\ell=\zo{\cento GY}_\ell$  and ${\rm A}_G(Y)=\{1\}$.}

(C) Type $\AA$

(C.1) Assume $G={\rm GL}_n$, $Y\subseteq G^F$.

The representation of $Y$ on the space
$\F{}^n$ is semi-simple, $\F{}^n=\oplus_EV_E$, where $E$ is a set of  \irr\ representations of $Y$, $V_E$ is isotypic with some multiplicity $m(E)$. With these notations one has an isomorphism  $\cent GY\cong \prod
_E{\rm GL}_{m(E)}$. So 
$\cent GY$ is connected with a  connected center, (a), (d) and (e) hold.

A  maximal  $F$-stable torus $T_Y$ of $\cent GY$ is  defined by a decomposition of
 $\F{}^n$ as an $F$-stable direct sum $\oplus_{i\in I}V_i$ of \irr\ $\F{} Y$-modules (Schur's lemma).  A maximal torus
 of $G$ such that $T_Y\subseteq T$ and $Y\subseteq\nor GT$ is defined by
a choice in each  $\F{}$-space $V_i$ of an $Y$-stable family of generating and linearly independant lines, and $T$ is $F$-stable if the family is $F$-stable. Such a family exists because any representation of a nilpotent finite group is monomial. As $F$ acts by permutation on the set  $I$ an $F$-stable family exists and defines $T$. Hence (b) holds. 

One knows that $G={\rm GL}_n$    is isomorphic to $G^*$. Assertion (c) is clear in that case.

(C.2) Assume $\z G$ is connected.

There exists a regular covering $\pi\colon K:=\zo G.H\to G$ where $H$ is a covering of $[G,G]$ and a direct product of linear groups,  and $\zo G.H$ a central product. Furthermore the restriction of $\pi$, $K^F\to G^F$ is onto. There is some $\ell$-subgroup $X$ of $K^F$ such that $\pi(X)=Y$. 

 By (C.1) our Proposition is satisfied in $H$. As well it is satisfied in $K$. For any finite $\ell$-subgroup $Y$ of $K$, there exists a finite $\ell$-subgroup $Y_1$ of $H$ with $Y_1.\zo G=Y.\zo G$.  As for (a) (b) and (d) for $Y$ it follows from the equalities $K=\zo G.H$, $\zo K=\zo G.\zo H$  and (a) (b) and (d) for $Y_1$. As for (c), $K^*$ is a central product $\zo G^*.H^*$ and a couple of dual Levi subgroups in duality in $H$ and $H^*$  give a couple of dual Levi subgroups in $K$, $K^*$. The sequence ($\cal K$) in (A.v) gives ${\rm A}_K(Y)={\rm A}_H(Y_1)=1$. 

A special case is $\zo G=1$, that is $G$ adjoint. By Proposition~1.2.3 $\pi(\cento H X)=\cento GY$. By Lemma~1.2.1 $\cent
GY/\pi(\cent HX)$ is an $\ell$-group. Thus ${\rm A}_G(Y)$ is an $\ell$-group. The exact sequence ($\cal K$) implies that $\cent {[H,H]}X/\cento{[H,H]}X$ is an $\ell$-group for any finite $\ell$-subgroup $X$  of $H$. This apply to any morphism $G\to G_{\rm ad}$ hence 

{\sl in the sequence ($\cal K$) the group $\cent{[G,G]}X/\cento{[G,G]}X$ is an $\ell$-group. So is ${\rm A}_G(Y)$ for any $G$ of type $\AA$.}

Applying once more Proposition~1.2.3 and Lemma~1.2.1 to $\pi\colon K\to G$  with $\pi(X)=Y$ and to $K^F\to G^F$ if $Y\subseteq H^F$, one obtains (a),  (b), (e) (recall that $\Ker \pi\cap[K,K]$ is isomorphic to $\cF(G)$) and $\z{\cento GY}=\zo{\cento G Y}$. On dual side one has a regular embedding $\pi^*\colon G^*\to K^*$. A dual $\cento {K^*}{X'}$ of $\cento KX$ may be written $\zo{K^*}.\cento{\pi(G^*)}{\pi^*(Y')}$ for some finite $\ell$-subgroup $Y'\subseteq  G^*$ and one obtains easily (c) for $Y$ in $G$ from (c) for $X$ in $K$ : $\cento {G^*}{Y'}$ is isomorphic to $\cento{\pi^*(G^*)}{\pi^*(Y')}$ in duality with $\cento G Y$.

(C.3) End of the proof in type $\AA$

Let us consider a regular imbedding  $G\subseteq H$ defined over $\F q$ with $F$-actions. One has $\cF(G)=\cF(H)$, $\z H=\zo H$. By (C.2), (a) to (e) hold in $H$. For any $\ell$-subgroup $Y$ of $G$ one has $\cento HY=\zo H.\cento GY$. (a) and (b) for $Y$ in $H$ imply (a) and (b) for $Y$ in $G$, with $T_Y=\zo H.(T_Y\cap G)$, $T_Y$ maximal in $\cento HY$, $T_Y\cap G$ maximal in $\cento GY$, $T=\zo H.(T\cap G)$....

 As for (d), one knows that $\zo H.\z{\cento GY}=\z {\cento HY}=\zo{\cento HY}=\zo H.\zo{\cent GY}$, hence $\zo{\cento GY}/\zo H\cap \zo{\cento GY}$ and $\z{\cento GY}/\zo H\cap \z{\cento GY}$ are isomorphic. But $\zo H\cap G=\z G$ and $\zo G\subseteq \cento GY$, thus  
$(\zo H\cap \z{\cento GY})/(\zo H\cap \zo{\cento GY})$ is isomorphic to a section of $\z G/\zo G$.

By a covering dual of the embedding, say  $\pi^*\colon G^*\to H^*$, then $\pi^*(T^*)$ is a dual of $T\cap G$. If a group $Y'$ is given by (c) for $Y$ in $H$, then the duality between $(T_Y\subseteq \cento HY)$ and $(T^*_{Y'}\subseteq \cento{H^*}{Y'})$ induces a duality between $(T_Y\cap G\subseteq  \cento G Y)$ and $(\pi^*(T_{Y'})\subseteq  \cento {G^*}{\pi^*(Y')})$. If $\cent H{T_Y}$ and $\cent{G^*}{T^*_{Y'}}$ are dual Levi subgroups, so are $\cent G{T_Y\cap G}$ and $(\cent{G^*}{\pi^*(T^*_{Y'})}$. That is (c) for $Y$ in $G$.

To compare ${\rm A}_G(Y)$ and ${\rm A}_H(Y)$  use the sequence ($\cal K$). The two groups are quotient of the same $\ell$-group with included kernels : with $C^\circ=\cento{[G,G]}Y$,  $\zo G\cap[G,G]/\zo G\cap C^\circ$ injects in $\zo H\cap [G,G]/\zo H\cap C^\circ$, hence ${\rm A}_H(Y)$ is a quotient of ${\rm A}_G(Y)$ whose cokernel is $F$-isomorphic to a section of $\z{[G,G]}/\zo G\cap C^\circ$, hence of $\z G/\zo G$. If $\z G^F_\ell=\zo G^F_\ell$ and $\cF(G)^F_\ell=1$, then ${\rm A}_H(Y)^F$ and ${\rm A}_G(Y)^F$ are isomorphic and ${\rm A}_H(Y)^F=1$ by (C.2), that is the ``rational part" of (e).

(D) End of the proof.

Assuming $Y\subseteq G^F$, we use the \dec\ in central product defined in section 1.1.5 : put $\b G'=\zo G.\b G$. Then $G=\a G.\b G'$ and $\z G/\zo G$ is isomorphic to a direct product  $\z
{\a G}/\zo{ G}\times \z{\b G'}/\zo G$.

Any $y\in G_\ell^F$ writes in a unique way in $(\a G)_\ell^F\times (\b G)_\ell^F$. 
 Let  $\a Y$ and $\b Y$ be the
projections of $Y$ on $(\a G)_\ell^F$ and $(\b G)_\ell^F$. One has
$\cent GY=\cent{\a G}{\a Y}.\cent{\b G'}{\b Y}$, $\cento GY=\cento{\a G}{\a
Y}.\cento{\b G'}{\b Y}$, so that  ${\rm A}_G(Y)$ is isomorphic to ${\rm A}_{\a G}(\a Y)\times
{\rm A}_{\b G'}(\b Y)$. By definition of
$\b G$, $\z G^F_\ell/\zo G^F_\ell$ is isomorphic
to $\z {\a G}^F_\ell/\zo
{\a G}^F_\ell$ and ${\cF}(G)^F_\ell$ is isomorphic \`a ${\cF}(\a G)^F_\ell$. The
 assertions (a), (d) and (e) are true for $\a Y$ in $\a G$ and $\b Y$ in $\b G'$, there are
 true for  $Y$ in  $G$. 

The properties required in (b) and (c) go up from $\a G$ and $\b G'$ to $G$ by straightforward constructions.   

With clear notations $T_Y=(T_Y\cap \a G).(T_Y\cap \b G')=T_{\a Y}.T_{\b Y}$, $T=(T\cap \a G).(T\cap \b G')=\a T.\b T$ (central products on $\zo G$), with $T_{\a Y}=(T^{\a Y})^\circ$, $T_{\b Y}=(T^{\b Y})^\circ$, 
hence $(T^Y)^\circ =T_Y$, $W(\cento{G}Y,T_Y)\cong W(\cento {\a G}{\a Y}, T_{\a Y})\times W(\cento {\b G}{\b Y}, T_{\b Y})$

Dualities between 
$( T^*_{\a Y'}\subseteq \cento{{\a G}^*}{\a Y'})$ and $(T_{\a Y},\cento{\a G}{\a Y}$, between $( T^*_{\b Y'}\subseteq \cento{({\b G'})^*}{\b Y'})$ and $(T_{\b Y},\cento{\b G'}{\b Y})$ define, through the morphism $\rho\colon G^*\to {\a G}^*\times (\b G')^*$, whose kernel is a central torus,  a duality between  $(T^*_{Y'}\subseteq \cento{G^*}{Y'})$ and $(T_Y\subset\cento GY)$ where $Y'=\rho(\a Y'\times \b Y')$ and $T^*_{Y'}=\rho(T^*_{\a Y'}\times T^*_{\b Y'})$. 

Similarly dualities between $(\a T\subseteq \cent{\a G}{T_{\a Y}})$ and $(\a T^*\subseteq \cent{{\a G}^*}{T^*_{\a Y'}})$, between $(\b T\subseteq \cent{\b G'}{T_{\b Y}})$ and $(\b T^*\subseteq \cent{(\b G')^*}{T^*_{\b Y'}})$ give a duality between $(T\subseteq \cent G{T_Y})$ and $(T^*\subseteq \cent{G^*}{T^*_{Y'}})$ where $T^*=\rho(\a T^*\times \b T^*)$.
\bull

  \medskip\noindent{\bf 1.2.7. Proposition. }{\sl Let $(G,F)$ be of classical type in odd caracteristic, let $Y$ be a $2$-subgroup de
$G^F$.

(a) $\cento GY$ is reductive.

(b)  For any (eventually $F$-stable if $Y\subseteq G^F$) maximal torus $T_Y$ of $\cento GY$, there exists
a maximal  (eventually $F$-stable)  torus
$T$ of
$G$ such that $T_Y\subseteq T$ and $Y\subseteq\nor GT$. Thus
${\rm W}(\cento GY,T_Y)$ is isomorphic
to a subgroup of ${\rm W}_{G}(\cent G{T_Y})\cong\nor G{\cent G{T_Y}}/\cent G{T_Y}$. 

(c) If $G$ has no component of type $\BB$ or $\CC$, then assertion (c) with $\ell=2$ of Proposition~1.2.6 is true.

(d) The group $\z{\cento GY}/\zo
{\cento GY}$ is isomorphic to a quotient of $\z G/\zo G$ and $\cF(\cento GY)$ is
isomorphic to a subgroup of
$\cF(G)$, with $F$-action when $Y\subseteq G^F$. 

(e) ${\rm A}_G(Y)$ is an $2$-group. If ${\cF}(G)_2=\{1\}$ and $(\z G/\zo G)_2=\{1\}$ then
${\rm A}_G(Y)=\{1\}$.  If $Y\subseteq G^F$,
${\cF}(G)^F_2=\{1\}$ and $(\z G/\zo G)^F_2=\{1\}$ then
${\rm A}_G(Y)^F=\{1\}$. 

(f) If $L$ is a  Levi subgroup of $G$ and $Y\subseteq L$, ${\rm A}_L(Y)$ is
 isomorphic to a subgroup of ${\rm A}_G(Y)$. If furthermore $Y\subseteq G^F$, ${\rm A}_L(Y)^F$ is
 isomorphic to a subgroup of ${\rm A}_G(Y)^F$.
}

 \preuve (a sketch of) We have seen in the proof of Proposition~1.2.6, using central products and regular isotypic morphisms, how to reduce to rationally irreducible types. Details are left to the reader.
 
In type $\AA$ $2$ is good,  Proposition~1.2.6 applies.

 Let $f$ be a non degenerate bilinear form on a space $V$ on $\F{}$, defined on $\F q$, assume $G={\rm SO}(f)$. Thus
$G^F$ is a symplectic or orthogonal group on  $\F q$.

 Under $Y$-action $V$ is a direct and orthogonal sum of isotypic $\F {} [Y]$-modules 
 $$V=V_+\oplus V_-\oplus(\oplus_{\{E,E'\}}(V_E\oplus V_{E'}))$$ Here we have denoted $V_+$ (resp. $V_-$) the space of fixed points (resp. antifixed points : $y.v=-v$ for all $y\in Y$) and $\{E,E'\}$
belongs to a set of isomorphism classes of pairs of contragredient \irr\ representations. Let $f_+$ (resp. $f_-$) be the restrictions of $f$ to $V_+$ (resp. $V_-$). Then $\cent GY$ is
the intersection with $G$ (inside $ {\rm GL}(V)$) of the direct product $$O(f_+)\times
O(f_-)\times (\times_{\{E,E'\}} {\rm GL}_{m(E)})$$  Any factor
${\rm GL}_{m(E)}$ is contained in $G$, it acts on  $V_E\oplus V_{E'}$ where $V_E$ and $V_{E'}$ are totally isotropic spaces with dual basis. 
One sees that
$\cento GY$ is reductive, that gives (a) for all $G$. 

Clearly (b) holds in $G$ and
${\rm A}_G(Y)$ is a $2$-group. These properties are preserved by isotypic morphisms.
There exists a group $( G_0,F)$ of same type such that $\zo{ G_0}=\z{ G_0}$ and $\cF( G_0)=1$. One may verify  that, for any $2$-subgroup $ Y_0$ of $ G_0^F$, $\cent {G_0}{ Y_0}$ and $\z{\cent{G_0}{Y_0}}$ are connected, and $\cF(\cent{G_0}{Y_0})=1$. One deduce (d), (e) and (f)  as in Proposition~1.2.6. 

In types $\AA$ and $\DD$ there exists groups isomorphic to a dual, so (c) in Proposition~1.2.6 is true.
\bull

\bigskip \noindent{\bf 1.3.   Jordan \dec\ in $\II{G^F}$}

The properties we recall in that section are essentially due to G. Lusztig [26], [27]. In Proposition~1.3.2 we state what we need on Jordan decomposition of \irr\ characters when the center of $G$ is connected. As a consequence of ``non-multiplicity condition" in an isotypic morphism $H\to G$ (Proposition~1.3.3), Propositions~1.3.6, 1.3.7 describe the link between Jordan \dec s in an isotypic morphism. 

\medskip\noindent{\bf 1.3.1. Some facts and notations. } 

Let $L$ be an $F$-stable Levi complement of the unipotent radical of a parabolic subgroup $P$ of $G$, then is defined a virtual ${\cal O}[G^F\times L^F]$-bimodule, where ${\cal O}$ may be  the ring of an $\ell$-modular splitting  system $({\cal O},K,k)$ (see 1.1.2), and, as a consequence, " Deligne-Lusztig induction ", which maps any ${\cal O} L^F$-module to a virtual ${\cal O} G^F$-module,  see [20] chapter 11. By extension from ${\cal O}$ to $K$ and linearity is defined a linear map from $\CF K{L^F}$ to $\CF K{G^F}$, `` twisted induction" in [16] \S 8.3,  where it is denoted  $R_{L\subseteq P}^G$. 
Denote
$\Lu {L\subseteq P}G$ the linear map  
$\epsilon_G\epsilon_L R_ {L\subseteq P}^G$ where 
$\epsilon_G=(-1)^{r(G)}$, and $r(G)$ is the semi-simple $\F q $-rank
of
$G$,  we name it {\bf Lusztig induction}, the dual map with respect to standard scalar products on spaces of central functions is called {\bf Lusztig restriction} and denoted $\slu {L\subseteq P}G$.  In a isotypic morphism $\sigma\colon (G,F)\to (H,F)$, as defined in 1.1.4, one has,  by [6] Proposition~1.1 and with evident notations
$$\Res{}{\sigma\colon G^F\to H^F}\circ\Lu{M\subseteq Q}H=\Lu{\sigma^{-1}(M)\subseteq\sigma^{-1}(Q)}G\circ \Res{}{\sigma\colon \sigma^{1}(M)^F\to M^F}{}\leqno{(1.3.1.1)}$$
The character formula gives ([20] 12.17) 
$$\forall \chi\in\CF K{L^F}\quad(\Lu {L\subseteq P}G\chi)(1)={|G^F|_{p'}\over |L^F|_{p'}}\chi(1) \leqno{(1.3.1.2)}$$
Some properties, as formula (1.3.1.2), may be independant of the choice of $P$, given $L$, and the notation $\Lu
{L\subseteq P}G$ is frequently simplified in $\Lu LG$. The map $\Lu {L\subseteq P}G$ itself may be independant of the choice of $P$; it is the case if $L$ is a torus, see also Definition~2.2.1 and Proposition~2.2.4. The subspace of $\CF K{G^F}$ generated by  $\cup_T\Lu TG(\CF K{T^F})$ when $T$ ranges over $F$-stable maximal tori of $(G,F)$, is called the space of {\bf uniform functions on $G^F$}.  The orthogonal projection on the space of uniform functions of a central function $\chi$ on $G^F$ will be denoted $\piu ^G(\chi)$; one has ([20] 12.12)
$$\piu^G(\chi)=\sum_{\{(T,\theta)\}/G^F}|W(G,T)^F|^{-1}.\scal{\Lu TG\theta}\chi{G^F}\Lu TG\theta\leqno{(1.3.1.3)}$$
By Mackey decomposition formula ([20] 11.13) and (1.3.1.3) $\piu$ commute with Lusztig induction and restriction :$$\piu^L\circ \slu LG=\slu LG\circ\piu^G,\quad \piu^G\circ \Lu LG=\Lu LG\circ \piu^L \leqno{(1.3.1.4)} $$
Since the regular representation of $G^F$ is a uniform function ([20] 12.14) $\piu$ preserves the value on $1$.

Let
$(G,F)$ and 
$(G^*,F)$ be in duality and let
$s\in G^{*\,F}$ be semi-simple, the  {\bf rational Lusztig series} $\ser{G^F}s$ is a subset of $\II {G^F}$. The set $\ser{G^F}s$ is defined by the following property :  for any couple of $F$-stable maximal torus $(T\subseteq G,T^*\subseteq G^*)$ in dual conjugacy classes, if $\theta\in(T^F)^\wedge$ corresponds by duality to $s\in T^{*\,F}$, then $\Lu TG\theta$ writes in $\Z \ser{G^F}s$. One has a partition
$\II {G^F}=\cup_{(s)}\ser{G^F}s$ where $(s)$ ranges over $G^{*\,F}$-conjugacy class of semi-simple elements of $G^{*\,F}$, see [20] 14.41, [16] \S~8.4. Thus if, with evident notations, $(T_1,\theta_1)$ is not $G^F$-conjugate to $(T_2,\theta_2)$ ---and that is equivalent to non $G^{*\,F}$-conjugacy between $(T^*_1,s_1)$ and $(T^*_2,s_2)$  see [20] 13.13--- then ${\Lu {T_1}G{\theta_1}}$ and ${\Lu {T_2}G{\theta_2}}$ have no common \irr\ constituent. 

Elements of  $\ser{G^F}1$ are said to be {\bf unipotent}. If $\sigma\colon (G,F)\to (H,F)$ is an isotypic morphism, the restriction through $ \sigma$ restricts to a one-to-one map $\ser{H^F}1\to \ser{G^F}1$ [20] 13.20.

One has $\piu^G(K\ser{G^F} s)\subseteq K\ser{G^F}s$ and the orthogonal projection on the space $K\ser{G^F}s$ commute with $\piu^G$. Let us design by $\piu^{G,s}$ the  product of these projections :
$$\piu^{G,s}\colon K\II{G^F}\to \piu^G(K\ser{G^F}s)=\piu^G(K\II{G^F})\cap K\ser{G^F}s\leqno{(1.3.1.5)}$$  

Lusztig [26] has defined an orthogonal basis $\{\Lu fG\}_f$ of $\piu^G(K\ser{G^F}1)$ with indexation in some set of representations of $W.\gen F$.  There is a partition ${\frak F}$ of that set in so called {\it families}, defining a decomposition of  $\piu^G(K\ser{G^F}1)$ in an orthogonal sum of subspaces ${\cal FU}(G,1,{\frak f})=\sum_{f\in{\frak f}}K\Lu fG$, such that  for any $\chi\in\ser{G^F}1$ there exists a unique family ${\frak f}\in{\frak F}$  such that $\piu^G(\chi)\in {\cal FU}(G,1,{\frak f})$ :
$$\piu^G(K\ser{G^F}1)=\perp_{{\frak f}\in{\frak F}}{\cal FU}(G,1,{\frak f}), \quad \ser{G^F}1=\cup_{\frak f\in{\frak F}}\ser{G^F}1_{\frak f}, \quad \piu^G(\ser{G^F}1_{\frak f})\subseteq  {\cal FU}(G,1,{\frak f})\leqno{(1.3.1.6)}$$
 The decomposition of $\piu^G(\chi)$ on the basis $\{\Lu fG\}_{f\in{\frak f}}$  of ${\cal FU}(G,1,{\frak f})$ is known, see [26], Chapter 4, for details.

When  $s\in (G^*)^F_{\ell'}$ one defines $\lser{G^F}s:=\cup_{t\in \cento{G^*}s^F_\ell}\ser{G^F}{st}$.
We denote $\bl{G^F}{s}$ the set of blocks $b$ of $G^F$  such that $\II
b\cap\lser{G^F}s\neq \emptyset$. By [9] (see  [16]
Theorem~9.12)
$$b\in \bl{G^F}{s}\;{\rm if\; and \;only\;if\; } \II b\cap \ser{G^F}s\neq\emptyset \;{\rm if\; and \;only\;if\; }\II
b\cap\lser{G^F}s\neq \emptyset \leqno{(1.3.1.7)}$$
An element of  $\bl{G^F}1$ is said to be unipotent.

When $G$ is not connected
we define $\ser{G^F}1$ as the set of $\chi\in\II{G^F}$ that cover some $\chi^\circ\in\ser{G^{\circ\, F}}1$. Elements of $\ser{G^F}1$ may be defined as the  \irr\ components of the $\Lu TG 1$, with suitable definitions of so called {\it quasi-tori} $T$,  and of Lusztig induction [21]. Similarly an 
$\ell$-bloc $B$ of $G^F$ is said to be unipotent if it covers a unipotent block of $G^{\circ\,F}$, or equivalently if $\II B\cap \ser{G^F}1\neq \emptyset$.

In the following Proposition we introduce one of our first important tool, the Jordan \dec\ in $\II{G^F}$, see  [26] 4.23,  [20] 13.24.
 Recall that, when $\z G$ is connected and $s$ is semi-simple in $G^*$, $\cent{G^*}s$ is connected (Propositions 1.1.3 and 1.2.4). By general conventions we assume that the duality between $L$ and $L^*$, $G$ and $G^*$ are defined  around the same pair of dual maximal tori $(T\subseteq L, T^*\subseteq L^*)$ and so on for any ``dual" sets $\{L_j\}_j$, $\{L^*_j\}_j$ of Levi subgroups of $G$, $G^*$ with a common maximal torus $T\subseteq \cap_jL_j$, $T^*\subset\cap_jL^*_j$. Similar restrictions apply to dualities over $\F q$ between $F$-stable groups. It follows that the choice of duality between $L$ and $L^*$, hence the eventual different choices, if $s\in L^*$,  of $\Psi_{L,s}$ do not affect $\Lu LG\circ\Psi_{L,s}$.

\medskip\noindent{\bf 1.3.2. Proposition. }{\sl For any 
connected algebraic reductive group $(G,F)$ defined on $\F q$,  with  connected center, and dual $(G^*,F^*)$, there exist
one-to-one maps, named Jordan \dec s,
$$\Psi_{L,s}\colon \ser{\cent{L^*}s^F}1\to
\ser{L^F}s $$
defined for any couple of $F$-stable Levi subgroups $(L\subseteq G,L^*\subseteq G^*)$ in dual conjugacy classes and  any 
semi-simple element $s$ in $ L^{*\,F}$ and that satisfy the following properties :

 Extend $\Psi_{L,s}$  by linearity  $$\Psi_{L,s}\colon K\ser{\cent{L^*}s}1\to K \ser{L^F}s$$ 
 
(i) On orthogonal projections on the spaces of uniform functions : 

For any couple $(T\subseteq L,T^*\subseteq L^*)$ of  $F$-stable maximal tori in dual conjugacy classes with $s\in T^{*\,F}$ and any $\la\in\ser{\cent{L^*}s^F}1$, one has 
$$\scal{\Psi_{L,s}(\la)}{\Lu TL(\Psi_{T,s}(1_{T^{* \,F}})}{L^F}=\scal{\la}{\Lu{T^*}{\cent{L^*}s}(1_{T^{* F}})}{\cent{L^*}s^F}\leqno{\rm
(J.1)}$$
$$\Lu TL(\Psi_{T,s}(1_{T^{* F}})=\Psi_{L,s}(\Lu{T^*}{\cent{L^*}s}(1_{T^{* F}}))\leqno{\rm
(J.2)}$$

(ii) Assume $s$ is central in $L^*$, the duality between $(L,F)$ and $(L^*,F)$ defines    $\hat
s:=\Psi_{L,s}(1_{L^F})\in\II{L^F}$ such that  $\hat s(1)=1$. For any semi-simple $t\in L^{*\,F}$, one has  $\ser{L^F}{st}=\hat s\otimes
\ser{L^F}t$ and,   for
any
$\la\in \ser{\cent{L^*}t^F}1$, $$\Psi_{L,st}(\la)=\hat s\otimes \Psi_{L,t}(\la).$$

(iii) If $\cent{G^*}s\subseteq L^*$, then for any parabolic subgroup  $P$ of $G$ with Levi complement $L$ one has $$\Lu{L\subseteq
P}G\circ\Psi_{L,s}=\Psi_{G,s}.$$

(iv) Let $(\sigma,\sigma^*)$ be an isomorphism from $((L,F),(L_1^*,F))$ to
$((L_1,F),(L^*,F))$ and $s=\sigma^*(s_1)$. Then  $(\sigma,\sigma^*)$ induces
isomorphisms $\tau\colon L^F\to L_1^F$ et $\tau_s^*\colon
\cent{L_1^*}{s_1}^F\to
\cent{L^*}s^F$,  and there is a  commutative diagram with one-to-one maps
$${\def\normalbaselines
{\baselineskip20pt\lineskip3pt
\lineskiplimit3pt}\matrix
{{\ser{\cent{L^*_1}{s_1}^F}1}&{\mapleft{{\rm
Res}_{\tau^*_s}}}&{\ser{\cent{L^*}s^F}1}\cr
{\mapdown{\Psi_{L_1,s_1}}}&{}&{\mapdown{\Psi_{L,s}}}
\cr{\ser{L_1^F}{s_1}}&{\mapright{\Res{}\tau}}&{\ser{L^F}s}\cr}}$$

(v) $\Psi_{L,s}$ is ``functorial with respect to 
central products in $L$".

(vi) For any $\la\in\ser{\cent{L^*}s^F}1$, one has
$\Psi_{L,s}(\la)(1)={|L^F]_{p'}\over |\cent{L^*}s^F|_{p'}}\la(1)$.}

\medskip\noindent {\it Comments. } 

The case $s=1$ : $\Psi_{G,1}\colon \ser{G^{*\,F}}1\to
\ser{G^F}1$ is a bijection between unipotent series of groups in  duality. The hypothesis ``$\z G$ is connected" may be dropped (see Proposition~1.3.6 for a general statement in case $\z G$ is not connected). 

In case of  tori $T$, $T^*$ in dual conjugacy classes for $(L,L^*)$ ---and so for $(G, G^*)$ ---  $\Psi_{T,s}(1_{T^{* F}})$ is the image of $s$ by an isomorphism $T^{*\,F}\to (T^F)^\wedge$,  the construction of it assume choices of isomorphisms $(\Q/\Z)_{p'}\cong {\bar\Q_\ell}^\times$ and $(\Q/\Z)_{p'}\cong {\F{}}^\times$, that are fixed one for all, and depends of the duality between $T$ and $T^*$. Thanks to  [20] 13.13 $\Lu {T}L(\Psi_{T,s}(1_{T^{* F}}))$ is well defined as it is implicitely assumed in (J1), (J2). Recall the notation $$\Lu{T^*}L(s):= \Lu TL(\Psi_{T,s}(1_{T^{* F}}))\leqno{(1.3.2.1)}$$

The relative Weyl groups $W_G(T)^F_{\Psi_{T,s}(1_{T^{* F}})}$ and $W_{\cent{G^*}s}(T^*)^F$ are isomorphic. By [20] 11.16 we have $\scal{\Lu {T^*}Gs}{\Lu{T^*}Gs}{G^F}=\scal{\Lu {T^*}{\cent{G^*}s}1}{\Lu{T^*}{\cent{G^*}s}1}{\cent{G^*}s^F}$. By [20] 13.12, 12.12 and (J1) and (J2) above  $\Psi_{G,s}$ send an orthogonal basis of $\piu^{\cent{G^*}s}(K\ser{\cent{G^*}s^F}1$ to  an orthogonal  basis of  $\piu^G(K\ser{G^F}s)$,  "Jordan \dec\ commute with $\piu$"  and restricts to an isomorphism of metric spaces $$\piu^G\circ \Psi_{G,s}=\Psi_{G,s}\circ \piu^{\cent{G^*}s,1},\quad\Psi_{G,s}\colon \piu^{\cent{G^*}s}(K\ser{\cent{G^*}s^F}1)\,{\smash{\mathop{\relbar\joinrel\rightarrow}\limits^{\sim}} }\,\piu^G(K\ser{G^F}s)\leqno{(1.3.2.2)}$$ (here we have used notation (1.3.1.5)). More precisely by $\Psi_{G,s}$ and the decomposition (1.3.1.5) we may define  families in $\ser{G^F}s$ : $\ser{G^F}s_{\frak f}:=\Psi_{G,s}(\ser{\cent{G^*}s^F}1_{\frak f})$, ${\cal FU}(G,s,{\frak f}):=K\ser{G^F}s_{\frak f}$
$$\Psi_{G,s}({\cal FU}(\cent{G^*}s,1,{\frak f})={\cal FU}(G,s,{\frak f})\quad  \piu^G(K\ser{G^F}s)=\perp_{{\frak f}\in{\frak F}}{\cal FU}(G,s,{\frak f})\leqno{(1.3.2.3)}$$
In (1.3.2.3) the set $\frak F$ is defined from $(W(\cent{G^*}s),F)$.

Let $M$ and $M^*$ be Levi subgroups in duality of $L$ and $L^*$ respectively, such that $s\in M^*$, then $\cent{M^*}s$ is a Levi subgroup of $\cent{L^*}s$. Let $\al\in\ser{\cent{L^*}s^F}1$, $\la=\Psi_{L,s}(\al)$. Let be  any couple of $F$-stable maximal tori $(T\subseteq M,T^*\subseteq M^*)$ in dual conjugacy classes with $s\in T^*$, put $\theta=\Psi_{T,s}(1_{T^{* F}})$.  We have $\Lu TM \theta=\Psi_{M,s}(\Lu {T^*}{\cent{M^*}s} { 1_{T^{* F}}})$. By transitivity of Lusztig induction and (J.1) we have
$$\scal{\slu ML \la}{\Lu {T}M \theta}{M^F}=\scal{\slu {T}L \la}{\theta}{{T}^F}=\scal {\slu {T^*}{\cent{L^*}s}\al}{ 1_{T^{* F}}}{{T^*}^F}=\scal{\slu {\cent{M^*}s}{\cent {L^*}s} \al}{\Lu {T^*}{\cent{M^*}s} { 1_{T^{* F}}}}{\cent{M^*}s^F}$$
Thus we enforce (1.3.2.2) in another commutation formula with Lusztig induction, we may combine with (1.3.1.4).  The preceding equalities give 
$$(\piu^{M,s}\circ \slu ML\circ \Psi_{L,s})(\al)=(\Psi_{M,s}\circ \piu^{\cent{M^*}s,1}\circ \slu {\cent{M^*}s}{\cent{L^*}s})(\al), \quad \al\in\ser{\cent{L^*}s^F}1\leqno {(1.3.2.4)}$$ and by adjunction 
$$(\Psi_{L,s}\circ \Lu{\cent{M^*}s}{\cent{L^*}s}\circ \piu^{\cent{M^*}s}) (\al)=(\piu^L\circ\Lu ML\circ \Psi_{M,s})(\al), \quad \al\in\ser{\cent{M^*}s^F}1\leqno (1.3.2.5)$$

The functoriality of the Jordan \dec\ (assertion (iv)) apply in the following situation :

 Let $(L_j)_j$ and $(L^*_j)_j$ be two pairs ($j=1,2$) of $F$-stable Levi subgroups of $G$ and $G^*$ respectively, $L_j$ and $L^*_j$ in dual classes. Assume that for some $g\in G^F$ and $s\in (L^*_1\cap L^*_2)^F$, $(L_1,\ser{L^F}s)^g=(L_2,\ser{G^F}s)$. The duality between $L_j$ and $L^*_j$ is defined around $F$-stable couples (Torus$_j$ $\subset$  Borel$_j$) with root data in duality. As $g$ is defined modulo $\nor{G^F}{L,\ser{L^F}s}$ we may assume that $g$ sends such a couple $(T_1 \subseteq B_1)$ onto $(T_2\subseteq B_2)$, assuming by choice $s\in {T^*_j}^F$. Then there exists $g^*\in {G^*}^F$ that sends $(T^*_2\subseteq B^*_2)$ onto $(T^*_1\subseteq B^*_1)$ so that $g^*s=sg^*$. Then  $g^*$ and $g$ induce by interior automorphisms dual morphisms $L_1\to L_2$ and $L^*_2\to L^*_1$ ($g^*$ is defined modulo $T^*_1$) and (iv) applies : one has $\Psi_{L_1,s}(\al)^g=\Psi_{L_2,s}(\lexp{g^*}\al)$ for any $\al\in\ser{\cent{{L_1^*}^F}s}1$.

 As for (v) functoriality is as follows : assume $G=G_1.G_2$ a central product over a torus of $F$-stable groups $G_j$. Let $\pi\colon G_1\times G_2\to G$ so defined and let $\pi^*\colon G^*\to G^*_1\times G_2^*$ a dual morphism. Put $\pi^*(s)=(s_1,s_2)$, with $s_j\in (G_j^*)^F$, so that $\pi^*(\cent {G^*}s)=\cent {G^*_1}{s_1}\times \cent{G_2^*}{s_2}$. Similarly $L=\pi(L_1\times L_2)$ and $\pi^*(L^*)=L^*_1\times L^*_2$. By what we saw in case of unipotent series one has a one-to-one map $\sigma\colon \ser {\cent {L^*_1}{s_1}^F}1\times \ser {\cent {L^*_2}{s_2}^F}1\to \ser {\cent {L^*}{s}^F}1$.  As $(s_1,s_2)=\pi^*(s)$, elements of $\ser{L_1^F}{s_1}$ and $\ser{L_2^F}{s_2}$ have equal restriction on the kernel of $\pi$ that is the restriction of $\hat{s_j}$ in the notation of assertion (ii). Then, using the symbol $\otimes$  in two different senses, one may write $\Psi_{L,s}(\sigma(\al_1\otimes \al_2))=\Psi_{L_1,s_1}(\al_1)\otimes \Psi_{L_2,s_2}(\al_2)$.
  
\medskip\noindent{\bf 1.3.3. Remarks. } A generalization of (J.2) in the form $$\Lu {L\subseteq P}G\circ\Psi_{L,s}=\Psi_{G,s}\circ
\Lu{\cent{L^*}s}{\cent{G^*}s}\leqno{\rm
(J.3)}$$ with parabolic subgroups suitably defined or, better!, dropped, would be quite useful, as well as Mackey decomposition formula for Lusztig induction [20] chapter 11.

(a) In type {\bf A}  with connected center all central functions
are uniform. In that case (J.2) implies (J.3).
Asai and Shoji have shown that (J.3) is true in classical type with connected center for any $L$, see [29] and [24]
Appendix, and unicity of Jordan's decomposition follows in that case. In our Appendix, Proposition~5.3, we give an elementary proof of (J3) in classical types, assuming Mackey decomposition formula for $\slu LG\circ\Lu LG$ and  knowing Asai's formulas that give $\Lu{\cent{L^*}s}{\cent{G^*}s} \al$ for any $\al\in\ser{\cent{L^*}s}1$.

(b) If $s=1$, (J.3) is true  
because 
$\ser{G^F}1$ has a generic parametrization and Lusztig induction  is ``generic on unipotent functions" (see
[10] 1.33). If $s$ is central, by assertion (ii) in Proposition~1.3.2,
(J.3) goes from 1 to $s$.

More generally, if
$\cent{G^*}s\subseteq L^*$, a Levi subgroup of $G^*$ in the dual class of $L$,
$\Lu LG$ restricts to a one-to-one map
$\ser{L^F}s\to \ser{G^F}s$ independantly of choices of parabolics.
 It follows that 
(J.3) is satisfied  if $\cent{G^*}s$ is a Levi subgroup of $G^*$ :

Indeed let
$G(s)$ be a Levi subgroup of $G$ in duality with $\cent{G^*}s$, let $\hat s=\Psi_{L,s}(1_{L^F})$,
 and let $L(s)$ be a Levi subgroup of $G(s)$ in
duality with $\cent{L^*}s$. Let
$\al\in\ser{\cent{G^*}s^F}1$ and let  $\beta \in \ser{\cent{L^*}s^F}1$. By assertion (ii) we have 
$$\Psi_{G,s}(\al) =
\Lu{G(s)}G(\hat s \otimes
\Psi_{G(s),1}(\al)), \quad\Psi_{L,s}(\beta)=\Lu{L(s)}L((\Res{G(s)^F}{L(s)^F}\hat s)\otimes
\Psi_{L(s),1}(\beta))$$ hence
$$\eqalign{\Lu
LG\big(\Psi_{L,s}(\beta)\big)&=\Lu{L(s)}G\big((\Res{G(s)^F}{L(s)^F}\hat
s)\otimes
\Psi_{L(s),1}(\beta)\big)=\Lu{G(s)}G\big(\hat s \otimes
\Lu{L(s)}{G(s)}(\Psi_{L(s),1}(\beta))\big)\cr &= \Lu{G(s)}G\big(\hat
s
\otimes
\Psi_{G(s),1}(\Lu{\cent{L^*}s}{\cent{G^*}s}\beta)\big)=
\Psi_{G,s}(\Lu{\cent{L^*}s}{\cent{G^*}s}\beta)\,.\cr}$$

There are deeper and stronger results in that hypothesis, existence of a perfect isometry [8] and, better, of a Morita equivalence (see [5], or [16] Chapters 10--12 for details).

Now let two semi-simple elements
$s$, $t$  in $(G^*)^F$, with coprime order, such that
$st=ts$. Assume $\cent{G^*}t$ is a Levi subgroup of $G^*$, let $G(t)$ be a dual Levi subgroup in $G$. For any $\al\in\ser{\cent{G^*}{st}^F}1$ one has
$\Psi_{G,st}(\al)=\Lu{G(t)}G(\hat t\otimes
\Psi_{G(t),s}(\al)).$

(c) If $L$ is split, $\Lu LG$ is Harish-Chandra induction and (J.3) is satisfied. Indeed Jordan \dec\ is defined such that Harish-Chandra series correspond.
If $L_0$ is an $F$-stable Levi complement of an $F$- parabolic subgroup of  $G$ (split Levi subgroup), and $\la\in\ser{L_0^F}s$ is cuspidal, that is $(L_0,\la)$ is a cuspidal datum in $(G,F)$, to $(L_0,\la)$ there corresponds a cuspidal datum $(L_0^*(s),\al)$ in  $(\cent{G^*}s,F)$ where $L^*_0(s)=\cent{L^*}s$ for some Levi subgroup $L^*$ of $G^*$ in duality with $L$, $\Psi_{L,s}(\al)=\la$ and 
the Hecke algebras one obtain as endomorphism algebras of Harish-Chandra modules are built on isomorphic relative rational Weyl groups
$W:={\rm W}_{G^F}(L_0,\la)\cong {\rm W}_{\cent{G^*}s^F}(L^*_0(s),\al)$
(see [26],  especially
theorems~8.6 and 4.23). Then for all
$\beta\in\ser{\cent{G^*}s^F}1$,  one has 
$$\scal{\Lu
{L_0^*(s)}{\cent{G^*}s}\al}\beta{\cent{G^*}s^F}=\scal{\Lu
{L_0}G\la}{\Psi_{G,s}(\beta)}{G^F}$$
because that integer is the degree of the same element of $\II W$ that is associated to $\al$ and to $\Psi_{L,s}(\al)$ by the maps $\II {W}\to \ser{\cent{G^*}s^F}1$ and 
$\II {W}\to \ser{G^F}s$. Furthermore if one extends linearly the preceding maps to 
 $\eta_G\colon \Z\II{W}\to
\Z\ser{G^F}s$, $\eta_{G^*,s}\colon \Z\II{W}\to
\Z\ser{\cent{G^*}s^F}1$ and if  $L$ is a split Levi-subgroup of $G$ with $L_0\subseteq L$, one has  $\Lu
LG\circ\eta_L=\eta_G\circ\Ind{{\rm W}_{L^F}(L_0,\zeta)}{ W}$, and similarly $\Lu
{\cent{L^*}s}{\cent{G^*}s}\circ\eta_{L^*,s}
=\eta_{G^*,s}\circ\Ind{{\rm W}_{\cent{L^*}s^F}(L^*_0(s),\al)}{W}$.
But
$\Psi_{G,s}\circ\eta_{G^*,s}=\eta_G$ and 
$\Psi_{L,s}\circ \eta_{\cent{L^*}s,s}=\eta_L$,  (J.3) is satisfied.

From the construction of $\Psi_{G,s}$ with the Weyl group of $\cent{G^*}s$, it follows that,
if
$[\cent{L^*}s,\cent{L^*}s]=[\cent{G^*}s,\cent{G^*}s]$, (J.3) is true. Under that assumption one has, when $\al\in\ser{\cent{G^*}s^F}1$,
$$\Lu
LG(\Psi_{L,s}(\Res{\cent{G^*}s^F}{\cent{L^*}s^F}\al))
=\Psi_{G,s}(\al)$$
(here to simplify the formula we have identify  {\it via} restriction $\ser{\cent{G^*}s}1$ with $\ser{[\cent{G^*}s,\cent{G^*}s]^F}1$).


\medskip\noindent{\bf 1.3.4. Proposition. 
}{\sl Let
$\sigma \colon (G,F)\to (H,F)$ be an  isotypic morphism  between connected reductive 
algebraic  groups defined on $\F q$. Restriction from $G^F$ to $H^F$ has no multiplicity : 
$$\forall (\chi,\xi)\in\II{G^F}\times \II{H^F}\quad \scal{\Res{H^F}{G^F}\xi}{\chi}{G^F}\in\{0,1\}\, .$$}
{\it Proof. } Let $\eta\colon (K,F)\to (H,F)$ be a regular covering defined on
$\F q$, $[K,K]$ is simply connected. The adjoint groups of $H^*$,
$G^*$, $K^*$ are isomorphic and are in duality with $[K,K]$. Thus $[K,K]\to [H,H]$ factor through $[G,G]$. One has
commutative diagrams

$$\diag { [K,K]}{}{K}{\eta}{H}{\sigma}{G}{}{}\quad \diag { [K,K]^F}{}{K^F}{\eta}{H^F}{\sigma}{G^F}{}{}$$ where $K^F\to H^F$ is onto,  as the kernel of $\eta$ is a torus, and $\sigma(G^F)$ is a normal subgroup of $H^F$. By a theorem of  Lusztig, (see
[16] Theorem~15.11, chapter 16), the restriction from
$K^F$
to
$[K,K]^F$ has no multiplicity. Thus the composed restriction
from $H^F$ to $[K,K]^F$ has no multiplicity, hence there is no mutiplicity between $H^F$ and $G^F$.
\bull 

\medskip\noindent{\bf 1.3.5. Embeddings and Jordan \dec. } 

Let $G\subseteq H$ be a regular embedding. What happens to Jordan \dec\ in $\II H$, as described by Proposition~1.3.2, by restriction from $H^F$ to $G^F$ ? An answer is given by Lusztig in [27]. 

In type $\AA$  Bonnaf\'e [3] for non twisted type  and Cabanes [12] for twisted type have define an explicit one-to-one map from $\ser{\cent{G^*}s^F}1$ onto $\ser{G^F}s$, where $\ser{\cent{G^*}s^F}1$ has to be understood in the extended definition for non connected reductive groups (see 1.3.1). See our section 5.4 for more details.

Let $\sigma^*\colon H^*\to G^*$ be a dual morphism
of the inclusion of $G$ in $H$, it is a regular covering and  by restriction of $\sigma^*$ the map $H^{*\,F}\to G^{*\,F}$ is onto. Let $t\in H^{*\,F}$ and $\sigma^*(t)=s\in G^{*\,F}$. We have seen in Proposition~1.2.4 that ${\rm A}_{G^*}(s)$ is isomorphic to a subgroup of $\cF(G^*)=\Ker \sigma^*\cap [H^*,H^*]$. By duality  $(\Ker{\sigma^*})^F$ is isomorphic to $\II{H^F/G^F}$ [16] (15.2). Therefore there is an injective map  $$\sigma_{  H,s}\colon {\rm A}_{G^*}(s)^F
\mapright{}\,
\II{  H^F/G^F}\leqno{(1.3.5.1)}$$
If $  M$ is an $F$-stable Levi subgroup of $H$ and $L=M\cap G$ such that $s$ belongs to a dual Levi subgroup $L^*$ of $G^*$, then, through the injective morphism of ${\rm A}_{L^*}(s)^F$ in ${\rm A}_{G^*}(s)^F$ and isomorphism $H^F/G^F\cong M^F/L^F$,  $\sigma_{  M,s}$ is the restriction of $\sigma_{  H,s}$.

The existence of $\sigma_{H,s}$ allows us to define a subgroup $\tau_{  H,s}(A')$ 
of $H^F$ for any subgroup $A'$ of ${\rm A}_{G^*}(s)^F$ by the formula
$$\sigma_{  H,s}( A')=\II{  H^F/\tau_{  H,s}(A')}\leqno{(1.3.5.2)}$$
 
One has $G^F\subseteq \tau_{H,s}(A')$ and $A_1\subseteq A_2$ implies $\tau_{  H,s}(A_2)\subseteq  \tau_{  H,s}(A_1)$. 

By duality between an  abelian group and its group of characters $\sigma_{H,s}$ gives an isomorphism 
$$\sigma_{  H,s}^\vee\colon H^F/\tau_{H,s}({\rm A}_{G^*}(s)^F)
\mapright{}\,({\rm A}_{G^*}(s)^F)^\wedge
\leqno{(1.3.5.3)}$$

The group $(H^F/G^F)^\wedge$ acts by tensor product on $\II{H^F}$. It has been proved that  $\sigma_{  H,s}({\rm A}_{G^*}(s)^F)$ is the stabilizer of the subset $\ser{  H^F}t$ of $\II{H^F}$.  So  ${\rm A}_{G^*}(s)^F$ acts on $\ser{H^F}t$. The same group ${\rm A}_{G^*}(s)^F$, as quotient $\cent{G^*}s^F/\cento{G^F}s^F$,  acts on $\ser{\cento{G^*}s^F}1$. The set $\ser{\cento{G^*}s^F}1$ identifies with $\ser{\cento{H^*}t^F}1$ by restriction through  $\sigma^*$. A fundamental result [27] is that $\Psi_{H,t}$ is a morphism for these two dual operations of ${\rm A}_{G^*}(s)^F$ : 
$$\Psi_{H,t}(\lexp a \al)\otimes \sigma_{H,s}(a)=\Psi_{H,t}(\al),\quad (a\in{\rm A}_{G^*}(s)^F, \al\in\ser{\cent{H^*}t^F}1)\leqno (1.3.5.4)$$

For any $\chi\in\ser{H^F}t$, $\Res{H^F}{G^F}\chi$ belongs to $\Z \ser{G^F}s$ [16] Proposition~15.6. Thus $({\rm A}_{G^*}(s)^F)^\wedge$ acts on $\ser{G}s$ through $(\sigma^\vee_{H,s})^{-1}$ (1.3.5.3). That action is independant of the choice of $H$ in the regular embedding $G\subseteq H$ and of $t$ such that $\sigma^*(t)=s$. Using the non-multiplicity property in restriction from $H^F$ to $G^F$ (Proposition~1.3.4), one obtains [27] Proposition~8.1, [16] Corollary~15.14 :

\medskip\noindent{\bf 1.3.6. Proposition. }{\sl  Let $\sigma\colon (G,F)\to (H,F)$ be a regular embedding, let $\sigma^*\colon (H^*,F)\to (G^*,F)$ be a dual morphism, $t$ semi-simple in $H^{*\,F}$, $s=\sigma^*(t)$. Denote $A={\rm A}_{G^*}(s)^F$. There is a bijective map between sets of orbits 
$$\ser{G^F}s/{A}^\wedge\longleftrightarrow\ser{\cento{G^*}s^F}1/A\;.\leqno{
(1.3.6.1)}$$ 
Let $\al\in\ser{\cento{G^*}s^F}1$. To the orbit of $\al$ under $A$ in $\ser{\cento{G^*}s^F}1$ there corresponds the set of \irr\ components of $\Res{H^F}{G^F}(\Psi_{H,t}(\al))$, it is a regular orbit under $A^\wedge/(A_\al)^\perp\cong (A_\al)^\wedge$.}

 With
notations introduced in Appendix, 5.1  $${\rm if}\; \Psi_{H,t}(\al)\in\II{H^F\mid \chi}, {\rm then}\;\tau_{  H,s}({\rm
A}_{G^*}(s)^F_\al)={\rm I}^{H^F}_{G^F}(\Psi_{H,t}(\al))=H^F_\chi.\leqno{(1.3.6.2)}$$

 \medskip\noindent{\bf 1.3.7. Proposition. }{\sl Let
$\sigma \colon (G,F)\to (H,F)$ be an  isotypic morphism  between connected reductive 
algebraic  groups defined on $\F q$.  Let
$\sigma^*\colon H^*\to G^*$ be a dual morphism,   $K^*=\Ker{\sigma^*}$, and let $t$ be a semi-simple element of $H^{*\,F}$, $s=\sigma^*(t)\in G^*$,
$\zeta\in\ser{H^F}t$. 

(a) $\sigma^*$ defines by restriction from $H^{*\,F}$ to $G^{*\,F}$ a bijection
$\ser{\cento{H^*}t^F}1\to
\ser{\cento{G^*}s^F}1$ and an injective morphism $\tau \colon {\rm A}_{H^*}(t)\to  {\rm A}_{G^*}(s)$ that transforms the action of ${\rm A}_{H^*}(t)^F$ on $\ser{\cento{H^*}t^F}1$ in the action on $\ser{\cento{G^*}s^F}1$ of its image in ${\rm A}_{G^*}(s)^F$. The quotient ${\rm A}_{G^*}(s)/\tau({\rm A}_{H^*}(t))$ is isomorphic to a subgroup
of $K^*\cap [H^*,H^*]$. 

(b) Let $\al\in\ser{\cento{H^*}t^F}1 $ be in the orbit under ${\rm A}_{H^*}(t)^F$ that is associated to the orbit of $\zeta\in \ser{H^F}t$ under $({\rm A}_{H^*}(t)^F)^\wedge$ by (1.3.6.1).  

(b.1) If $\z H$ is connected, then
$\Psi_{H,s}(\al)=\zeta$.  
 
(b.2) $\Res{}{G^F\to H^F}\zeta$ is a sum of 
$|{\rm A}_{G^*}(s)^F_\al/\tau({\rm A}_{H^*}(t)^F_\al)|$ distinct elements of $\ser{G^F}s$ and there are $|{\rm A}_{G^*}(s^*)^F|.|{\rm A}_{H^*}(t)^F_\al |/|{\rm A}_{G^*}(s^*)^F_\al |.|{\rm A}_{H^*}(t)^F|$ elements of $\ser{H^F}t$ with equal restriction $\Res{}{G^F\to H^F}\zeta$.}

\preuve  
  We know by Proposition~1.2.3 that  $\sigma^*(\cento{H^*}t).\zo{G^*}=\cento{G^*}s$ hence $\sigma^*$ defines a map 
$\tau\colon {\rm A}_{H^*}(t)\to {\rm A}_{G^*}(s)$ that commute with $F$. By Lemma~1.2.1 the quotient ${\rm A}_{G^*}(s)/\tau({\rm A}_{H^*}(t))\cong \cent{G^*}s/\sigma^*(\cent{H^*}t)$ is isomorphic to a subgroup of $\Ker \sigma\cap[H^*,H^*]$.  The induced morphism $\cento{H^*}t\to  \cento{G^*}s$ is isotypic and identify the unipotent series of $\cento{H^*}t^F$ and $\cento{G^*}s^F$.  

Consider regular embeddings defined over $\F q$ as in the proof of 1.1.4 (d) , and dual morphisms :
$$\def\normalbaselines
{\baselineskip20pt\lineskip3pt
\lineskiplimit3pt}\matrix
{
G&\hookrightarrow{}&G_0\cr
\mapdown{\sigma}&&\mapdown{\sigma_0}\cr
H&\hookrightarrow{}&H_0\cr}  \quad\quad \matrix
{
G^*&\mapleft{}&G^*_0\cr
\mapup{\sigma^*}&&\mapup{\sigma^*_0}\cr
H^*&\mapleft{} &H^*_0\cr}  \leqno{(1.3.7.1)} $$
Horizontal maps in the right diagram are coverings with central tori as kernels.

(b.1) follows directly from Propositions 1.2.3 and 1.3.6.

(b.2) follows from Proposition~1.3.6 and Clifford theory applied in the commutative diagram
$$\def\normalbaselines
{\baselineskip20pt\lineskip3pt
\lineskiplimit3pt}\matrix
{
G^F&\hookrightarrow{}&G_0^F\cr
\mapdown{\sigma}&&\mapdown{\sigma_0}\cr
H^F&\hookrightarrow{}&H_0^F\cr}  \leqno{(1.3.7.2)}$$
where all maps have invariant images and abelian cokernels.

 Let $t_0$ be a semi-simple element in $ H_0^{*\,F}$ of image $t\in H^*$ and put  $s_0=\sigma^*_0(t_0)\in G^*_0$, so that $s_0$ maps on $s$. By restriction in (1.3.7.1) we obtain one-to-one maps between $\ser{\cento {H^*}t^F}1$, $\ser{\cento {H_0^*}{t_0^F}}1$, $\ser{\cento
{G^*}s^F}1$, $\ser{\cento {G^*_0}{s_0}^F}1$ and we identify theses sets. The groups ${\rm A}_{H^*}(t)^F$, ${\rm A}_{H^*_0}(t_0)^F$, ${\rm A}_{G^*}(s)^F$,
 and ${\rm A}_{G_0^*}(s_0)^F$ act respectively on these unipotent series.  As $\z{G_0}$ and $\z{H_0}$ are connected, ${\rm A}_{G_0^*}(s_0)=1={\rm A}_{H_0^*}(t_0)$. Furthermore $\tau$ commute with $F$ and actions on unipotent series. 

Let $\zeta_0\in\II{H^F_0\mid \zeta}\cap \ser{H_0^F}{t_0}$. Assume $\zeta_0=\Psi_{H_0,t_0}(\al)$ for a Jordan \dec\ $\Psi_{H_0,t_0}$ as in Proposition 1.3.5, where $\al \in\ser{\cent{H_0^*}{t_0}^F}1$. We have by  (1.3.6.2)
$$\Res{H_0^F}{H^F}(\zeta_0)=\sum_{h\in H_0^F/\tau_{H_0,t_0}({\rm A}_{H^*}(t)^F_\al)}\lexp h{\zeta}\leqno{(1.3.7.3)}$$
hence $\zeta_0(1)=\zeta(1).|{\rm A}_{H^*}(t)^F_\al|$. There are $|{\rm A}_{H^*}(t)^F/{\rm A}_{H^*}(t)^F_\al |$ elements in $\ser{H_0^F}{t_0}$ with equal restriction to $H^F$. We call ``multiplicative factor of $\Res{H^F_0}{H^F}$ above $\al$" the quotient of cardinals of the subsets of $\ser{{H}^F}{t}$ and $\ser{{H_0}^F}{t_0}$ that   correspond to $\al$, here it  is $|{\rm A}_{H^*}(t)_\al^F|^2/|{\rm A}_{H^*}(t)^F|$.

As the kernel of $\sigma_0$ is a torus, $\sigma_0$ restricts on a surjective morphism $G_0^F\to H_0^F$, so that $\Res{}{\sigma_0\colon G_0^F\to H_0^F}$ sends $\ser{H^F_0}{t_0}$ in $ \ser{G_0^F}{s_0}$. We may assume $$\chi_0:=\Res{}{\sigma_0\colon G_0^F\to H_0^F}(\zeta_0)=\Psi_{G_0,s_0}(\al)\leqno{(1.3.7.4)}$$ for suitable Jordan \dec\ $\Psi_{G_0,s_0}$ and identification. Then we have, if $\chi_0\in\II{G_0^F\mid \chi}$ :
$$\Res{G_0^F}{G^F}(\chi_0)=\sum_{g\in G_0^F/\tau_{G_0,s_0}({\rm A}_{G^*}(s^*)^F_\al)}\lexp g{\chi}\leqno{(1.3.7.5)}$$

By what we have seen in 1.3.5 $$ \tau_{G_0,s_0}({\rm A}_{G^*}(s^*)^F_\al))\subseteq\tau_{G_0,s_0}(\tau({\rm A}_{H^*}(t^*))^F_\al)=\sigma_0^{-1}(\tau_{H_0,t_0}({\rm A}_{H^*}(t^*)^F_\al))$$
with $\chi_0(1)=\chi(1).|{\rm A}_{G^*}(s^*)^F_\al |$ and the multiplicative factor of $\Res{G_0^F}{G^F}$ above $\al$ is $|{\rm A}_{G^*}(s)^F_\al |^2/|{\rm A}_{G^*}(s)^F|$.

In view of (1.3.7.2), (1.3.7.3), (1.3.7.4) and (1.3.7.5) we conclude by transitivity of restriction and that $\chi(1).|{\rm A}_{G^*}(s)^F_\al |=\zeta(1).|{\rm A}_{H^*}(t)^F_\al |$. 

The multiplicative factor of $\Res{H^F}{G^F}$ above $\al$ is $|{\rm A}_{G^*}(s^*)^F_\al |^2.|{\rm A}_{H^*}(t)^F |/|{\rm A}_{G^*}(s^*)^F |.|{\rm A}_{H^*}(t)_\al^F|^2$.
\bull

\bigskip\noindent{\bf 1.4. Theorem. }{\sl  Let $(G,F)$ be a reductive algebraic  group over an
algebraic closure of a prime field  $\F p$, defined over $\F q$, with Frobenius
endomorphism $F$, in duality over $\F q$ with ($G^*,F)$. Let $G^F$ be the subgroup of rational points.  Assume that Mackey decomposition formula holds for Lusztig induction inside closed $F$-stable subgroups of $G$ and of $G^*$.

Let $\ell$ be an odd
prime number, different from $p$. Assume that $\ell\geq 7$ if $G$ has a component of type $\EE_8$, $\ell\geq 5$ if $G$ has a component of non classical type. Let $(G,F)\to( H,F)$ be a regular emdedding with a dual morphism $H^*\to G^*$. Let $s\in G^{*\, F}_{\ell'}$, image of $t \in H^{*\, F}_{\ell'}$.

There exist a reductive group defined over $\F q$, we denote $(G(s),F)$, such that

(A) $(G(s)^\circ,F)$ and $(\cento{G^*}s,F)$ are in duality, $G(s)/\cento{G^*}s$ is isomorphic to ${\rm A}_{G^*}(s)$ (see (1.2.0)), and so acts on $G(s)$ coherently with the action  of ${\rm A}_{G^*}(s)$ on the root datum of $\cento{G^*}s$. 

(B) There exist a one-to-one map $${\cal B}_{G,s}\colon \bl{G(s)^F}1\to \bl {G^F}s$$ 
 from the set of unipotent $\ell$-blocks of $G(s)^F$ onto the set of $\ell$-blocks of $G^F$ in series $(s)$ such that
 
 (B.1) If $G=H$, then Jordan decomposition defines a one-to-one map from $\II b$ onto $\II {{\cal B}_{H,t}(b)}$ (where $b\in\bl{H(t)^F}1$) : if $h\in \cent{H^*}t^F_\ell$, so that $\cent{H^*} {th}=\cent{H(t)^*}h$, then 
 $$(\Psi_{H,th}\circ{\Psi_{H(t),h}}^{-1})(\II b\cap\ser{H(t)^F}h)=\II{{\cal B}_{H,t}(b)}\cap\ser{H^F}{th}$$
 
 (B.2) ${\cal B}_{G,s}$ and ${\cal B}_{H,t}$ respect Clifford theory in the following sense.
 
 The regular embedding $G\to H$ defines an embedding $G(s)^\circ \to H(t)$. Let $c$ in $\bl{H(t)^F}1$, $b$ in $\bl{G(s)^F}1$, $b_0\in\bl{G(s)^{\circ\,F}}1$ such that $c$ restricts to $b_0$. Then the block ${\cal B}_{H,t}(c)$ covers ${\cal B}_{G,s}(b)$ if and only if $b$ covers $b_0$.
 
 (B.3) There is a one-to-one height preserving map $\Psi_b$ from $\II b$ onto $\II{{\cal B}_{G,s}(b)}$ such that $$\Psi_b(\zeta)(1)|G(s)^F|_{p'}=\zeta(1)|G^F|_{p'}$$
  for all $\zeta\in\II b$ 
  
  (B.4) The defect groups of $b\in\bl{G(s)^F}1$ and of ${\cal B}_{G,s}(b)$ are isomorphic. The Brauer categories of $b$ and of ${\cal B}_{G,s}(b)$ are equivalent.}
  
  \preuve The all proof is contained in sections 2 to 5, we use here Propositions~2.1.4, 2.1.7, 2.1.10, 2.3.5, 2.3.6, 3.1.1, 3.4.1, 3.4.2, 4.1.2, 4.2.4.
  
  (A) The group $G(s)$ such that (A) holds is constructed in sections 3.1, 3.2, 3.3, see Proposition~3.1.1.
  
  (B.1) By  Proposition~2.1.7 a unipotent $\ell$-block $b$ of $H(t)^F$ is defined by a $H(t)^F$-conjugacy class of unipotent cuspidal data $(L(t),\al(t))$ in $(H(t),F)$. See definitions~2.1.1 : $L(t)$ is a $d$-split Levi subgroup of $H$, $\al(t)$ is $d$-cuspidal in $\ser{L(t)^F}1$), that gives the block $b=b_{H(t)^F}(L(t),\al(t))$.  Equivalently $b$ is defined by a $\cent{H^*}t^F$-conjugacy class of unipotent cuspidal data $(L(t)^*,\al)$ in $(\cent{H^*}t,F)$ where $\al(t)=\Psi_{L(t),1}(\al)$. 
  
  By Proposition~2.1.7 again, an $\ell$-block $B$ of $H^F$ in series $(t)$ is defined by a $H^F$-conjugacy class of cuspidal data $(L,\la)$ in series $(t)$ in $(H,F)$. 
  
  Thanks to Jordan decomposition inside $d$-split Levi subgroups of $H$, there is a one-to-one map between the two sets of classes of cuspidal data, $(L(t)^*,\al)/\cent{H^*}t^F\mapsto (L,\Psi_{L,t}(\al))/H^F$, where $L$ is in duality with $\cent{G^*}{\zo{L^*}_{\phi_d}}$, see Propositions~2.1.4 and 2.1.10. So is defined $B:={\cal B}_{H,t}$ in Proposition~3.4.1, by $B={\cal B}_{H,t}(b_{H(t)^F}(L(t),\Psi_{L(t),1}(\al)))=b_{H^F}(L,\Psi_{L,t}(\al))$.
  
  One knows (1.3.1.7) that $\II B\subseteq \cup_{h\in \cent{H^*}t^F_\ell}\ser{H^F}{th}$. Given $h\in\cent{H^*}t^F_\ell$,  replace $(G,s,s_0)$ by $(H,t,h)$ in Propositions~2.3.5, 2.3.6, where $\II{B}\cap\ser{H^F}{th}$ is computed as a generalized $d$-HC series. Let $H(h)$ be an $E$-split Levi subgroup of $H$ in duality with $\cent{H^*}h$ (notation in coherence with $G(s)$, but $h=h_\ell$)). Using (2.3.5.3), (2.3.6.2), (2.3.6.3) we obtain a one-to-one map, restriction of $\Psi_{H,th} $, from $\ser{\cent{\cent{H^*}h}t^F}{(H^*_{th},\al_h)}$ onto $\II B\cap\ser{H^F}{th}$, where $(H^*_{th},\al_h)$ is a unipotent $d$-cuspidal datum in $(\cent{\cent{H^*}h}t^F,F)$, $(H^*_{th},\al_h)$ being defined from $(L(t)^*,\al)$ as $(\cent{L^*_0}s,\al_0)$ from $(L^*_s,\al)$ in Proposition~2.3.6 :
  $$  \Psi_{H,th}\big( \ser{\cent{\cent{H^*}h}t^F}{(H^*_{th},\al_h)}\big)=\II B\cap\ser{H^F}{th}$$
    
  One has  $\cent{\cent{H^*}h}t=\cent{\cent{H^*}t}h$ and  $(H^*_{th},\al_h)$ defines by Proposition~2.1.7 a $H(t)^F$-conjugacy class of $d$-cuspidal data in series $(h)$ in $(H(t),F)$, say $(L_h(t),\la_h(t))$, where $\la_h(t)=\Psi_{L_h(t),h}(\al_h)$. Propositions 2.3.5, 2.3.6 apply again in that special case, replace $(G,s,s_0)$ by $(H(t),1,h)$ or by $(H(t),t,h)$ ($t$ is central in $H(t)^*$), $B$ by $B(t):=b_{H(t)^F}(L(t),\al(t))$. We have $\II {B(t)}\cap \ser{H(t)^F}h=\ser{H(t)^F}{(L_h(t),\al_h(t)}$ hence 
  $$\Psi_{H(t),h}\big(\ser{\cent{\cent {H^*}t}h^F}{(H^*_{th},\al_h)}\big)=\II b\cap\ser{H(t)^F}{h}$$
  
   That is (B.1). If $G=H$, $\Psi_b$ is the restriction of various $ \Psi_{H,th}\circ {\Psi_{H(t),h}}^{-1}$, for $h\in\cent{H^*}t^F_\ell$.
   
   The equality on degrees in (B.3) when $G=H$ is an immediate consequence of the definition of $\Psi_b$ and degree formula (vi) in Proposition~1.3.2, applied in $H$ and in $H(t)$ (we have yet used  $\cent{H(t)^*}h=\cent{H^*}{th}$). 
   
  There is a kind of duality between the two Clifford situations for $G(s)^{\circ\,F}\subseteq G(s)^F$ and $G^F\subseteq H^F$, when restricted to the sets $\II b$ or $\II{{\cal B}_{G,s}(b)}$ respectively, under the action of $A:={\rm A}_{G^*}(s)^F\cong G(s)^F/G(s)^{\circ\,F}$ on $\bl{G(s)^{\circ\,F}}1$ and of $A^\wedge$ on $\bl{H^F}t$. See Proposition
  3.4.1, where is proved (B.2). 
  
  (B.3) The  existence of $\Psi_b$ is proved in Proposition 3.4.2. The degree formula in (B.3) for $G$ is deduced from the degree formula for $H$.
  
  As $G(s)^F/G(s)^{\circ\,F}$ is prime to $\ell$, if a block $b$ of $G(s)^F$ covers  a block $b_0$ of $G(s)^{\circ\,F}$, then a defect group of $b_0$ is a defect group of $b$. The defect groups of $b$ and ${\cal B}_{G,s}(b)$ are computed in Proposition~4.1.2, they are isomorphic. 
  
  Then the degree formule in (B.3) imply that $\Psi_b$ preserves height.
  
 The Brauer category of a block $B$ of a finite group $X$ is a small category whose objects are $\ell$-subpairs that contain $(\{1\},B)$ and morphisms are produced by restriction of interior automorphisms of $X$. By fusion theorems (see [32], \S47 and (48.3)), once the defect groups $D$ of a block is  given, the Brauer category is entirely defined by the  groups of automorphisms 
 $$\nor X {Y,b_Y}/Y.\cent XY, \quad Y\subseteq D,\quad (\{1\},B)\subset (Y,b_Y)\subset (D,b_D)$$
 Thus (B.4) follows from Proposition 4.2.4. \bull

  \bigskip\noindent{\bf 1.5. $\ell=2$ for classical types in odd characteristic }

In this section we obtain Jordan \dec\ of blocks for $\ell=2$ in classical types. We assume Mackey's decomposition formula.

Two facts help us  : there is only one unipotent $2$-block and if $s$ is a semi-simple element of odd order in $G^{*\, F}$, then $\cento{G^*}s$ is a Levi subgroup of $G^*$. Short of type $\AA$, where $2$ is good, $\cent{G^*}s$ is connected.

\medskip\noindent{\bf 1.5.1. Hypotheses and notations. } {\sl $G$ is a reductive group defined over $\F q$, $q$ odd. All components of $G$ have classical type. 

Let $s\in G^{*\,F}$ be semi-simple of odd order. Let $G(s)^\circ$ be a Levi subgroup of $G$ in the dual $G^F$-conjugacy class of the $G^{*\,F}$-class of  $\cento{G^*}s$, duality around dual maximal torii $T$, $T^*$. Let $G(s)\subseteq \nor G{G(s)^\circ}$ such that $G(s)^\circ\subseteq G(s)$, and $G(s)/G(s)^\circ\cong \nor{G(s)}T/\nor{G(s)^\circ}T$ has image ${\rm A}_{G^*}(s)$ by the antimorphism between Weyl groups $ \nor GT/T$ and $\nor {G^*}{T^*}/T^*$.}

\medskip\noindent{\bf 1.5.2. Proposition. }{\sl Assume 1.5.1. 

(a) 
$G^F$ has only one unipotent $2$-block, the  principal block. 

(b) If
${\rm A}_G(s)^F=\{1\}$, and that is the case in types $\BB$, $\CC$ and $\DD$, 
$G^F
$ has only one $2$-block $b(s)$ in  series $(s)$, such that $\II{b(s)}={\cal E}_2(
{G^F},s)$. 

(c)  A $2$-Sylow subgroup of $G(s)^F$ is  a defect group of any
 $2$-block in series $(s)$. Such a block has central defect  group if and only if  $\cento{G^*}s=T^*$  and $T^F_2\subseteq \z G$.
 }

 The condition
$T^F_2\subseteq \z G$ is satisfied when  $G^F={\rm SL}_{2m}(q)$, $m$ is odd, $q\equiv 3\pmod 4$ and $T$ is a   Coxeter torus in $G$.
There is a similar example in other classical types, $G$ simply connected.

\preuve   The notation ${\cal E}_\ell( {G^F},s)$  has been introduced in 1.3.1, see (1.3.1.7). 

All properties reduce easily to the rationally \irr\ case. 

Assertion~(a) is [16], Theorem~21.14.

(b) Assuming $G(s)=G(s)^\circ$, the virtual bimodule defining $\Lu {G(s)}G$ allows to construct a perfect isometry [8] and  Morita equivalence [5] between the principal block $b_1$ of $G(s)^F$, unique in series $(1)$ in $G(s)^F$ by (a), and $b(s)$, therefore unique in series $(s)$ in $G^F$. The application $(\la\mapsto \Lu {G(s)}G(\Psi_{G(s),s}(1)\otimes \la))$ is a one-to-one map from $\II{b_1}={\cal
E}_2(G(s)^F,1)$ onto $\II{b(s)}={\cal
E}_2(G^F,s)$. Furthermore a $2$-Sylow subgroup of $G(s)^F$ is a common defect group of the two blocks and the map $(\la\mapsto \Lu {G(s)}G(\Psi_{G(s),s}(1)\otimes \la))$ preserve height by degrees formulas (1.3.1.2) and (vi) in Proposition~1.3.2. 

Clearly (c) is true if ${\rm A}_{G^*}(s)^F=\{1\}$.

That happens if
$G$ has only types $\BB$, $\CC$, $\DD$, because ${\cal F}(G)$ is then a $2$-group (Proposition~1.2.6  (e)).

 (c)  We have to consider rational types $\AA$, $\lexp 2\AA$. 
 
 Let
$G\to H$ be a regular embedding, let
$t$ semi-simple in $H^{*\,F}_{2'}$ of image 
$s$ by a dual map $\sigma\colon H^*\to G^*$. There exists a Levi subgroup $H(t)$ of $H$ such that $G(s)^\circ=H(t)\cap G$ and $H(t)$ is in duality with $\cent{H^*}t$. Assertion (b) applies to $(H,t)$. Let $B(t)$ be the unique $2$-block in series $(t)$ of $H^F$.

Any element of $\ser{H^F}t$ (resp. ${\cal E}_2(H^F,t)$) restricts on $G^F$ on a sum of elements of $\ser{G^F}s$ (resp. ${\cal E}_2(G^F,s)$) and all $\ser{G^F}s$ (resp. ${\cal E}_2(G^F,s)$) appear this way. Thus $B(t)$ covers all $2$-blocks in series $(s)$ of $G^F$ and these blocks are $H^F$-conjugate. 

We have seen that a $2$-Sylow subgroup $E$ of $H(t)^F$ is a defect group of $B(t)$. One has $G^F.H(t)^F=H^F$, $G^F.E/G^F$ is a $2$-Sylow subgroup of $H^F/G^F$ and $E\cap G^F$ is a $2$-Sylow subgroup of $G(s)^F$. By  [28] Chapter 5, 5.16, $E\cap G^F$ is a defect group of any block of $G^F$ covered by $B(t)$.

 \bull

We can now describe Brauer subpairs by standard arguments of local theory of blocks.

\medskip\noindent{\bf 1.5.3. Proposition. }{\sl Assume 1.5.1. Let $Y$ be a $2$-subgroup of $G(s)^F$, $T_Y$ be a maximal $F$-stable torus of $\cento GY$, $T$ be a maximal $F$-stable torus of $G$ such that $Y\subseteq \nor GT$ and $T_Y\subseteq T$. Let $T^*$ be a maximal torus in $\cento{G^*}s$ in the $G^{*\, F}$-conjugacy class of the $G^F$-class of $T$. Let $T^*\to T^*_Y$ be a dual morphism of the inclusion of $T_Y$ in $T$, it defines a Lusztig series $\ser{\cento GY^F}{s_Y}$. 

If $b$ is a $2$-block of  $\cent GY^F$ that covers a block in series $(s_Y)$ of $\cento GY^F$ and $(\{1\},B)\subset (Y,b)$ is an inclusion of subpairs in $G^F$, then $B$ is  in series $(s)$.}

\preuve We know that $\cento GY$ is reductiveand the existence of $T$ has been proved in Proposition~1.2.7. 

Given a group $(\cento GY^*,F)$ in duality with $\cento GY$ around torii $T_Y$, $T_Y^*$, $s_Y\in T_Y^{*\, F}$ may be defined by $\Psi_{T_Y,s_Y}(1)=\Res{T^F}{T_Y^F}(\Psi_{T,s}(1))$, so that $\ser{\cento GY^F}{s_Y}$ is defined. In types $\AA$ or $\DD$ (c) of Proposition~1.2.7 applies hence we may see $\cento GY^*$ as an algebraic  subgroup $C^*$ of $G^*$ such that $s\in C^*$. Then we may assume $s_Y=s$.

(a) If $s=1$, our claim is Brauer's First Main Theorem, the inclusion between subpairs formed with principal blocks $B_1$ and $b$ : thanks to (a) in Proposition~1.5.2, $\cent GY^F/\cento GY^F$ is a $2$-group, the principal block of $\cent GY^F$ is the unique block that covers the principal block $b_1$ of $\cento GY^F$. 

If $s$ is central in $G$, one has blocks $B=\Psi_{G,s}(1)\otimes B_1$, $b_0=\Psi_{\cento GY,s_Y}(1)\otimes b_1$, $\Psi_{\cento GY,s_Y}(1)$ extends to $\cent GY^F$, so is defined $b$. The inclusion of subpairs follows. 

(b) Assume first $\z G$ connected and $Y$ cyclic, $Y=\gen y$.   

By Proposition~1.5.2 is given a $2$-block $B(s)$ of $G^F$.
We have  $T_Y=T$. Let $L_y:=\cento {G(s)}y=\cento{\cento Gy}s$, $L_y$ is a Levi subgroup of $\cento G y$.  Let $\hat s:=\Psi_{G(s),s}(1)$, $\hat s_y=\Psi_{L_y,s_y}(1)$, hence $\hat s_y=\Res{G(s)^F}{L_y^F}\hat s$. Let  $\xi$ in ${\cal E}_2(G(s)^{ F},1)$, then $\chi=\Lu LG(\hat s\otimes \xi)\in\II{B(s)}$. If $b(s)$ is the $2$-block of $G(s)^F$ in series $(s)$ ($\z{G(s)}$ is connected) , $b(s).(\hat s\otimes \xi)\in\II{b(s)}$. Let $\xi_y\in {\cal E}_2(L_y^{F},1)$, then $\chi_y:=\Lu{L_y}{\cento Gy}(\hat s_y\otimes \xi_y)\in \II{b_y}$ for some $2$-block $b_y$ of $\cento G y^F$ in series $(s_y)$. 

Using the second Brauer Main Theorem, to prove an inclusion of subpairs $(1,B(s))\subset (\gen y,b)$ in $G^F$, where $b$ covers $b_y$, we may consider only {\it connected subpairs} [16] 21.1, so  it is sufficient to prove $b_y.d^y(\chi)\neq 0$. As $b_y.\chi_y=\chi_y$, consider $\scal{b_y.d^y(\chi)}{\chi_y}{\cento Gy^F}=\scal{d^y(\chi)}{\chi_y}{\cento Gy^F}$. By commutation formula [16], Theorem 21.4,  adjunction and part (a) of the proof applied in in $G(s)$, we have 
$$\eqalign{\scal{d^y(\chi)}{\chi_y}{\cento Gy^F}&=\scal{d^y(\Lu {G(s)}G(\hat s\otimes\xi))}{\Lu{L_y}{\cento Gy}(\hat s_y\otimes\xi_y}{\cento Gy^F}\cr &=\scal{d^y(\slu {G(s)}G(\Lu {G(s)}G(\hat s\otimes  \xi)))}{\hat s_y\otimes \xi_y}{L_y^F}\cr &=\scal{d^y\big(b(s).\slu {G(s)}G(\Lu {G(s)}G(\hat s\otimes  \xi))\big)}{\hat s_y\otimes \xi_y}{L_y^F}\cr}$$
To compute $\slu {G(s)}G(\Lu {G(s)}G(\hat s\otimes  \xi))$ we use Mackey formula. As $G(s)$ is in duality with $\cent{G^*}s$, if $g\in G^F$ and $T^g\subseteq G(s)$, then $(T^g,(\Res{G(s)^F}{T^F}\hat s)^g)$ is conjugate to $(T,\Psi_{T,s}(1))$ only if $g\in G(s)$. So $$b(s).\slu {G(s)}G(\Lu {G(s)}G(\hat s\otimes  \xi))=\hat s \otimes \xi$$
As $\hat s_y$ is restriction of $\hat s$ we obtain
$$\scal{d^y(\chi)}{\chi_y}{\cento Gy^F}=\scal{d^y(\xi)}{\xi_y}{L_y^F}$$
Given $b_y$, the existence of $\xi_y\in\II{b_y}$ and $\xi\in {\cal E}_2(G(s)^F,1)$ such that $\scal{d^y(\xi)}{\xi_y}{L_y^F}\neq 0$ is just our claim for unipotent blocks inside $G(s)$, we proved it in (a).

(c1) To easy induction, thanks to Proposition~1.2.7 ---we use freely--- assume now $\z G$ connected and $[G,G]$ simply connected, with no restriction on the $2$-subgroup $Y$ of $G^F$. 

By (b) and induction we may assume that our claim is true for smaller groups $Y_1$ ($|Y_1|<|Y|$) and for groups $G_1$ with strictly smaller semi-simple rank.

We may assume $\z G^F_2.\z{\cent GY}^F_2\subseteq Y$, with no change on $\cent GY$. Let $L_Y=\cent{G(s)}Y=\cent{\cent GY}s$, $L_Y$ is a Levi subgroup of $\cent GY$. 

If $\z Y\neq \z G^F_2$, let $y\in\z Y\setminus \z G$. One has $Y\subseteq \cent Gy$. If $s_Y$, $s_y$ are given from $s$ as above, we have, by part (a) of the proof, three $2$-blocks, $b(s_y)$, $b(s_Y)$, $b(s)$ of respectively $\cent Gy^F$, $\cent GY^F$, $G^F$. 

By induction on the semi-simple rank of $G$ we have an inclusion of Brauer subpairs in $\cent Gy^F$ : $(1,b(s_y))\subset (Y, b(s_Y))$. That inclusion may be seen in $G^F$ : $(\gen y,b(s_y))\subset (Y,b(s_Y))$. By (b) we have an inclusion  in $G^F$ : $(1,b(s))\subset (\gen y,b(s_y))$. That imply by transitivity  $(1,b(s))\subset (Y,b(s_Y))$ in $G^F$.

It may happens that $\z Y=\z G^F_2$. If $\z Y$ is not a $2$-Sylow subgroup of $\cent{G(s)}Y^F$, let $x\in (T_Y)^F_2\setminus \z Y$ ($x$ exists thanks to a good choice of $T_Y$). Consider $Y_1=\gen{Y,x}$. We may take $T_{Y_1}=T_Y$, then $s_{Y_1}=s_Y$ and apply the preceding result on $Y_1$ to obtain an inclusion $(1,b(s))\subset (Y_1,b(s_{Y_1}))$ of subpairs in $G^F$. Inductively we have an inclusion in $\cent GY^F$ : $(1,b(s_Y))\subset (\gen x,b(s_{Y_1})$ which may be read in $G^F$ as $(Y,b(s_Y))\subset (Y_1,b(s_{Y_1}))$. By unicity in Brauer's correspondance we have $(1,b(s))\subset 
(Y,b(s_Y))$.

It may happens that $\z Y=\z G^F_2$ is a $2$-Sylow subgroup of $\cent{G(s)}Y^F$ ! Let $Y_1$ be a subgroup of $Y$ with index 2. If $\cent GY=\cent G{Y_1}$  we have our claim by inductive hypothesis. One sees that if $\cent G{Y_1}\neq \cent GY$, then $\cent G{Y_1}$ is not rationnally irreducible and that imply $\z{\cent G{Y_1}}^F_2\neq \z G^F_2$. Take $y\in Y\setminus Y_1$. We argue as above : we have an inclusion $(1,b(s))\subset (Y_1,b(s_{Y_1}))$ in $G^F$, as well as $(1,b(s_{Y_1}))\subset (\gen y,b(s_Y)$ in $\cent G{Y_1}^F$.

(c2) Assuming only $\z G$ connected, there exists an isotypic epimorphism $H\to G$, with central kernel, such that $\z H$ is connected and $[H,H]$ simply connected (see Proposition~1.1.4 (d)). Then $G^F$ is a quotient of $H^F$ with central kernel. We  consider $\II{G^F}$ as a subset of $\II{H^F}$. If $s$ has image $t\in H^{*\, F}$ by a dual morphism, $t$ defines a $2$-block $b(t)$ of $H^F$, and  we have $\II{b(s)}\subseteq \II{b(t)}$. Similarly, if $X$ is a $2$-subgroup of $H^F$ of image $Y$ in $G^F$, $\cento GY$
is image of $\cent HX$, $T_Y$ is image of some maximal torus $T_X$ in $ \cent HX$. There is an element of odd order  $t_X$  in $T_X^{*\, F}\subseteq \cent HX^*$, defining a $2$-block $b(t_X)$ of $\cent HX^F$. We may assume that  $t_X$ is image of $s_Y$. Then by construction $\II{b(s_Y)}\subseteq \II{b(t_X)}$ where $b(s_Y)$ is a $2$-block of $\cento GY^F$. 
The inclusion of subpairs $(1,b(t))\subset (X,b(t_X))$ in $H^F$ given by (c1) implies an inclusion of ``connected subpairs" $(1,b(s))\subset (Y,b(s_Y))$ in $G^F$, so an inclusion $(1,b(s))\subset (Y,b_Y)$ in $G^F$ where $b_Y$ is the $2$-block of $\cent GY^F$ that covers $b(s_Y)$.

(c3) To reach the more general case we have to consider a regular embedding $G\subseteq H$.

Let $H^*\to G^*$ be a dual map, $t$ semi-simple in $H^{*\,F}_{2'}$ of image $s$, $H(t)$ a Levi subgroup of $H$ in duality with $\cent{H^*}t$ such that $G(s)^\circ=H(t)\cap G$. Let $B(t)$ be the unique $2$-block idempotent of $H^F$ in series $(t)$, $B(s)$ be the sum of $2$-blocks idempotents of $G^F$ in series $(s)$. As $H^F$ stabilizes $\ser{G^F}s$ and ${\cal E}_2(H^F,t)$ restricts in $\Z {\cal E}_2(G^F,s)$, $B(t)$ covers the $2$-blocks in series $(s)$ of $G^F$ and these are $H^F$-conjugate, components of $B(s)$. Let $X=G^F.H^F_{2'}\subseteq H^F$, $H^F/X$ is an abelian  $2$-group. Thus $B(t) $ covers a unique $2$-block $B'(t)$ of $X$ and $B'(t)$ covers the $2$-blocks in series $(s)$ of $G^F$. The embedding $\cento GY\subseteq \cento HY$ is regular and we may define $t_Y$,  $s_Y$ semi-simple of odd order in $\cent HY^F$, $\cento GY^F$ respectively, such that $(1,B(t))\subset (Y,B(t_Y))$, where $B(t_Y)$ covers the $2$-blocks in series $(s_Y)$ of $\cento GY^F$. As $\cent GY^F/\cento GY^F$ is a $2$-group, we may speak of ``$2$-blocks in series $(s_Y)$ of $\cent GY^F$". If $b$ is such a block, it is covered by $B(t_Y)$ in $\cent HY^F$, and by a well-defined $2$-block $B'(t_Y)$ in $\cent XY$. 

Proposition~5.1.4, (a) (i) applies with $(1,Y,G^F,X)$ instead of $(U,V,Y,X)$. Assume $(1,B)\subset (Y,b)$ in $G^F$. As $B'(t_Y)$ covers $b$ and $(1,B'(t))\subset (Y,B'(t_Y))$ in $X$, $B'(t)$ covers $B$, hence $B$ is in series $(s)$.
\bull
 
\medskip\noindent{\bf 1.5.4. Proposition. }{\sl  Assume 1.5.1. 
Let $D$ be a  $2$-Sylow subgroup of $G(s)^{\circ\,F}$, $s_D\in\cento GD^*$ as in Proposition~1.5.3. The set of $2$-blocs in series $(s)$ of $G^F$ is a regular orbit under the action of $({\rm A}_{\cento
GD^*}(s_D)^F)^\wedge$. }
 
 \preuve In types $\BB$, $\CC$, $\DD$, as $s$ has odd order,  we have ${\rm A}_{G^*}s)=1$, as well as ${\rm A}_{\cento
GD^*}(s_D)=1$ (see Propositions 1.2.4, 1.1.3 and 1.2.7 (d)). There is only one $2$-block in series (s) by Proposition~1.5.2. So we may assume that $G$ has type $\AA$.

Then, as said at the beginning of the proof of Proposition~1.5.3,  for any $2$-subgroup $Y$ of $G(s)^{\circ\,F}$,  we may see $\cento{G}Y^*$ as an $F$-stable reductive subgroup $C^*$ of
$G^*$ with $s\in C^*$ (Proposition~1.2.7 (c)). As $\cento{G^*}s$ is a Levi subgroup of $G^*$ and  $\cento{C^*}s$ is a Levi subgroup of $C^*$, for some torus $S$, $\cento{G^*}s=\cent{G^*}S$, hence $\cento{C^*}s= \cent{C^*}S=C^*\cap\cento{G^*}s$. By Proposition~1.2.2 (b), with $(G^*, C^*, \gen s)$ instead of $(G,H,X)$,  ${\rm A}_{C^*}(s)$ is a subgroup of ${\rm A}_{
G^*}(s)$. If $\z G $ is connected ${\rm A}_{G^*}(s)=1$ and our claim is clear.

So we consider as usual a regular embedding
$G\subseteq H$, with dual map $H^*\to G^*$,
$t\in H^*$ semi-simple of odd order and of image $s^*$, 
$H(t)\subseteq H$, a Levi subgroup in duality with $\cento{H^*}t\subseteq H^*$ and such that $G(s)^\circ=H(t)\cap G$. We have $H=\zo H.G$, $H(t)=\zo H.G(s)^\circ$, $\cento HD=\zo H.\cento GD$, hence 

\noindent
$H^F=G^F.H(t)^F=G^F.\cent{H^F}D$, $\cento HD^F=\zo H^F.\cento GD^F$ and $\cent{H^F}D=\cento{H}D^F.\cent{G^F}D$. 

The unique $2$-block $B(t)$ in series $(t)$ of $H^F$ covers any $2$-block $b$ in series $(s)$ of $G^F$. If $(D,b_D)$ is a maximal pair of   $b$, $b_D$ is covered by a block $B_D$ of $\cent{H^F}D$, such that $(1,B(t))\subset (D,B_D)$ in $H^F$. Any $2$-block with central defect $b_D$ of $\cent{G^F}D$ covers one block with central defect $b^\circ _D$ of $\cento {G^F}D$. By conjugacy of maximal subpairs we see that $(H^F)_b=G^F.\nor{H^F}D_{b_D^\circ}$.

 The embedding $\cento GD\subset
\cento HD$ is regular and  $b_D^\circ$ has central defect.
By Proposition~1.5.2 (c), $\cento{G(s)^\circ}D$ is a torus, say $S_D$. As well
 $T_D:=\zo H.S_D=\cento {H(t)}D$ is a torus. There exists $t_D\in {T_D}^{*\, F}$, with image $s_D$, related to $t\in ({\z H.S})^{*\, F}$ as $s_D\in {S_D}^*$ is related to $s\in S^*$. Then $b_D^\circ$ is covered by the unique block $B_D$ in series $(t_D)$ of $\cent HD^F$. The canonical character of $B_D$ is $\hat {t_D}:=\Res{\cent Ht^F}{T_D^F}(\Psi_{\cent Ht,t}(1))$.
 By Proposition~1.3.6, $\Res{\cent HD^F}{\cento GD^F}\hat{t_D}$  is a sum of elements of $\ser{\cento GD^F}{s_D}$, a regular orbit under  $({\rm A}_{\cento
GD^*}(s_D)^F)^\wedge$. Each one is the canonical character of one of the blocks of $\cento GD^F$ covered by $B_D$. 
In other words, with notations of (1.3.5.2) $H_b=\tau_{H,s}({\rm A}_{\cento
GD^*}(s_D)^F)$ where ${\rm A}_{\cento
GD^*}(s_D)$ is viewed as a subgroup of
de ${\rm A}_{G^*}(s)$. 
\bull

\medskip\noindent{\bf 1.5.5. Proposition. }{\sl  Assume 1.5.1. There is a one-to-one map ${\cal B}_{G,s}$ from the set of unipotent $2$-blocks of $G(s)^F$ to the set of $2$-blocks in series $(s)$ of $G^F$ such that,  if $b\in{\bl{G(s)^F}1}$,

(1)  $b$ and  ${\cal B}_{G,s}(b)$ have a common defect group.

(2) There is a one-to-one map $\Psi_{b}$ from  $\II b$ onto $\II{{\cal B}_{G,s}(b)}$ that preserves height.

(3) The Brauer's category of $b$ and ${\cal B}_{G,s}(b)$ are isomorphic.
}

\medskip\noindent{\it On the proof. }

When ${\rm A}_{G^*}(s)^F=1$, as we have seen in the proof of Proposition~1.5.2, $b$ is the principal and unique unipotent $2$-block of $G(s)^F$,  $\II b={\cal E}_2(G(s)^F,1)$,  ${\cal B}_{G,s}(b)$ is the unique $2$-block in series $(s)$ of $G^F$, the map  $\Psi_{b}\colon \II b \to \II{{\cal B}_{G,s}(b)}$ is given using  Lusztig induction $(\la\mapsto \Lu {G(s)}G (\Psi_{G(s),s}(1)\otimes \la))$. The map $\Psi_{b}$ preserves height thanks to degree's formulas given in section~1.3.

Assuming ${\rm A}_{G^*}(s)^F\neq 1$, $G$ has type $\AA$. Consider now a regular embedding $(G,F)\subseteq (H,F)$ and $t$ semi-simple of odd order in $H^{*\, F}$ with image $s$ by a dual morphism. Using Clifford theory of \irr\ representations and of blocks and knowing that there is no multiplicity in restrictions of \irr\ representations from $H^F$ to $G^F$ and from $G(s)^F$ to $G(s)^{\circ\,F}$, the crucial step is to verify the following combinatorial facts (see the proofs of Propositions 3.4.1 and 3.4.2 for more detailed arguments): 

(a) If the unique $2$-block in series $(t)$ of $H^F$ covers exactly $m$ blocks of $G^F$ (therefore in series $(s)$), then the  unique unipotent $2$-block of $G(s)^{\circ\,F}$ is covered by exactly $m$ blocks of $G(s)^F$. Here $m$ is given in Proposition~1.5.4.

(b) Let $\chi\in{\cal E}_2(H^F,t)$, hence $\chi=\Lu {H(t)}H(\Psi_{H(t),t}(1)\otimes \la)$, where $\la\in{\cal E}_2(H(t)^F,1)$. Here we have $\la\in \ser{H(t)^{F}}{t_1}$, where $t_1\in \cent{H^*}t^F_{2}$. Let $s_1$ be the image of $t_1$ in $\cento{G^*}s^F_2$. There is a Levi subgroup $H(tt_1)\subseteq H(t)$ in the dual class of $\cent{T^*}{tt_1}$ and $G(ss_1)^{\circ\,F}:=H(tt_1)\cap G$ is a Levi subgroup of  $G(s)^\circ$  in the dual class of $\cento{G^*}{ss_1}$ (recall that $G$ has type $\AA$). Let $\al\in\ser{H(tt_1)^F}1$ such that $\la=\Psi_{H(t),t_1}(\al)$.  The inclusion $G(ss_1)^\circ\subseteq H(tt_1)$ is a regular embedding. We may identify $\ser{H(tt_1)^F}1$ and $\ser{G(ss_1)^{\circ\,F}}1$, so that $\la$ defines $\mu:=\Psi_{G(s)^\circ,s_1}(\al)\in\ser{G(s)^{\circ\,F}}{s_1}\subseteq \II {b(s)}$. 

We have to show that if $\chi$ covers exactly $n$ elements of ${\cal E}_2(G^F,s)$, then $\mu$ is covered by $n$ elements of $\II{G(s)^F}$. 

The set of \irr\ components of $\Res{H^F}{G^F}\chi$ is a regular orbit under $({\rm A}_G(ss_1)_\al^F)^\wedge$
and ${\rm A}_G(ss_1)$ is isomorphic to ${\rm A}_{G(s)^\circ}(s_1)$. Thus 
$n=|{\rm A}_G(ss_1)_\al^F|=|{\rm A}_{G(s)^\circ}(s_1)^F_\al| $ is the number of elements of $\II{G(s)^F}$ that cover $\mu$.

The proof of (3) is left to the reader. It is similar and simpler that the proof of property (B.4) in Theorem~1.4, see section 4 below. 
\bull

 \vfill\eject

 \noindent{\bf 2. Cuspidal data, Generalized $d$-Harish-Chandra theory and blocks}
 
 \bigskip In all this section are given $(G,F)$, $(G^*,F)$, $\ell > 2$, $d$ as in 1.1.2 and 1.1.5, see also Assumption~2.1.2.
 
 \bigskip\noindent{\bf 2.1. Facts on cuspidal data and blocks}
 
 \medskip\noindent{\bf 2.1.1. Definitions. }{\sl An element $\chi$ of $\II {G^F}$ is said to be {\bf $d$-cuspidal} when for any proper
$d$-split Levi subgroup $L$ of $G$ and any parabolic subgroup $P$ of $G$ admitting $L$ as a Levi complement, one has $\slu {L\subseteq P}G\chi=0$.  

A {\bf $d$-cuspidal datum in $(G,F)$} is a couple
$(L,\la)$, where
$L$ is a $d$-split Levi subgroup of $G$ and
  $\la$ is a $d$-cuspidal element of $\II{L^F}$. It is said {\bf in series (s)}, if $\Lu LG\zeta$ writes in Lusztig series $\ser{G^F}s$, $s$ a semi-simple element of $G^{* F}$.
  
  Let us denote in this hypothesis  $\ser{G^F}{(L,\la)}$ the set of $\chi\in\II{G^F}$ such that $\scal{\Lu {L\subseteq P}G\la}\chi{G^F}\neq 0$ for some $P$.}

Unipotent $d$-cuspidal data are described in [10]. Once $d$-split Levi subgroups are known, classification of $d$-cuspidal data in series $(s)$ reduces to the existence of $d$-cuspidal elements in $\ser{G^F}s$. Then one may rely $d$-cuspidality with Jordan decomposition. A first result is the following, to compare with Propositions~2.2.3, 2.1.4. 

\medskip\noindent{\bf 2.1.1.1. }{\sl Assume $\cent{G^*}s$ is connected and that commutation formula $\Lu {L\subseteq P}G\circ\Psi_{L,s}=\Psi_{G,s}\circ
\Lu{\cent{L^*}s}{\cent{G^*}s}$ holds for any $d$-split Levi subgroup in $(G,F)$. Then $\chi\in\ser{G^F}s$ is $d$-cuspidal if and only if  

$\zo{\cent{G^*}s}_{\phi_d}=\zo{G^*}_{\phi_d}$ and $\chi=\Psi_{G,s}(\al)$, where $\al$ is $d$-cuspidal in $\ser{\cent{G^*}s^F}1$. }

 \preuve By Propositions~1.3.2, 1.3.6, $\Psi_{G,s}$ is well defined. The commutation formula is (J3), introduced in 1.3.3. Let $\al\in\ser{\cent{G^*}s^F}1$.
  Let $L$ be any proper $d$-split Levi subgroup in $(G,F)$, and $L^*$ in the dual $G^{*\,F}$-conjugacy class such that $s\in L^*$. As $L^*\neq G^*$, $\zo{L^*}_{\phi_d}\neq \zo{G^*}_{\phi_d}$. 
 
 If $\zo{\cent{G^*}s}_{\phi_d}=\zo{G^*}_{\phi_d}$ and $\al$ is $d$-cuspidal, then $\cent{L^*}s$ is a proper $d$-split Levi subgroup of $\cent{G^*}s$ so that $\slu {\cent{G^*}s}{\cent{L^*}s}\al = 0$. By commutation formula (J3) $\slu GL \Psi_{G,s}(\al)=0$ : $\Psi_{G,s}(\al)$ is $d$-cuspidal.
 
 If $\Psi_{G,s}(\al)$ is $d$-cuspidal, by (J3) $\slu {\cent{G^*}s}{\cent{L^*}s}\al = 0$. That imply ${\cent{L^*}s}\neq \cent{G^*}s$ for any proper $d$-split Levi subgroup of $G^*$, hence $\cent{G^*}{\zo{\cent{G^*}s}_{\phi_d}}=G^*$, equivalently $\zo{\cent{G^*}s}_{\phi_d}=\zo{G^*}_{\phi_d}$. Then if $L^*_s$ is any proper $d$-split Levi subgroup of $\cent{G^*}s$, $L^*:=\cent{G^*}{\zo{L^*_s}_{\phi_d}}$ is a proper Levi subgroup of $G^*$ and one has 
 $\slu {\cent{G^*}s}{L^*_s}\al = 0$ by (J3): $\al$ is $d$-cuspidal.  \bull
 
 To use  [15], we have to make some assumption on $\ell$.

 
\medskip\noindent{\bf 2.1.2. Assumption on $(G,F,\ell,d)$. }{\sl $G$ is defined on $\F q$ by $F$, $\ell$ is odd, $d$ is the order of $q$ modulo $\ell$, ``{block}" means ``$\ell$-block". If some component of $G$ has exceptional type, or if some rational component of $G$ has type $\lexp 3\DD_4$, then $\ell\geq 5$. If some component of $G$ has type $\EE_8$, then $\ell\geq 7$.}

Assumption~2.1.2 on $(G,F,\ell,d)$ implies 2.1.2 on $(H,F,\ell,d)$ for any $F$-stable connected algebraic subgroup $H$ of $G$. Among consequences of Assumption~2.1.2, recall that $\ell$ is good for $G$ and that any $d$-cuspidal unipotent \irr\ representation of $G^F$ is the canonical \irr\ representation of a block of $G^F$ with central defect [16] Proposition~22.16. 

From properties of Lusztig induction with respect to isotypic morphisms and of Jordan \dec\ (see Proposition~1.3.2),  transitivity of Lusztig restriction ($\slu TL\circ\slu {L\subseteq P}G=\slu TG$) and the {\it uniform criterion}  on $d$-cuspidal unipotent data [10] 3.13  one deduces  the following equivalence.  
 
\medskip\noindent{\bf 2.1.3. Proposition. }[15] Proposition~1.10, (i). {\sl Let $(G,F)$ and $(G^*,F)$ be in duality and $s\in(G^*)^F_{\ell'}$ semi-simple, $\chi\in\ser{G^F}s$. Assume 2.1.2 on $(G,F,\ell,d)$. The following assertions (i) and (ii)  are equivalent

(i) For every $F$-stable maximal torus $T$ of $G$ such that $T_{\phi_d}\not\subseteq \z G$, one has $\slu TG\chi=0$.

(ii) $\zo{\cento{G^*}s}_{\phi_d}\subseteq \z{G^*}$ and $\chi$ corresponds by  Proposition 1.3.6 to a $\cent{G^*}s^F$-orbit of $d$-cuspidal unipotent \irr\ characters of $\cento{G^*}s^F$ by Lusztig's parametrization.
 }

\medskip\noindent{\bf 2.1.4. Proposition. }[15] Theorem~4.2.  {\sl Assume~2.1.2 on $(G,F,\ell,q)$. Let $s$ be a semi-simple $\ell'$-element of ${G^*}^F$. Any  $d$-cuspidal datum in series $(s)$ of $(G,F)$ may be defined as follows 

Let $(L^*_s,\al)$ be a unipotent $d$-cuspidal datum of $(\cento{G^*}s,F)$, let $L^*:=\cent{G^*}{\zo {L^*_s}_{\phi_d}}$, a $d$-split Levi subgroup of $G^*$, let $L$ be a Levi subgroup of $G$ in the dual $G^F$-conjugacy class of the $(G^*)^F$-class of $L^*$. To the orbit of $\al$ under ${\rm A}_{L^*}(s)^F$ there corresponds, by (1.3.6.1), an orbit $\Lambda$ under $({\rm A}_{L^*}(s)^F)^\wedge$ in $\ser{L^F}s$. For any $\la\in\Lambda$, $(L,\la)$ is a $d$-cuspidal datum of $(G,F)$.

If ${\rm A}_{G^*}(s)^F=\{1\}$, the map $(L^*_s,\la)\mapsto (L,\zeta)$ so defined induces a bijection from the set of $\cent{G^*}s^F$-conjugacy classes of unipotent $d$-cuspidal data in  $(\cent{G^*}s,F)$ to the set of $G^F$-conjugacy classes of  $d$-cuspidal data in series $(s)$ of $(G,F)$.}

  \medskip\noindent{\it Comments. }We comment here the last assertion of Proposition~2.1.4.
 
 When $A_{G^*}(s)^F=1$ in the definition above one may write $\la=\Psi_{L,s}(\al)$. In this case $(L,\la)$ is defined up to $G^F$-conjugacy : $s$ is defined up to ${G^*}^F$-conjugacy, $(L^*_s,\al)$ is defined up to $\cent{G^*}s^F$-conjugacy, hence $(L^*,s)$ is defined up to ${G^*}^F$-conjugacy, finally $(L,\ser{L^F}s)$ is defined up to $G^F$-conjugacy, so is $(L,\la)$ thanks to  (iv) of Proposition~1.3.2. 

Assume now that two unipotent $d$-cuspidal data $(L^*_{s,j},\al_j)$ ($j=1,2$) define $G^F$-conjugate $(L_2,\la_2)$, $(L_1,\la_1)$. Thus $(L_j,\ser{L^F}s)$ ($j=1,2$) are $G^F$-conjugate and, by our comments on (iv) of  Proposition~1.3.2, there exists some $g^*\in\cent{G^*}s^F$  inducing a dual morphism $L^*_1=\lexp{g^*}L^*_2$, so that $\la_2=\Psi_{L_2,s}(\al_2)={\la_1}^g=\Psi_{L_2,s}(\lexp{g^*}{\al_1})$, hence $\al_2=\lexp{g^*}{\al_1}$.

\medskip\noindent{\bf 2.1.5.  Proposition. }{\sl Assume 2.1.2 on $(G,F,\ell,d)$. 


(a) Let $\sigma\colon (G,F)\to (H,F)$ be an isotypic morphism.

Let $M$ be a $d$-split Levi subgroup of $H$, $\mu\in\ser{M^F}t$, $t\in(M^*)^F_{\ell'}\subseteq H^*$, $L=\sigma^{-1}(M)$ and $\la$ be an \irr\ component of $\Res{} {L^F\to M^F} \mu$. Then $(M,\mu)$ is a $d$-cuspidal datum in $(H,F)$ if and only if $(L,\la)$ is a $d$-cuspidal datum in $(G,F)$.

(b) Let $(G,F)=(G_1,F).(G_2,F)$ be a central product of connected reductive algebraic groups defined on $\F q$. Let $L=L_1.L_2$ be an $F$-stable Levi subgroup of $G$ where $L_j=L\cap  G_j$, $\la\in\II{L^F}$, $\la_j\in\II{L_j^F}$ such that $\la\in\II{L^F\mid \la_1\otimes\la_2}$. Then $(L,\la)$ is a $d$-cuspidal datum in $(G,F)$ if and only if, for $j=1,2$, $(L_j,\la_j)$ is a $d$-cuspidal datum in $(G_j,F)$.

(c) Assume $G=\a G$ (see \S~1.5.1). There is only one $G^F$-conjugacy class of $d$-cuspidal unipotent data in $(G,F)$, such as $(T,1_{T^F})$ where $T$ is a diagonal torus.}

\medskip\noindent {\it On proofs. } In an isotypic morphism $\sigma\colon G\to H$, the sets of Levi subgroups $L$ of $G$ and $M$ of $H$ correspond bijectively by $L\mapsto M=\sigma(L).\zo H$ and $M\mapsto L=\sigma^{-1}(M)$. Then $L$ is $d$-split if and only if $M$ is $d$-split. To verify (a) we may assume $L=H$ and $M=G$. If $\mu$ is unipotent, then $\Res{H^F}{G^F}\mu\in\ser{G^F}1$. By (1.3.1.1) $d$-cuspidality of $\mu$ is equivalent to $d$-cuspidality of $\Res{H^F}{G^F}\mu$. If $\mu\in\ser{H^F}t$ and $\sigma^*(t)=s$, $\la\in\ser{G^F}s$ and there is some $\al\in\ser{\cento{H^*}t^F}1$ that corresponds to $\la$ and $\mu$ (once $\ser{\cento{H^*}t^F}1$ and $\ser{\cento{G^*}s^F}1$ are identified, see Proposition~1.3.7). By Proposition~2.1.4, $\mu$ is $d$-cuspidal if and only if $\al$ is $d$-cuspidal, if and only if $\la$ is $d$-cuspidal.

(b) is clear if the product is direct. The central quotient morphism is isotypic.

(c) If $G=\a G$, the $d$-split Levi subgroups of $G$ are the diagonal Levi subgroups. For any diagonal Levi subgroups $M\subseteq L$ and  $\chi\in\ser{L^F}1$,   $\slu ML\chi\neq 0$. If $(L,\al)$ is a unipotent $d$-cuspidal datum, $L$ is a diagonal torus, therefore $\al=1_{L^F}$. 
\bull

Our classification of blocks of $G^F$ is given by the two following propositions (see also [16] 21.7, Chapter~4,
Exercise 4).

\medskip\noindent{\bf 2.1.6. Proposition. }[15] Theorem~2.5. 
{\sl  Assume 2.1.2 on $(G,F,\ell,d)$. Let $K$ be an $E$-split Levi subgroup of $G$. Let $b_K$ be a block of $K^F$ in series $(s)$, where $s\in (K^*)^F_{\ell'}\subseteq G^*$. There exists an block $B$ of
$G^F$ in series $(s)$, we denote $\Lu KG b_K$, such that, for all $\xi\in\II{b_K}\cap
\ser{K^F}{s}$, any \irr\ component of  $\Lu {K\subseteq P}G\xi$
belongs to $\II B$, whatever be the parabolic $P$ with Levi complement $K$. Furthermore if $K^F=\cent{G^F}{\z K^F_\ell}$, then
$(\{1\},\Lu KG b_K)\subseteq (\z
K^F_\ell,b_K)$, an inclusion of Brauer subpairs in $G^F$.}

\medskip\noindent{\bf 2.1.7. Proposition. }[15] Theorem~4.1.
{\sl 
 Assume 2.1.2 on $(G,F,\ell,d)$. Let $s$ be a semi-simple element of $(G^*)^F_{\ell'}$ and let $(L,\la)$ be a $d$-cuspidal datum  in series $(s)$ in $(G,F)$. Then is defined an block $b_{G^F}(L,\la)$ of
$G^F$ with  the following  properties :

 (a) $\la\in \II{b_{L^F}(L,\la)}$.
 
 (b) If $M$ is a $d$-split Levi subgroup of $G$ and
$(L,\la)$  is a $d$-cuspidal datum in  $(M,F)$, then, with notations of Proposition~2.1.6, $\Lu MG(b_{M^F}(L,\la))=b_{G^F}(L,\la)$.
  
  The map $(L,\la)\mapsto b_{G^F}(L,\la)$ is one-to-one from the  set of $G^F$-conjugacy classes of  $d$-cuspidal data in $\ell'$-series in $(G,F)$ onto the set of blocks of $G^F$.}
  
\smallskip\noindent{\it Some comments on Proposition~2.1.7. }

(a) For bad primes, see [23] who suggests that blocks correspond to $d$-cuspidal data $(L,\la)$ where
$\la$ has central defect group, that is not always the case for all $d$-cuspidal data if $\ell$ is bad for $G$.

(b)  The group $(G^F/[G,G]^F)^\wedge$ acts on the set of blocks of $G^F$ ``by tensor product"
 that is by the equality (see Appendix, 5.1)
$$\theta\otimes b_{G^F}(L,\la)=b_{G^F}(L,(\Res{G^F}{L^F}\theta)\otimes
\la),\quad \theta\in(G^F/[G,G]^F)^\wedge$$
thanks to the equality  $\Lu LG((\Res{G^F}{L^F}\theta)\otimes
\la)=\Lu LG\theta\otimes\la $ [20] 12.6. In section 2.4 we compute the stabilizer of $b_{G^F}(L,\la)$ in terms of the unipotent $d$-cuspidal datum corresponding to
$(L,\la)$  by Proposition~2.1.4. 

(c) More generally in a central product $G=G_1.G_2$, by Proposition~2.1.5 a $d$-cuspidal datum $(L,\la)$ in $(G,F)$ covers a product of $d$-cuspidal data $(L_1\times L_2,\la_1\otimes\la_2)$ in $(G_1\times G_2,F)$. Furthermore $G^F/G_1^F.G_2^F$ is  isomorphic to $L^F/L_1^F.L_2^F$ by Lemma~1.2.1. By (1.3.1.1) applied to $\sigma\colon G_1\times G_2\to G$, $(G^F/G_1^F.G_2^F)^\wedge$ acts on the set of blocks of $G^F$ defined by $d$-cuspidal data with support $L_1.L_2$ as $(L^F/L_1^F.L_2^F)^\wedge$ acts on $d$-cuspidal elements  in $\II{L^F}$.

\medskip\noindent{\bf 2.1.8.  Corollary. }{\sl  Let $(G,F,\ell,d,s,L,\la)$ as in Proposition~2.1.7 and $\mu\in\II{L^F}$. Then $\Lu {L\subseteq P}G\la$ and $\Lu {L\subseteq P}G\mu$ have a common \irr\ component if and only if $\la$ and $\mu$ are conjugate under $\nor GL^F$.}

\preuve If $\la$ and $\mu$ are $\nor LG^F$-conjugate, then $\Lu {L\subseteq P}G\la=\Lu {L\subseteq P}G\mu$.

Assume that $\Lu {L\subseteq P}G\mu$ and $\Lu {L\subseteq P}G\la$ have a common component  $\chi\in \ser {G^F}s$. There exist some $d$-cuspidal datum $(L_0,\la_0)$ in $(L,F)$ such that $\mu\in\II{b_{L^F}(L_0,\la_0)}$. Then $\chi\in\II{b_{G^F}(L_0,\la_0)}\cap\II{b_{G^F}(L,\la)}$ by Proposition~2.1.6. Hence  
$b_{G^F}(L,\la)=b_{G^F}(L_0,\la_0)$, $(L,\la)$ and $(L_0,\la_0)$ are 
$G^F$-conjugate and $\mu=\la_0$.
\bull

In the following Proposition we precise the induction on blocks (Proposition~2.1.6) with respect to the parametrization by conjugacy classes of $d$-cuspidal data (Proposition~2.1.7).

\medskip\noindent{\bf 2.1.9. Proposition}{ \sl Assume 2.1.2 on $(G,F,\ell,d)$.  Let $K$ be an $E$-split Levi subgroup of $G$ and $(L_K,\la_K)$ be a $d$-cuspidal datum in $(K,F)$ in an $\ell'$-series, defining the block $b_{K^F}(L_K,\la_K)$ (Proposition~2.1.7). Let $L:=\cent G{\zo {L_K}_{\phi_d}}$. Then $\Lu{L_K}L\la_K$ is $d$-cuspidal and one has, with notations of Proposition~2.1.6,
$$\Lu KG(b_{K^F}(L_K,\la_K))=b_{G^F}(L,\Lu{L_K}L\la_K)$$}

\preuve The important fact is that $\Lu{L_K}L\la_K$ is $d$-cuspidal and we first prove it. By definition $L$ is the smallest $d$-split Levi subgroup of $G$ such that $L_K\subseteq L$.

If $G$ is a central product over $F$, $G=G_1.G_2$, then $K_i:=K\cap G_i$ is $E$-split in $G_i$ and   $K=K_1.K_2$. Let  $L_{K,i}=L_K\cap G_i$, hence  $L_K=L_{K,1}.L_{K,2}$, $L_{K,i}$ is $d$-split in $K_i$. By (b) in Proposition~2.1.5,  $\Res{L_K^F}{L_{K,1}^F.L_{K,2}^F}\la_K$ writes $\la_{K,1}\otimes \la_{K,2}$ where each $\la_{K,i}$ is $d$-cuspidal. We have $L_i:=L\cap G_i=\cent{G_i}{\zo{L_{K,i}}_{\phi_d}}$. If $\Lu{L_{K,i}}{L_i}(\la_{K,i})$ is $d$ cuspidal for $i=1,2$,  as $\Res{L^F}{L_1^K.L_2^F}(\Lu {L_K}L\la_K)=\Lu{L_{K,1}}{L_1}(\la_{K,1})\otimes \Lu{L_{K,2}}{L_2}(\la_{K,2})$, by Proposition 2.1.5 again, $\Lu {L_K}L\la_K$ is $d$-cuspidal. 

So we may assume $G$ rationally \irr.

(a) In a first step we assume that $K$ is a maximal proper $E$-split Levi subgroup of  $G$.

If $K$ is $d$-split, $L_K$ is $d$-split and $L=L_K$, we are done.  
If $G=\zo G.\b G$, then $K$ is $d$-split by 1.1.5.3.
 
Assume now that $G\neq \zo G.\b G$ and $K$ is not $d$-split. Let $\AA_n((\epsilon q)^{da})$ be the rational type of $G$. Then  $K$ has type $\AA_{n'}((\epsilon q)^{da})\times \AA_m((\epsilon q)^{da\ell})$, where $(n+1)= (m+1)\ell+n'+1$ and the type of $L_K$ has  the form $\big(\times_j\AA_{n'_j}((\epsilon q)^{da})\big)\times \big (\times_i\AA_{m_i}((\epsilon q)^{da\ell}))\big)$.  

Then the rational type of $L$ is $\big(\times_j\AA_{n'_j}((\epsilon q)^{da})\big)\times \big (\times_i\AA_{\ell m_i}((\epsilon q)^{da}))\big)$ and $L^*$ differs from $(L_K)^*$ only on the right side of that product.

We see that $L$ and $L_K$  have a common Coxeter torus. With a coherent choice of dualities $L^*$ and 
$(L_K)^*$  have a common Coxeter torus.  

If $\mu_K\in\ser{L_K^F}s$, as $\mu_{K}$ is $d$-cuspidal, $\cento{(L_K)^*}{s}$ is a Coxeter torus $T^*$ of $(L_K)^*$. By definition of $L^*$, $L^*=\cent{G^*}{\zo{(L_K)^*}_{\phi_d}}$ and  we have $\cent{L^*}s\cap(L_K)^*=\cent{(L_K)^*}s=T^*$. For each value of the index $i$ in the decomposition above, the components of rational types $\AA_{ m_i}((\epsilon q)^{da\ell})$ and $\AA_{\ell m_i}((\epsilon q)^{da})$ of $L^*$ and $(L_K)^*$ respectively are represented  on the same vector space $V_i$ of dimension $da\ell(m_i+1)$. As $\cent{(L_K)^*}s$ is a Coxeter torus, product of Coxeter torii of each component,  the semisimple element $s$ has $da\ell(m_i+1)$ distinct eigenvalues with multiplicity one on $V_i$. Hence the $i$-component of $\cent{L^*}s$ is a maximal torus of the $i$-component of $L^*$, that is $\cent{L^*}s=T^*$, a Coxeter torus of $L^*$.  If $T$ is a Coxeter torus of $L$, it is in the dual class  of $T^*$. In case $\z G$ is connected, we have $\la_{K}=\Lu {T}{L_K}(\Psi_{T,s}(1))$, hence $\Lu{L_K}{L}\la_K=\Lu TL (\Psi_{T,s}(1))$; it is  a cuspidal element in $\ser{L^F}t$. In general case $\la_K$ corresponds to $1_{T{^* F}}$ as well as $\Lu{L_K}{L}\la_K$ by Propositions 2.1.4 and 1.3.6, and  we have our claim.

(b) We use induction on the semi-simple rank of $G$.

If $K$ is not proper maximal as $E$-split Levi subgroup of $G$, let $K_1$ be a maximal proper $E$-split Levi subgroup of $G$ such that $K\subset K_1$. Denote $L_1:=\cent{K_1}{\zo {L_K}_{\phi_d}}$, $L_1$ is an $E$-split Levi subgroup of $G$. We have $\zo {L_1}_{\phi_d}=\zo L_{\phi_d}=\zo {L_K}_{\phi_d}$ hence $L=\cent G {\zo {L_1}_{\phi_d}} $. By induction hypothesis we know that
$\Lu{K}{L_1}\la_K$ is $d$-cuspidal. By transitivity of Lusztig induction and (a) $\Lu  {L_1}L(\Lu{K}{L_1}\la_K)=\Lu KL\la_K$ is $d$-cuspidal. 

(c) On blocks :

By Proposition~2.1.7 any \irr\ component $\xi$ of $\Lu{L_K}K\la_K$ belongs to $\II{b_{K^F}(L_K,\la_K)}$. Then by Proposition~2.1.6 any \irr\ component of $\Lu KG\xi$ belongs to $\II{\Lu KG(b_{K^F}(L_K,\la_K))}$. We see that any \irr\ component of $\Lu {L_K}G\la_K$ belongs to $\II{\Lu KG(b_{K^F}(L_K,\la_K))}$. But $\Lu {L_K}G\la_K=\Lu {L}G(\Lu{L_K}L\la_K)$ writes in $\Z \II{b_{G^F}(L,\Lu {L_K}L\la_k)}$ hence $\Lu KG(b_{K^F}(L_K,\la_K))=b_{G^F}(L,\Lu {L_K}L\la_K)$.
\bull

Propositions 2.1.4, 2.1.6, 2.1.7 suggest to compare blocks in series $(s)$ of $G^F$ and unipotent blocks of $\cent{G^*}s^F$ using $d$-cuspidal data. By Jordan decomposition, the relation is quite simple if the center of $G$ is connected :

\medskip\noindent{\bf 2.1.10. Proposition. }{\sl Assume 2.1.2 on $(G,F,\ell,q)$, $\z G$ connected and $s$ semi-simple in $(G^*)^F_{\ell'}$. By Proposition~2.1.4 there is a one-to-one map from the set 
 of $G^F$-conjugacy  classes of
$d$-cuspidal data in series $(s)$ in
$(G,F)$  and the set 
 of $\cent{G^*}s^F$-conjugacy  classes of
 unipotent
$d$-cuspidal data in  $(\cent{G^*}s,F)$. 

As a consequence, by Proposition~2.1.7, there is a one-to-one map from the set 
 of  blocks $B$ of $G^F$ such that
$\II B\subseteq \lser{G^F}s$ onto the set of unipotent blocks of $\cent{G^*}s^F$.}

\preuve In Proposition~2.1.4 the group
$L$ is defined up to $G^F$-conjugacy. The  $G^F$-conjugacy class of
$(L,\ser{L^F}s)$ depends on the $\cent{G^*}s^F$-conjugacy class of
$L^*_s$. 

Assume now that two couples $(L_1,\ser{L_1^F}s)$ et $(L_2,\ser{L_2^F}s)$ are conjugate by some $g\in G^F$. Then there exists $h\in{G^*}^F$ such that $(L_1^*,s)=(L^*_2,s)^h$. Indeed dualities are defined around couple of torii in duality
$(T_j\subset L_j,T^*_j\subset L^*_j)$ ($j=1,2)$). Up to $L^F_1$-conjugacy, one may assume $\lexp g{T_1}=T_2$. Then there is $h\in(G^*)^F$ such that $\lexp hT^*_1=T^*_2$. By $(g,h)$-conjugacy the root datum of $(L_1,L_1^*)$ with respect to $(T_1,T^*_1)$ is sent on the root datum of $(L_2,L^*_2)$ with respect to $(T_2,T^*_2)$. As $g$ sends $\ser{L_1^F}s$ on $\ser{L_2^F}s$, up to $(L^*_2)^F$-conjugacy we may assume that $h\in\cent{G^*}s^F$. Therefore the $d$-split Levi subgroups $L^*_{1,s}:=\cent{L^*_1}s$ and $L^*_{2,s}:=\cent{L^*_2}s$ are $\cent{G^*}s^F$-conjugate. 

Denote
${\cal S}(L,s)$ the set of $d$-cuspidal
$\la\in\ser{L^F}s$.

Once $s$ and the groups $L$, $L^*$ are given, we are reduced to consider the set of orbits under $\nor{G}L^F$ on ${\cal S}(L,s)$. On the other side is the set of orbits under  $\nor {\cento{G^*}s}{L^*_s}^F)$ on ${\cal S}(L^*_s,1)$. One knows that
 $\nor
GL^F$ acts on Lusztig series in $L^F$ as ${\rm W}_G(L)^F={\rm W}_{G^*}(L^*)^F$ acts on classes of semi-simple elements of  $(L^*)^F$. We obtain a one-to-one map $$\big(\cup_{n\in \nor GL^F}\lexp
n{{\cal S}(L,s)}\big)/\nor GL^F\mapright\cong
{\cal S}(L,s)/\nor{G}{L,\ser{L^F}s}^F$$  

As $\z L$ and $\z G$ are connected $\Psi_{G,s}$ et $\Psi_{L,s}$
satisfy to Proposition~1.3.2, specially assertion (iv). The group
$\nor{G^F}{L,\ser{L^F}s}/L^F$ is isomorphic to the relative Weyl group 
${\rm W}_{\cent{G^*}s}({L^*_s})^F$, isomorphic to $
\nor{{\rm W}_G(s)^F}{{\rm W}_L(s)}/{\rm W}_L(s)^F
$ (see 2.2.5 below). It  acts on the root datum of
$L$ with respect to $T$ and on the root datum of $L^*_s:=\cent{L^*}s$ with respect to $T^*$. {\it Via}
these isomorphisms and Jordan decomposition, the
actions on
$\ser{L^F}s$ and on $\ser{{L^*_s}^F}1$ are exchanged. Let $\la\in {\cal S}(L,s)$, $\la=\Psi_{L,s}(\al)$, $\al\in {\cal S}(L^*_s,1)$ (Proposition~2.1.4). We have one-to-one maps on quotients  $$\eqalign{\nor{G^F}{L,\ser{L^F}s}/\nor{G^F}
{L}_{\la} &\cong
{\rm W}_{\cent{G^*}s}(L^*_s)^F/
{\rm W}_{\cent{G^*}s}({L^*_s)^F_{\al}}\cr
&\cong\nor{\cent{G^*}s}
{L^*_s}^F /\nor{\cent{G^*}s}
{L^*_s}_{\al}^F.\cr}$$

  Let
${\cal U}={\cal S}(L^*_s,1)$ be the set of  $d$-cuspidal in $\ser{(L^*_s)^F}1$.  The one-to-one map ${\cal S}(L,s)\longleftrightarrow{\cal U}$ induces a one-to-one map ${\cal S}(L,s)/\nor{G^F}{L,\ser{L^F}s}\longleftrightarrow{\cal U}/\nor{\cent{G^*}s}
{L^*_s}^F$.  \bull

\medskip\noindent{\bf 2.1.11. Proposition. On central defect. }{\sl Assume 2.1.2 on $(G,F,\ell,d)$.

(a) Let
$\sigma \colon (G,F)\to (H,F)$ be an  isotypic morphism defined on $\F q$.  Let
$\chi\in\II{H^F}$ with central $\ell$-defect. Any \irr\ component of  $\Res{H^F}{G^F}\chi$
has central $\ell$-defect.

 (b) Let $G=\a G.\b G$ the decomposition defined in 1.1.5.2. Let $\chi\in\II{G^F\mid{\a\chi\otimes\b \chi}}$, where $\a \chi\in \II{\a G^F}
$ and $\b \chi\in \II{\b G^F}
$. Then $\chi$ has central $\ell$-defect if and only if $\a\chi$ and $\b \chi$ have central $\ell$-defect.

(c) Let
$(L,\la)$ be a $d$-cuspidal datum in series $(s)$. 

(i) Assume $G=\b G$. The block $b_{G^F}(L,\la)$ has central defect group if and only if
$\z{\cento {G^*}s}^F_\ell=\z{G^*}^F_\ell$ and $L=G$. Then $\la$ is the canonical character of 
$b_{G^F}(L,\la)$. 

(ii) Assume $G=\a G$. The block $b_{G^F}(L,\la)$ has central defect group if and only if
$\cento {G^*}s$ is a torus with a dual $T$ in $G$ such that  $T^F_\ell= \z{G^F}_\ell$. Let $\theta=\Psi_{T,s}(1_{T^{* F}})\in(T^F)^\wedge$.
Then $T$ is a Coxeter torus of $(G,F)$ and $L=T$. The canonical character of $b_{G^F}(T,\la)$  is a component of $\Lu TG(\theta)$.}

\preuve 
(a) is a consequence of non-multiplicity  in $\Res{}\sigma$ (Proposition 1.3.4). Let
$\xi\in\II{\sigma(G^F)}\subseteq
\II{G^F}$ such that $\chi\in\II{H^F\mid \xi}$ and let $X$ be the stabilizer of $\xi$ in $H^F$. One has  $\chi(1)=|H^F/X|.\xi(1)$ and the assumption on $\chi$ writes $\chi(1)_\ell=|H^F/\z {H^F}|_\ell$. As $\sigma(G^F).\z {H^F}\subseteq X$ and by isomorphism theorem $\sigma(G^F).\z{H^F}/\z{H^F}\cong  \sigma(G^F)/\sigma(G^F)\cap \z{H^F}$, we have $$\xi(1)_\ell=|X/\sigma(G^F).\z{H^F}|_\ell.|\sigma(G^F)/\sigma(G^F)\cap \z{H^F}|_\ell\leqno{(2.1.11.1)}$$
But $\z{H^F}=\z H^F$, hence $ \sigma(G^F)\cap \z{H^F}\subseteq \z{\sigma(G)}$. Moreover $\sigma(\z G)=\z{\sigma(G)}$ and, since the kernel $K$ of $\sigma$ is contained in $\z G$, $\sigma(\z G)\cap\sigma(G^F)=\sigma(\z G\cap G^F.K)=\sigma(\z G^F)=\sigma(\z{G^F})$  so that $\sigma(G^F)\cap \z{H^F}=\sigma(\z G^F)$  and by report in (2.1.11.1)
$$\xi(1)_\ell=|X/\sigma(G^F).\z{H^F}|_\ell.|\sigma(G^F)/\sigma(\z G^F)|_\ell\leqno{(2.1.11.2)}$$

As the restriction of $\sigma$ to $G^F$ has kernel $K^F\subseteq \z{G^F}$, we have an isomorphism 
$$G^F/\z{G^F}\cong \sigma(G^F)/\sigma(\z{G^F})$$

As $\xi(1)$ divides $|G^F/\z{G^F}|$, (2.1.11.2) shows that $\xi$ has central defect and the order of the group $X/\sigma(H^F).\z{H^F}$ is prime to $\ell$.
Indeed the order of $s$ is prime to $\ell$. 

(b) The morphism $\a G\times\b G\to G$ is isotypic, hence (a) gives half of the equivalence.

 As $\a G\times \b G\to G$ is onto, $|G^F|=|\a G^F|.|\b G^F]$. By definition of $\b G$, $\z{\b G}^F$ has order prime to $\ell$, hence the kernel and cokernel of $\a G^F\times \b G^F\to G^F$ are $\ell'$-groups, by 1.1.5.2, Proposition~1.1.5.3 and Lemma~1.2.1. It follows that $\z G^F_\ell$ is isomorphic to $\z{\a G}^F_\ell$, and, by Clifford theory, that  $\chi(1)_\ell=\a \chi(1)_\ell.\b\chi(1)_\ell$. Thus $\chi(1)_\ell.|\z {G^F}|_\ell=|G^F|_\ell$ is equivalent to ``  $\a\chi(1)_\ell.|\z {{\a G}^F}|_\ell=|{\a G}^F|_\ell$ and $\b \chi(1)_\ell.|\z {{\b G}^F}|_\ell=|{\b G}^F|_\ell$ ".

(c) By [16], Proposition~22.16,
if $\al\in\ser{G^F}1$ is $d$-cuspidal,
$b_{G^F}(\al)$ has  central defect. Reciprocally if  $(L,\al)$ is a
$d$-cuspidal unipotent datum, we have $L=\cento G{\z L^F_\ell}$ and $\z L^F_\ell$ is
contained in a defect group of $b_{G^F}(L,\al)$ (Proposition~2.1.6).
Hence $b_{G^F}(L,\al)$ has central defect group if and only if $L=G$.  When $G=\a G$ and $G$ is not a torus there is no unipotent block with central defect group (Proposition~2.1.5).

When $s=s_{\ell'}$, $\z G^F_\ell$ is contained in the kernel of any $\chi\in\ser{G^F}s$ 
 because the caracteristic function of $\z
G^F_\ell$ is uniform [DM], 12.21. If $\chi$ is in a block of $G^F$ with central defect group, $\chi$ is the canonical character of its block. So  let
$\al\in\ser{\cento{G^*}s^F}1$ in the orbit under ${\rm A}_{G^*}(s)^F$
associated to
$\chi$ by Proposition~1.3.6. We know that ${\rm A}_{G^*}(s)^F$ is an
$\ell'$-group, so that by (vi) in Proposition~1.3.2, Propositions~1.3.4, 1.3.6, 
$$\chi(1)_\ell=\al(1)_\ell|(G^*)^F/\cento{G^*}s^F|_\ell$$
  But
$|(G^*)^F|=|G^F|$ hence $b_{G^F}(\chi)$ has central defect  group if and only if  $$|\z{G}^F|_\ell.\al(1)_\ell=|\cento{G^*}s^F|_\ell$$

(i) Assume $\a G$ is a torus. We have $\z
G^F_\ell=\zo G^F_\ell$, as well for $G^*$, hence $|\z
G^F|_\ell=|\z{G^*}^F|_\ell$. Clearly $\z {G^*}\subseteq \z{\cento{G^*}s}$. We see that $b_{G^F}(\chi)$ has central defect  only if 
$b_{\cento{G^*}s^F}(\al)$ has central defect and $\z {G^*}^F_\ell=
\z{\cento{G^*}s}^F_\ell$. 
By [13], Proposition 2.2, as $L^*$ is $d$-split in $G=\b G$, one has $L^*=\cento{G^*}{\z {L^*}^F_\ell}$ (In fact (ii) in [13], 2.2, is true under the hypothesis ``$\ell\in\gamma(G^*,F)$", that is satisfied if $G=\b G$). But
$\z{L^*}^F_\ell\subseteq\z{\cento{G^*}s}^F_\ell=\z{G^*}^F_\ell$. It follows that $L^*=G^*$,
$L=G$,  $\al$ is $d$-cuspidal.

ii) Assume  $G=\a G$, non abelian. Let $T$ be a maximal $F$-stable torus in $L$ in duality with a diagonal
torus $T_s^*$ of  $(\cento{G^*}s,F)$. A defect group of $b_{G^F}(L,\la)$ is extension of
$T^F_\ell$ by an $\ell$-Sylow subgroup of ${\rm W}_G(s)^F$. It is central in $G^F$ if and only if
$T^F_\ell\subseteq
\z G$ and
${\rm W}_G(s)^F_\ell=\{1\}$. Then $T$ is a 
Coxeter torus of $(G,F)$. Thus a diagonal torus $T^*_s$ of $(\cento{G^*}s,F)$ is a Coxeter torus of  $(G^*,F)$. As a Coxeter torus it is a maximal $F$-stable proper Levi subgroup in $G^*$, hence $\cento{G^*}s=T^*_s$ and $L=T$. In a regular embedding $G\to H$, with dual $H^*\to G^*$, $T$ is send in a Coxeter torus $S$ of $H$, $T^*$ is image of a Coxeter torus $S^*$ of $H^*$, $s$ is  image of $t\in S^{*\,F}$ such that $\cent{H^*}t=S^*$. We see that $ \Lu SH(\Psi_{H,t}(1_{S^{* F}}))$ is the unique element of $\ser{H^F}t$. We have  $\Lu TG\theta=\Res{H^F}{G^F}(\Lu SH(\Psi_{H,t}(1_{S^{* F}}))$ by (1.3.1.1). The Lusztig series $\ser{G^F}s$ is the set of \irr\ components of $\Lu TG \theta $. Such a component has central defect if and only if $T^F_\ell\subseteq \z G$. \bull

We note that (b) cannot be generalized to any isotypic morphism. Consider $G={\rm SL}_\ell\subseteq H:={\rm GL}_\ell$ on $\F q$ and assume that $(q-1)_\ell=\ell$. The two  groups $G^F$ and $H^F$ have isomorphic $\ell$-centers. If $T$ is a coxeter torus of $H$, $|T^F|_\ell=\ell^2$ and $|(T\cap {\rm SL}_\ell)^F|_\ell=\ell$. By  (c) (ii) in Proposition~2.1.11, with $T^*=\cent{H^*}s$, $\chi:=\Psi_{H,s}(1_{T^{* F}})$ has not central defect, but the components of $\Res{H^F}{G^F}\chi$ have central defect.

\bigskip\noindent{\bf 2.2. Generalized Harish-Chandra theory and blocks}

\medskip\noindent{\bf 2.2.1. Definition. }{\sl Let $d$ be an integer and $s$ a semi-simple element of $G^{*\,F}$. We say that  {\bf Generalized $d$-Harish-Chandra's theory holds in $\ser{G^F}s$} (or shortly {\bf G.$d$-HC holds in $\ser{G^F}s$}) if and only if :

 for any $\chi\in\ser{G^F}s$ there exists a $d$-cuspidal datum $(L,\la)$ in $(G,F)$, uniquely defined up to $G^F$-conjugacy, and $a\neq 0$, such that 
 $$\slu {L\subseteq P}G\chi=a.\big(\sum_{g\in \nor GL^F/\nor GL^F_\la}\la^g\big),\leqno{(2.2.1.1)}$$ independantly of the choice of the parabolic subgroup $P$. 
 
 Then  $\ser{G^F}s=\cup_{\{(L,\la)\}/G^F}\ser{G^F}{(L,\la)}$ and it is a partition.
  }
  
  Note that if $\scal \chi{\Lu LG\la}{G^F}\neq 0$ and $\chi\in\ser{G^F}s$, then there is a $d$-split Levi subgroup $L^*$ in $G^*$, in the dual conjugacy class of $L$,  such that $s\in L^*$ and $\la\in\ser{L^F}s$. That is why the property may be considered inside each rational series.

 If $d=1$, Generalized $d$-Harish-Chandra's theory is just classical Harish-Chandra's theory and holds in any type, any series [20] Chapter~6.
 
 If $G=\a G$, $d$-split Levi subgroups are diagonal Levi subgroups. That imply that G. $d$-HC theory reduces in non twisted type to classical Harish-Chandra's theory. Then a decisive fact is that if $\chi\in\II{G^F}$ and $L$ is a split Levi subgroup of $G$, then $\slu LG$ is an effective representation. In twisted type $\lexp 2\AA$ with $G=\a G$ in view of {\it  Ennola-duality} we obtain a similar property at least when $\z G$ is connected.  Then there exist a function $$\epsilon^G\colon  \II{G^F}\to \{-1,1\}\leqno{(2.2.1.2)}$$ such that for any diagonal Levi subgroup $L$ and any $(\la,\chi)\in\II{L^F}\times\II{G^F}$ one has (see [20] 15.4 and reference)$$\epsilon^G(\chi)\epsilon^L(\la)\scal{\Lu LG\la}\chi{G^F}\geq 0\leqno{(2.2.1.3)}$$
 Furthermore if $\chi(1)=\chi'(1)$, then $\epsilon^G(\chi)=\epsilon^G(\chi')$. In non twisted type we define $\epsilon$ as a constant application with value $1$. We may conjecture that (2.2.1.3) holds for $G=\a G$ as well when $\z G$ is not connected. Indeed under a conjecture relying Gelfand-Graev characters an Lusztig induction, one may extend Jordan decomposition for any series and (J3) [3], [12]. In that case, as we verify it in our Appendix, section 5.4, G.$d$-HC theory  holds in  $\ser {\cent{G^*}s^F}1$ and transfers to $\ser{G^F}s$.

\smallskip 
If Proposition~2.1.7 applies then $\cup_{\{(L,\la)\}/G^F}\II{b_{G^F}(L,\la)}$ is known to be a partition, hence
 
 \medskip\noindent{\bf 2.2.2. Proposition. }{\sl  Assume 2.1.2 on $(G,F,\ell,d)$ and $s$ semi-simple in $(G^*)_{\ell'}^F$.   Generalized $d$-Harish-Chandra's theory holds in $\ser{G^F}s$, if and only if, for any $d$-cuspidal datum $(L,\la)$ in series $(s)$  
  in $(G,F)$ 
 $$\ser{G^F}{(L,\la)}=\II{b_{G^F}(L,\la)}\cap\ser{G^F}s .$$}

 From the important paper of M. Brou\'e, G. Malle J. Michel [10] it is known that Generalized $d$-Harish-Chandra theory  holds in $\ser{G^F}1$. The fundamental theorem 3.2 in [10] gives precise results on degrees, value of scalar products, as $a$ in (2.2.1.1), and relation with Lusztig's map. We have retain the partition in ``Generalized Harish-Chandra series" $\ser{G^F}{(L,\la)}$ and formula (2.2.1.1), a consequence of transitivity theorem [10] 3.11.
 
 When $\z G$ is connected  and commutation between Lusztig induction and Jordan \dec\ (J3) in 1.3.3 is true, as it is the case in  classical types (see [27], [24] Appendix), and  if Proposition~2.1.4 applies, {\it via} $\Psi_{G,s}$, G.$d$-HC holds in $\ser{G^F}s$ because it holds in $\ser{\cent{G^*}s^F}1$.

 Let us consider now exceptional  types. The proof of G.$d$-HC in unipotent series for exceptional types uses two properties : 
 
 (i) Mackey decomposition formula for restriction of  Lusztig induction [20] 11.13;  
 
 (ii) orthogonal projection on the space of uniform functions (1.3.1.3). 
 
From now on we enforce assumption 2.1.2 by Mackey decomposition formula. Mackey formula is proved for classical types and any $q$ and for exceptional types if $q>2$ [4].
 It gives us the norm of $\Lu {L\subseteq P}G \la$ for any $d$-cuspidal unipotent datum $(L,\la)$ in $(G,F)$, independantly of  $P$.  To prove G.$d$-HC in $\ser{G^F}s$ we need Mackey formula  for and inside $d$-split Levi subgroups of $G$ and of $\cent{G^*}s$.

\medskip\noindent{\bf 2.2.3. Assumption. }{\sl Assume 2.1.2 on $(G,F,\ell,d)$ and that Mackey decomposition formula for Lusztig induction holds inside closed $F$-stable subgroups of $G$ or $G^*$. }
 
 \medskip\noindent{\bf 2.2.4. Proposition. }{\sl Assume 2.2.3 on $(G,F,\ell,d)$ and $\z G$ is connected. Let $s$ be semi-simple in $(G^*)^F_{\ell'}$. Then Generalized $d$-Harish-Chandra theory, as defined in 2.2.1,  holds in $\ser{G^F}s$.
 
 Furthermore if $(L,\la)$ is a $d$-cuspidal datum in series $(s)$ in $(G,F)$ associated to a $d$-cuspidal unipotent datum $(L^*_s,\al)$   in $(\cent{G^*}s,F)$ by Proposition~2.1.4, there exist a one-to-one map  
 $$\tilde\Psi_{G,s}(L,\la)\colon \ser{\cent{G^*}s^F}{(L^*_s,\al)}\to \ser{G^F}{(L,\la)}$$ such that $$(\tilde\Psi_{G,s}(L,\la)(\beta))(1)={|G^F]_{p'}\over |\cent{G^*}s^F|_{p'}}\beta(1)\leqno{(2.2.4.1)}$$ and $$\Lu LG\la=\sum_{\beta\in\ser{\cent{G^*}s^F}{(L^*_s,\al)}}\scal{\beta}{\Lu {L^*_s}{\cent{G^*}s}\al}{\cent{G^*}s^F}\tilde\Psi_{G,s}(L,\la)(\beta)\leqno{(2.2.4.2)}$$}
 
 The proof of Proposition~2.2.4 is given in 2.2.8.

 Clearly a candidate for $\tilde\Psi_{G,s}(L,\la)$ is the restriction of $\Psi_{G,s}$ to $\ser{\cent{G^*}s^F}{(L^*_s,\al)}$, but there is a doubt on the definition of $\Psi_{G,s}$ in some exceptional cases, and we show only that 
 $$\piu^G\circ \tilde\Psi_{G,s}(L,\la)=\piu^G\circ\Psi_{G,s}(L,\la)$$
 and that imply (2.2.4.1) by  (vi) in Proposition~1.3.2.
 
 To transfer G.$d$-HC from unipotent series to any Lusztig series, we have to use some classical properties on "relative Weyl groups" as the following
 
 \medskip\noindent{\bf 2.2.5. $W$-argument. } 
{\sl Let $L$ be a  connected  reductive $F$-stable subgroup   of $G$,  and $H$ a subgroup of $G^F$ such that $L^F\subseteq H\subseteq  \nor G L$. Let $T$ be a maximally split torus of 
$L$. Put ${\rm W}(HL,T)=\nor {HL}T.T/T$, a subgroup of ${\rm W}(G,T)$. There is a split short exact sequence
$$1\to {\rm W}(L,T)^F\to  {\rm W}(HL,T)^F\to H/L^F\to 1$$ }
\noindent{\it Proof. }  
By $L^F$-conjugacy  of maximally split tori in $L$ we have 
$(HL)^F=H=\nor{HL}T^F.L^F$ hence by an isomorphism theorem  $H/L^F\cong \nor{HL}T^F/\nor{L}T^F$. 
 As $T$ is connected one has  ${\rm W}(HL,T)^F\cong  \nor{HL}T^F/T^F$ and ${\rm W}(L,T)^F\cong \nor LT^F/T^F$, thus our claim.
The extension is split : the  stabilizer in $H$ of an $F$-stable couple (torus
$\subseteq $ Borel)  in $L$ is a complement of $\nor{L}T^F$ in $\nor{HL}T^F$.
\bull

Note that Frattini's argument on maximally split tori inside $F$-stable Borel subgroups may be used in type $\lexp 2\AA$ for diagonal tori and Borel subgroups.

Anti-isomorphisms between $W(G,T)$ and $W(G^*,T^*)$ for groups in duality extends easily to some relative Weyl groups.

As a first example assume that $L$ is an $F$-stable Levi subgroup of $G$, with $(T\subseteq L\subseteq G,F)$ in duality with $(T^*\subseteq L^*\subseteq G^*,F)$.

Then $\nor GL/L\cong \nor{W(G,T)}{W(L,T)}/W(L,T)= W_G(L)\cong W_{G^*}(L^*)$,  $ W_G(L)^F\cong W_{G^*}(L^*)^F$ by 2.2.5. If
$s\in (T^*)^F$, put $\theta=\Psi_{T,s}(1_{T^{* F}})$ and let
$\nor{G^F}{L,\ser{L^F}s}$ be the stabilizer of
$\ser{L^F}s$ in $\nor{G}
 L^F$.  Then, as $(T,\theta)$ is defined by the series $(s)$ modulo $L^F$-conjugacy, $\nor{G^F}{L,\ser{L^F}s}=\nor{G^F}{L,T,\theta}.L^F$, so that 
$\nor{G^F}{L,\ser{L^F}s}/L^F$ is isomorphic to $
\nor{G^F}{T,\theta}/\nor{L^F}{T,\theta}$. One have an
isomorphism
$$\nor{G^F}{L,\ser{L^F}s}/L^F\cong
(\nor{G^{* F}}{L^*}\cap\cent{G^*}s)/\cent{L^{* F}}s\leqno{(2.2.5.1)}$$
We need a restriction of that isomorphism in case of cuspidal data. To simplify notations, as $T$ and $T^*$ are fixed we write $W(G)$ instead of $W(G,T)$ and so on. Given $s\in T^*\subseteq G^*$ and anti-isomorphism between $W(G,T)$ and $W(G^*,T^*)$ as above, we abbreviate the image of $W(\cento{G^*}s,T^*)$ and of $\nor{\cent{G^*}s}{T^*}/T^*$ in $W^\circ_G(s)$ and $W_G(s)$ respectively, two subgroups of $W(G,T)$, so that the last isomorphism becomes $$\nor{G^F}{L,\ser{L^F}s}/L^F\cong
\nor{W_G(s)}{W_L(s)}^F/W_L(s)^F\leqno{(2.2.5.2)}$$ With these notations we have

\medskip\noindent{\bf 2.2.6.  Proposition. }{\sl Let $(L,\la)$ be a $d$-cuspidal datum in series $(s)$ of an $F$-stable Levi subgroup $M$ of $G$. Assuming  Weyl groups are defined over dual $F$-stable torii $T$, $T^*$ with $ T^*$ maximally split in $\cento{L^*}s$, one has an
isomorphism 
$$\big[\nor G M\cap
\nor{G^F}{L,\ser{L^F}s}\big]/L^F\cong \big[
\nor{{\rm W}(G)}{{\rm W}(M)}\cap \nor{{\rm W}_G(s)^F}{W^\circ_L(s)}\big]/{\rm W}_L(s)^F.$$}
\preuve    From the hypotheses we use the equalities $L=M\cap \cent{G}{\zo{L}_{\phi_d}}$, $L^*=M^*\cap
\cent {G^*}{\zo {L^*}_{\phi_d}}$ and 
$\zo{L^*}_{\phi_d}=\zo{\cento{L^*}s}_{\phi_d}$ (Proposition~1.4.2).
Therefore $\nor {G}{M}^F\cap \nor {G}{L}=\nor {G}{M}^F\cap \nor
{G}{\zo {L}_{\phi_d}}$, $\nor {G^*}{M^*}^F\cap \nor {G^*}{L^*}=\nor {G^*}{M^*}^F\cap \nor
{G^*}{\zo {L^*}_{\phi_d}}$. As $s\in L^*$, $\nor{G^*}{L^*}^F\cap
\cent{G^*}s=\nor{\cent{G^*}s}{\zo{L^*}_{\phi_d}}^F$. Now
$\zo{L^*}_{\phi_d}=\zo{\cento{L^*}s}_{\phi_d}$ implies
$\nor{G^*}{M^*}\cap \nor{G^*}{L^*}\cap
\cent{G^*}s^F=\nor{G^*}{M^*}\cap\nor{\cent{G^*}s}{ \cento{L^*}s}^F$.
  Thus
 $$\nor {G^*}{M^*}\cap\nor{G^{* F}}{L^*,s}=(L^*)^F.[
\nor{G^*}{M^*}\cap\nor{\cent{G^*}s}{ \cento{L^*}s}]^F.$$
Using 2.2.5, (2.2.5.1) and (2.2.5.2)  that last equality gives our claim.
 \bull
 
 Thanks to commutation of $\piu$ with Lusztig induction (1.3.1.4) $\piu^G\circ \Lu {L\subseteq P}G$ is known and independant of $P$. In [10] p.55 we read --- translation in $G^F$ of a generic result and using our notations ---
 
 {\it It will turn out that there exists an essentially unique element $\gamma \in \Z \ser{G^F}1$ with $\piu^G(\gamma)=\piu^G(\Lu LG \al)$ of minimal norm, and that this norm coincides with the norm of $\Lu LG\al$ calculated from the Mackey formula.}

\medskip\noindent{\bf 2.2.7. Lemma. }{\sl Assume 2.1.2 on $(G,F,\ell, d)$ and $\z G$ is connected. Let $s$ be a semi-simple $\ell'$-element of $G^{*\,F}$. 
 Let $(L^*_s,\al)$ be a $d$-cuspidal unipotent datum in $(\cent{G^*}s,F)$ and let $(L,\la)$ be an associated $d$-cuspidal datum in series $(s)$ in $(G,F)$, hence $\la=\Psi_{L,s}(\al)$ (see Proposition 2.1.4). 
 
 (a) The relation
$$\Lu LG(\Psi_{L,s}(\al))=\Psi_{G,s}(\Lu {L^*_s}{\cent{G^*}s}\al)\leqno{(2.2.7.1)}$$
is true when one of the following conditions is satisfied : 

(i) $\cent{G^*}s$ is a Levi subgroup of $G^*$.

(ii)   $L^*_s$ has type $\AA$.

(iii) Mackey decomposition formula holds between $E$-split Levi subgroups in $G$ and in $\cent{G^*}s$ and 
 $\Lu{L^*_s}{\cent{G^*}s}\al$ is  the unique element $\xi$ of $\Z\ser{\cent{G^*}s^F}1$ of minimal square norm such that  $$\piu^{\cent{G^*}s}(\xi)=\Lu{L^*_s}{\cent{G^*}s}(\piu^{L^*_s}(\al))$$

  When  (iii) is satisfied $\Lu LG\la$ is  the unique element $\chi$ of $\Z\ser{G^F}s$ of minimal square norm such that  $\piu^{G}(\chi)=\Lu LG(\piu^L(\la))$.
   
 (b) If $\tilde\Psi_{G,s}(L,\la)$ exists as in Proposition~2.2.4 for any $(L,\la)$ defined from a proper $d$-cuspidal unipotent datum $(L^*_s,\al)$ in $(\cent{G^*}s,F)$, one has the partition
$$\ser{G^F}s=\cup_{(L,\la)/G^F}\ser{G^F}{(L,\la)}\leqno{(2.2.7.2)}$$ }

\preuve Recall that $L$ is in duality with $L^*:=\cent{G^*}{\zo{L^*_s}_{\phi_d}}$. In formula (2.2.7.1) the parabolic subgroup is omitted because the result of induction is independant of it.

(i) A special case of (i) is $s=1$. Then $\Psi_{G,1}$ is a bijection $\ser{(G^*)^F}1\to \ser{G^F}1$ which commute with Lusztig induction, see [10]. 

If $s$ is central in $G^*$, then, by Proposition~1.3.2, (ii), tensor product by $\Psi_{G,s}(1)$ induces a bijection $\ser{G^F}1\to \ser{G^F}s$ hence   $\Psi_{L,s}(\al)=\Psi_{L,s}(1)\otimes \Psi_{L,1}(\al)$. (2.2.7.1) follows from the cases $s=1$ and  equalities $\Lu LG(\Psi_{L,s}(\al))=\Psi_{G,s}(1)\otimes \Lu LG(\Psi_{L,1}(\al))$, $\Lu LG(\Psi_{L,1}(\al))=\Psi_{G,1}(\Lu {L^*_s}{G^*}(\al))$.  

Assume now that $\cent{G^*}s$ is a Levi subgroup of $G^*$ and use  (iii) in Proposition~1.3.2 : let $G(s)$ be a Levi subgroup of $G$ in duality with $\cent{G^*}s$ such that $L^*_s$, a Levi subgroup of $G^*$, is in duality with $L(s):=G(s)\cap L$. Then $\Psi_{L,s}(\al)=\Lu{L(s)}L(\Psi_{L(s),s}(1_{L^{* F}_s})\otimes \Psi_{L(s),1}(\al))$. We have
$$\Lu LG(\Psi_{L,s}(\al))=\Lu{G(s)}G(\Lu{L(s)}{G(s)}(\Psi_{L(s),s}(\al))),\quad \Lu{L(s)}{G(s)}(\Psi_{L(s),s}(\al))=\Psi_{G(s),s}(\Lu{L^*_s}{\cent{G^*}s}\al)$$ (the case $s$ central above). Furthermore $\Lu{G(s)}G$ restricts to a bijection from $\ser{G(s)^F}s$ to $\ser{G^F}s$ that exchanges $\Psi_{G(s),s}$ and $\Psi_{G,s}$ : $\Lu{G(s)}G(\Psi_{G(s),s}(\Lu{L^*_s}{\cent{G^*}s}\al))=\Psi_{G,s}(\Lu{L^*_s}{\cent{G^*}s}\al)$, (2.2.7.1) follows.

(ii) If $L^*_s$ has type $\AA$, any central function on ${L^*_s}^F$ is uniform, hence $\Psi_{L,s}(\al)\in\piu^L(\Z \ser{G^F}s)$ by (1.3.2.2), and $\Lu LG (\Psi_{L,s}(\al))$ is uniform by (1.3.1.4). Then (1.3.2.5) applies and gives (2.2.7.1).

(iii) The square norms of $\Lu{L^*_s}{\cent{G^*}s}\al$ and $\Lu LG\la$ are given by Mackey decomposition formula : 

$||\Lu LG\la||^2=\scal{\Lu LG \la}{\Lu LG\la}{G^F}=\scal{\la}{\slu LG(\Lu LG \la)}{G^F}$... But as  $\la$ is $d$-cuspidal and  $L$ is $d$-split the formula for $\slu LG(\Lu LG \la)$ reduces to $\slu LG(\Lu LG\la)=\sum_{g\in \nor G L^F/L^F}\lexp g\la$, so that $||\Lu LG\la||^2=|\nor GL^F_\la/L^F|$.

Similarly, $||\Lu {L^*_s}{\cent{G^*}s}\al||^2=|\nor {\cent{G^*}s}{L^*_s}^F_\al/{L^*_s}^F|$. 

By Proposition~2.2.6 the quotients $\nor {G^F}{L,\ser{L^F}s}/L^F$ and $\nor {\cent{G^*}s}{L^*_s}^F/{L^*_s}^F$ are isomorphic. By assertion (iv) of Proposition~1.3.2,  $\nor {\cent{G^*}s}{{L^*_s}}^F_\al/{L^*_s}^F$ is isomorphic to $\nor GL^F_\la/L^F$, so that 
$$||\Lu LG\la||^2=||\Lu {L^*_s}{\cent{G^*}s}\al||^2\;.$$ 

By (1.3.2.5) and (1.3.1.4) $\piu^G(\Lu LG\la)=\Psi_{G,s}(\piu^{\cent{G^*}s}(\Lu{{L^*_s}}{\cent{G^*}s}\al))$. 

More generally "$\Psi$ commute with $\piu$" by formulas (1.3.2.2). Thus there exists a unique $\chi\in\Z\ser{G^F}s$ such that $\piu^G(\chi)=\Psi_{G,s}(\piu^{\cent{G^*}s}(\Lu{{L^*_s}}{\cent{G^*}s}\al))$ and $||\chi||^2=||\Lu {{L^*}(s)}{\cent{G^*}s}\al||^2$ and it is $\chi=\Psi_{G,s}(\Lu {{L^*_s}}{\cent{G^*}s}\al)$. We have $\chi=\Lu LG\la$.

(b) Assume  the existence of $\tilde\Psi_{G,s}(L,\la)$  for all proper unipotent $d$-cuspidal data $({L^*_s},\al)$ in $(\cent{G^*}s,F)$. Summing on the set of $\cent{G^*}s^F$-conjugacy classes on such data denote
$N:=|\cup\ser{\cent{G^*}s^F}{({L^*_s},\al)}|$. Thanks to Propositions  2.1.7, 2.2.2 and hypothesis one has $$N=\sum|\ser{\cent{G^*}s^F}{(L^*_s,\al)}|=\sum|\ser{G^F}{(L,\Psi_{L,s}(\al))}|\leq \sum|\II{b_{G^F}(L,\psi_{L,s}(\al))}\cap\ser{G^F}s|\leqno{(2.2.7.3)}$$ 
where $L$ corresponds to ${L^*_s}$ as in Proposition~2.1.4 so that $L^*_s=\cent{L^*}s$ and $(L,\Psi_{L,s}(\al))$ runs in the set of $G^F$-conjugacy classes of $d$-cuspidal data in series $(s)$ such that $\cent{G^*}s\neq L^*_s$. 

Let $N_0$ be the number of $d$-cuspidal elements of $\ser{\cent{G^*}s^F}1$. 
By G. $d$-HC in $\ser{\cent{G^*}s^F}1$, $N_0+N$ is $|\ser{\cent{G^*}s^F}1|$. By the existence of $\Psi_{G,s}$, $|\ser{G^F}s|=|\ser{\cent{G^*}s^F}1|$ and $N_0$ is greater that the number of elements of the complement of $ \cup_{L\neq L_0}\II{b_{G^F}(L,\la)}\cap\ser{G^F}s$ in $\ser{G^F}s$ where  $L_0$ is a Levi subgroup of $G$ in duality with $\cent{G^*}{\zo{\cent{G^*}s}_{\phi_d}}$. There are $N_0$ non $G^F$-conjugate $d$-cuspidal data in series $(s)$ of the form $(L_0,\Psi_{L_0,s}(\al))$ in $G^F$, in different $G^F$-conjugacy classes (see Proposition 2.2.6) and $L_{0,s}^*=\cent{G^*}s$ implies $||\Lu {L_0}G(\Psi_{L_0,s}(\al))||^2=1$, so that $\II{b_{G^F}(L_0,\Psi_{L_0,s}(\al)}\cap \ser{G^F}s=\{\Lu {L_0}G(\Psi_{L_0,s}(\al))\}$.  In (2.2.7.3) above we have $N=\sum|\II{b_{G^F}(L,\psi_{L,s}(\al))}\cap\ser{G^F}s|$, hence  the partition of $\ser{G^F}s$ by blocks is the partition by $d$-series $\ser{G^F}{(L,\Psi_{L,s}(\al))}$, that is (2.2.7.2). \bull

\medskip\noindent {\bf 2.2.8} {\it Proof of Proposition~2.2.4. } 

Note  that, once (2.2.7.2) is proved for all $G$ under assumption~2.2.3, then, for any $\chi\in\ser{G^F}{(L,\la)}$, (2.2.1.2), equivalently (2.2.4.2), is true, with $a=\scal\chi{\Lu LG \la}{G^F}$. Indeed let   $\xi$ be some \irr\ component of $\slu LG\chi$. By  (2.2.7.2) applied to $(L,F)$, $\xi$ belongs to  $\ser{L^F}{(M,\mu)}$ for some $d$-cuspidal datum $(M,\mu)$ in $(L,\la)$. Then $(M,\mu)$ is a $d$-cuspidal datum in $(G,F)$ and, by Propositions~2.1.7 and 2.1.6, we have 
$$\xi\in\II{b_{L^F}(M,\mu)},\quad \Lu LG(b_{L^F}(M,\mu))=b_{G^F}(M,\mu),\quad\chi\in\II{b_{G^F}(M,\mu)}\cap \II{b_{G^F} (L,\la)}$$ hence $(L,\la)$ and $(M,\mu)$ are $G^F$-conjugate.

Using Lemma~2.2.7 (b), we have to prove the existence of $\tilde\Psi_{G,s}(L,\la)$ for any proper unipotent $d$-cuspidal data $(L^*_s,\al)$ in $\cent{G^*}s,F)$.

If $G$ is a central product over a torus $G_1.G_2$ of $F$-stable reductive subgroups with connected centers, $G^F={G_1}^F.{G_2}^F$, $s$ has image $(s_1,s_2)$ in ${G_1}^*\times {G_2}^*$, the space $K\ser{\cent{G^*}s^F}1$ is an orthogonal product of $K\ser{\cent{{G_1}^*}{s_1}^F}1$ and $ K\ser{\cent{{G_2}^*}{s_2}^F}1$, $\Psi_{G,s}$ is defined by  $\Psi_{G_1,s_1}\times\Psi_{G_2,s_2}$. That's why we may assume $(G,F)$ irreducible. 

When (2.2.7.1) holds the restriction of $\Psi_{G,s}$ to $\ser{\cent{G^*}s^F}{(L^*_s,\al)}$ as $\Psi_{G,s}(L,\la)$ is a good choice. In classical types (2.2.7.2) is a consequence of G.$d$-HC in $\ser{\cent{G^*}s^F}1$, Proposition~2.1.10 and commutation between Lusztig induction and Jordan \dec. 

We have to consider groups $(G,F)$ of exceptional types when eventually Lemma~2.2.7 (a) don't apply for some $d$-cuspidal datum.
 
Unfortunately there are unipotent $d$-cuspidal data which don't satisfy any of the assumptions  (i), (ii), (iii) in  Lemma~2.2.7 on $({L^*_s},\al)$. It happens when two $d$-cuspidal \irr\ characters  have equal projections on the space of uniform functions. The equality  $\piu^L(\al)=\piu^L(\al')$ for distinct elements of $\ser{L^F}1$ occurrs only for algebraically conjugate representations in exceptional types. 

We recall in Table 1 all these cases as described in [26] Chapter~4 and Appendix.

  \medskip
  {\centerline{\bf Table 1}}
  \medskip
  
\centerline{\vbox{\offinterlineskip
\hrule
\halign{&\vrule#&\strut\quad#\hfil\quad&\vrule#&\strut\quad#\hfil\quad&\vrule#&\strut\quad#\hfil\quad\cr
height3pt&\omit&&\omit&&\omit&&\omit&\cr
&Type&&$\al$'s&&$d$&& central defect&\cr
height3pt&\omit&&\omit&&\omit&&\omit&\cr
\noalign{\hrule}
height2pt&\omit&&\omit&&\omit&&\omit&\cr
&{\bf G}$_2$&&$\al(1)={1\over 3}q\phi_1(q)^2\phi_2(q)^2$&&$1,2$&&$\ell\geq5$&\cr
height2pt&\omit&&\omit&&\omit&&\omit&\cr
&{\bf F}$_4$&&F$_4[\theta^j]$, $j=1,2$&&$1,2,4,8$&&$\ell\geq3$&\cr
height1pt&\omit&&\omit&&\omit&&\omit&\cr
&&&F$_4[\pm i]$&&$1,2,3,6$&&$\ell\geq5$&\cr
height2pt&\omit&&\omit&&\omit&&\omit&\cr
&{\bf E}$_6$&&E$_6[\theta^j]$, $j=1,2$&&$1,2,4,5,8$&&$\ell\geq5$&\cr
height2pt&\omit&&\omit&&\omit&&\omit&\cr
&$\lexp 2{\bf E}_6$&&$\lexp 2{\rm E}_6[\theta^j]$, $j=1,2$&&$1,2,4,8,10$&&$\ell\geq5$&\cr
height2pt&\omit&&\omit&&\omit&&\omit&\cr
&{\bf E}$_7$&&E$_6[\theta^j,1]$, $j=1,2$&&$4,5,7,8,10,14$&&$\ell\geq5$&\cr
height1pt&\omit&&\omit&&\omit&&\omit&\cr
&&&E$_6[\theta^j,\epsilon]$, $j=1,2$&&$4,5,7,8,10,14$&&$\ell\geq5$&\cr
height2pt&\omit&&\omit&&\omit&&\omit&\cr
&{\bf E}$_8$&&E$_6[\theta^j,1]$, $j=1,2$&&$4,5,7,8,10,12,14,20,24$&&$\ell\geq5$&\cr
height1pt&\omit&&\omit&&\omit&&\omit&\cr
&&&E$_6[\theta^j,\epsilon]$, $j=1,2$&&$4,5,7,8,10,12,14,20,24$&&$\ell\geq5$&\cr
height1pt&\omit&&\omit&&\omit&&\omit&\cr
&&&E$_6[\theta^j,\epsilon']$, $j=1,2$&&$4,5,7,8,10,14,15,20,30$&&$\ell\geq5$&\cr
height1pt&\omit&&\omit&&\omit&&\omit&\cr
&&&E$_6[\theta^j,\epsilon'']$, $j=1,2$&&$4,5,7,8,10,14,15,20,30$&&$\ell\geq5$&\cr
height1pt&\omit&&\omit&&\omit&&\omit&\cr
&&&E$_6[\theta^j,r]$, $j=1,2$&&$2,4,5,7,8,10,14,15,18,20,24$&&$\ell\geq5$&\cr
height1pt&\omit&&\omit&&\omit&&\omit&\cr
&&&E$_6[\theta^j,r']$, $j=1,2$&&$2,4,5,7,8,10,14,18,20,30$&&$\ell\geq5$&\cr
height1pt&\omit&&\omit&&\omit&&\omit&\cr
&&&E$_8[\la^j]$, $j=1,2,3,4$&&$1,2,3,4,6,7,8,9,12,14,18,24$&&$\ell = 3$, $\ell\geq7$&\cr
height1pt&\omit&&\omit&&\omit&&\omit&\cr
&&&E$_8[-\theta^j]$, $j=1,2$&&$1,4,5,7,8,9,10,12,14,15,20$&&$\ell\geq5$&\cr
height1pt&\omit&&\omit&&\omit&&\omit&\cr
&&&E$_8[\theta^j]$, $j=1,2$&&$1,4,5,7,8,9,10,14,20,24,30$&&$\ell\geq5$&\cr
height1pt&\omit&&\omit&&\omit&&\omit&\cr
&&&E$_8[\pm i]$&&$1,2,3,5,6,7,9,10,14,15,18,30$&&$\ell\geq 3$&\cr
height1pt&\omit&&\omit&&\omit&&\omit&\cr}
\hrule}}

We have used notations of [26] to design elements of $\ser{L^F}1$ in the second column. In the third one are the $d\in\N$ such that $\al$ is $d$-cuspidal and $\phi_d$ divides the polynomial order $P_{L,F}$ of $(L,F)$ (see 1.1.5; note that if $\phi_d$ does not divides $P_{L,F}$ and $L$ is a $d$-split Levi subgroup of $G$, then $L=G$ and any $\chi\in\II{G^F}$ is $d$-cuspidal). In the last column are the odd prime numbers $\ell$ such that $\al$ has central $\ell$-defect group, assuming that $d$ is the order of $q$ mod $\ell$, so that $d$ divides $(\ell-1)$.

Inspection of the matrices of decomposition in [26] 4.14---4.16 shows

\medskip\noindent {(2.2.8.2)  {\sl Assume $G$ simple of any type. If $\al,\beta\in\ser{G^F}1$ and $\piu^G(\al)$, $\piu^G(\beta)$ are proportional, then $\al=\beta$ short of the cases $\piu^G(\al)=\piu^G(\beta)$ listed in Table 1.}

Assume that $\cent{G^*}s$ is contained in a proper $d$-split Levi subgroup $M^*$ of $G^*$, with dual $M$ in $G$. By (iii) in Proposition~1.3.2, $\Lu LG$ induces a bijection $\ser{M^F}s\to \ser{G^F}s$. $M^F$-conjugacy classes of $d$-cuspidal data in series $(s)$ and $G^F$-conjugacy classes of $d$-cuspidal data in series $(s)$ are in natural bijection. If our claim on G.$d$-HC holds in $M$, it holds in $G$. If (2.2.7.1) holds in $M$, it holds in $G$.

So we assume now that
 $s$ is  isolated in $(G^*)^F_{\ell'}$, that is $\cent{G^*}s$ is not contained in any proper Levi subgroup of $G^*$. Assume further the existence in $(\cent{G^*}s,F)$ of a $d$-cuspidal unipotent datum $({L^*_s},\al)$ where ${L^*_s}$ is not of type $\AA$, so that (a) in Lemma~2.2.7 don't apply. There are few cases to consider, we list  in Table 2 in ascending rank, and use parametrizations of $\ser{{L^*_s}^F}1$ given in [26], short of $\ser{\lexp 3{\DD}_4(q)}1$, notations of [10] Table 1.

 \medskip
  {\centerline{\bf Table 2}}
  \medskip
  
\centerline{\vbox{\offinterlineskip
\hrule
\halign{&\vrule#&\strut\quad#\hfil\quad&\vrule#&\strut\quad#\hfil\quad&\vrule#&\strut\quad#\hfil\quad&\vrule#&\strut\quad#\hfil\quad\cr
height3pt&\omit&&\omit&&\omit&&\omit&&\omit&\cr
&Case&&Type of $G$&&Type of $\cent{G^*}s$&& d&& $({L^*_s},\al)$&\cr
height3pt&\omit&&\omit&&\omit&&\omit&&\omit&\cr
\noalign{\hrule}
&1&&{\bf F}$_4$&&$\BB_4$&&2&&$(\BB_2,{{0,2}\choose{ 1 \;}})$&\cr
height2pt&\omit&&\omit&&\omit&&\omit&&\omit&\cr
&2&&&&&&4&&$(\BB_2,{{0,1}\choose{ 2 \;}})$&\cr
height2pt&\omit&&\omit&&\omit&&\omit&&\omit&\cr
&3&&&&$\AA_1\times\CC_3$&&2&&$(\CC_2,{{0,2}\choose{ 1 \;}})$&\cr
height2pt&\omit&&\omit&&\omit&&\omit&&\omit&\cr
&4&&$\EE_7$&&$\DD_6\times\AA_1$&&2&&$(\DD_4,{{1,3}\choose{ 0,2}})$&\cr
height2pt&\omit&&\omit&&\omit&&\omit&&\omit&\cr
&5&&$\EE_8$&&$\DD_8$&&2&&$(\DD_4,{{1,3}\choose{0,2}})$&\cr
height2pt&\omit&&\omit&&\omit&&\omit&&\omit&\cr
&6&&&&&&3&&$ \lexp 2{\DD}_5$&\cr
height2pt&\omit&&\omit&&\omit&&\omit&&\omit&\cr
&7&&&&&&4&&$ {\DD}_4$&\cr
height2pt&\omit&&\omit&&\omit&&\omit&&\omit&\cr
&8&&&&&&6&&$ \lexp 2{\DD}_5$&\cr
height2pt&\omit&&\omit&&\omit&&\omit&&\omit&\cr
&9&&&&&&8&&$\lexp 2{\DD}_4$&\cr
height2pt&\omit&&\omit&&\omit&&\omit&&\omit&\cr
&10&&&&$\EE_7\times\AA_1$&&2&&$(\EE_7,[512_a], [512'_a])$&\cr
height2pt&\omit&&\omit&&\omit&&\omit&&\omit&\cr
&11&&&&&&2&&$(\lexp 2{\EE}_6,\lexp 2E_6[\theta^j])$ $ j=1,2$&\cr
height2pt&\omit&&\omit&&\omit&&\omit&&\omit&\cr
&12&&&&&&2&&$(\DD_4,{{1,3}\choose{0,2}})$&\cr
height2pt&\omit&&\omit&&\omit&&\omit&&\omit&\cr
&13&&&&&&3&&$(\lexp 3{\DD}_4\times\AA_1,\lexp 3{D}_4[-1]\otimes \al')$&\cr
height2pt&\omit&&\omit&&\omit&&\omit&&\omit&\cr
&14&&&&&&6&&$(\lexp 3{\DD}_4\times\AA_1,\phi_{2,1}\otimes\al')$&\cr
height2pt&\omit&&\omit&&\omit&&\omit&&\omit&\cr
&15&&&&$\EE_6\times\AA_2$&&3&&$(\lexp 3{\DD}_4,\lexp 3{D}_4[-1])$&\cr
height2pt&\omit&&\omit&&\omit&&\omit&&\omit&\cr
&16&&&&$\lexp 2\EE_6\times\lexp 2\AA_2$&&2&&$({\DD}_4,{{1,3}\choose{0,2}})$&\cr
height2pt&\omit&&\omit&&\omit&&\omit&&\omit&\cr
&17&&&&&&6&&$(\lexp 3{\DD}_4,\phi_{2,1})$&\cr
height3pt&\omit&&\omit&&\omit&&\omit&&\omit&\cr}
\hrule}}

(a) When $\cent{G^*}s$ has a type $\AA$ component, the space $\piu^{\cent{G^*}s}(K\ser{\cent{G^*}s^F}1)$ of uniform unipotent functions on $\cent{G^*}s^F$ is an orthogonal product of the corresponding spaces for the two components and differ from the all space of unipotent functions only on the other side, we have only to prove our claim on non type $\AA$  part. Thus the case 10 is clear.

(b) It happens, in cases 2, 3, 6, 8, 9,  that ${L^*_s}$ is a maximal proper $d$-split Levi subgroup of $\cent{G^*}s$ and $\cent{G^*}s$ have a classical type. Then $\Lu{{L^*_s}}{\cent{G^*}s}\al$ is given by Asai's $d$-hook formula ($d$ odd) [24], [10] (3.5) or a $d/2$-cohook formula [10] (3.9) in "one step". The formulas show that $\Lu{L^*_1}{\cent{G^*}s}\al$ is an algebraic sum of irreducible $\chi_j$ all in distinct families, families defined from $(W(\cent{G^*}s),F)$. In fact (iii) of (a) in Lemma~2.2.7 is satisfied, hence the commutation formula (2.2.7.1). 

(c) In some cases there exists between ${L^*_s}$ and $\cent{G^*}s$ a maximal proper $d$-split Levi subgroup $L^*_1$ of $\cent{G^*}s$ and $L^*_1$ is also a $d$-split Levi subgroup of $G^*$. Let $L_1$ be a ($d$-split) Levi subgroup of $G$ in the dual class such that $L\subseteq L_1$. As $s$ is central in $L^*_1$  (2.2.7.1) applies between $L$ and $L_1$ :  $\Lu L{L_1}(\Psi_{L,s}(\al))=\Psi_{L_1,s}(\Lu{{L^*_s}}{L^*_1}(\al))$. 

Assume first that $\cent{G^*}s$ has classical type. For any $\beta\in\ser{(L^*_1)^F}{(L^*_s,\al)}$ the hook or cohook formula gives, as in (a), $\Psi_{G,s}(\Lu{L^*_1}{\cent{G^*}s}\beta)=\Lu{L_1}G(\Psi_{L_1,s}(\beta))$. Finally we have (2.6.7.1) again.
This applies to cases 1, 4, 5, 7 where $L_1^F$ is $\BB_3(q).(q+1)$, $\DD_5(q).(q+1)^2$, $\lexp 2\DD_7(q).(q+1)$, $\lexp 2\DD_6(q).(q^2+1)$ respectively.

(d) In cases 12 to 17, ${L^*_s}$ is itself a Levi subgroup of $G^*$, $L^F$ is $\DD_4(q).(q+1)^4$, (resp. $\lexp 3\DD_4(q).(q^3-1).\AA_1(q)$, $\lexp 3\DD_4(q).(q^3+1).\AA_1(q)$, $\lexp 3\DD_4(q).(q^2+q+1)^2$, $\DD_4(q).(q+1)^4$, $\lexp 3\DD_4(q).(q^2-q+1)^2$)  and $\Lu{{L^*_s}}{\cent{G^*}s}\al$ is given on $\EE_7$ (resp. ...$\EE_6$, $\lexp 2\EE_6$)-side by the series 14 (resp. 17, 27, 5, 9, 11) in [10] Table 2 and (a) (iii) of Lemma~2.2.7 applies.

In double case 11, with ${L^*_s}=L^*_1$, $L_1^F$ is $\lexp 2 \EE_6(q).(q+1)^2$. Instead of Asai's formula we need the series 15 and 16 in [10] Table 2 : it is a case where  the two generalized characters $\Lu{{L^*_s}}{\cent{G^*}s}(\lexp 2E_6[\theta^j])$ cannot be distinguished by their projection on the space of uniform functions, but the known common norm and value on 1 imply their form. Thus  (a) (iii) of Lemma~2.2.7 does not applies nevertherless  $$\{\ser{G^F}{(L,\Psi_{L,s}(\lexp 2E_6[\theta^j])}\}_{j=1,2}=\{\Psi_{G,s}(\ser{\cent{G^*}s^F}{(L^*_s,\lexp 2E_6[\theta^j])}\}_{j=1,2}$$ 
It is the only \irr\ case where (2.2.7.1) is not proved with our combinatorial arguments and where $\tilde\Psi_{G,s}(L,\la)$ may differ from the  restriction of $\Psi_{G,s}$. 
\bull

 \bigskip\noindent{\bf 2.3. On \irr\ characters in a block of $G^F$, $\z G $ connected }

We want to describe $\II {b_{G^F}(L,\la)}$ for any $d$-cuspidal datum $(L,\la)$ in an $\ell'$-series in $(G,F)$. In [15] Theorem~2.8 it is shown, using induction on blocks (Proposition~2.1.6), $\z G$ connected or not, how to recover $\II b\cap\ser{G^F}{s}$ knowing $\II b\cap\ser{G^F}{s_{\ell'}}$ for any bloc $b$ of $G^F$.  To obtain our main theorem we refer to the case of connected center and rely the above connection to relation between unipotent $d$-cuspidal data in $(\cento{G^*}{s},F)$ and unipotent $d$-cuspidal data in $(\cento{G^*}{s_{\ell'}},F)$, in Proposition~2.3.5. A similar way is to refer to $d$-cuspidal data in any series, as in Proposition~2.3.6.

\medskip\noindent{\bf 2.3.1. Notation. }
{\sl Let $L$ and $M$ be two $F$-stable connected  reductive subgroups of $(G,F)$, $\la\in\ser{L^F}1$, $\mu\in\ser{M^F}1$.  Denote
$$(L,\la)\simex {G^F}(M,\mu)$$
the following equivalence relation $${\rm there \; exists\;} g\in G^F \;{\rm such \; that }\quad[L,L]=[M,M]^g\;{\rm
and \;}\Res{{L}^F}{[L,L]^F}\la=(\Res{M^F}{[M,M]^F}\mu)^g.$$}

The relation $\simex {G^F}$ is the extension by $G^F$-conjugacy of the relation $\sim$  introduced in [15]  (see Definition 23.1 and Proposition 23.2 in [16]). 

\smallskip Assertion (c) of Proposition~2.1.5 suggests that there are few $d$-cuspidal \irr\ characters and few $d$-cuspidal data. The following Proposition is a consequence of that fact, and will be used to simplify the description  of $\II{b_{G^F}(L,\la)}\cap
\ser{G^F}{st}$ (see [16] Theorem~23.2 in case $s=1$). 

 \medskip\noindent{\bf  2.3.2. Proposition. }{\sl  Let $H$ be an $E$-split Levi subgroup of
$G$. 

(a) Let  $(L_H,\la_H)$ be a unipotent $d$-cuspidal datum in $ (H,F)$. There exists a unipotent $d$-cuspidal datum $(L,\la)$ in $(G,F)$ such that $[L,L]=[L_H,L_H]$, $L_H=L\cap H$, and $\Res {L^F}{[L,L]^F}\la=\Res{L_H^F}{[L_H,L_H]^F}\la_H$.

(b)  Let $(L,\la)$ et $(M,\mu)$ be two unipotent
$d$-cuspidal data in  $(G,F)$, and $(L_H,\la_H)$, $(M_H,\mu_H)$ be two unipotent
$d$-cuspidal data in  $(H,F)$ such that $[L,L]=[L_H,L_H]$ and $[M,M]=[M_H,M_H]$. Then
$L$ and $M$ are $G^F$-conjugate if and only if $L_H$ and $M_H$ are $H^F$-conjugate.

(c) The relation
$(L,\la)\simex {G^F}(L_H,\la_H)$ define a one-to-one map between the set of $G^F$-conjugacy classes of $d$-cuspidal unipotent data $(L,\la)$
in $(G,F)$ such that  $[L,L]\subseteq gHg^{-1}$ for some $g\in G^F$ and the set of $H^F$-conjugacy classes of $d$-cuspidal unipotent data $(L_H,\la_H)$
in $(H,F)$.}

\preuve  Assertion~(a) follows from (i) in  [16] Proposition~23.3 where the relation 

\centerline{$[L,L]=[L_H,L_H] \quad{\rm and}\quad \Res {L^F}{[L,L]^F}\la=\Res{L_H^F}{[L_H,L_H]^F}\la_H$}

\noindent is denoted $(L,\la)\sim(L_H,\la_H)$. We have added the equality $L_H=H\cap L$. 

Assume $[L,L]=[L_H,L_H]$. By [13] Proposition~1.7, $\zo L$ is a maximal torus in $\cento G{[L,L]}$ and $\zo{L_H}$ is a maximal torus in $\cento H{[L_H,L_H]}=H\cap\cento G{[L,L]}$. There exist some $z\in \cento G{[L,L]}$ such that $\zo {L_H}\subseteq \zo L^z$, hence $\zo{L_H}= H\cap \zo L^z$. Thus $L_H=\zo {L_H}.[L_H,L_H]=(H\cap \zo L^z).[L,L]=H\cap L^z$. We may replace $(L,\la)$ by $(L^z,\la^z)$, preserving the above two equalities.

 Clearly (b) implies (c).

Proof of (b) :

If $G=\a G$, then $H=\a H$ and there is only one conjugacy class of $d$-cuspidal data in $G$ and in $H$. So we assume $\b G\neq \{1\}$ and $H\neq G$ and  prove the Proposition inductively on the semi-simple rank  of $G$. 

The restriction along an isotypic morphism is a bijection on unipotent series and Proposition~2.1.5 allows us to assume $(G,F)$ rationally \irr.
By properties of scalar descent we assume that $G$ is \irr.

(b.i) Assume that  $H$ is  $d$-split in $G$, so are $L_H$ and $M_H$.

 Assume first that $L_H=(M_H)^h$ for some $h\in H^F$. Then $[L,L]=[M,M]^h$ hence, by [16] Proposition~22.8, $L$ and $M^h$ are $\cento G{[L,L]}^F$-conjugate, so that $L$ and $M$ are $G^F$-conjugate.
 
 Assume now that $L$ and $M$ are $G^F$-conjugate.

In classical types ${\bf A}$, $\BB$, $\CC$, $\DD$ with $d\geq 1$  [10,
\S~3],  and in exceptional  types with  $d=1$ [26] Appendix one sees on tables that

\smallskip {\sl (i.a)  $H=\a H$ is a torus or $\b H$ is rationally  \irr. 

(i.b) If $L$  is not a torus, then 
$\b L$ is rationally \irr. 

(i.c) If $\b L$ and $\b M$ have same type, then
$L$ and $M$ are $G^F$-conjugate.}

From these facts $\b L=[L,L]=[L_H,L_H]= (L_H\b)$ (eventually $=1$). In $(G,F)$ (resp. in $(H,F)$) the conjugacy class of $L$ (resp. $M$) is defined by $\b L$ (resp. $\b{(L_H)}$, hence our claim.

In exceptional type with  $d>1$ Table 1 in  [10]  show that (i.a) and (i.b) are no more true, but (i.c) is  true for each $d$.
One may verify, using properties (i.a), (i.b), (i.c) in classical \irr\ type and (i.c) in exceptional \irr\ type, that a proper maximal $d$-split Levi subgroup of 
$G$ cannot contain two non conjugate unipotent $d$-cuspidal  data as $(L_H,\la_H)$ with $(L_H\b )$ of same type as $\b L$ and one conclude by induction.  

(b.ii) In the general case with  $G=\b G$ and $H$ a proper $E$-split Levi subgroup of $G$, let $K$ be a proper $d$-split Levi subgroup of $G$ that contains $H$ (1.1.5.3). By (a) there exist unipotent $d$-cuspidal  data $(L_K,\la_K)$ and $(M_K,\mu_K)$ in $(M,F)$ such that $(L_H,\la_H)\simex {K^F}(L_K,\la_K)$ and $(M_H,\mu_H)\simex {K^F}(M_K,\mu_K)$.
 Induction applies in  $K$ : $L_K$ and $M_K$ are $K^F$-conjugate. By part (b.i) of our proof $L$ and $M$ are $G^F$-conjugate.
\bull

 In the following proposition we make a connection between  induction of  blocks (Proposition~2.1.7) and cuspidal data, using [15] and [16].
 
 \medskip\noindent{\bf 2.3.3. Proposition. }{\sl Assume 2.1.2 on $(G,F,\ell,d)$ and $\z G$ is connected. Let $H$ be an $E$-split Levi subgroup of $G$ with a dual $H^*$ in $G^*$,
such that $s\in (H^*)^F_{\ell'}$. Let
$(L_H,\la_H)$ be a $d$-cuspidal datum in  $(H,F)$ in series $(s)$, let $(L^*_{s,H},\al_H)$
be a $d$-cuspidal unipotent datum in $(\cent{H^*}s,F)$
such that  $L^*=\cent{H^*}{\zo{L^*_{s,H}}_{\phi_d}}$ belongs to the dual $(H^*)^F$-conjugacy class of the $H^F$-conjugacy class of $L$, and
$\la_H=\Psi_{L_H,s}(\al_H)$ (see Proposition~2.1.4). 

Let $(L^*_s,\al)$ be a  $d$-cuspidal
unipotent datum in $(\cent{G^*}s,F)$ such that
$(L^*_{s,H},\al_H)
\simex{\cent{G^*}s^F} (L^*_s,\al)$. Then $(L^*_s,\la)$ defines a $d$-cuspidal datum $(L,\la)=(L,\Psi_{L,s}(\al))$ in series $(s)$ in $(G,F)$, (Proposition~2.1.4). 

We have $$\Lu HG(b_{H^F}(L_H,\la_H))=b_{G^F}(L,\la)$$ and, if
$\xi\in
\II{b_{H^F}(L_H,\la_H)}$ and  $s_0\in\z{H^*}^F_\ell$, then $\Lu HG
(\Psi_{H,s_0}(1) \otimes\xi)\in \Z\II{b_{G^F}(L,\la)}$.  

}

\preuve We know by [16] Proposition~22.8, (iii) that  $\zo {L^*_s}_{\phi_d}$ (resp. $\zo
{L^*_{H,s}}_{\phi_d}$) is a maximal
$\phi_d$-subgroup in
$\cento {\cent{G^*}s}{[L^*_s,L^*_s]}$ (resp.
$\cento {\cent{H^*}s}{[L^*_{s},L^*_s]}$). By Sylow's theorem on
$\phi_d$-subgroups
([16] Theorem~13.18) there exist
$c\in\cento{\cent{G^*}s}{[L^*_s,L^*_s]}^F$ such that
$\zo{L^*_{H,s}}_{\phi_d}\subseteq \lexp c {\zo{L^*_s}_{\phi_d}}$. Up to conjugacy by $c$ we may assume
$\zo{L^*_{H,s}}_{\phi_d}\subseteq  {\zo{L^*_s}_{\phi_d}}$. By construction and up to $H^F$-conjugacy  we may assume 
$\zo{L_{H}}_{\phi_d}\subseteq  {\zo{L}_{\phi_d}}$ (Proposition~2.1.4). 

By Proposition~2.1.6, for any
$\xi\in
\II{b_{H^F}(L_H,\la_H)}\cap \ser{H^F}s$, any suitable parabolic $P$, $\Lu {H\subseteq P}G\xi$
belongs to  $\Z \II B$ where 
$B=\Lu HG(b_{H^F}(L_H,\la_H))$.  We prove by induction on the rank of $G$ that  $B=b_{G^F}(L,\la)$ : our claim for $s_0 =
1$.  

Let be $G=\a G.\b G$ the \dec\ in central product we defined in 1.1.5.2
and use the standard dichotomy whereas $\b G\subseteq H$ or not (recall functorial properties of $d$-cuspidal data, Proposition~2.1.5). 

If
$G=\a G$ there is only one block $B_{G^F}(s)$ and  one $G^F$-conjugacy class of cuspidal data in $(L,\la)$ series $(s)$ :  $(L,\la)$ is such that $L^*_s$
is a diagonal torus in  $\cent{G^*}s$ and $\al=1$ (Proposition~2.1.5, (c)). In that case $H=\a H$ and the equality $\Lu HG(B_{G^F}(s))=B_{H^F}(s)$ is evident. Assume more generally that  
$\b G\subseteq H $. Then $H=(H\cap \a G).\b G$ and there is a natural bijection from the set of $\ell$-blocks of $H^F$ in series $(s)$  to  the set of $\ell$-blocks of ${\b G}^F$ in some series $(\b s)$, (or from the set of $H^F$-conjugacy classes of $d$-cuspidal data $(L_H,\la_H)$ in series $(s)$ and the set of ${\b G}^F$-conjugacy classes of $d$-cuspidal data in series $(\b s)$), and as well to the set of $G^F$-conjugacy class of $d$-cuspidal data of ${G}^F$ in  series $(s)$ : $(L_H,\la_H)\to (L,\la)$, with $L_H\cap\b G=L\cap \b G$ and $\Res{{L_H}^F}{(L_H\cap \b G)^F}\la_H=\Res{L^F}{(L\cap \b G)^F}\la.$ Now
$\Lu HG$ reduces to
$\Lu{H\cap
\a G}{\a G}$, so that $\Lu HG$ commutes with that bijection.  

If $\b G$ is not contained in $H$, let $M$ be a proper $d$-split  Levi subgroup of $G$ that contains  $H$  (1.1.5.3). We have
$\zo M_{\phi_d}\subseteq \zo H_{\phi_d}\subseteq \zo
{L_H}_{\phi_d}\subseteq
\zo L_{\phi_d}$ hence $L\subseteq M$. In the dual class there exist a $d$-split Levi subgroup $M^*$ of $G^*$ such that $L^*\subseteq M^*$. Inductive hypotheses says that  $\Lu HM\xi\in\Z\II{b_{M^F}(L,\la)}$ for any  $\xi\in
\II{b_{H^F}(L_H,\la_H)}\cap \ser{H^F}s$. By Proposition~2.1.7, (ii)
 $\Lu HG\xi=\Lu MG(\Lu HM\xi)\in\Z\II{b_{G^F}(L,\la)}$, that is $\Lu HG(b_{H^F}(L_H,\la_H))=b_{G^F}(L,\la)$.

Let $s_0\in\z{H^*}^F_\ell$, $\cent{G^*}{s_0}$ is a Levi subgroup of $G^*$. Let $G(s_0)$ be a Levi subgroup of $G$  in the dual class of $\cent{G^*}{s_0}$ and such that $H\subseteq G(s_0)$. Denote $\theta_0=\Psi_{G(s_0),s_0}(1)$. As $s_0$ is central in $(G(s_0)^*)^F_{\ell}$, $\theta_0(1)=1$  (see Proposition~1.3.2, (ii)) and $\theta_0$ has order a power of $\ell$. Furthermore $\Lu H{G(s_0)}(\Psi_{H,s_0}(1)\otimes \xi)=\theta_0\otimes \Lu H{G(s_0)}\xi$. By [15]
Theorem~2.8, for any component $\xi_0$ of $\Lu H{G(s_0)}\xi$,
$\Lu{G(s_0)}G\xi_0$ and $\theta_0\otimes\Lu{G(s_0)}G\xi_0$ belongs to the same block of  $G^F$ that is  $b_{G^F}(L,\la)$.
 \bull

From [9] we know that if $(L,\la)$ is in series $(s)$ then $\II{b_{G^F}(L,\la)}\subseteq \lser{G^F}s$. In Proposition~2.3.3 above appear some non $\ell'$-series $(ss_0)$. So it's time to describe entirely $\II {b_{G^F}(L,\la)}$.

\medskip\noindent{\bf 2.3.4.  Notation } {\sl When $s,s_0$ are semi-simple in $(G^*)^F$ and $s_0=(ss_0)_\ell$, define $\bl{G^F}{s,s_0}$ as a set of blocks  of $G^F$ by $$B\in \bl{G^F}{s,s_0}\;{\rm if\;and\;only\;if}\;\II B\cap\ser{G^F}{ss_0}\neq \emptyset$$
Thus $\bl{G^F}{s,1}=\bl{G^F}s$ and $\bl{G^F}{s,s_0}\subseteq \bl{G^F}{s}$. 
 }

\medskip\noindent{\bf 2.3.5. Proposition. }{\sl Assume 2.2.3 on $(G,F,\ell,d)$. Let $(L,\la)$ be a $d$-cuspidal datum in series $(s)$ in $(G,F)$ ($s$ semi-simple, $s\in (G^*)^F_{\ell'}$) and $(L^*_s,\al)$ an associated $d$-cuspidal unipotent datum in $(\cent{G^*}s,F)$ by Proposition~2.1.4. Assume   $\z G$  connected if $s\neq 1$. Let $s_0\in\cent{G^*}s^F_\ell$. Then
$\II{b_{G^F}(L,\la)}\cap\ser{G^F}{ss_0}\neq \emptyset$ if and only if there exist a $d$-split Levi subgroup $L^*_{ss_0}$ of $\cent{G^*}{ss_0}$ such that $[L^*_{ss_0},L^*_{ss_0}]$ is $\cent{G^*}s^F$-conjugate to $[L^*_s,L^*_s]$. For such $L^*_{ss_0}$ there exist a unique $d$-cuspidal $\al_0\in\ser{(L^*_{ss_0})^F}1$ such that $$(L^*_{ss_0},\al_0)\simex{\cent{G^*}s^F}(L^*_s,\al)\leqno{(2.3.5.1)}$$
When (2.3.5.1) holds, let $G(s_0)$ be a Levi subgroup of $G$ in the  dual $G^F$-conjugacy class of the $G^{*\,F}$-conjugacy class of $\cent{G^*}{s_0}$, and let $(L_0,\la_0)$ be a $d$-cuspidal datum in series $(s)$ in $(G(s_0),F)$ associated to $(L^*_{ss_0},\al_0)$. One has 
$$\II{b_{G^F}(L,\la)}\cap\ser{G^F}{ss_0}=\{\Lu {G(s_0)}G(\Psi_{G(s_0),s_0}(1)\otimes \xi)\mid \xi\in\ser{G(s_0)^F}{(L_0,\la_0)}\}\leqno{(2.3.5.2)}$$
}

\preuve
Note that (2.3.5.2) has been proved  for $s_0=1$ (Propositions~2.2.2 and 2.2.4). If $\z G$ is connected, $\z{G(s_0)}$ is connected (Proposition~1.1.3 (b)), so  this case gives also $$\ser{G(s_0)^F}{(L_0,\la_0)}=\Psi_{G(s_0),s}(\ser{\cent{\cent{G^*}{s_0}}s^F}{(L^*_{ss_0},\al_0)}\leqno{(2.3.5.3)} $$

 If $s=1$, Proposition~2.3.5 is contained in [13] Theorem~4.4. 
 So we assume $\z G$ connected. Then G.$d$-HC holds in $\ser{G(s_0)^F}s$ by Proposition~2.2.4.
 
 Assuming $[L^*_{ss_0},L^*_{ss_0}]=[L^*_s,L^*_s]$,  the equality  $\Res{(L^*_s)^F}{[L^*_s,L^*_s]^F} \al=\Res{(L^*_{ss_0})^F}{[L^*_{ss_0},L^*_{ss_0}]^F} \al_0$ defines a unique $\al_0$ in $\ser{(L^*_{ss_0})^F}1$. Then if $L^*_{ss_0}$ is $d$-split, $(L^*_{ss_0},\al_0)$ is a $d$-cuspidal unipotent datum in $(\cent{G^*}{ss_0},F)$ and we have (2.3.5.1).

We may write $G$ as a central product of rationally \irr\ components $G_j$ with connected center $\z {G_j}=\z G$.  Then by the usual process $G^F$ is a central quotient of $\times_j{G_j}^F$, $ss_0$ has image $(s_js_{j,0})_j$ in $\times_j{G_j}^*$, $\ser{G^F}{ss_0}$ is in bijection with $\times_j\ser{{G_j}^F}{s_js_{j,0}}$, $(L,\la)$ is image of $\times _j(L_j,\la_j)$, and $\II{b_{G^F}(L,\la)}$ is the set of elements of $\times_j\II{B_{G_j^F}(L_j,\la_j)}$ whose kernel contains the kernel of $\times_j{G_j}^F\to G^F$, and so on...

So we assume $(G,F)$ rationally \irr. 

If $G=\a G$, then there is only one conjugacy class of $d$-cuspidal data in series $(s)$ in $(G,F)$ and in $(G(s_0),F)$ ($L^*_s$ and $L^*_{ss_0}$ are diagonal torii in $\cent{G^*}s$ and $\cent{G^*}{ss_0}$ respectively) and only one $\ell$-block $b$ of $G^F$ such that $\II{b}\cap \ser{G^F}{ss_0}\neq \emptyset$. The condition (2.3.5.1) is satisfied. One has $\II{b}=\cup_{s_0}\ser{G^F}{ss_0}$. With notation of the Proposition, when $\xi$ runs in $\ser{G(s_0)^F}s$, $\Lu {G(s_0)}G(\Psi_{G(s_0),s_0}(1)\otimes \xi)$ runs in $\ser{G^F}{ss_0}$ by Proposition~1.3.2, (ii) and (iii). That gives (2.3.5.2).

If $G=\zo G.\b G$, $\cent{G^*}{s_0}$ is an $E$-split Levi subgroup of $G^*$  :  $\cF (G)^F$ being prime to $\ell$, by Proposition~1.1.3 $(\z {\cent{G^*}{s_0}}/\zo{\cent{G^*}{s_0}})^F$ is prime to $\ell$, hence $s_0\in \zo{\cent{G^*}{s_0}}$, hence $\zo{\cent{G^*}{s_0}}_{\phi_E}\neq 1$ by definition of $E$. We may apply Proposition~2.3.3 with $H=G(s_0)$. It gives us an inclusion $$\Lu {G(s_0)}G\big(\Psi_{G(s_0),s_0}(1)\otimes \II{b_{G(s_0)^F}(L_0,\la_0)}\cap \ser{G(s_0)^F}s\big)\subseteq \II{b_{G^F}(L,\la))}$$
Assume (2.3.5.1). 
By (2.3.5.2) for $1$ in $G(s_0)$ we know that 
$$\ser{G(s_0)^F}s=\cup_{(L_0,\la_0)}\II{b_{G(s_0)^F}(L_0,\la_0)}\cap \ser{G(s_0)^F}s$$
a disjoint union when $d$-cuspidal data $(L_0,\la_0)$ are considered modulo $G(s_0)^F$-conjugacy. 
Now $G(s_0)^F$-conjugacy on $d$-cuspidal data $(L_0,\la_0)$ in series $(s)$ in $(G(s_0),F)$ corresponds to $\cent{G^*}{ss_0}^F$ conjugacy of unipotent $d$-cuspidal data $(L^*_{ss_0},\al_0)$ in $(\cent{G^*}{ss_0},F)$ (Proposition~2.1.4), that corresponds, {\it via} (2.3.5.1)  to $\cent{G^*}s^F$-conjugacy of some unipotent $d$-cuspidal data $(L^*_s,\al_s)$ in $(\cent{G^*}s,F)$ (Proposition~2.3.2). Applying Proposition 2.1.4 again and Proposition~1.3.2, (ii) and (iii), we see that \irr\ in series (s) and in different blocks $b_{G(s_0)^F}(L_0,\la_0)$ are sent by $\xi\mapsto \Lu {G(s_0)}G\big(\Psi_{G(s_0),s_0}(1)\otimes \xi)$ in series $(ss_0)$ and in different blocks of $G^F$. Equality (2.3.5.2) is proved.
\bull

\medskip\noindent{\bf Proposition~2.3.6. }{\sl Assume 2.2.3 on $(G,F,\ell,d)$ and $\z G$ connected. Let $s$, $s_0$ in $G^{*\,F}$ such that $(ss_0)_\ell=s_0$. Let $G(s_0)$ be a Levi subgroup of in dual $G^F$-conjugacy class of the $G^{*\,F}$-conjugacy class of $\cent{G^*}{s_0}$.

(a) One defines a one-to-one map from the set of $G(s_0)^F$-conjugacy class of $d$-cuspidal data in series $(s)$ in $(G(s_0),F)$ onto the set of $G^F$-conjugacy class of $d$-cuspidal data in series $(ss_0)$ in $(G,F)$ as follows :

Let $\al_0\in\ser{\cent{L^*_0}s^F}1$, $\la_0=\Psi_{L_0,s}(\al_0)$ such that $(L_0,\la_0)$ be a $d$-cuspidal datum in series $(s)$ in $(G(s_0),F)$. Here $L^*_0$ belongs to  the dual $\cent{G^*}{s_0}^F$-conjugacy class of $L_0$, $s\in L^*_0$. Define $M^*$ by
$$M^*:=\cent{G^*}{\zo{L_0^*}_{\phi_d}}\leqno{(2.3.6.1)}$$
Let $M$ be in the dual $G^F$-conjugacy class of the $(G^*)^F$-conjugacy class of $M^*$ and $\mu$ be defined by 
$$\mu=\Lu{L_0}M(\Psi_{L_0,s_0}(1)\otimes\la_0)$$
Then $(M,\mu)$ is a $d$-cuspidal datum in series $(ss_0)$ in $(G,F)$. 

Generalized $d$-Harish-Chandra theory holds in $\ser{G^F}{ss_0}$.

(b) Induction on blocks $\Lu{G(s_0)}G$, as defined in Proposition~2.1.6, restricts to a one-to-one map from $\bl{G(s_0)^F} s$ onto $ \bl{G^F}{s,s_0}$. 

Let $L^*_0$, $\al_0$, $(L_0,\la_0)$ and $ (M,\mu)$ as  in (a).  Let $(L,\la)$ be a $d$-cuspidal datum  in series $(s)$ in $(G,F)$, defined, thanks to Proposition~2.1.4, by some $d$-cuspidal unipotent datum $(L^*_s,\al)$ in $(\cent{G^*}s,F)$ such that $$(L^*_s,\al)\simex{\cent{G^*}s^F}(\cent{L^*_0}s,\al_0)\,.$$ One has 
$$\Lu{G(s_0)}G(b_{G(s_0)^F}(L_0,\la_0))=b_{G^F}(L,\la), \leqno{(2.3.6.2)}$$
$$\eqalign{\II{b_{G^F}(L,\la)}\cap \ser{G^F}{ss_0}&=\ser{G^F}{(M,\mu)}\cr&=\Lu{G(s_0)} G\big(\Psi_{G(s_0),s_0}(1)\otimes\ser{G(s_0)^F}{(L_0,\la_0)}\big)\cr
&=\Psi_{G,ss_0}(\ser{\cent{G^*}{ss_0}^F}{(\cent{L^*_0}s^F,\al_0)}\cr}\leqno{(2.3.6.3)}$$
}

\preuve Let  $S:=\zo{L_0^*}_{\phi_d}$. From (2.3.6.1) we have
$$L^*_0=\cent{M^*}{s_0},\quad S=\zo{M^*}_{\phi_d}=\zo{\cent{L^*_0}s}_{\phi_d}\leqno(2.3.6.4)$$
By Proposition~2.1.4 $(\cent{L^*_0}s,\al_0)$ is a $d$-cuspidal unipotent datum in $(\cent{G^*}{ss_0},F)$. As $s_0$ is central in the dual of $L_0$, we have (Proposition~1.3.2, (ii))
$$\Psi_{L_0,s_0}(1)\otimes \la_0=\Psi_{L_0,ss_0}(\al_0)$$
By Proposition~1.3.2, (iii) which apply to $\Psi_{M,ss_0}$ with $\cent{G^*}{ss_0}\subseteq L^*_0$ we have 
$$\mu=\Psi_{M,ss_0}(\al_0)$$
Clearly $(M,\mu)$ is defined up to $G^F$-conjugacy.

To prove that $(M,\mu)$ is $d$-cuspidal, consider a $d$-split Levi subgroup $L_1$ in $G$ such that $L_1\subseteq M$ and $\scal{\slu M{L_1}\mu}{\la_1}{L_1^F}\neq 0$ for some $\la_1\in\II{L_1^F}$. There is a $d$-split Levi subgroup $L^*_1$ in $G^*$, in duality with $L_1$, such that $ss_0\in L^*_1$ and $L^*_1\subseteq M^*$. We may assume that $\la_1\in\ser{L_1^F}{ss_0}$ by a good choice of $s_0$ (defined up to $L_1^{*\,F}$-conjugacy). Define $\al_1$ by $\la_1=\Psi_{L_1,ss_0}(\al_1)$. There is a Levi subgroup $L_{1,0}$ of $L_1\cap L_0$ in duality with $\cent{L^*_1}{s_0}$ so that (Proposition~1.3.2, (ii) again)
$$\Psi_{L_1,ss_0}(\al_1)=\Lu{L_{1,0}}{L_1}(\Psi_{L_{1,0},ss_0}(\al_1)),\quad \scal{\Lu {L_1}M(\Psi_{L_1,ss_0}(\al_1))}\mu{L^F}\neq 0$$
 By Proposition~1.3.2, (ii), $\Psi_{L_{1,0},ss_0}(\al_1)=\Psi_{L_{1,0},s_0}(1)\otimes\Psi_{L_{1,0},s}(\al_1))$ hence
 $$\Lu {L_1}M\la_1=\Lu{L_{1,0}}M(\Psi_{L_{1,0},s_0}(1)\otimes\Psi_{L_{1,0},s}(\al_1))=\Lu{L_0}M(\Psi_{L_0,s_0}(1)\otimes \Lu{L_{1,0}}{L_0}(\Psi_{L_{1,0},s}(\al_1))\leqno{(2.3.6.5)}$$
 As $\Lu{L_0}M$ is isometric on $\ser{L_0^F}{s_0}$ and on $\ser{L_0^F}{ss_0}$, using transitivity of Lusztig induction we have
 $$\eqalign{\scal{\Lu{L_1}M\la_1}\mu{M^F}&=\scal{\Lu{L_{1,0}}{L_0}(\Psi_{L_{1,0},ss_0}(\al_1))}{\Psi_{L_{0},ss_0}(\al_0)}{L_0^F}\cr&= \scal{\Lu{L_{1,0}}{L_0}(\Psi_{L_{1,0},s}(\al_1))}{\Psi_{L_{0},s}(\al_0)}{L_0^F}\cr&=\scal{\Lu{\cent{L_1^*}{ss_0}}{\cent{M^*}{ss_0}}\al_1}{\al_0}{\cent{M^*}{ss_0}^F}\neq 0\cr}\leqno{(2.3.6.6)}$$
 the second equality thanks to Proposition~1.3.2, (ii), the third one by Proposition~2.2.4, formula (2.2.4.2) applied to $L_{1,0}\subseteq L_0$. As $\al_0$ is assumed to be $d$-cuspidal, we obtain $\cent{L_1^*}{ss_0}=\cent{M^*}{ss_0}$ and $\al_0=\al_1$. That imply $\zo{L^*_1}_{\phi_d}=\zo{\cent{L_1^*}{s_0}}_{\phi_d}=\zo{\cent{L_1^*}{ss_0}}_{\phi_d}=\zo{M^*}_{\phi_d}$ by (2.3.6.4), hence $L^*_1=M^*$ and $(L_1,\la_1)=(M,\mu)$.
 
 We have proved that our construction of $(M,\mu)$ from $(L_0,\la_0)$ in (a) provides a $d$-cuspidal datum in series $(ss_0)$ in $(G,F)$. By Proposition~2.2.4 we know that the set of $G(s_0)^F$-conjugacy classes of $d$-cuspidal data in series $(s)$ in $(G(s_0),F)$ define a partition of $\ser{G(s_0)^F}s$.
 Using the one-to-one map $(\xi\mapsto \Lu{G(s_0)}G(\Psi_{G(s_0),s_0}(1)\otimes \xi)$ from $\ser{G(s_0)^F}s$ onto  $\ser{G^F}{ss_0}$ we obtain a partition of $\ser{G^F}{ss_0}$. By transitivity formula $\Lu MG\circ \Lu{L_0}M=\Lu{G(s_0)}G\circ \Lu{L_0}{G(s_0)}$ it is sent by $\Lu{G(s_0)}G$ on the partition of $\ser{G^F}{ss_0}$. Non $G(s_0)^F$-conjugate data $(L_0,\la_0)$ in series $(s)$ in $(G(s_0),F)$ define by our construction in (a) non $G^F$-conjugate $d$-cuspidal data in series $(ss_0)$ in $(G,F)$. The definitions of $\mu$ and $\la_0$ from $\al_0$ imply
 $$\ser{G^F}{(M,\mu)}=\Lu{G(s_0)}G(\Psi_{G(s_0),s_0}(1)\otimes \ser{G(s_0)^F}{(L_0,\la_0)}\leqno{(2.3.6.7)}$$

 Furthermore G.$d$-HC holds in $\ser{G^F}{ss_0}$  as it holds in $\ser{G(s_0)^F}s$.
 
 One can recover $(L_0,\la_0)$ from $(M,\mu)$ without using the partition (2.3.6.7) given by G.$d$-HC theory in $\ser{G(s_0)^F}s$ : consider $M^*\subseteq G^*$ in duality with $M$ and with $ss_0\in M^*$, so that $\mu=\Psi_{M,ss_0}(\al)$, $\al_0\in\ser{\cent{M^*}{ss_0}^F}1$. Then $\cent{M^*}{s_0}$ is $d$-split in $\cent{G^*}{s_0}$. Take $L_0$ in the dual $G(s_0)^F$-conjugacy class of $\cent{M^*}{s_0}$, $\la_0=\Psi_{L_0,s}(\al_0)$, one has $\mu=\Lu{L_0}M(\Psi_{L_0,s_0}(1)\otimes \la_0)$, so that by construction (a) $(L_0,\la_0)$ gives $(M,\mu)$. Using (2.3.6.6)  one  verifies easily that $\mu$ is $d$-cuspidal.
 
 (b) By Proposition~1.2.5, $G(s_0)$ is an $E$-split Levi subgroup of $G$. As $\al_0$ is $d$-cuspidal, $(L_0,\Psi_{L_0,s}(\al_0))$ is a $d$-cuspidal datum in series $(s)$ in $(G(s_0),F)$. By Propositions~2.1.6, 2.1.9, $b_{G(s_0)^F}(L_0,\Psi_{L_0,s}(\al_0))$ and $\Lu {G(s_0)}G(b_{G(s_0)^F}(L_0,\Psi_{L_0,s}(\al_0)))$ are well defined.
 
 By Proposition~2.3.2, with $(\cent{G^*}s,\cent{G^*}{ss_0}, (\cent{L_0^*}s,\al_0))$ instead of $(G,H,(L_H,\la_H))$, we have 
 
\noindent $(L^*_s,\al)\simex{\cent{G^*}s^F}(\cent{L^*_0}s,\al_0)$ for some $d$-cuspidal unipotent datum $(L^*_s,\al)$ in $(\cent{G^*}s,F)$ and may assume $\cent{G^*}{s_0}\cap L^*_s=\cent{L^*_0}s$. Then $(L^*_s,\al)$ defines by Proposition~2.1.4 a $d$-cuspidal datum $(L,\la)$ in series $(s)$ in $(G,F)$, such that, for some  dual $L^*$ of $L$,   $s\in L^*\subseteq G^*$, $L^*_s=L^*\cap \cent{G^*}s$ and $\la=\Psi_{L,s}(\al)$. 

By Propositions~2.3.3, 2.3.5,  $\Lu{G(s_0)}G\xi\in\II{b_{G^F}(L,\la)}$  and $\Lu{G(s_0)}G(\Psi_{G(s_0),s_0}(1)\otimes \xi)\in\II{b_{G^F}(L,\la)}$ for any $\xi\in\ser{G(s_0)^F}{(L_0,\la_0)}$, hence $b_{G^F}(L,\la)\in\bl{G^F}{s,s_0}$. In other words and with notations of Proposition~Ê2.1.7 we have  (2.3.6.2)
$$\Lu {G(s_0)}G(b_{G(s_0)^F}(L_0,\la_0))=b_{G^F}(L,\la)\in\bl{G^F}{s,s_0}\leqno{(2.3.6.8)}$$
 
 Given $(\cent{L^*_0}s,\al_0)$ as above, $(L^*_s,\al)$ is defined up to $\cent{G^*}s^F$-conjugacy, the map given in (2.3.6.8) $\bl{G(s_0)^F}s\to \bl{G^F}{s,s_0}$ is well defined. By  Proposition~2.3.2,  $\cent{G^*}s^F$-conjugacy on $(L^*_s,\al)$ implies $\cent{G^*}{ss_0}^F$-conjugacy on $(\cent{L^*_0}s,\al_0)$, hence $\Lu{G(s_0)}G$ is one-to-one from $\bl{G(s_0)^F}s$ onto $ \bl{G^F}{s,s_0}$.
 
 From (2.3.6.2), (2.3.6.8) and Proposition 2.1.6 we deduce $\ser{G^F}{(M,\mu)}\subseteq \II{b_{G^F}(L,\la)}$. Now $G^F$-conjugacy on $(L,\la)$ implies $\cent{G^*}s^F$-conjugacy on $(L^*_s,\al)$, hence $\cent{G^*}{ss_0}^F$-conjugacy on $(\cent{L^*_0}s,\al_0)$, hence $G(s_0)^F$-conjugacy on $(L_0,\la_0)$, finally $G^F$-conjugacy on $(M,\mu)$, that is equality on sets $\ser{G^F}{(M,\mu)}$ inside $ \ser{G^F}{ss_0}$. As we have seen, these last one form a partition of $\ser{G^F}{ss_0}$. Therefore
 $$ \II{b_{G^F}(L,\la)}\cap \ser{G^F}{ss_0}=\ser{G^F}{(M,\mu)}$$
 hence the second equality of (2.3.6.3), thanks to (2.3.6.7).
 
 By definition of $(L_0,\la_0)$, $\ser{G(s_0)^F}{(L_0,\la_0)}=\Psi_{G(s_0),s}(\ser{\cent{G^*}{ss_0}^F}{(\cent{L^*_0}s^F,\al_0)})$. As $s_0$ is central in $G(s_0)^*$, $\Psi_{G(s_0),s_0}(1)\otimes \Psi_{G(s_0),s}=\Psi_{G(s_0),ss_0}$. By Proposition~1.3.2 (iii), $\Lu{G(s_0)}G\circ\Psi_{G(s_0),ss_0}=\Psi_{G,ss_0}$. The last equality of (2.3.6.3) follows.
 \bull

\medskip\noindent{\bf 2.4. On blocks and their \irr\ characters in $G^F$ when $\z G$ is not connected }

As usual we consider a regular embedding $$(G,F)\to (H,F)$$ (see 1.1.4 (b)) and a dual morphism $H^*\to G^*$ to obtain by Clifford theory results on $d$-cuspidal data and blocks when the center of $G$ is not connected. Clearly  assumption 2.2.3 on $(H,F,\ell,d)$ is equivalent to 2.2.3 on $(G,F,\ell,d)$. 

We have only to consider types for which the fundamental group is non trivial. If the fundamental group of $G$ is trivial $\z{[G,G]}=\{1\}$ and $H$ is a direct product $[G,G]\times \zo H$. Any $G$ is a direct product $G=G_1\times G_2$ where $G_2=[G_2,G_2]$ and contains the components of type $\EE_8$ or $\FF_4$ or $\GG_2$, so that $\z{[G_2,G_2]}=1$. Therefore we assume in that section that $G=G_1$. 
The proof 2.2.8 of Proposition~2.2.4 shows that the commutation formula (2.2.7.1) holds.

From Proposition~2.2.4, specially (2.2.4.1), (2.2.4.2), and  section~5.2 in Appendix, we have

\bigskip\noindent{\bf 2.4.1.  }{\sl Assume 2.2.3 on $( H,F,\ell,d)$ and  $\z H$ is connected.  Let $( M, \mu)$ be a $d$-cuspidal datum in an $\ell'$-series  in
$( {H},F)$. 
There exist 
$\xi\in\II{ {H}^F}$ such that  $\scal{\Lu
{ M} {H} \mu} \xi{ {H}^F}\in\{-1,1\}$ and, for any $\chi\in\ser{H^F}{(M,\mu)}$ different from $\xi$,  $\xi(1)\neq \chi(1)$.}

 In the proof of the following Proposition we use 2.4.1  as  Geck and Bonnaf\'e on $1$-cuspidality in a regular embedding (see [3], 12.C).
 
\medskip\noindent{\bf
Proposition~2.4.2. }{\sl  Let $(G,F)\to (H,F)$ be a regular embedding. Assume 2.2.3 on $( H,F,\ell,d)$.  Let $( M, \mu)$ be a $d$-cuspidal datum in an $\ell'$-series in
$( {H},F)$. 
Let 
 $L:= M\cap G$, $ \nu:=\Res{ M^F}{L^F} \mu$,
$\la$ an \irr\ component  of $ \nu$. Then
 $(L,\la)$
is a $d$-cuspidal datum in $(G,F)$. 

(a) One has 
$\nor{G^F}{L, \nu}=\nor{G^F}{L,\la}$,

(b) If $\chi\in\ser{G^F}{(L,\la)}$,
then $(H^F)_\chi\subseteq G^F.{ (M^F)}_{\la}$.}

\preuve
$(L,\la)$ is a $d$-cuspidal datum by Proposition~2.1.5.

(a)  The restriction from $ M^F$ to $L^F$ has no
multiplicity, we have
$ \nu=\sum_{g\in
 M^F/( M^F)_{\la}}\lexp g{\la}$.  
 From 2.4.1 we deduce
\smallskip\noindent  (2.4.2.1) $\quad$ {\sl There exist some $\chi\in\II{G^F}$ such that
$|\scal{\Lu {L}G \nu}{\chi}{L^F}|=1$. }

\smallskip 
 The restriction from $H^F$ to $G^F$ has no multiplicity and $H^F/G^F$ is abelian. Given two elements of $\II{ {H}^F}$ their restrictions to 
$G^F$ are disjoint or equal, if equal their degrees are equal. Let 
$ \xi\in\ser{H^F}{(L,\mu)}$ as in  2.4.1. Then if $\chi\in\ser{H^F}{(M,\mu)}$ and $\chi\neq \xi$, $\scal{\Res{ {H}^F}{G^F} \xi}{\Res{ {H}^F}{G^F} \chi}{G^F}=0$. Thanks to  the equality $\Res{ {H}^F}{G^F}(\Lu {M} {H} \mu)=\Lu {L}G \nu$, we see that components of $\Res{ {H}^F}{G^F} \xi$  occurr with multiplicity  $\pm 1$ in $\Lu {L}G \nu$, our claim (2.4.2.1).

The actions of $ M^F/L^F$ and
$\nor{G^F}L/L^F$ on  $\II{L^F}$ commute : as
$[ M^F,G^F]\subseteq L^F$, we see that for any
$\la'\in\II{L^F}$ and any
$x\in M^F$, we have $\nor{G^F}{L,\la'}=\nor{G^F}{L,\lexp x{\la'}}$,
hence
$\nor{G^F}{L,\la}\subseteq\nor{G^F}{L, \nu}$.

 The group 
$\nor{G^F}{L, \nu}$ acts on the set of \irr\ components of  $ \nu$, a regular orbit under $ M^F/ (M^F)_{\la}$. The set $X$ of $x\in M^F$ such that $\lexp x \la=\lexp n\la$ for some $n\in \nor{G^F}{L,\nu}$ is  a subgroup  of
$ M^F$ such that  $M_\la\subseteq X$ and $|X/(M^F)\la|=|\nor{G^F}{L, \nu}/\nor{G^F}{L,\la}|$. 

We have 
$ \nu=\sum_{x\in X/(M^F)_\la}\sum_{h\in M^F/X}\lexp{xh}{\la}$  and $\scal {\Lu LG\lexp {xh}\la}\chi{G^F}=\scal {\Lu LG\lexp {h}\la}\chi{G^F}$, hence
 $$\scal{\Lu {L}G \nu}\chi{G^F}=
|\nor{G^F}{L, \nu}/\nor{G^F}{L,\la}|\Big(\sum_{h\in M^F/X}\scal
{\Lu {L}G{\lexp h{\la}}}\chi{G^F}\Big)$$ By assertion (2.4.2.1) $|\nor{G^F}{L, \nu}/\nor{G^F}{L,\la}|=1$ .

(b) As  $\chi\in\ser{G^F}{(L,\la)}$, if $h\in H^F$ and $\lexp h\chi=\chi$, then $\chi\in \II{b_{G^F}(L,\la)}\cap \II{b_{G^F}(\lexp hL,\lexp h\la)}$
hence the $d$-cuspidal data $(L,\la)$ and $(\lexp hL,\lexp h\la)$ are $G^F$-conjugate. We have  $H^F=G^F.M^F$
 hence $g\in G^F.({ M^F})_{\la}$.
\bull

\medskip\noindent {\bf  2.4.2.2. }  We may translate assertion (b) of Proposition~2.4.2 as follows. If $\chi\in\ser{G^F}s$ and $\chi$ corresponds to the orbit under ${\rm A}_{G^*}(s)^F$  of $\beta\in\ser{\cento{G^*}s^F}1$ by (1.3.6.1), then $H_\chi=\tau_{H,s}({\rm A}_{G^*}(s)^F_\beta)$. Similarly, $M_\la=\tau_{M,s}({\rm A}_{L^*}(s)^F_\al)$, where $\al$ is unipotent $d$-cuspidal in $\II{\cento{L^*}s^F}$ such that $\beta\in\ser{\cento{G^*}s^F}{(\cento{L^*}s,\al)}$. There are a one-to-one morphism $\pi\colon {\rm A}_{L^*}(s)^F\to {\rm A}_{G^*}(s)^F$ (Proposition~1.2.6, (f)) and an isomorphism $M^F/L^F\cong H^F/G^F$ so that $G^F.\tau_{M,s}({\rm A}_{L^*}(s)^F)=\tau_{H,s}(\pi({\rm A}_{L^*}(s)^F))$. Assertion (b) writes 
$$\pi({\rm A}_{L^*}(s)^F_\al)\subseteq {\rm A}_{G^*}(s)^F_\beta$$
By G.$d$-HC in $\ser{\cento{G^*}s^F}1$ we have $ {\rm A}_{G^*}(s)^F_\beta\cap \pi({\rm A}_{L^*}(s)^F)\subseteq \pi({\rm A}_{L^*}(s)^F_\al)$. Finally $$ {\rm A}_{G^*}(s)^F_\beta\cap \pi({\rm A}_{L^*}(s)^F)= \pi({\rm A}_{L^*}(s)^F_\al)\, .$$

In the following, assuming a choice of suitable dual $F$-stable tori $T\subseteq M\subseteq H$, $T^*\subseteq M^*\subseteq H^*$, one has  isomorphisms between Weyl groups  ${\rm W}(G,T\cap G)\cong {\rm W}(H,T)$, ${\rm W}(M\cap G,T\cap G)\cong {\rm W}(M,T)$, ... to symplify notations we omit  reference to the torii. 
Consider   an action of  ${\rm W}_{ {H}}( M)^F\times
\II{ {H}^F/G^F}$ on $\II{ M^F}$ as follows :
$$\eqalign{{\rm W}_{G^F\subset  {H}^F}( M)={\rm W}_{ {H}}( M)^F\times
\II{ {H}^F/G^F}\colon\;
\II{ M^F}&\to
\II{ M^F}\cr  (w,\theta)\;\colon \quad\quad\quad\mu&\mapsto \lexp w \mu\otimes
(\Res{ {H}^F}{M^F}\theta^{-1}).\cr}$$ 

 \medskip\noindent{\bf Proposition~2.4.3. }{\sl Let $G,H,F,\ell,d,M,\mu,L,\la,\nu$ as in Proposition 2.4.2. Define $M^F_G( \mu)$ by $$M^F_G( \mu)=M^F\cap(\cap_\theta \Ker \theta)\quad (\theta\in (H^F/G^F)^\wedge,\;\Lu MH\mu\otimes\theta=\Lu MH\mu)$$ 

(a) Let $ t\in M^{*\,F}$ such that
$\mu\in\ser{M^F} t$, and $\al\in\ser{\cent{M^*}t^F}1$ such that
 $\Psi_{M, t}(\al)=\mu$. Dualities being defined around  maximal $F$-stable torii $T$, $T^*$ in $M$, $M^*$, with $t\in T^*$, let   $s\in L^{*\,F}$ be the image of $t$ and ${\rm A}_{G^*}(s,L^*)^F$ be the image of
$\nor{{\rm W}_G(s)}{{\rm W}^\circ_L(s)}^F$ in ${\rm A}_{G^*}(s)^F$. One has
$$
G^F.M^F_G(\mu)=\tau_{ {H},s}({\rm A}_{G^*}(s,L^*)^F_\al).$$

(b) The isomorphism of Proposition~2.2.6 sends  ${\rm W}_{ {H}}(M)^F_{\mu}$ onto 
${\rm W}_{\cent{ {H}} t}(\cent{M}t)^F_\al$.

One has three short exact sequences :
$${\rm W}_{ {H}}(M)^F_{\mu}\longrightarrow {\rm W}_{G^F\subset
 {H}^F}(M)_\mu\longrightarrow (M^F/M^F_G(\mu))^\wedge,$$
$$\big(M^F/\tau_{M,s}({\rm A}_{L^*}(s)^F_\al)\big)^\wedge\longrightarrow {\rm W}_{G^F\subset
 {H}^F}(M)_\mu\longrightarrow {\rm W}_{G}(L)^F_\la,$$
$${\rm W}_{ {H}}(M)^F_{\mu}\longrightarrow {\rm W}_{G}(L)^F_\la\longrightarrow
{\rm A}_{G^*}(s,L^*)_\al/{\rm A}_{L^*}(s)^F_\al.$$

(c) If
$\chi\in\ser{G^F}{(L,\la)}$, then $G^F.M^F_G(\mu)\subseteq ({ {H}^F})_\chi$.} 

\preuve 
(a) The inclusion
$M^F_G(\mu)\subseteq \Ker \theta$ is equivalent to
$G^F.M^F_G(\mu)\subseteq \Ker \theta$. 

We have ${\rm A}_{L^*}(s)\cong {\rm W}_{L^*}(s)/{\rm W}^\circ_{L^*}(s)$. In Proposition~1.2.6 is defined an injective morphism from ${\rm A}_{L^*}(s)$ to ${\rm A}_{G^*}(s)$ by $w{\rm W}^\circ_{L^*}(s)\mapsto w{\rm W}^\circ_{G^*}(s)$ (notations of (2.2.5.2) and Proposition~2.2.6). With same definition that map extends to $\nor{{\rm W}_{G^*}(s)}{{\rm W}^\circ_{L^*}(s)}/{\rm W}^\circ_{L^*}(s)$, so is defined 
${\rm A}_{G^*}(s,L^*)^F$. As well  ${\rm A}_{G^*}(s,L^*)^F$ may be defined as the stabilizer of the $\cent{G^*}s^F$-conjugacy class of $L^*_s$ in ${\rm A}_{G^*}(s)^F$.

 (a.1) Let $\theta\in(H^F/G^F)^\wedge$, we claim
 
 \centerline{\sl $(\Lu
 {M} {H}\mu)\otimes \theta=\Lu
 {M} {H}\mu$ if and only if there exist $w\in W_H(M)^F$ such that $(w,\theta)\in {\rm W}_{G^F\subset
 {H}^F}(M)_\mu$.}
 
 We have $(\Lu
 {M} {H}\mu)\otimes \theta=\Lu
 {M} {H}(\mu\otimes \Res{ {H}^F}{M^F}\theta)$.
 If
$(\Lu
 {M} {H}\mu)\otimes \theta=\Lu
 {M} {H}\mu$,  $\theta$
stabilises the series $\ser{ {H}^F} t$, hence $\tau_{ {H},s}({\rm A}_{G^*}(s)^F)\subseteq \Ker \theta$ (see section 1.3.5) and the order of $\theta$ is  prime
to $\ell$.  
By Corollary~2.1.8,
$\mu$ and
$\mu\otimes(\Res{ {H}^F}{M^F}\theta)$
are conjugate under $\nor{ {H}}M^F$.  As $\mu$ is fixed under ${\rm W}(M)^F$,  the equality $(\Lu
 {M} {H}\mu)\otimes \theta=\Lu
 {M} {H}\mu$ 
is equivalent to
`` {\sl  there exist $w\in {\rm W}_H(M)^F$  such that $\mu\otimes\Res{ {H}^F}{M^F} 
\theta=\lexp w\mu$} ",
that is $(w,\theta)\in {\rm W}_{G^F\subset
 {H}^F}(M)_\mu$.

 
 (a.2) Claim :
 
 \centerline {\sl $(w,\theta)\in {\rm W}_{G^F\subset
 {H}^F}(M)_\mu$ imply $\theta\in\sigma_{H,s}({\rm A}_{G^*}(s,L^*)^F_\al).$ }
 
 We have $\lexp
w{(\Res{M^F}{L^F}\mu)}=\nu\in\Z\ser{L^F}s$. 

By Proposition~2.2.6, $w\in
{\rm W}(M)^F.\nor{{\rm W}_G(s)}{{\rm W}^\circ_L(s)}^F/{\rm W}(M)^F$. As ${\rm W}(M^*)\cap W_{G^*}(s)=W^\circ _{L^*}(s)$ by antiisomorphism $W(H)\to W(H^*)$ and isomorphism theorem, $w$ has  image  $\tilde w^*\in\nor{{\rm W}_{G^*}(s)}{{\rm W}^\circ_{L^*}(s)}^F/{\rm W}^\circ_{L^*}(s)^F$ hence in ${\rm A}_{G^*}(s,L^*)^F$ :
$$a:=\tilde w^*{\rm W}_{G^*}^\circ(s)\in
{\rm A}_{G^*}(s,L^*)^F$$ 

As seen in 1.3.5,  (1.3.5.1), $\sigma_{ {H},s}(a)\in(H^F/G^F)^\wedge$, $\Res{ {H}^F}{M^F}(\sigma_{ {H},s}(a))\in(M^F/L^F)^\wedge$. To  $\sigma_{ {H},s}(a)$ there corresponds by duality $z\in\z{G^*}^F$  such that  $\lexp at=tz$. By Proposition~1.3.2, (ii),  $\sigma_{ {H},s}(a)=\Psi_{H,z}(1)$ and $\Res{ {H}^F}{M^F}(\sigma_{ {H},s}(a))=\Psi_{M,z}(1)$ ($\Psi_{H,z}(1)$ is a uniform function). 

Any element of $n\in \nor{H}{T}^F$ with image  $w\in \nor{{\rm W}( {H})}{{\rm W}(M)}^F/{\rm W}(M)^F$ induces an automorphism  of $(M,F)$ that stabilizses $\ser{M^F}s$. To $nT\in W(H)^F$ there corresponds $n^*T^*$, where $n^*\in\nor{H^*}{T^*}^F$ maps onto $\tilde w^*$,  and $n$, $n^*$ induce dual automorphisms of $(M,F)$, $(M^*,F)$ respectively. Furthermore $n^*$ acts on $\cent{M^*}s$ and on $\ser{\cent{M^*}s^F}1$ as $a$. 
By (iv) of Proposition~1.3.2, $$\Psi_{M, t z}(\lexp a\al)=\lexp w\mu\in\ser{M^F}{tz}\leqno{(2.4.3.2)}$$ 
 By Proposition~1.3.2, (ii) again we obtain $$\Psi_{M, t }(\lexp a\al)=\lexp w\mu\otimes
\Res{ {H}^F}{M^F}( 
\sigma_{ {H},s}(a)^{-1})\leqno{(2.4.3.3)}.$$ 
From the hypotheses $\mu=\Psi_{M,t}(\al)$ and $\mu\otimes\Res{ {H}^F}{M^F} 
\theta=\lexp w\mu$ and (2.4.3.3)  follows $$\Psi_{M, t }(\lexp a\al)\otimes
\Res{ {H}^F}{M^F}( 
\sigma_{ {H},s}(a))=\Psi_{M,t}(\al) \otimes\Res{ {H}^F}{M^F} 
\theta.$$   By section 1.3.5, (1.3.5.4)  $\Res{ {H}^F}{M^F}( 
\sigma_{ {H},s}(a)\otimes \theta^{-1})=\sigma_{M,s}(b)$ where $b\in {\rm A}_{L^*}(s)^F$ and $\lexp b\al=\lexp a\al$. Then $\theta=\sigma_{H,s}(ab^{-1})$ and $ab^{-1}\in 
{\rm A}_{G^*}(s,L^*)^F_\al$.

So (2.4.3.1) implies 
$$\theta=\sigma_{ {H},s}(c)^{-1}\; {\rm where}\;  c\in
{\rm A}_{G^*}(s,L^*)^F_\al.\leqno{(2.4.3.4)}$$

(a.3) Assume  (2.4.3.4). Let $z\in\z{G^*}^F$ such that $\lexp c t=tz$. Let $\tilde w^*\in\nor{{\rm W}_{G^*}(s)}{{\rm W}^\circ_{L^*}(s)}^F/{\rm W}^\circ_{L^*}(s)^F$ with image  $c=\tilde w^* {\rm W}^\circ_{G^*}(s)$ in ${\rm A}_{G^*}(s,L^*)^F_\al$. As above to prove (2.3.4.2) we see that there exists $\bar w$ in $
{\rm W}(M)^F.\nor{{\rm W}_G(s)}{{\rm W}^\circ_L(s)}^F/{\rm W}(M)^F\subseteq W_H(M)^F$ with image $\tilde w^*$ such that (iv) of Proposition~1.3.2 apply: $\Psi_{M, t z}(\al)=\lexp w\mu$ hence $\Psi_{M, t }(\al)=\lexp w\mu\otimes \Res{H^F}{M^F}(\sigma_{H,s}(c)^{-1})$. By definition of the action of ${\rm W}_{G^F\subset H^F}(M)$,  $(w,\theta)\in{\rm W}_{G^F\subset H^F}(M)_\mu$.

We have proved that $\Lu MH\mu=\Lu MH\mu\otimes \theta$ is equivalent to (2.4.3.4). By definition of $\tau_{H,s}$ we have  $\cap_{a\in{\rm A}_{G^*}(L^*,s)^F_\al}\Ker (\sigma_{H,s}(a))=\tau_{H,s}({\rm A}_{G^*}(L^*,s)^F_\al)$, hence $G^F.M^F_G(\mu)=\tau_{H,s}({\rm A}_{G^*}(s,L^*)^F_\al)$.

(b) The first assertion is a consequence of (iv) in Proposition~1.3.2.

(b.1) We proved in (a.2), (a.3) that  ${\rm W}_{G^F\subset H^F}(M)_\mu$ maps onto  $(H^F/\tau_{H,s}({\rm A}_{G^*}(s,L^*)^F_\al))^\wedge$ by  projection on the right side. But $H^F/\tau_{H,s}({\rm A}_{G^*}(s,L^*)^F_\al)\cong M^F/\tau_{H,s}({\rm A}_{G^*}(s,L^*)^F_\al)\cap M^F=M^F/M^F_G(\mu)$. That gives the first exact sequence.

(b.2) Clifford theory shows that the image of ${\rm W}_{G^F\subset H^F}(M)_\mu$ by projection on left side is ${\rm W}_H(M)^F_\nu$. We have ${\rm W}_H(M)^F_\nu={\rm W}_G(L)^F_\nu$. By (a) of Proposition~2.4.2, ${\rm W}_G(L)^F_\nu={\rm W}_G(L)^F_\la$. 

By (1.3.5.4) $\mu=\mu\otimes \Res{H^F}{M^F}\theta^{-1}$ is equivalent to $\Res{H^F}{M^F}\theta^{-1}=\sigma_{M,s}(a)$ where $a\in {\rm A}_{L^*}(s)^F_\al$. The kernel of the left projection is $\{1\}\times (H^F/G^F.\tau_{M,s}({\rm A}_{L^*}(s)^F_\al))^\wedge$, isomorphic to 
$(M^F/\tau_{M,s}({\rm A}_{L^*}(s)^F_\al))^\wedge$.

(b.3) Clifford theory without multiplicities between $L^F$ and $M^F$ implies $(M^F/L^F)^\wedge_\mu=(M^F/M^F_\la)^\wedge$ (see Appendix B), that is $\tau_{M,s}({\rm A}_{L^*}(s)^F_\al)=M^F_\la$.   The projection of ${\rm W}_{G^F\subseteq H^F}(M)_\mu$ on right side contains the kernel of the projection on the left side, that is $M^F_G(\mu)\subseteq M^F_\la$. As well  ${\rm W}_{ {H}}(M)^F_{\mu}\subseteq  {\rm W}_{G}(L)^F_\la$. The preceding exact sequences imply an isomorphism ${\rm W}_{G}(L)^F_\la/{\rm W}_{ {H}}(M)^F_{\mu}\cong M^F_\la/M^F_G(\mu)$. The last quotient is in duality, {\it via} functions $\tau$, with ${\rm A}_{G^*}(s,L^*)^F_\al/{\rm A}_{L^*}(s)^F_\al$, we obtain the third exact sequence.

(c) Let $\chi\in\ser{G^F}{(L,\la)}$. Note that $H^F_\chi=\cap_{H^F_\chi\subseteq \Ker\theta} \Ker \theta$. So let $\theta\in(H^F/G^F)^\wedge$ such that $H^F_\chi\subseteq \Ker \theta$. By (1.3.1.1) and Proposition~2.1.7 there exist $\xi\in\II{H\mid \chi}\cap \ser{H^F}{(M,\mu)}$. We know that $\xi\otimes\theta=\xi$, hence  $\Lu MG\mu=\Lu MG\mu\otimes \theta$. By definition of $M^F_G(\mu)$, $M^F_G(\mu)\subseteq H^F_\chi$.
\bull

\medskip\noindent{\bf 2.4.4. Proposition. }{\sl  Let $\sigma \colon (G,F)\to (H,F)$ be a regular embedding and $\sigma^*$ a dual morphism. Assume 2.2.3 on $( H,F,\ell,d)$.  Let $( M, \mu)$ be a $d$-cuspidal datum in series $(t)$ in
$( {H},F)$, a dual $M^*$ of $M$,  given as a Levi subgroup of  $H^*$, $\mu=\Psi_{M,t}(\al)$,  $\al\in\ser{\cento{ M^*}t^F}1$. 
Let 
 $L:= \sigma^{-1}(M)$, $L^*=\sigma^*(M)$, $s=\sigma^*(t)$, $\la\in \ser{L^F}s$ such that $\mu$ covers $\la$.

(a) The set of blocks $B$ of $ H^F$ that cover $b_{G^F}(L,\la)$ and with $\II
B\subset\lser{ H^F}t$ is a regular orbit under ${\rm A}_{G^*}(s)^F/{\rm A}_{G^*}(s,L^*)^F_\al$.

(b) The set  of blocks of  $G^F$ that are covered by $b_{ H^F}( M, \mu)$ is a regular orbit under
$ H^F/\tau_{ H,s}({\rm A}_{L^*}(s)^F_\al)$.

(c) If $\cent{H^*}t=\a {\cent{H^*}t}$, hence if $G=\a G$, then the blocks in series (s) of $G^F$ are conjugate under $H^F$.} 

\preuve 
 We know that
$\Res{}{ L^F\to M^F } \mu$ is a sum of  $M^F$-conjugate of 
$\la$, and that $\la$ is $d$-cuspidal.  By Propositions~2.1.6, 2.1.7 and formula (1.3.1.1) any $\xi\in\ser{H^F}{(M,\mu)}$ covers some element of $\ser{G^F}{(L,\lexp m\la)}$ ($m\in M^F$). Hence $b_{ H^F}( M, \mu)$
covers 
$b_{G^F}(L,\lexp m\la)$ and only these blocks of $G^F$  [28] Chapter 7,
Lemmas 5.3, 5.7. From Proposition~2.1.7 again we have (see [15] Remark~2.7) 

\smallskip \noindent{\bf 2.4.5. }{\sl Domination between  blocks of
$ H^F$ and $G^F$ is equivalent to domination between
  conjugacy classes of  $d$-cuspidal data.
}

 (a)  Given a block $B=b_{H^F}(M,\mu)$ in series $(t)$ of $ H^F$, modulo 
$(G^*)^F$-conjugacy on
$s$, 
$( H^*)^F$-conjugacy on $ t$ and
$G^F$-conjugacy on $d$-cuspidal data in
series $(s)$ in $(G,F)$, we may fix $ t$, $s$,
$ M^*$, $L^*$, $ M$, $L$, $\mu$, $\la$. 

The blocks  of $ H^F$ that cover $b_{G^F}(L,\la)$  
are the $ b_{ H^F}( M, \mu)\otimes \theta=b_{H^F}(M,\mu\otimes\Res{H^F}{M^F}\theta)$    where $\theta\in (H^F/G^F)^\wedge$. These blocks are in series $(t)$ when
$\theta\in \sigma_{ H,s}({\rm A}_{G^*}(s)^F)$ (1.3.5, (1.3.5.1)). The number of such  blocks  is the number of orbits on the set of $d$-cuspidal $\mu'\in\ser{M^F}t$ under $ {\rm W}_H(M)^F$ that are fused under
${\rm W}_s:={\rm W}_{ H}( M)^F\times
\Res{ H^F}{ M^F}(\sigma_{ H,s}({\rm A}_{G^*}(s)^F))$. By
Proposition~2.4.3, (b), ${\rm W}_s$ contains $
{\rm W}_{G^F\subset H^F}( M)_ \mu$. As ${\rm W}_{ H}( M)^F\subseteq W_s$, 
the orbit of $W_s$ on  the set of  blocks of $H^F$ we consider is regular under
${\rm W}_s/{\rm W}_{G^F\subset H^F}( M)_ \mu$ isomorphic to ${\rm A}_{G^*}(s)^F/{\rm A}_{G^*}(s,L^*)^F_\al$, that is  (a) of the Proposition -- recall that the action of ${\rm A}_{G^*}(s)^F$ on blocks of $H^F$ is by isomorphism with $(H^F/\tau_{H,s}({\rm A}_{G^*}(s)^F)^\wedge$.

(b)    On $G$-side, the set of $G^F$-conjugacy classes
of
$d$-cuspidal data in series $(s)$ with support conjugate to $L$  is in bijection with the $N:={\rm W}_G(L)^F=\nor{G}L^F/L^F$-orbits of $d$-cuspidal elements of $\ser{L^F}s$. The  $d$-cuspidal data  that are covered by the $H^F$-class of $(M,\mu)$ are in the orbit  of $\la$ under ${\rm W}_H(M)^F\times  (M^F/L^F)$, where $W_H(M)^F$ acts as $W_G(L)^F$. The stabilizer of the ${\rm W}_G(L)^F$-orbit of $\la$ is just ${\rm W}_G(L)^F\times M^F_\la$. 
By (1.3.6.2), $M^F/M^F_\la\cong ({\rm A}_{L^*}(s)^F_\al)^\wedge$, hence (b) is proved.

(c) We have $\cent{H^*}t=\a {\cent{H^*}t}$ if and only if $\cento{G^*}s=\a {\cento{G^*}s}$. In that case the only $d$-cuspidal unipotent data of $(\cento{H^*}t,F)$ are the $(T^*,1_{(T^*)^F})$ where $T^*$ is a diagonal torus of $\cent{H^*}t$. Let $T_t$ (resp. $T_s:=T_t\cap G$) be is in the dual $H^F$-conjugacy class of $T^*$ (resp. in the dual $G^F$-conjugacy class of $\sigma^*(T^*)$),   and $M:=\cent H{(T_t)_{\phi_d}}$ (resp. $L:=\cent G{(T_s)_{\phi_d}}=M\cap G$).  By Proposition~2.1.4 $(M,\mu)$ is, up to $H^F$-conjugacy the unique $d$-cuspidal datum in series $(t)$.  Any element of $\cent{G^*}s^F$ stabilizes the $\cento{G^*}s^F$-conjugacy class of diagonal tori. We have ${\rm A}_{G^*}(s)^F={\rm A}_{G^*}(s,\sigma^*(T^*))^F_1$, a coherent result with the fact that there is only one block in series $(t)$ of $H^F$ !  The group ${\rm A}_{L^*}(s)^F $ may be different from 1.
 \bull

\vfill\eject

\medskip\noindent{\bf 3. The group $G(s)$,  ``in duality with"  $\cent{G^*}s$. }

In this chapter we describe a good candidate to be named " a dual of $(\cent{G^*}s,F)$ " in case $\cent{G^*}s$ is not connected, and denote it as $(G(s),F)$, as announced in Theorem~1.4. The construction is based on dual root data with $F$-action, and we assume properties of functorial type, see Proposition 3.1.1, (A). Of course it may be generalized to construct a dual of a non-connected algebraic reductive group, eventually defined over $\F q$. We  describe the set $\II{G(s)^F}$ and its partition in blocks.  
To obtain  the set of blocks of $G(s)^F$, we precise $G(s)^F$-conjugacy classes of maximal Brauer pairs in Proposition~3.1.2. In the two studies we have to verify  a non-multiplicity property in a Clifford  theory with quotient a subgroup of $A_{G^*}(s)^F$ ((B) in Proposition~3.1.1, (C) in Proposition~3.1.2). 
All proofs in section 3.2 are inductive, the "minimal case" being when $G$ is \irr\ and simply connected, and $s$ is rationally quasi-isolated in $G^*$. 

Then, in section 3.3, we may rely ``unipotent ''blocks of $G(s)^F$ and their \irr\ representations to blocks in series (s) of $G^F$ and their \irr\ representations, the main result of this paper.

\medskip\noindent{\bf 3.1. Propositions } 

\medskip\noindent{\bf 3.1.1.  Proposition. }{\sl  When $(G,F)$ and
$(G^*,F)$ are dual
algebraic reductive groups defined on  $\F q$ and      $s\in G^{*\,F}$ is  semi-simple, denote ${{\frak E}(G^*,s)}$ the short exact sequence 
$$1\to \cento{G^*}s\to \cent{G^*}s\to {\rm A}_{G^*}(s)\to 1\leqno{{\frak E}(G^*,s)}$$

(A)  Let
$(G (s)^\circ,
 F)$ be in  duality with $(\cento{G^*}s,F)$ around a maximally split
root datum of  $(\cento{G^*}s,F)$. There exists
 an extension  $$1\to G(s)^\circ\to G(s)\to {{\rm A}_{G^*}(s)}\to 1 \leqno{{\frak
E}(G,s)}$$
 where ${{\rm A}_{G^*}(s)}$ and $F$
act on the root datum of $G(s)^\circ$ by transposition of their action on the dual one, and the action of  $F$ on $G(s)^\circ$ and on
${\rm A}_{G^*}(s)$ extends in an   action on $G(s)$.
 
 The  various exact sequences ${\frak E}(G,s)$  may be defined satisfying the following properties :
 
(A.1) If
$(G,F)$ is a direct product 
$\times_i(G_i,F_i)$,  and
$s=(s_i)_i\in G^*$,
$(s_i\in  G^*_i$), then $({\frak E}(G,s),F)$ is isomorphic to the direct product of
extensions $({\frak E}(G_i,s_i),F)$.

(A.2) If
$\cento{G^*}s$ is a Levi subgroup of $G^*$, let  $L(s)$ be a Levi subgroup in the dual $G^F$-conjugacy class. Then
$G(s)$ is isomorphic to the subgroup $N$ of $\nor G {L(s)}$ such that $N/{L(s)}$ is the image of $\cent{G^*}s/\cento{G^*}s$ by the isomorphism of relative Weyl groups ${\rm W}_{G^*}(\cento{G^*}s)\cong {\rm W}_G(L(s))$.

(A.3)  If $\sigma\colon (G,F)\to (H,F)$ is an isotypic morphism and $\sigma^*\colon H^*\to G^* $
a dual one, assume $s=\sigma^*(t)$ with $t\in H^{*\, F}$. One has $\cento{G^*}s=\sigma^*(\cento {H^*}t)$ and $\sigma^*(\cent {H^*}t)\subseteq \cent{G^*}s$ so is defined 
 $\alpha\colon {\rm A}_{H^*}(t)\to{\rm A}_{G^*}(s)$. Let  $G(t)$ be the inverse image of 
$\alpha({\rm A}_{H^*}(t))$ in
$G(s)$ by ${\frak E}(G,s)$, and let ${\cal E}(G,t)$ be the restriction of ${\cal E}(G,s)$ to $\alpha({\rm A}_{H^*}(t))$. There is a morphism of short exact sequences
$${\def\normalbaselines{\baselineskip20pt\lineskip3pt
\lineskiplimit3pt}\matrix{{{\cal E}(G,t)}&{}&1&\to&G(s)^\circ&\mapright{}&G(t)&\mapright{}&\alpha({\rm A}_{H^*}(t))&\to&1\cr
{\mapdown{}}&{}&{}&{}&\mapdown{\sigma_t}&&\mapdown{}&&\mapdown{\cong}&{}&{}\cr 
{{\frak E}(H,t)}&\quad&1&\to&H(t)^\circ&\mapright{}&H(t)&
\mapright{ }&{{\rm A}_{H^*}(t)}&\to&1\cr} }$$
 where $\sigma_t$ is a dual isotypic morphism of the restriction of $\sigma^*\colon  \cento{H^*}t\to
\cento{G^*}s$ and all maps
commute with
$F$. 

(B) "Non-multiplicity condition" :

For any semi-simple $s_1\in \cento{G^*}s^F$ with order prime to the order of $s$,$$\forall (\chi_0,\chi)\in\ser{G(s)^{\circ\, F}}{s_1}\times \II{G(s)^F}\quad \scal{\Res{G(s)^F}{G(s)^{\circ F}}\chi}{\chi_0}{G(s)^{\circ F}}\in\{0,1\}\, .$$
}

It is clear that, when $(G,G^*,F,s)$ is given, ${\frak E}(G,s)$ is defined up to an $F$-isomorphism. If $s$ is central in $G^*$, we may assume $G(s)^\circ=G=G(s)$.

The following Proposition and its Corollary describe Clifford theory of unipotent blocks between $G(s)^{\circ\, F}$ and $G(s)^F$. We fix a bijection $  \Psi_{G(s)^\circ, 1}\colon \ser{G(s)^{\circ F}}1\to\ser{\cento{G^*}s^F} 1$, $(\al\mapsto\al)$, and so for $d$-split Levi subgroups. By Proposition~1.3.1, these bijections preserves $d$-cuspidality and relative Weyl groups.  Thus there is a bijection between the set of conjugacy classes of $d$-cuspidal unipotent data in $(G(s)^\circ ,F)$ and the analogous set in $(\cento{G^*}s,F)$. 
Note that, as $G(s)^\circ$ is connected, the exact sequence ${\frak E}(G,s)$ gives by restriction an isomorphism $$G(s)^F/G(s)^{\circ F}\cong {\rm A}_{G^*}(s)^F\, .$$

\medskip\noindent{\bf 3.1.2.  Proposition. }{\sl Assumption 2.2.3 on $(G,F,\ell,d)$. Let  $(G,F)$, $(G^*,F)$, $s\in G^{*\,F}$, ${{\frak E}(G^*,s)}$, ${{\frak E}(G,s)}$ as in Proposition~3.1.1. Let $$\rho \colon G(s)^F\to  {\rm A}_{G^*}(s)^F$$ be given by restriction to $F$-fixed points in ${{\frak E}(G,s)}$. 

Let $(L(s),\al)$ be a $d$-cuspidal unipotent datum of $(G(s)^\circ, F)$ and $b:=b_{G(s)^{\circ F}}(L(s),\al)$ the unipotent block so defined by Proposition~2.1.7. Let $L^*_s$ be a $d$-cuspidal Levi subgroup of $\cento{G^*}s$ in the dual class of $L(s)$ and $L^*:=\cent{G^*}{\z {L^*_s}_{\phi_d}}$. Let $D$ be a defect group of $b$, and $\al_D\in\II{\cent{G(s)^{\circ F}}D}$ be the canonical \irr\ in a block $b_D$ such that $(1,b)\subseteq (D,b_D)$. Let $\beta_D\in\II{\cent{G(s)^F}D\mid \al_D}$. One has

(A) $\rho(\cent{G(s)^F}D)={\rm A}_{L^*}(s)^F$ ;

(B) $\nor{G(s)^F}D_{\al_D}\subseteq \nor{G(s)^F}D_{\beta_D}$ ;

(C) "Non-multiplicity condition" : 
 $$ \forall (\chi,\chi_0)\in\II{\cent{G(s)^F}D}\times \II{b_D}, \quad\scal{\Res{\cent{G(s)^F}D}{\cent{G(s)^{\circ F}}D}\chi}{\chi_0}{\cent{G(s)^{\circ F}}D}\in\{0,1\}\,.$$}
 
 \medskip\noindent{\bf 3.1.3. Corollary. }{\sl Notations and assumptions of Proposition~3.1.2 on $(G,F,\ell,d,\rho,(L(s),\al),b, L^*_s$).
  Let $B$ be a block of $G(s)^F$ that covers $b$. One has 
  $$\rho(G(s)^F_b)=\rho({\rm A}_{G^*}(s,L^*)^F_\al),\quad \rho(I^{G(s)^F}_{G(s)^{\circ F}}(B))={\rm A}_{L^*}(s)^F_\al\,.$$}

\medskip\noindent{\it Proof of Corollary 3.1.3. }  Let $A:={\rm A}_{G^*}(s)^F$. The notation ${\rm A}_{G^*}(s,L^*)^F_\al $ was introduced in Proposition  2.4.3. As $\zo{L^*_s}_{\phi_d}=\zo{L^*}_{\phi_d}$, $L^*_s=L^*\cap \cento{G^*}s$ and ${\rm A}_{G^*}(s,L^*)^F_\al $ is the stabilizer of $(L^*_s,\al)$ in $A$.

The $d$-cuspidal unipotent data of $(G (s)^\circ, F)$ that define a fixed block $b$ of $G(s)^{\circ F}$ are conjugate hence  $(G(s)^F)_b=G(s)^{\circ\, F}.\nor{G(s)^F}{L(s),\al}$. We have $$\rho(\nor{G(s)^F}{L(s)})=A_{L^*_s},\quad A_b={A}_{L^*_s,\al} ={\rm A}_{G^*}(s,L^*)^F_\al\,.\leqno{(3.1.3.1)}\,.$$

Let
$C:=\cento{G(s)^\circ}{[L(s),L(s)]}$. An
$\ell$-Sylow subgroup of $C^F$
  is a
defect group of $b$ ([13] Theorem~4.4 or [16] Theorem~22.9). The normalizer $\nor{C}{\zo{L(s)}}$
contains such an 
$\ell$-Sylow  because $\zo{L(s)}$ is a maximal $F$-stable torus of  $C$ and contains a 
maximal
$\phi_d$-subgroup  of
$C$ [16] Proposition~22.8 or [13] Proposition~1.6. So we assume $D\subseteq
\nor{G(s)^\circ}{\zo {L(s)}}$ and $[L(s),L(s)]^F\subseteq \cent{G(s)^{\circ F}}D\subseteq L(s)$. We have $(1,b)\subseteq (D,b_D)$ for some block $b_D$ of $\cent{G(s)^{\circ F}}D$ with central defect and the canonical \irr\ $\al_D$ in $b_D$ is the only element in $\II{b_D}$ such that $\Res{L(s)^F}{\cent{G (s)^{\circ F}}D}\al=\al_D$ [16] Proposition~15.9, Lemma~22.18. Assertion (A) in Proposition~3.1.2 implies
$$\rho(\cent{G(s)^F}D_{\al_D})=A_{L^*}(s)_\al\leqno {(3.1.3.2)}$$
By conjugacy of maximal Brauer subpairs in a fixed block, one has $G(s)^F_b=G(s)^{\circ\, F}.\nor{G(s)^F}{D,b_D}$. But $\nor{G(s)^F}{D,b_D}=\nor{G(s)^F}{D,\al_D}$. Using (3.1.3.1),
$$\rho(\nor{G(s)^F}{D,\al_D})={A}_{L^*_s,\al} \,.\leqno{(3.1.3.3)}$$

It is known that $D$ is a split extension of
$Z:=\z{L(s)}^F_\ell$ by an
$\ell$-Sylow subgroup of ${\rm W}(C,\zo{L(s)} )^F$ and that 
 $Z$ is a maximal
 normal abelian subgroup of
$D$ [15] Lemma~4.16. As $L(s)=\cento{G(s)^\circ}Z$ ([13] Proposition~3.3), $\nor{G(s)}D^F\subseteq \nor{G(s)}Z^F\subseteq \nor{G(s)}{L(s)}^F\subseteq \nor{G(s)}C^F$. 

That  implies  $\rho(\nor{G( s)^F}{D})\subseteq A_{L(s)}=A_{L^*_s}$. 

By Frattini's argument we have $\nor{G(s)}C^F=C^F.\nor{G(s)}D^F$, hence $\rho(\nor{G(s)}{L(s)}^F)\subseteq \rho(\nor{G(s}D^F)$.
We get $$\rho(\nor{G( s)^F}{D})= A_{L^*_s}\leqno{(3.1.3.4)}$$
 With notations of section~5.1 in Appendix, let $I(B):=I^{G(s)^F}_{G(s)^{\circ F}}(B)$. We have $\rho(I(B))\subseteq \rho(G(s)^F_b)=A_{L^*_s,\al}$ by (B.0) and (3.1.3.1). Thanks to Proposition 3.1.2, (C), Proposition 5.1.4   applies, with $(G(s)^F,G(s)^{\circ\, F})$
instead of
$(G,H)$. By  Proposition~3.1.2, (B), the last equality in (c) of Proposition 5.1.4 simplify in $$\rho(I^{G(s)^F}_{G(s)^{\circ F}}(B))=\rho(\cent{G(s)^F}D_{\al_D})$$ The equality we claim follows from (3.1.3.2).
\bull 

\medskip\noindent{\bf 3.2. Proofs of Propositions 3.1.1 and 3.1.2. Minimal cases }
  
  \medskip\noindent{\bf 3.2.1. Preliminary remarks }

In 3.2 we assume $G$ \irr, simply connected, $s$  "rationally quasi-isolated" in $G^*$ and ${\rm A}_{G^*}(s)^F\neq\{1\}$. By Proposition 1.2.4, the types $\GG_2$, ${\bf F}_4$ and $\EE_8$ are excluded. A semi-simple element $s$ of $G^*$ is said  quasi-isolated in $G^*$ if $\cent{G^*}s$ is not contained in a Levi subgroup of a proper parabolic subgroup of $G^*$ or equivalenly if $\zo{\cent{G^*}s}=\zo{G^*}$. We say that a semi-simple element $s$ of $(G^*)^F$ is  rationally quasi-isolated in $G^*$ if $\cent{G^*}s^F$ is not contained in an $F$-stable Levi subgroup of a proper parabolic subgroup of $G^*$. A classification of quasi-isolated elements in reductive groups is given in [2]. If $s\in G^{*\,F}$, then $\cent{G^*}{\zo{\cent{G^*}s}}$ is an $F$-stable Levi subgroup of $G^*$ : a rationally quasi-isolated element of $G^*$ is quasi-isolated.

In types $\BB$, $\CC$ and $\DD$, there is no non central quasi-isolated semi-simple elements if $\F{}$ has characteristic 2 [2].

Let $s$, $\rho$, $(L(s),\al)$, $L^*_s$, $L^*$, $L$, $D$, $\al_D$, $\beta_D$, $s_1$ be defined as in the Propositions to prove. Recall the properties to verify : 
two non-multiplicity conditions, (B) in Proposition~3.1.1  and (C) in Proposition~3.1.2, and, with notations of Proposition~3.1.2, 

(A) $\rho(\cent{G(s)^F}D)={\rm A}_{L^*}(s)^F$,

(B) $\nor{G(s)^F}D_{\al_D}\subseteq \nor{G(s)^F}D_{\mu_D}$.

 Put $$A:=A_{G^*}(s)^F\, .$$We note first several simple facts we use freely in sectionÊ 3.2 :

 {\sl As $A$ is abelian, non-multiplicity condition is equivalent to maximal extensiblity :  $\chi_0\in\ser{G(s)^{\circ \,F}}{s_1}$ extends to its stabilizer in $G(s)^F$.  If $A$ is cyclic then the maximal extensibility condition is satisfied.

If $L(s)=G(s)^\circ$, then (C) in Proposition~3.1.2 follows from (B) in Proposition~3.1.1; if furthermore  $\z{G(s)^\circ}^F_\ell\subseteq \z{G(s)}^F$ then (A) is satisfied.

If ${\rm A}_{L^*}(s)^F=\{1\}$ then (A) implies (B).

If $|A|$ is a prime then (B) is satisfied. }

Indeed, if $|A|$ is prime, then 

$\bullet$ $\rho(\cent{G(s)^F}D)=\{1\}$ and then $\beta_D=\al_D$ in (b), or 

$\bullet$
$\rho(\cent{G(s)^F}D)=A$ hence $ \nor {G(s)^F}D=\nor {G(s)^{\circ F}}D.\cent{G(s)^F}D$.

 In any case (B) is true.
 
 To verify (A) and (B) in types $\AA$, $\DD$ and $\EE$ we use the fact that $\cento{G^*}s$ and $G(s)^\circ$ have the same type. As $G^*$ is adjoint there exists an isotypic morphism $G(s)^\circ\to \cento{G^*}s$ that extends to $$\pi_s\colon G(s)\to \cent{G^*}s\;{\rm with}\; \pi_s(G(s)^F)=\cent{G^*}s^F\,.$$ Then to prove (A) we shall verify that
 $${\rm If}\; {\rm A}_{L^*}(s) \neq \{1\},\; {\rm then}\; [L^*,L^*]\subseteq \cento{G^*}{\pi_s(D)}\,.\leqno{(3.2.1.1) }$$
 
\medskip\noindent{\bf 3.2.1.2. Lemma. }{\sl Assume all roots of $G$ have same length. One has $\rho(\cent{G(s)^F}D)\subseteq {\rm A}_{L^*}(s)^F$. If furthermore (3.2.1.1) holds, then  $\rho(\cent{G(s)^F}D)={\rm A}_{L^*}(s)^F$.}

\noindent\medskip{\it Proof. } When all roots have same length, $G$ and $G^*$ have same type. Consider, as in the proof of Corollary~3.1.3, a maximally split torus of $S$ of $L(s)$, in duality with a maximally split torus $S^*$ of $L^*_s$. The duality between $G(s)^\circ$ and $\cento{G^*}s$ may be defined by dual root data with respect to $(S,S^*)$. One know that $D$ is an $\ell$-Sylow subgroup of $C^F$ where $C:=\cento{G(s)^\circ}{[L(s),L(s)]}$ [16] Theorem~22.9.  By [16] Proposition~22.7, as $\zo {L(s)}$ is a maximal torus in $C$ and contains a maximal $\phi_d$-subgroup of $C$, one may assume that $D\subseteq \nor{G(s)^\circ}{\zo {L(s)}}$, so that $D$ is a split extension of $Z:=\z {L(s)}^F_\ell$ by an $\ell$-Sylow subgroup of $W(C,\zo{{L(s)}})^F$. Furthermore, $Z$ is caracteristic in $D$  by [15] Lemma~4.16 and ${L(s)}=\cento{G(s)^\circ}Z$  by [13]  Proposition~3.3. 

If $\ell$ divides the order of the kernel $Z_0$  of $G(s)^F\to \cent{G^*}s^F$ then $d=1$ in type $\AA$, $G(s)^\circ$ is a Levi subgroup of $G$, $Z_0$ is the kernel of $G^F\to G^*$, hence $G=\ga$,  and $D$ is a Sylow subgroup of $G(s)^{\circ\, F}$, ${L(s)}$ is a diagonal torus of $G(s)^\circ$  and we'll see in 3.2.2 (i) that $\cent{G(s)^F}D\subseteq G(s)^\circ$ and ${\rm A}_{L^*}(s)=\{1\}$. 

Assume the kernel of the restriction of $\pi_s$ on groups of rational points is prime to $\ell$. Then $\z {L^*_s}_\ell^F=\pi_s(\z {L(s)}^F_\ell)$  and $\cent{\cent{G^*}s^F}{\pi_s(D)}=\pi_s(\cent{G(s)^F}D)$. By Proposition 1.2.5, we have $L^*_s=\cento{\cento{G^*}s}{\z {L^*_s}^F_\ell}$ and $L^*=\cento{G^*}{\z{L^*}^F_\ell}$. As $L^*_s\subseteq L^*$, we have $\z {L^*}^F_\ell\subseteq\z {L^*_s}^F_\ell$, $\cento{G^*}{\z{L^*_s}^F_\ell}\subseteq L^*$. 

Hence $\pi_s(\cent{G(s)^F}D)\subseteq \cent{\cent{G^*}s}{\pi_s(Z)}^F\subseteq \cent{L^*}s$. The last inclusion implies $\rho(\cent{G(s)^F}D)\subseteq {\rm A}_{L^*}(s)$.

 We have $\cent {L^*}s=\zo {L^*}.\cent{[L^*,L^*]}s$ hence $\cent {L^*}s=\cento {L^*}s.\cent{[L^*,L^*]}s$. If (3.2.1.1) holds, then  $\cent{\cent{L^*}s}{\pi_s(D)}.\cento{L^*}s^F=\cent{L^*}s^F$ so that $\rho(\cent{G(s)^F}D)={\rm A}_{L^*}(s)^F$. 
\bull

\medskip\noindent{\bf 3.2.2. Type A}

Here we assume that $G={\rm SL}_{n+1}$, a subgroup of ${\rm GL}_{n+1}$,  and $G^F={\rm SL}_{n+1}(r)$ where $r=\epsilon q$ (by convention ${\rm SL}_a(-q)$ is ${\rm SU}_a(q^2/q)$)  so that  $G^*={\rm PGL}_{n+1}$. Let $\tG={\rm GL}_{n+1}$ acting on $\F{}^{n+1}$. To a natural $F$-epimorphism $$\pi \colon\tG\to G^*$$ there corresponds $G\to {\rm GL}_{n+1}$ between duals.

There exists $\ts\in\tG^F$ with the following properties :

(a) $\pi(\ts)=s$, $\ts$  is semi-simple and of order prime to $\ell$, so that $\pi(\cento{\tG}\ts)=\cento {G^*} s$.

(b) Let $\Gamma$ be the set of eigenvalues of $\ts$ and, for $\gamma\in\Gamma$, let $V_\gamma$ be the corresponding
 eigenspace. All spaces $V_\gamma$ have the same dimension $m$, so that $n+1=|\Gamma|.m$. 
  
 (c) The group ${\rm A}_{G^*}(s)$ acts regularly on $\Gamma$ by translation in $\F{}^\times$ and on the set $\{V_\gamma\}_{\gamma\in\Gamma}$ {\it via} a morphism $v\colon {\rm A}_{G^*}(s)\to \F{}^\times$ so that for any $\tg\in \tG$ with $\pi(\tg)\in\cent {G^*}s$
and any $\gamma\in \Gamma$ then $\tg(V_\gamma)= V_\gamma'$ where $\gamma'=v(\rho(\pi(\tg))).\gamma$. Let $\zeta$ be of order $|{\rm A}_{G^*}(s)|$ in $\F{}^\times $, $\zeta$ is a generator of $v({\rm A}_{G^*}(s))$. 

Let $c_0$ be the order of the orbit of $\zeta$ under the map $(\Phi\colon \F{}\to\F{},\quad \gamma\to \gamma^r)$. Then $A:={\rm A}_{G^*}(s)^F$ has order $c:=|\Gamma|/c_0$ and $\zeta^{c_0}$ is a generator of $v(A)$. One sees easily that the orbit $\omega$ of $\gamma\in\Gamma$ under $\gen{\zeta^{c_0},\Phi}$ acting on $\Gamma$ has order $c\delta(\omega)$, where $\delta(\omega)$ is the order of the orbit of $\gamma$ under $\Phi$. To $\omega$ there corresponds a rational component of $\cent \tG \ts$ of type $[{\rm GL}_m(r^{\delta(\omega)})]^{c}$. One has $|\Gamma|=(\sum_\omega c\delta(\omega))$, hence  $c_0=\sum_\omega \delta(\omega)$.

If there is more than one orbit in $\Gamma$ under $\Phi$, there exists an $F$-stable proper Levi subgroup $K^*$ of $G^*$ such that $\cent{G^*}s^F\subseteq K^*$, hence $A={\rm A}_{K^*}(s)^F$. Assuming $s$ rationally quasi-isolated, 
 there is only one orbit $\omega$ and write $\delta=\delta(\omega)=c_0$. Then $\cento{G^*}s$ has  $c$ rational components, corresponding to a decomposition of $\F{}^{n+1}$ in a direct sum $\oplus_{j\in[1,c]}V_j$ and  $A$ acts on $[1,c]$. 

Let $f:=d/(d,\delta)$. A $d$-cuspidal unipotent datum of $(\cento{G^*}s,F)$ or of $(G (s)^\circ, F)$ is defined by a set of $c$ partitions without any $f$-hook, that is to say $c$ $f$-cores : for $j=1,\dots, c$, $\kappa(j)$ is a partition of $k(j)$, and $g(j)=(m-k(j))/f\in \N$. We have $L^*_s=\pi(\tL^*)$, where $(\tL^*)^F\cong \times_j\big([{\rm GL}_1(r^{f\delta})]^{g(j)}\times {\rm GL}_{k(j)}(r^\delta)\big)$ (if $f=1$, then $k(j)=0$ and $g(j)=m$ for all $j$). Then $L^*=\pi(M^*)$ where $(M^*)^F\cong \big(\times_j[{\rm GL}_{f\delta/d}(r^d)]^{g(j)}\big)\times{\rm GL}_N(r)$ with $N=\delta(\sum_jk(j))$. One sees that
$${\rm A}_{L^*}(s)^F=\{a\in A\mid {\rm if} \;g(j)\neq 0)\;{\rm then}\; a(j)=j\}, \quad A_{L^*_s}=\{a\in A\mid g\circ a =g\}, \leqno{(3.2.2.1)}$$ 
$${\rm A}_{L^*}(s)^F_\al= \{a\in {\rm A}_{L^*}(s)^F\mid \kappa\circ a =\kappa\},\quad A_{L^*_s,\al}=\{a\in A\mid \kappa\circ a =\kappa\}\,.\leqno{(3.2.2.2)}$$
 Here we consider ${\rm A}_{L^*}(s)^F$ as a subgroup of $A$. Indeed  if $g(j)\neq 0$, any element of $\pi^{-1}((L^*)^F)$ stabilizes a non null subspace of $V_j$. By its component $\pi({\rm GL}_N(r))\cap\cent{G^*}s^F$, $\cent{L^*}s^F$ acts on the set of $j$ with $g(j)=0$ as freely as $\cent{G^*}s^F$. That gives ${\rm A}_{L^*}(s)^F$. The three others equality are clear.
 
 As $A$ is cyclic we have only to verify (A) and (B).
 
 As said in  condition (A2) in Proposition~3.1.3, we may assume that $G(s)^\circ$ is a Levi subgroup of $G$, with same rational type that $\cento{G^*}s$. With evident notations, ${L(s)}=L_1\dots L_j\dots L_c$. A defect group $D$ of $b$ is an $\ell$-Sylow subgroup of $\cento{G(s)^\circ}{[{L(s)},{L(s)}]}^F$ and $\pi^{-1}(\cento{G(s)^\circ}{[{L(s)},{L(s)}]}^F)\cong \times_j\big({\rm GL}_{fg(j)}(r^\delta)\times  {\rm GL}_{k(j)}(r^\delta)\big)$ : $D$ is a  central product of the ${\pi(D_j)}$, where $D_j$ is an $\ell$-Sylow subgroup of a wreath product  ${\rm GL}_1(r^{f\delta})\wr {\cal S}_{g(j)}$ (because the torus with rational-points group ${\rm GL}_1(r^{f\delta}))$ is a minimal $d$-split Levi subgroup of the component  ${\rm GL}_{fg(j)}(r^\delta)$). 
 
 (o) Let us consider the case $d=1$, that is  $G=\ga$, to complete the proof of Lemma~3.2.1.2.  
 
 As ${L(s)}$ is a diagonal torus of $G(s)^\circ$ and $M=\cent G{\zo {L(s)}_{\phi_1}}$, ${L(s)}^F={\rm SL}_{n+1}\cap[{\rm GL}_1(r^\delta)]^{mc}$ and $M^F=   {\rm SL}_{n+1}\cap [{\rm GL}_\delta(r)]^{mc}$ and $L^*=\pi( [{\rm GL}_\delta(r)]^{mc})$. One sees that ${\rm A}_{L^*}(s)=\{1\}$. One has $(r^\delta-1)_\ell=(r-1)_\ell$ because $\delta$ is a divisor of $|A|$, prime to $\ell$. It follows that $\z {L(s)}^F_\ell=\z M^F_\ell$ and, by Proposition~1.1.6, that $M^F=\cent{G^F}{\z {L(s)}^F_\ell}$. Then $\cent{G(s)^F}D\subseteq \cent G{\z {L(s)}^F_\ell}\subseteq M$.
 
 (i) Assume  that $m=1$.  Then $G(s)^\circ$ and $\cento{G^*}s$ are tori, ${L(s)}=G(s)^\circ$, $\al=1_{{L(s)}^F}$, $D=G(s)^{\circ\, F}_\ell$, $b_D$ is the principal block. In case $f>1$,   $D=\{1\}$ and $\zo{\cento{G^*}s}_{\phi_d}=\{1\}$. It follows that $L^*=G^*$ and $\cent{G(s)^F}D=G(s)^F$.  Hence (3.2.1.1), (A) and (B) are true in that case. In case $f=1$, $d$ divides $\delta$. As $\delta$ is prime to $\ell$, $(r^\delta-1)_\ell=(r^d-1)_\ell$.  In other words $\cento {G^*}s^F_\ell\subseteq \cento{G^*}s_{\phi_d}=\zo{L^*}_{\phi_d}$, so that $D\subseteq \z M_{\phi_d}$. This implies $\pi_s(D)\subseteq {\z L^*}$ hence (3.2.1.1) and (A) by Lemma~3.2.1.2.
The condition (B) follows, $b_D$ beeing the principal block, $\al_D$ and $\beta_D$ are the unit characters.
 
 (ii) Assume now that $g(j)=0$ for any $j\in[1,c]$, but $m>1$. Then $f>1$ and again $D=\{1\}$, $A_{L^*}(s)=A_{G^*}(s)$, $\cent{G(s)^F}D=\nor{G(s)^F}D=G(s)^F$, (3.2.1.1), (A) and (B) are clearly true.
 
 (iii) If $g(j)\neq 0$ for any $j$, by (3.2.2.1) ${\rm A}_{L^*}(s)^F=\{1\}$. But $\pi(D_j)\neq \{1\}$ for any $j$ and $A$ acts regularly on $[1,c]$, hence $\rho(\cent{G(s)^F}D)=\{1\}$, $\al_D=\beta_D$ and (A), (B) are true.
 
 (iv) In the general case consider $J_0=\{j\in[1,c]\mid g(j)=0\}$ and $J_1=[1,c]\setminus J_0$. That partition defines a decomposition $\F{}^{n+1}=V_0\oplus V_1$ and dual maximal proper $F$-stable Levi subgroups $H=H_0.H_1$, $H^*=H^*_0.H^*_1$ (central products, $H^*_j$ and $H_j$ are  in duality "up to an isotypic morphism") of $G$  and $G^*$. One has $G(s)^\circ=G(s)^\circ_0.G(s)^\circ_1$ with $G(s)^\circ_i=H_i\cap G(s)^\circ$, $\cento{G^*}s=(H_0^*\cap \cento{G^*}s).(H_1^*\cap \cento{G^*}s)$. One may write $s=s_0.s_1\in (H^*_0)^F_{\ell'}.(H^*_1)^F_{\ell'}$ and ${\rm A}_{H^*_0}(s_0)^F\times {\rm A}_{H^*_1}(s_1)^F$ may be viewed as a subgroup of ${\rm A}_{G^*}(s)^F$ : it is the stabilizer in ${\rm A}_{G^*}(s)$ of the decomposition $G(s)^\circ_0.G(s)^\circ_1$. By  (3.2.2.1) $A_{L^*_s}\subseteq {\rm A}_{H^*_0}(s_0)^F\times {\rm A}_{H^*_1}(s_1)^F$. 
 As seen in the extreme cases (ii) and (iii) above, where one $J_i$ is $[1,c]$, 
 $$H^*_0\subseteq L^*, \quad {L(s)}=(H_0\cap G(s)^{\circ F}).({L(s)}\cap H_1),\quad D\subseteq H_1$$
  hence (3.2.1.1). By (3.2.2.1) again, ${\rm A}_{L^*}(s)^F={\rm A}_{H^*_0}(s_0)^F$. Now $\cent{G(s)^{\circ F}}D=G(s)_0^{\circ F}.\cent{ G(s)^{\circ F}_1}D$ and  $\al_D=\al_0\otimes \al_1$, where $\al_0$ has null defect. Furthermore $\rho(\cent{G(s)^F}D)={\rm }A_{H^*}(s_0)^F$ (see (iii)) that is (A). Then $\beta_D=\beta_1\otimes \al_1$. Similarly $\nor{G(s)^F}D$ stabilizes the \dec\ of $G(s)^\circ$ and, by (ii) and (iii), 
  $$\rho(\nor{G(s)^F}D)=A_{H^*_0}(s_0)^F.\rho(\nor{G(s)^{\circ F}_1}D)\subseteq \rho(\cent{G(s)^F}D).\rho(\nor{G(s)^{\circ F}_1}D).$$  This proves (B).

 \smallskip\noindent{\bf Type \BB }

Let $\pi\colon {\rm Sp}_{2n}\to G^*={\rm PSp}_{2n}$. 
Some $\hat s\in \pi^{-1}(s)$ has order $2$ or $4$ and $\hat s$ is conjugate to $-\hat s$ under ${\rm A}_{G^*}(s)$. In the standard action on $\F{}^{2n}$, $\hat s$ has eigenvalue $1$ and $-1$ with multiplicity $2m$, and primitive $4$-roots of $1$ with multiplicity $(n-2m)$. So
$\cento{G^*}s$ has type $\CC_m\times\AA_{n-2m-1}\times \CC_m$ with $0\leq m\leq n/2$, $G(s)^\circ$ has type $\BB_m\times\AA_{n-2m-1}\times \BB_m$
(with usual conventions on $\AA_0$, $\BB_0$, $\BB_1$).

A generator of $A_{G^*}(s)$ exchange the factors of type $\CC_m$, if there are some, and acts as diagram automorphism on the factor of type $\AA$, if there is some.  Then $L^*_s$ has type $\CC_{m_1}\times \AA_{m_3}\times \CC_{m_2}$  for some $m_j$ with $m_1, m_2\leq m$, $m_3\leq n-2m-1 $ and $L^*$ has type $\CC_{m_1+m_2+m_3+1}$, short of the case $m_1=m_2=0$, where $L^*$ has type $\AA_{m_3}$.

If $L^*_s\neq \cento{G^*}s$, that is ${L(s)}\neq G(s)^\circ$, then ${\rm A}_{L^*}(s)=\{1\}$ and $D$ has a non trivial intersection with one of the \irr\ components of $G(s)^\circ$, so that $\cent{G(s)^F}D\subseteq G(s)^\circ$. 

If $L^*_s=\cento{G^*}s$, then the defect group $D$ of $b_{G(s)^{\circ F}}({L(s)},\al)$ is central in $G(s)^{\circ\, F}$ and contained in the type $\AA$ factor. One has $D=\z{G(s)^\circ}^F_\ell$. The polynomial order of $\zo{G(s)^\circ}$ has only cyclotomic factors $\phi_1$ and $\phi_2$ and $\ell>2$, hence $\z{G(s)^\circ}^F_\ell\subseteq\z{G(s)^\circ}_{\phi_d}$. But $L^*=\cent{G^*}{\z {\cento{G^*}s}_{\phi_d}}$ and ${\rm A}_{L^*}(s)$  acts trivially on $\z{\cento{G^*}s}_{\phi_d}$. Dually ${\rm A}_{L^*}(s)$ acts trivially on $ \z{G(s)^\circ}_{\phi_d}$, therefore on $D$, that is (A).

\smallskip\noindent{\bf Type C }

Take $G={\rm Sp}_{2n}$, $G^*={\rm SO}_{2n+1}$ and let $\hat G$ be the spin group of same type with $\pi\colon \hat G\to G^*$ a natural quotient. a semi-simple $s\in G^*$ has non connected centralizer in $G^*$ if and only if it has eigenvalues $1$  and $-1$ in the standard action on $\F{2n+1}$ [16] Proposition~16.25.  We may consider two cases for quasi-isolated semi-simple elements in $G^*$ :

{\it Case 1. } {\sl $s$ is not isolated, has d'ordre 2,  the fixed-point space of $s$ in $\F{}^{2n+1}$ has dimension $(2n-1)$, and $\cento{G^*}s$ is a maximal proper Levi subgroup of  type $ \BB_{n-1}$.}

 $G(s)^\circ$ is a Levi subgroup of $G$ and has type $\CC_{n-1}$. There exists an involution $x\in{\rm SO}_{2n+1}(q)$ acting on $\cento{G^*}s$ such that $xyx= y^{-1}$ if $y\in \zo{\cento{G^*}s} $ and $xyx=y$ if $y\in[\cento{G^*}s,\cento{G^*}s]$. Then $\rho(x)$ generates  $A$ and $x$ acts on $G(s)^\circ$ on the same way : $xyx= y^{-1}$ if $y\in \zo{G(s)^\circ} $ and $xyx=y$ if $y\in[G(s)^\circ,G(s)^\circ]$. 
 
If $\zo{G(s)^\circ}^F_\ell=\{1\}$ and ${L(s)}\neq
G(s)^\circ$, then   $\{1\}\neq D\subseteq [G(s)^\circ,G(s)^\circ]$, hence $\rho(\cent{G(s)^F}D)=A$.  
$\zo{\cento{G^*}s}\subseteq L^*=\cent{G^*}{\zo {L^*_s}}_{\phi_d}$ $s$ is quasi-isolated in $L^*$ so that
${\rm A}_{L^*_s}(s)={\rm A}_{G^*}(s)$. 

If $\zo{G(s)^\circ}^F_\ell\neq\{1\}$, then $d=1$, $G(s)^\circ$ is
a $d$-split Levi subgroup of $G$, and $G(s)^\circ=\cento{G}{\zo{G(s)^\circ}^F_\ell}$, hence $M\subseteq G(s)^\circ$, $L^*\subseteq \cento{G^*}s$ and  ${\rm A}_{L^*}(s)=\{1\}$. Furthermore $\zo{G(s)^\circ}^F_\ell\subseteq D$, $G(s)^\circ=\cento{G}{\zo{G(s)^\circ}^F_\ell}$ [16] Lemma~13.17 and $\cento{G}D=\cent GD$ (Proposition~1.2.4). Thus $\cent{G(s)^F}D\subseteq G(s)^\circ$ and (A) is true.

{\it Case 2. } {\sl   $s$ is isolated, has order $ 2$,  the fixed-point space of $s$ in $\F{}^{2n+1}$ has dimension $2n+1-2m$ and 
$\cento{G^*}s$ has type
$\DD_m\times \BB_{n-m}$ ($1<m\leq n$).}

Then $\cento{G^*}s$ is isomorphic to
SO$_{2m}\times{\rm SO}_{2n-2m+1}$, $s$ belongs to the first component and $\cent{G^*}s$ is isomorphic to $  ({\rm O}_{2m}\times{\rm O}_{2n-2m+1})\cap {\rm SO}_{2n+1}$. But ${\rm O}_{2k+1}={\rm SO}_{2k+1}\times\gen{-{\rm Id}}$. On $G(s)^\circ\cong{\rm  SO}_{2m}\times{\rm Sp}_{2n-2m}$ a generator of $A$ acts by a diagonal automorphism on the first component and
$G(s)\cong {\rm O}_{2m}\times {\rm Sp}_{2n-2m}$.

If ${L(s)}=G(s)^\circ$, then $D=\{1\} $ and  $L^*=G^*$, hence (A). 

In general $s$ is quasi-isolated in $L^*$ with ${\rm A}_{L^*}(s)\neq \{1\}$ if and only if $s\notin\z{L^*}$. That condition is equivalent to  "$[L^*_s,L^*_s]$ is not contained in the second component of $\cento{G^*}s$" or to  "$[{L(s)},{L(s)}]$ is not contained in the component $ {\rm Sp}_{2n-2m}$ of $G(s)^\circ$".  One sees easily that if $X$ is an $\ell$-subgroup of ${\rm SO}_{2k}(q)$ with a non null space of fixed points on $\F{}^{2k}$, then $\cent{{\rm O}_{2k}}X.{\rm SO}_{2k}={\rm O}_{2k}$. This applies to $D\cap{\rm  SO}_{2m}$ in the first factor of  $G(s)^\circ$ when  ${\rm A}_{L^*}(s)\neq \{1\}$ because $D$ centralizes $[{L(s)},{L(s)}]$.

 If ${\rm A}_{L^*}(s)=\{1\}$, then $d$ divides $m$ and $D$ contains an $\ell$-Sylow subgroup of ${\rm SO}_{2m}(q)$. One can verify that if $S$ is the minimal $d$-split Levi subgroup, it is a maximal torus in  ${\rm SO}_{2m}(q)$ and $\cent{{\rm O}_{2m}}{S^F_\ell}\subseteq {\rm SO}_{2m}$. That shows that $\cent{G(s)^F}D\subseteq G(s)^\circ$.

\smallskip\noindent{\bf 3.2.4. Type D.  } 

We assume
that $G$ is the spin group of  a non degenerate symmetric bilinear  form on a space of dimension $2n$ on $\F{}$ and  $G^*={\rm PSO}_{2n}$, all defined on $\F q$ by $F$. There are isotypic morphisms $\pi\colon G\to G^*$,   $G(s)^\circ\to \cento{G^*}s$ and a restriction $\pi_s\colon G(s)^{\circ\, F}\to \cento{G^*}s^F$ (see section 3.2.1).

(a) Assume $s$ a semi-simple quasi-isolated element of
$G^{*\,F}$  such that $|{\rm A}_{G^*}s^F|=4$.

 Then  $s^4=1$ and $\cento{G^*}s$ has type $\DD_m\times \AA_{n-2m-1}\times \DD_m$ where $2\leq m\leq n/2$  (and with some conventional notations : if $m=n/2$, then the component of ``type $\AA_{-1}$" is $\{1\}$ and $s$ is isolated) or $\cento{G^*}s$ has type $\AA_{n-3}$ (special case $m=1$ above). There exists $\hat s\in{\rm SO}_{2n}$, which has  image $s$ in PSO$_{2n}$ and four eigenvalues $1,-1, i,-i$ ($i^2=-1$) with respective multiplicities $m, m, n-2m, n-2m$ in the standard action of SO$_{2n}$ on $\F{}^{2n}$. 
 
 So $\cento{G^*}s$ is the image in ${\rm PSO}_{2n}$ of a subgroup of  ${\rm SO}_{2n}$ isomorphic to ${\rm SO}_{2m}\times {\rm GL}_{n-2m}\times {\rm SO}_{2m}$, or to ${\rm GL}_1\times{\rm GL}_{2n-2}\times {\rm GL}_1$ if $m=1$.
The dual $G(s)^\circ=G_1.G_0.G_2$ of $\cento{G^*}s$ maps dually in  ${\rm SO}_{2m}\times {\rm GL}_{n-2m}\times {\rm SO}_{2m}$ (or if $m=1$, $G(s)^\circ$ is a Levi subgroup of $G$, $G_1$ and $G_2$ are torii) and so is defined $\pi_s$. 

The groups of
diagram automorphisms of $\cento{G^*}s$ and $G(s)^\circ$ are isomorphic. Let $\delta_1$ be a diagram automorphism of order 2 of the component $G_1$ of type $\DD_m$ of $G(s)^{\circ F}$ (if $m=1$, $\delta_j$ acts by inversion on $G_j$) let $\tau$ be a diagram automorphism  defined on $\F q$ and of order 2 that exchange the two components $G_1$ and $G_2$, $\delta_2:=\tau\delta_1$. If $G_0\neq\{1\}$ let  $\gamma$ be a diagram automorphism of order 2  and defined on $\F q$ of $G_0$. If $G_0=\{1\}$ let $\gamma=1$ . One has $|\cF(G(s)^\circ)|=|\z{\cento{G^*}s}/\zo{\cento{G^*}s}|=4$.

The group $A:={\rm A}_{G^*}(s)^F$ is non cyclic if and only if $n$ is even. Then $A$ is generated by diagram automorphisms : $A=\gen{\tau.\gamma, \delta_1\delta_2}$. $\delta_1\delta_2$ is induced by the image in ${\rm PSO}_{2n}$ of $({\rm O}_{2m}\times {\rm GL}_{n-2m}\times {\rm
O}_{2m})\cap {\rm SO}_{2n}$. The extension $\cento{G^*}s\to \cent{G^*}s\to {\rm A}_{G^*}(s)$ is split. We may assume that $\tau$ acts on $({\rm O}_{2m}\times {\rm GL}_{n-2m}\times {\rm
O}_{2m})$ by exchange of the components ${\rm O}_{2m}$.
Thus  $G(s)^F$  is a quotient by a finite  central $2$-group of a subgroup of the direct product $({\rm Spin}_{2m}(q).\gen \delta)\wr{\gen \tau}\times {\rm GL}_{n-2m}(q).\gen \gamma$. 

Let $s_1\in \cento{G^*}s^F$ with odd order and $\xi\in\ser{G(s)^{\circ F}}{s_1}$. The kernel of $\xi$ contains $\z{G(s)^\circ}^F_2$. If $G(s)^F$ fixes $\xi$, $\xi$ has the form $\xi_1\otimes \xi_0\otimes \xi_2$ where $\gamma(\xi_0)=\xi_0$, $\delta_j(\xi_j)=\xi_j$ and 
$\tau(\xi_1)=\xi_2$. Then $\xi$ is fixed under the all group $\gen{\delta,\tau,\gamma}$. So $\xi_1$ extends to $\zeta_1\in\II{{\rm Spin}_{2m}(q).\gen\delta}$, $\xi_2$ extends to $\tau(\zeta_1)$, $\xi_0$ extends to $\la_0\in\II{{\rm GL}_{n-2m}(q).\gen \gamma}$. Then $\zeta_1\otimes \zeta_2$ extends to $({\rm Spin}_{2m}(q).\gen \delta)\wr{\gen \tau}$. It follows that $\xi$ extends to $G(s)^{\circ \,F}.\gen{\delta,\tau,\gamma}$. The condition (B) in Proposition~3.1.1 is satisfied.

If $n$ is odd, $A$ is generated by $\tau\delta_1\gamma$.

One has  ${L(s)}=L_1.L_0.L_2$ a central product, and $\al=\al_1\otimes\al_0\otimes \al_2$, where $(L_j,\al_j)$ is a $d$-cuspidal unipotent datum of $(G_j,F)$, eventually $L_0=\{1\}$  (Proposition~1.3.1.(ii), if $m=1$, $\al_j=1_{G_j^F}$ for $j=1,2$).

(a.1) If ${L(s)}=G(s)^\circ$, then $D=\z{G(s)^\circ}^F_\ell$. 

If $\z{G(s)^\circ}^F_\ell= \{1\}$ then  $D=\{1\}$ and
$L^*=G^*$  there is nothing to prove.

If  $\z{G(s)^\circ}^F_\ell\neq \{1\}$ then $d=1$.

If  $m=1$, then  $\cento{G^*}s$ is a $d$-split Levi subgroup of $G^*$ and ${L(s)}$ is a maximally split maximal torus in $G(s)^\circ$, $L^*_s$ is a d-split Levi subgroup of $G^*$, $L^*= L^*_s$, ${\rm A}_{L^*}(s)=\{1\}$, $D$ is an $\ell$-Sylow subgroup of $G(s)^F$, $\cent{G(s)^F}D\subseteq \cent{G(s)^F}{\z {L(s)}^F_\ell}={L(s)}$ (Proposition~1.1.6)), $\al_D=\beta_D$, there is nothing to prove.

 If $m>1$ : 

Then $D\subseteq G_0$, $G_0$ is a torus  of rank $1$, $n=2m+1$, $\gamma$ acts non trivially on  $D$ by $(x\mapsto x^{-1}$), hence  $\rho(\cent{G(s)^F}D)=\gen {\delta_1\delta_2} $. As $\ell$ is good and $\cF(G)$ is prime to $\ell$, $\cent{G}D=\cent{G}{(G_0)_{\phi_d}}=M$ (Proposition~1.1.6); $L^*$ is a maximal Levi subgroup of type $\DD_{n-1}$ of $G^*$ and ${\rm A}_{L^*}(s)=\gen{\delta_1\delta_2}$, that is (3.2.5.1) and (A). 

By definition  $\beta_D$ is an extension of $\al=\al_D$ to $\cent{G(s)^F}D$; $\beta_D$ is also the restriction  of some $\beta_1\otimes \al_0\otimes \beta_2$ where $\beta_j$ extends $\al_j$ to ${\rm O}_{2m}(q)$ for $j=1, 2$. If $\rho(\nor{G(s)^F}D_{\al_D})=A$, then $\al_2=\tau(\al_1)$ and $\gamma(\al_0)=\al_0$. There are four extensions of $\al_1\otimes\al_2$ to ${\rm O}_{2m}(q)\times {\rm O}_{2m}(q)$ : $(\theta_1\beta_1)\otimes (\theta_2\beta_2) $ where $\theta_j$ is linear with square 1. Then $\gen\tau$ has two fixed points on that set, such as $\beta'_1\otimes\beta'_2$ and $\theta_1\beta'_1\otimes \theta_2\beta'_2$ where $\theta_j\neq 1$ (and with a good choice of notations) and an orbit $\{\theta_1\beta'_1\otimes \beta'_2,\beta'_1\otimes \theta_2\beta'_2\}$. But the first two one's have same restriction to ${\rm O}_{2m}(q)\times {\rm O}_{2m}(q)\cap{\rm SO}_{2m}(q)$, as well as the two last one. In other words $\tau$ fixes any extension of $\al_1\otimes \al_2$ to ${\rm O}_{2m}(q)\times {\rm O}_{2m}(q)\cap{\rm SO}_{2m}(q)$. This shows that $G(s)^F$ stabilizes $\beta_D$, that is (B).

(a.2) Assume ${L(s)}\neq G(s)^\circ$.

If $m>1$,  then ${L(s)}$ has type $\DD_{m_1}\times \AA_{m_0}\times \DD_{m_2}$ ($m_0>0$) or  $\DD_{m_1}\times \DD_{m_2}$ and $M$ has type $\DD_{m_1+m_0+m_2+1}$ or $\DD_{m_1+m_2}$, with special cases where $m_j=1$ for some $j\in\{1,2\}$, SO$_2$ being a torus. If $m=1$, then $L^*$ has type $\DD_t$ with $t\geq 2$.
Going back  from $G^*$ to SO$_{2n}$ one sees that $[L^*_s,L^*_s]$ and  $[L^*,L^*]$ have equal spaces of fixed points  on $\F{}^{2n}$. As $[{L(s)},{L(s)}]\subseteq \cent {G(s)^\circ} D$,  $[L^*,L^*]\subset\cent{G^*}{\pi_s( D)}$, hence $\cent{L^*}s\subseteq \cent{G^*}{\pi_s(D)}$. This implies (A).  

The description of $A$ above shows that there exists $x\in\zo{L^*}$ and  $a\in A$ such that $a(x)\neq x$. Thus ${\rm A}_{L^*}(s)\neq A$,  $|{\rm A}_{L^*}(s)|\leq 2$ and the non-multiplicity condition (C) of Proposition~3.1.2 is satisfied.

Assume $|{\rm A}_{L^*}(s)|=2$ and $\rho(\nor{G(s)^F}D_{\al_D}=A$. The all group $A$ stabilizes the block $b_{G(s)^{\circ F}}({L(s)},\al)$.

 In case $m>1$ this implies $m_1=m_2$, $\tau(\al_1,\al_2)=(\la_2,\al_1)$, $\gamma(\al_0)=\al_0$, $\delta_j(\al_j)=\al_j$. As in case (a.1) the all group $\gen{\tau,\gamma,\delta_j}$ acts on $\cent{G(s)^\circ}D$ and fixes $\al_D$ and (B) follows as in (a.1).

If $m=1$ then $d>1$, $\tau\in{\rm A}_{L^*}(s)$ and ${\rm A}_{L^*}(s)\neq A$ is equivalent to $L_0\neq G_0$. Then $\rho(\nor{G(s)^F}D_{\al_D}=A$ is equivalent to  $\gamma(\al_0)=\al_0$. Once again the all group $\gen{\tau,\gamma,\delta_j}$ acts on $\cent{G(s)^\circ}D$ and fixes $\al_D$ : (B) follows as in (a.1).

(b) Assume   $s^2=1$ and
$\cento{G^*}s$ has type
$\DD_m\times \DD_{n-m}$, ($2\leq m<n/2$) or $\DD_{n-1}$ (special case $m=1$). 

Then $|{\rm A}_{G^*}(s)|=2$ and $A$ is generated by the image in $G^*$
of an 
element of SO$_{2n}\cap ({\rm O}_{2m}\times {\rm O}_{2(n-m)})$. When $m=1$ $\cento{G^*}s$ and $G(s)^\circ$ are Levi subgroups of $G^*$ and $G$ respectively. 

If ${L(s)}=G(s)^\circ$, then $D=\{1\}$ or $m=1$ and $d=1$. In case $m=1$, $G(s)^{\circ F}=\cento G D^F=M^F$ by Proposition~1.1.6. So (A) is satisfied.

Assume ${L(s)}\neq G(s)^\circ$. As in (a.2) above one sees that $[L^*,L^*]\subseteq\cent{G^*}{\pi_s( D)}$, hence (A).

(c) Assume  $s^2=1$ and $\cento{G^*}s$ is a Levi subgroup of type ${\bf A}_{n-1}$. 

Then $n$ is even, $s$ is the image  of $\hat s\in{\rm SO}_{2n}$ such that $\hat s^4=1$ and $\hat s$ has two eigenvalues with equal multiplicities $n/2$. The group $A$ has order 2  and is generated by a diagram automorphism. Then $D=\{1\}$ or ${\rm A}_{L^*}(s)=\{1\}$.

\medskip\noindent{\bf 3.2.5. Type E$_6$. } 

In that type if $s$ is a quasi isolated semi-simple element of $G^*$ and $|{\rm A}_{G^*}(s)|>1$, then $|{\rm A}_{G^*}(s)|=3$,

$s^3=1$ and 
$\cento{G^*}s$ has type
$(\AA_2)^3$ ($s$ isolated) or 
$\DD_4$ or $s^6=1$ and $\cento{G^*}s$ has type $({\bf A}_1)^4$.

One has only to verify that  
${\rm A}_{L^*}(s)^F\subseteq \rho(\cent{G(s)^F}D)$ and a sufficient condition is that  ${L^*}^F\subseteq \cent {G^*}{\pi_s(D)}$ (see 3.2.1). 
One sees that 
$\phi_{d\ell}$ does not divides the generic degree of $\cento{G^*}s$ and the order of the Weyl group of $\cento{G^*}s$ is prime to $\ell$ (recall that $\ell\geq 5$ in our assumption). It follows that any
$\ell$-subgroup of
$G(s)^{\circ F}$ (resp. $\cento{G^*}s$) is contained in a $\phi_d$-subgroup of $G(s)^\circ$ (resp. $\cento{G^*}s^F$). Recall that $\zo{{L(s)}}_{\phi_d}$ is a maximal 
$\phi_d$-subgroup of $\cento{G(s)^\circ}{[{L(s)},{L(s)}]}$ and $D$ an $\ell$-Sylow subgroup of $\cento{G(s)^\circ}{[{L(s)},{L(s)}]}$. So $D\subseteq\zo{{L(s)}}_{\phi_d}$ and  $\pi_s(D)\subseteq\zo{L^*_s}_{\phi_d}$. Hence $L^*=\cent
{G^*}{\zo{L^*_s}_{\phi_d}}\subseteq\cent {G^*}{\pi_s( D)}$.

\smallskip\noindent{\bf 3.2.6. Type E$_7$. }

In all cases we have to consider $|{\rm A}_{G^*}(s)|=2$ :

$s^4=1$, $\cento{G^*}s$ has
type
${\bf A}_3\times {\bf A}_3\times {\bf A}_1$ or

$s^4=1$, $\cento{G^*}s$ has type $\DD_4\times\AA_1\times \AA_1$ or

$s^6=1$, $\cento{G^*}s$ has type $\AA_2\times\AA_2\times\AA_2$ or

$s^2=1$, $\cento{G^*}s$ has type $\EE_6$ or

$s^2=1$, $\cento{G^*}s$ has type ${\bf A}_7$.

The generic polynomial order of $\EE_7(q)$ is, in symbolic notations, $0^{63}.1^7.2^7.3^3.4^2.5.6^3.7.8.9.10.12.14.18$

The proof we gave in 3.2.5 is available with a special attention in the last two cases when $\phi_5$, $\phi_7$  divide the polynomial order of $\cento{G^*}s$.

Assume $\cento{G^*}s$ has type $\EE_6$, $\cento{G^*}s$ is a Levi subgroup of $G^*$. If ${L(s)}=G(s)^\circ$ and $\z{{L(s)}}^F_\ell\neq \{1\}$, then $d=1$, $\zo{L^*_s}=\zo{L^*_s}_{\phi_1}$ and $L^*=L^*_s=\cent{G^*}{\z{L^*_s}^F_\ell}$, so that $L^*\subseteq\cent{G^*}{\pi_s(D)}$. That proves (A). More generally one sees that for any $d$ and any prime $\ell\geq 5$, denoting $E:=\{d\ell^b\mid b\in \N, b>1\}$, if  ${L(s)}$ is a $d$-split Levi subgroup of $G(s)^\circ$ then  $\zo {L(s)}_{\phi_{E}}=\zo {L(s)}_{\phi_d}$.  This imply $\zo {L^*_s}_{\phi_{E}}=\zo {L^*_s}_{\phi_d}$ and $L^*=\cento{G^*}{\z{L^*_s}^F_\ell}$, hence $L^*\subseteq \cent{G^*}{\pi_s(D)}$.

\medskip\noindent{\bf 3.3. From minimal cases to general one}

We prove Propositions~3.1.1 and 3.1.2 by induction on the semi-simple rank of $G$ and begin with the more evident result :

\medskip\noindent{\bf 3.3.1. Direct product }

If $(G,F)$ is a direct product $\times_j (G_j,F)$,  so is $(G^*,F)$, $s=(s_j)_j$ ($s_j\in G^{*\,F}_j$), ${\rm A}_{G^*}(s)=\times_j{\rm A}_{G^*_j}(s_j)$ and the condition (A.1) define the sequence ${\cal E}(G,s)$ as a direct product of the ${\cal E}(G_j,s_j)$. Any of the groups ${L(s)}$, $L$, $L^*$, $D$, $\cent{G(s)^F}D$, $\nor{G(s)^F}D$ decomposes in a direct product and $\al_D$, $\beta_D$ in a tensor product. Assertions (B) of Proposition 3.1.1 and  (A), (B) and (C) of Proposition 3.1.2 are true for $(G,s)$ if and only if they are true for each direct component.

\medskip\noindent{\bf 3.3.2. Isotypic morphisms }

The non-contradiction between  conditions (A.1), (A.2) and (A.3) of Proposition~3.1.1 are easily verified. We precise here condition (A.3).

Let $\sigma\colon (G,F)\to (H,F)$ be an isotypic morphism and $\sigma^*\colon H^*\to G^* $
a dual one, let $s=\sigma^*(t)$, as in (3.1.1.(A.3)). If $T^*$ is a maximally split torus in $\cento{G^*}s$, $(\sigma^*)^{-1}(T^*)$ is a maximally split torus in $\cento{H^*}t$. 
The restriction  of $\sigma^*$ to $\cento{H^*}t$ is isotypic and defines a  morphism of sequences

$${\def\normalbaselines{\baselineskip20pt\lineskip3pt
\lineskiplimit3pt}\matrix
{{\frak E}(H^*,t)&
\quad\quad&\cento{H^*}t&\mapright{}&\cent{H^*}t&\mapright{}&{\rm A}_{H^*}(t)
\cr
\mapdown{}&&\mapdown{\sigma^*}&&\mapdown{\sigma^*}&&\mapdown{\alpha}\cr
{\frak E}(G^*,s)&&\cento{G^*}s&\mapright{}&\cent{G^*}s&
\mapright{}&{{\rm A}_{G^*}(s)}\cr} }$$

There are tori $T$ in $G$ and $\sigma(T)$ in $H$ in duality with respectively $T^*$ and $(\sigma^*)^{-1}(T^*)$.
So are defined  the connected reductive groups
$(H(t)^\circ, F)$, $(G(s)^\circ, F)$ in duality with
respectivement
$(\cento{H^*}t,F)$ and $(\cento{G^*}s,F)$. There is a morphism 
$$\sigma_t \colon (G(s)^\circ,F)\to (H(t)^\circ,F)$$
defined by duality from the restriction of $\sigma^*$, and such that $\sigma_t$ and $\sigma$ have equal restriction on 
$\nor {G(s)^\circ} T$. 

Consider a given  sequence $$1\to G(s)^\circ\to
G(s)\to {\rm A}_{G^*}(s)\to 1\leqno{{\frak E}(G,s)}$$ The kernel $K$ of $\sigma_t$ is central and contained in $T$, so it is stable under the group
$\alpha({\rm A}_{H^*}(t))$. Let
$G(t)$ be the inverse image of $\alpha({\rm A}_{H^*}(t))$ in $G(s)$. By condition (3.1.1.(A.3)) 
${\frak E}(H,t)$ is given, isomorphic to
$$1\to H(t)^\circ=G(s)^\circ/K\to H(t)=G(t)/K\to \alpha({\rm A}_{H^*}(t))\to 1$$

 Note that $G(t)^F/G(s)^{\circ \, F}$ is isomorphic to $H(t)^F/H(t)^{\circ\,F}$. If conditions (3.1.1.(A.1)) and (3.1.1.(A.2)) are satisfied in the construction of ${\frak E}(G,s)$, they are as well satisfied for ${\frak E}(H,t)$. 
 
 Isotypic morphisms may be composed. Using functoriality of fiber product, one verify easily the compatibility of the constructions defined above. An interesting and  special case is when $H^*$ is adjoint, the proof of the following is left to the reader.
 
 \medskip\noindent{\sl Let $\sigma\colon (G,F)\to (H,F)$ be an isotypic morphism and $\sigma^*\colon H^*\to G^* $
a dual one, let $s=\sigma^*(t)$.  Let $\bar t$ be the image of $t$ and $s$ in a common adjoint group of $G^*$ and $H^*$ with  isomorphic simply connected duals $H_{\rm sc}$, $G_{\rm sc}$.  One has a commutative diagram of sequences  $$\diag{{\frak E}(G_{\rm sc},\bar t)}{\cong}{{\frak E}(H_{\rm sc},\bar t)}{}{{\frak E}(H,t)}{}{{\frak E}(G,s)}{}{}$$}

\medskip\noindent{\bf 3.3.3. Lemma. } {\sl Assume $G$ quasi-simple and simply connected. Then (B) in Proposition~3.1.1 and (A), (B), (C) in Proposition~3.1.2 are satisfied.}

\preuve Thanks to the studies in in section 3.2, we may assume that $s$ is not rationally quasi-isolated in $G^*$. Let $K^*=\cent{G^*}{\zo{\cento{G^*}s}}$. Let $K$ be a Levi subgroup of $G$ in the dual $G^F$-conjugacy class. One has ${\rm A}_{K^*}(s)^F={\rm A}_{G^*}(s)^F$. From $(K,K^*,s)$ one define an exact short sequence $1\to K(s)^\circ\to K(s)\to {\rm A}_{K^*}(s)\to 1$ and the subsequence $$1\to K(s)^\circ\to K(s)^\circ K(s)^F\to {\rm A}_{K^*}(s)^F\to 1$$ is isomorphic to $$1\to G(s)^\circ\to G(s)^\circ.G(s)^F\to {\rm A}_{G^*}(s)^F\to 1$$ For any $d$-split Levi subgroup $L^*_s$ of $\cento{G^*}s$ one has $\cent{G^*}{\zo{L^*_s}_{\phi_d}}\cap \cent{G^*}s^F=\cent{K^*}{\zo{L^*_s}_{\phi_d}}\cap \cent{G^*}s^F$. The conditions (B) in Proposition~3.1.1 and (A), (B), (C) in Proposition~3.1.2 are equivalent for the two exact sequences and satisfied by 3.2.
\bull

\medskip\noindent{\bf 3.3.4. Lemma. } {\sl Assume $G$  rationally \irr, $G=[G,G]$ and $\cF(G)=\{1\}$. Then (B) in Proposition~3.1.1 and (A), (B), (C) in Proposition~3.1.2 are satisfied.}

\preuve $G$ is a direct product $\prod_{1\leq j\leq k}G_j$ where $G_j$ is \irr\ and $F(G_j)=G_{j+1}$ ($j<k$), $F(G_k)=G_1$. On dual side $G^*=\prod_{1\leq j\leq k}G^*_j$ and $F(G^*_j)=G^*_{j+1}$, $F(G^*_k)=G^*_1$, $s=(s_j)_{1\leq j\leq k}$ and $F(s_j)=s_{j+1}$, $F(s_k)=s_1$. It follows that the short exact sequence 
$$1\to \cento{G^*}s\to \cent{G^*}s\to {\rm A}_{G^*}(s)\to 1\leqno{{\frak E}(G^*,s)}$$
is isomorphic to 
$$1\to \prod_{1\leq j\leq k}\cento{G^*_j}{s_j}\to \prod_{1\leq j\leq k}\cent{G^*_j}{s_j}\to \prod_{1\leq j\leq k}{\rm A}_{G^*_j}(s_j)\to 1\leqno{\prod_{1\leq j\leq k}{\frak E}(G_j^*,s_j)}$$
 and the short exact sequence 
$$1\to {G}(s)^\circ\to G(s)\to {\rm A}_{G^*}(s)\to 1\leqno{{\frak E}(G,s)}$$
is isomorphic to (see (A.1) in Proposition~3.1.1)
$$1\to \prod_{1\leq j\leq k}{G_j}(s_j)^\circ\to \prod_{1\leq j\leq k}{G_j}(s_j)\to \prod_{1\leq j\leq k}{\rm A}_{G^*_j}(s_j)\to 1\leqno{\prod_{1\leq j\leq k}{\frak E}(G_j,s_j)}$$

About groups of rational points one have isomorphisms $G^{F}\cong G_1^{F^k}$, $G^{*\,F}\cong G_1^{*\,{F^k}}$, ${\rm A}_{G^*}(s)^F\cong {\rm A}_{G^*_1}(s_1)^{F^k}$ and  isomorphic extensions $$\big[1\to G(s)^{\circ\, F}\to G(s)^F\to {\rm A}_{G^*}(s)^F\to 1\big]\;\cong\;\big[ 1\to G_1(s_1)^{\circ\, F^k}\to G_1(s_1)^{F^k}\to {\rm A}_{G^*_1}(s_1)^{F^k}\to 1\big]\, ,$$ 
$$\big[1\to \cento{G^*}s^F\to \cent{G^*}s^F\to {\rm A}_{G^*}(s)^F\to 1\big]\;\cong\;\big[1\to \cento{G^*_1}{s_1}^{F^k}\to \cent{G^*_1}{s_1}^{F^k}\to {\rm A}_{G^*_1}(s_1)^{F^k}\to 1\big]\, .$$

Let $d':=d/(k,d)$. Then  $G(s)^{\circ\, F}$-conjugacy classes of  $d$-cuspidal data of $(G(s)^\circ, F)$ correspond to $G_1(s_1)^{F^k}$-conjugacy classes of  $d'$-cuspidal data of $(G_1(s_1)^\circ, F^k)$. The assertions (B) in Proposition~3.1.1 and (A), (B), (C) relative to $(G,F,s)$ and to $(G_1,F^k,s_1)$ are equivalent. 
\bull

\medskip\noindent{\bf 3.3.5. Lemma. } {\sl Assume $\cF(G)=\{1\}$ and $G$ rationally \irr. Then (B) in Proposition~3.1.1 and (A), (B), (C) in Proposition~3.1.2 are satisfied.}

\preuve (i) Assume $\zo G\cap[G,G]=\{1\}$. 

Then $G$ is a direct product : $G=\zo G\times [G,G]$ and $G^*=\zo G^*\times [G,G]^*$. Any object we consider in Propositions~3.1.1 and 3.1.2 decomposes in a direct product. For a torus there is nothing to prove, the result follows from the preceding Lemma.

(ii) Assume $\zo G\cap[G,G]\neq \{1\}$. 

One has $\z G=\zo G.\z{[G,G]}$ hence $\z G/\zo G$ is cyclic, so is ${\rm A}_{G^*}(s)$ and the assertions (A) of Proposition 3.1.1 and (C) of Proposition 3.1.2 are satisfied. In types $\BB$, $\CC$ and exceptionnel types, $\z G=\zo G$ and ${\rm A}_{G^*}(s)=\{1\}$, there is nothing to prove. 

Denote $H=[G,G]$. A morphism in duality with the inclusion $H\subseteq G$ is $i^*\colon G^*\to (G^*)_{\rm ad}$. Then ${\rm A}_{G^*}(s)$ is a subgroup of ${\rm A}_{H^*}(i^*(s))$. (A.3) in Proposition~3.1.1 describes the relation between ${\cal E}(G,s)$ and ${\cal E}(H,t)$ :
$${\def\normalbaselines{\baselineskip20pt\lineskip3pt
\lineskiplimit3pt}\matrix{{{\cal E}(H,s)}&{}&1&\to&H(t)^\circ&\mapright{}&H(s)&\mapright{}&\alpha({\rm A}_{G^*}(s))&\to&1\cr
{\mapdown{}}&{}&{}&{}&\mapdown{}&&\mapdown{\sigma_s}&&\mapdown{\cong}&{}&{}\cr 
{{\frak E}(G,s)}&\quad&1&\to&G(s)^\circ&\mapright{}&G(s)&
\mapright{ }&{{\rm A}_{G^*}(s)}&\to&1\cr} }$$
As $\cento{G^*}s\to\cento{H^*}t$ is onto, $\sigma_s$ is injective.  Let $t$, $\rho_H$, $(L_H(t),\al_H)$, $L^*_t$, $L^*_H$, $L_H$, ... be defined as $s$, $\rho$, $(L(s),\al)$, $L^*_s$, $L^*$, $L$, $D$, $\al_D$, $\beta_D$, $s_1$,  one has $$G=\zo G.H,\;{L(s)}=\zo G.\sigma_s(L_H(t)), \; \al_H=\Res{{L(s)}^F}{L_H(t)^F}\al,\;L_t^*=i^*(L^*_s),\;L_H^*=i^*(L^*)\,.$$  
As $D$ is an $\ell$-Sylow subgroup of $\cento{G(s)^\circ}{[{L(s)},{L(s)}]}^F$,  $D\subseteq \zo G.\sigma_s(D_H)$, $\sigma_s(D_H)=D\cap [G,G]$. With notations of section 3.2.1, $\pi_{t,H}=i^*\circ \pi_s\circ \sigma_s$. 

From (3.2.1.1) $[L^*_H,L^*_H]\subseteq \cento{H^*}{\pi_{s,H}(D_H)}$  one deduces $[L^*,L^*]\subseteq\cento{G^*}{\pi_s(D)}$ (Proposition~1.2.3). That gives (A) for $(G,s)$ and $\cent{G(s)^F}D=\sigma_s(\cent {H(s)^F}{D_H})$ and we may assume that $\beta_{D,H}$ is the restriction of  $\beta_D$. One has $\nor{G(s)}D=\zo G.\sigma_s(\nor{H(s)}{D_H})$ so that condition (B) for $H$ implies condition (B) for $G$.
\bull

 \medskip\noindent{\bf 3.3.6. Lemma }{\sl Let $\sigma\colon (H,F)\to (G,F)$ be an isotypic morphism whose kernel is a central torus and such that $\sigma(H)=G$, let $\sigma^*\colon G^*\to H^* $ be
a dual one, and let $t=\sigma^*(s)$. If Condition (B) in Proposition~3.1.1 and Proposition 3.1.2 are true with $(H,t)$ instead of $(G,s)$, then they are true for  $(G,s)$.}

 \medskip\noindent{\it Proof . } Denote $\rho_G$, $\rho_H$ the two morphisms defined as $\rho$ in Proposition~3.1.2. One has  $H^*=\zo {H^*}.\sigma^*(G^*)$, where $\sigma^*$ is an embedding, hence $\alpha\colon {\rm A}_{G^*}(s)\to {\rm A}_{H^*}(t)$ is an isomorphism. Let $K$ be the kernel of $\sigma$. By construction  $K\subseteq H(t)$, $G(s)$ is isomorphic to $H(t)/K$ and $G(s)^F$ to $H( t)^F/K^F$. 
 
 Let $s_1$ be semi-simple in $ \cento{G^*}s^F$ and $t_1:=\sigma^*(s_1)\in\cento{H^*}t^F$. If $s$ and $s_1$ have co-prime orders, $t$ and $t_1$ have co-prime orders. One has 
 $$\cento{\cento{H^*}t}{t_1}=\zo{H^*}.\sigma^*(\cento{\cento{G^*}s}{s_1}), \quad\cent{\cento{H^*}t}{t_1}=\zo{H^*}.\sigma^*(\cent{\cento{G^*}s}{s_1})\, .$$ 
 So, {\it via} the restriction of $\sigma^*$ to $\cento{G^*}s^F$, ${\rm A}_{\cento{G^*}s}(s_1)$ is isomorphic to ${\rm A}_{\cento{H^*}t}(t_1)$. 
 
By Proposition~1.3.6, a one-to-one map from $\ser{H(t)^{\circ \,F}}{t_1}$ onto $ \ser{G(s)^{\circ \,F}}{s_1}$ is induced by the  restriction from $H(t)^{\circ \,F}$ to $G(s)^{\circ \,F}$. As $G(s)^F=H(t)^F/K^F$, if condition (B) of Proposition~3.1.1 is true for $(H,t)$ it is true for $(G,s)$ (note that for any $\xi\in\ser{H(t)^{\circ\, F}}{\sigma^*(s_1)}$, the kernel of $\xi$ contains  $K^F$). Thus
 if it is satisfied for $(H,F,t,\sigma^*(s_1))$ it is satisfied for $(G,F,s,s_1)$.

 A unipotent $d$-cuspidal datum $(L(s),\al)$ in $(G(s)^\circ,F)$ is "image" by $\sigma$ of a unipotent $d$-cuspidal datum $(M(t),\beta)$ in $(H(t),F)$. Here $M(t)=\sigma^{-1}(L(s))$, so  that $\sigma$ induces a bijection $\ser{M(t)^F}1\to \ser{L(s)^F}1$ and we may write $\al=\beta$. 
 
 A defect group $D$ of $b_{G(s)^{\circ F}}(L(s),\al)$ is an $\ell$-Sylow subgroup of $\cento{G(s)}{[L(s),L(s)]}^F$. We have $ D=\sigma(E)$  where $E$ is an $\ell$-Sylow subgroup of $\cento{H(t)}{[L(t),L(t)]}^F$, that is a defect group of $b_{H^(t)^{\circ F}}(L(t),\al)$. Then  $\sigma(\nor{H(t)^F}E)\subseteq \nor{G(s)^F}D$ and  $\sigma(\cent{H(t)^F}E)\subseteq \cent{G(s)^F}D$, hence $\rho_H
(\nor{H(t)^F}E)\subseteq\al(\rho_G(\nor{G(s)^F}D) $ and $\rho_H
(\cent{H(t)^F}E)\subseteq \al(\rho_G(\cent {G(s)^F}D)$.

Let $a\in\rho_G(\nor{G(s)^F}D)$. There exists $g\in \nor{G(s)^F}D$ of image $a$ and order prime to $\ell$ (because the order of $a$ is prime to $\ell$). Hence there exists $g\in H(s)^F$ of order prime to $\ell$ and such that  $\sigma(h)=g$, $[h,E]\subseteq K^F.E$. But $E\cap K^F=K^F_\ell$ and $K^F$ is central in $H(s)^F$, so $[h,E]\subseteq E$. If $a\in\rho_G(\nor{G(s)^F}D)$ then we obtain $[h,E]\subseteq K^F_\ell$. But  $E$ is an $\ell$-group and $h$ an $\ell'$-element, so this implies $[h,E]=\{1\}$. We get 
 $$\nor {G(s)^F}D=\nor{H(t)^F}E/K^F,\quad \cent {G(s)^F}D=\cent{H(t)^F}E/K^F\leqno{(3.3.6.1)}$$
 $$\alpha(\rho_G(\nor {G(s)^F}D)= \rho_H
(\nor{H(t)^F}E),\quad \alpha(\rho_G(\cent {G(s)^F}D)= \rho_H
(\cent{H(t)^F}E)\leqno{(3.3.6.2)}$$
 On the other hand if $L^*_H:=\cent{H^*}{\z {L_H^*}_{\phi_d}}$ and $L^*_G:=\cent{G^*}{\z {L_G^*}_{\phi_d}}$ the isomorphism $\alpha$ gives $$\alpha({\rm A}_{L_G^*}(s))={\rm A}_{L^*_H}(t)\leqno{(3.3.6.3)}$$
 By (3.3.6.2) and (3.3.6.3) the assertions (A) in Proposition~3.1.2 for $(G,s)$ and for $(H,t)$ are equivalent.
 
 By (3.3.6.1) we have $\cent {G(s)^\circ}D^F=\cent{H(t)^\circ}E^F/K^F$. As  $\al_D$ is defined by restriction from $\al$,  the canonical $\beta_E\in \II{\cent{H(t)^F}E}$ is defined by $\beta=\al$. This shows that $\al_D$ is just the restriction {\it via} $\sigma$ of $\beta_E$. Thus any $\chi_0\in\II{b_D}$ is restriction of some $\xi_0\in\II{b_E}$ (evident notations). The assertions (B) of Proposition 3.1.2 for $(G,s)$ and for $(H,t)$ are equivalent.
 If the non-multiplicity condition (C) in Proposition~3.1.2. is true for $(H,t)$ it is true for $(G,s)$. 
 \bull

 \medskip\noindent{\bf 3.3.7. End of proof }
\preuve (a) Assume $\cF(G)=\{1\}$. Then  the \dec\ of $[G,G]$ in rationally \irr\ components is a direct product $[G,G]=\prod_jG_j$. Let $H_j=\zo G.G_j.$ and $H=\prod_jH_j$. There is a natural morphism 
$\sigma\colon H\to G$ whose kernel is a torus and such that  $\sigma(H)=G$. One concludes by 3.3.1 and Lemmas 3.3.5, 3.3.6.

(b)
Given $G$ without special assumption, let $\sigma\colon H\to G$ a dual morphism of a regular embedding $G^*\to H^*$. One has $\cF(H)=\{1\}$ and the kernel of $\sigma$ is a torus. One concludes by (a) and Lemma~3.3.6.
\bull 

\bigskip
\noindent{\bf 3.4. Jordan \dec\ on blocks. }

In this section we prove (B.1) of our main theorem 1.4. We have defined in 1.3.1 $\bl{G^F}s$ as the set of blocks $B$ of $G^F$ such that $\II b\cap\ser{G^F}s\neq \emptyset$ ($G$ connected), with the special case $\bl{G^F}1$ when $G$ is not connected.

\medskip\noindent{\bf 3.4.1. Proposition. }
{\sl Let $(G,F,\ell,d)$ with assumption 2.2.3, $\sigma\colon (G,F)\to (H,F)$ a regular embedding, $\sigma^*\colon H^*\to G^*$ a dual morphism. Let $t$ be a semi-simple element of $(H^*)^F_{\ell'}$ and $s=\sigma^*(t)$. 
Let  $$1\to G(s)^\circ\to G(s)\to {{\rm A}_{G^*}(s)}\to 1 \leqno{{\frak
E}(G,s)}$$ as constructed in Proposition~3.1.1.

 By Propositions~2.1.4 and 2.1.7 is defined a one-to-one map  $$ {\cal B}_{H,t}\colon \bl{H(t)^{F}}1\to\bl{H^F}t$$ ${\cal B}_{H,1}$ is identity.

 By compatibility with Clifford theory of blocks there exist  a one-to-one map $${\cal B}_{G,s}\colon\bl{G(s)^F}1\to \bl{G^F}s$$  such that 

 \noindent (3.4.1.1))$\quad$ If a block $b$ of $G(s)^F$ covers a block $b_0$ of $G(s)^{\circ\, F}$ and $B$ is the unique unipotent block of $H(t)^{F}$ that covers $b_0$ through the morphism $G(s)^{\circ \,F}\to H(t)^{F}$, then ${\cal B}_{G,s}(b)$ is covered by ${\cal B}_{H,t}(B)$.  

}

\preuve  ${\cal B}_{H,t}$ is defined by   composition of the following one-to-one maps 

$B=b_{H(t)^{F}}(M(t),\al_t)$, $(M(t),\al_t)$, a unipotent $d$-cuspidal datum in $(H(t),F)$ 

$\mapsto$ $H(t)^{F}$-conjugacy class of  $(M(t),\al_t)$ (Proposition~2.1.7)

$\mapsto$ $\cent{H^*}t^F$-conjugacy class of $(M^*_t,\al)$, a unipotent $d$-cuspidal datum in $(\cent{H^*}t,F)$,  (Proposition~2.1.4, $\al_t=\Psi_{M(t),1}(\al)$, $M^*_t$ in duality with $M(t)$)

$\mapsto$ $H^F$-congugacy class of $d$-cuspidal data $(M,\mu)$ in series $(t)$ in $(H,F)$ (Proposition~2.1.4, $\mu=\Psi_{M,t}(\al)$, $M^*\cap\cent{H^*}t=M^*_t$)

$\mapsto $ $b_{H^F}(M,\mu)$ (Proposition~2.1.7).

The isotypic morphism $\sigma_t\colon (G(s)^\circ,F)\to (H(t),F)$ is an embedding, because the kernel of the dual morphism $\cent{H^*}t\to \cento{G^*}s$ is a torus, and it induces by restriction a one-to-one map $$\Res{}{\sigma_t}\colon \ser{H(t)^{F}}1\to\ser{G(s)^{\circ\, F}}1$$ hence a one-to-one map between conjugacy classes of $d$-cuspidal unipotent data (Proposition~2.1.2) 
$$(M(t),\al_t)\mapsto (\sigma^{-1}(M(t)),\Res{}{\sigma_t}(\al_t))$$ hence a bijection between unipotent blocks (Proposition 2.1.7) $b_0\mapsto B$ such that 
$$\Res{}{\sigma_t}(\II{B}\cap \ser{H(t)^{F}}1)=\II{b_0}\cap \ser{G(s)^{\circ\, F}}1$$.

The existence of ${\cal B}_{G,s}$ for any $G$ under condition (3.4.1) follows from a combinatorial fact : given corresponding unipotent blocks $b_0$ and $B$ of $G(s)^{\circ \,F}$ and  $H(t)^F$ respectively, let

 $\bullet$ $m_0$ be the number of blocks of $G(s)^F$ that  cover  $b_0$,
 
 $\bullet$ $n_0$ be the number of $G(s)^F$-conjugate of $b_0$,
 
$\bullet$  $m$ be the number of blocks in series $(t)$ of $H^F$ that are in the $\II{H^F/G^F}$-orbit of ${\cal B}_{H,t}(B)$,
 
 $\bullet$ $n$ be the number of blocks in series $(s)$ of $G^F$ that are covered by ${\cal B}_{H,t}(B)$.

Indeed one has $m=n_0$ and $n=m_0$.

Assume $b_0=b_{G(s)^{\circ F}}(L(s),\al)$, hence, if $B=b_{H(t)^{F}}(M(t),\al)$ as above, $L(s)$ is the inverse image of $M(t)$ by the morphism $G(s)^\circ\to H(t)$ defined from $\sigma$. By Propositions~3.1.3 and 2.4.4 we have

$$m_0=| {\rm A}_{L^*}(s)^F_\al |=n, \quad n_0= | {\rm A}_{G^*}(s)^F/{\rm A}_{G^*}(s,L^*)_\al | = m$$
To the orbit of $\al$ under ${\rm A}_{G^*}(s)^F$ there correspond a set $\bl{H^F}t_{(\al)}$ of $m$ blocks of $H^F$ in series $(t)$ covering $n$ blocks in series $(s)$ of $G^F$. By the first part of the proof, $\bl{H^F}t_{(\al)}={\cal B}_{H,t}(\bl{H(t)^{F}}1_{(\al)})$, where $\bl{H(s)^{F}}1_{(\al)} $ is a set of $m$ unipotent blocks of $H(s)^{F}$, hence in one-to-one map with a set $\bl{G(s)^{\circ F}}s_{(\al)}$ of $m=n_0$ unipotent blocks of $G(s)^{\circ \,F}$, covered by $n=m_0$ (unipotent) blocks of $G(s)^F$.
\bull

\medskip\noindent{\bf 3.4.2. Proposition. }{\sl Hypothesis and notations of Proposition~3.4.1 on $(G,F,\ell,d,\sigma,H,t,s)$. Let $t_0$ in $\cent{H^*}t^F_\ell$ and $\sigma^*(t_0)=s_0$ in $\cento{G^*}s^F_\ell$. Let $B\in\bl{H(t)^F}1$, $b\in\bl{G(s)^{F}}1$. Define   $\ser{G(s)^F}{(s_0)}$ by
$$\ser{G(s)^F}{(s_0)}=\big\{\mu\mid\mu\in \II{G(s)^F\mid \la} \;{\sl for}\;{\sl some}\;\la\in\ser{G(s)^{\circ\, F}}{s_0}\big\}$$ 

There exist  one-to-one maps $\Psi_{H,B}$ from $\II{B}$ onto $\II{{\cal B}_{H,t}(B)}$ and $\Psi_{G,b}$ from $\II{b}$ onto $\II{{\cal B}_{G,s}(b)}$ such that,  
$${\sl for \;any\;} t_0\in\cent{H^*}t^F_\ell,\quad \Psi_{H,B}\big(\II B\cap \ser{H(t)^F}{t_0}\big)=\II{{\cal B}_{H,t}(B)}\cap\ser{H^F}{tt_0}\leqno{(3.4.2.1)}$$
$$\II{{\cal B}_{G,s}(b)}\cap \ser{G^F}{ss_0} \subseteq {\Psi_{G,b}}\big(\II b\cap \ser{G(s)^F}{(s_0)}\big) \leqno{(3.4.2.2)}$$
$${\sl for \;any\;} \mu\in\II b, \quad \Psi_{G,b}(\mu)(1).|G(s)^F|_{p'}=\mu(1).|G^F|_{p'}\leqno{(3.4.2.3)}$$}

\preuve There may be fusion  under $G(s)^F$ of Lusztig series in $G(s)^{\circ\, F}$, that's why we use a different notation  in (3.4.2.2), see part (B) of the proof.

Let be $b_0$, $b$, $B=b_{H(t)^{F}}(M(t),\al_t)$, $\al_t=\Psi_{M(t),1}(\al)$  as in Proposition 3.4.1 and its proof. 

(A) We first consider the case $G=H$.

We have  partitions [9], [16] Theorem 9.12 
$$\II{B}=\cup_{(t_0)}\big(\II{B}\cap \ser{H(t)^{F}}{t_0}\big),\quad \II{{\cal B}_{H,t}(B)}=\cup_{(h)}\big(\II{{\cal B}_{H,t}(B)}\cap \ser{H^F}{h}\big)$$
In the left equality $(t_0)$ runs on the set of $\cent{H^*}t^F$-conjugacy classes of $\ell$-elements of $\cent{H^*}t^F$. In the right one $(h)$ runs on the set of $H^{* \,F}$-conjugacy classes of semi-simple elements $h$ of $H^{* F}$ whose $\ell'$-component $h_{\ell'}$ is $H^{* \, F}$-conjugate of $t$. The map $(t_0)\mapsto (tt_0)$ is a one-to-one map between these two sets of conjugacy classes. 

To prove (3.4.2.1) we  define the restriction of $\Psi_{H,B}$  $$ \II{B}\cap \ser{H(t)^{F}}{t_0}\to \II{{\cal B}_{H,t}(B)}\cap \ser{H^F}{tt_0}$$ These two sets are described in Proposition~2.3.5. They are in bijection with  sets of components of two Generalized $d$-H.C. series : replace in 2.3.5 $(G,s,s_0,L,\la,\al,\al_0)$  by $(H(t),1,t_0,M(t),\Psi_{M(t_0),t}(\al),\al,\al(t_0))$ and then, with $\mu=\Psi_{M,t}(\al)$,  by $(H,t,t_0,M,\mu,\al,\al(t_0))$ .

On left side, unipotent block side, the source of $\Psi_{H,B}$, $\II{B}\cap \ser{H(t)^{F}}{t_0}$, if not empty, is in bijection with $\ser{H(t)(t_0)^F}{(M(t)(t_0),\al({t_0}))}$, where 
 $H(t)(t_0)$ is a Levi subgroup of $H(t)$ in the dual $H(t)^{F}$-conjugacy class of $\cent{\cent{H^*}t}{t_0}=\cent{H*}{tt_0}$  and 
 $(M(t)(t_0),\al({t_0}))$ is a $d$-cuspidal unipotent datum in $(H(t)(t_0),F)$,  associated by duality to a $d$-cuspidal unipotent datum $(M^*_{t t_0},\al(t_0))$ in $(\cent{\cent{H^*}t}{t_0},F)$ such that, $M^*_t$ being in the dual $\cent{G^*}t^F$-conjugacy class of the $H(t)^F$-conjugacy class of $M(t)$,  $(M^*_{tt_0},\al(t_0))\simex{\cent{H^*}t^F}(M^*_t,\al)$.  The map is $\Lu{H(t)(t_0)}{H(t)}(\Psi_{H(t)(t_0),t_0}(1)\otimes -)$
 $$ \nu\mapsto \xi=\Lu{H(t)(t_0)}{H(t)}(\Psi_{H(t)(t_0),t_0}(1)\otimes \nu), \quad \nu(1).|\cent{H^*}t^F|_{p'}=\xi(1).|\cent{H^*}{tt_0}^F|_{p'}\leqno{(3.4.2.4)}$$   
the equality thanks to  (1.3.1.2) and knowing that by duality   $|H(t)^{F}|=|\cent{H^*}t^F|$, $|H(t)(t_0)^F|=|\cent{H^*}{tt_0}^F|$.

On right side,   $\II{{\cal B}_{H,t}(B)}\cap \ser{H^F}{tt_0}$, if not empty, is in bijection with $\ser{H(t_0)^F}{(M(t_0),\mu(t_0))}$, where $H(t_0)$ is a Levi subgroup of $H$ in the dual conjugacy class of $\cent{H^*}{t_0}$, $(M(t_0),\mu(t_0))$ is a $d$-cuspidal datum in series $(t)$ in $(H(t_0),F)$ associated by Proposition~2.1.4 to a $d$-cuspidal unipotent datum $(M^*_{t,t_0},\al_t(t_0))$ in $(\cent{H^*}{tt_0},F)$ such that $(M^*_{t,t_0},\al_t(t_0))\simex{\cent{H^*}t^F}(M^*_t,\al)$. Clearly, as $\cent{\cent{H^*}t}{t_0}=\cent{H^*}{tt_0}$, the existence of such a $d$-cuspidal datum $(M^*_{tt_0},\al_t(t_0))$ is equivalent to the preceding one $(M^*_{t,t_0},\al(t_0))$ and we may assume 
$(M^*_{t,t_0},\al_t(t_0))=(M^*_{tt_0},\al(t_0))$. The map is $\Lu{H(t_0)}H(\Psi_{H(t_0),t_0}(1)\otimes -)$ 
$$ \zeta\mapsto \eta=\Lu{H(t_0)}{H}(\Psi_{H(t_0),t_0}(1)\otimes \zeta) ,\quad \zeta(1).|H^F|_{p'}=\eta(1).|\cent{H^*}{t_0}^F|_{p'}\leqno{(3.4.2.5)}$$   
by (1.3.1.2) and equality $|H(t_0)^F|=|\cent{H^*}{t_0}^F|$.

By Proposition~2.2.4 there exist  one-to-one maps 
$$\tilde\Psi_{H(t)(t_0),1}(M(t)(t_0),\al(t_0))\colon \ser{\cent{H^*}{tt_0}^F}{(M^*_{tt_0},\al(t_0))}\to\ser{H(t)(t_0)^F}{(M(t)(t_0),\al(t_0))}$$
$$\tilde\Psi_{H(t_0),t}(M(t_0),\mu(t_0))\colon\ser{\cent{H^*}{tt_0}^F}{(M^*_{tt_0},\al(t_0))}\to\ser{H(t_0)^F}{(M(t_0),\mu(t_0))}$$ 
The restriction of $\Psi_{H,B}$ we are looking for  is defined by  the commutation formula

$ \Lu{H(t_0)}{H}(\Psi_{H(t_0),t_0}(1)\otimes -) \circ \tilde\Psi_{H(t_0),t}(M(t_0),\mu(t_0))$\hfill

\hfill
 $=\Psi_{H,B}\circ \Lu{H(t)(t_0)}{H(t)}(\Psi_{H(t)(t_0),t_0}(1)\otimes -)\circ \tilde\Psi_{H(t)(t_0),1}(M(t)(t_0),\al(t_0))$.

From (2.2.4.1) we have, if $\nu_0\in\ser{\cent{H^*}{tt_0}^F}{(M^*_{tt_0},\al(t_0))}$, $\nu=\tilde\Psi_{H(t)(t_0),1}(M(t)(t_0),\al(t_0))(\xi_0)$, 

\noindent $\nu(1)=\nu_0(1)$ (trivial case, $t$ is central in the dual of $H(t)(t_0)$) and if $\zeta= \tilde\Psi_{H(t_0),t}(M(t_0),\mu(t_0))(\nu_0)$, 
$$\zeta(1).|\cent{H^*}{tt_0}^F|_{p'}=\nu(1).|\cent{H^*}{t_0}^F]_{p'}.$$
With notations of (3.4.2.4) and (3.4.2.5) we obtain $\eta(1).|\cent{H^*}t^F|_{p'}=\xi(1).|H^F|_{p'}$ whereas $\Psi_{H,B_0}(\xi)=\eta$ and $|\cent{H^*}t^F|=|H(t)^F|$, that is (3.4.2.3) for $(H,t)$. 

(B) Assuming now $\z G$ non-connected we may assume that type $\EE_8$ don't appear in $G$ (see first lines of 2.4 or 3.2.1). Then by the proof of  Proposition~2.2.4 the one-to-one map $\tilde\Psi_{H(t)(t_0),1}(M(t)(t_0),\al(t_0))$ (resp. $\tilde\Psi_{H(t_0),t}(M(t_0),\mu(t_0))$) is just the restriction of Jordan decomposition $\Psi_{H(t),t_0}$ (resp. $\Psi_{H(t_0),t_0}$). Hence $\xi=\Psi_{H(t),t_0}(\beta)$ for some $\beta\in \ser{\cent{H^*}{tt_0}^F}{(M^*_{tt_0},\al(t_0))}$ and $\Psi_{H,B_0}(\xi)=\Psi_{H,tt_0}(\beta)$.

To prove the existence of $\Psi_{G,b}$ we proceed as in the proof of the preceding Proposition, using properties of the functions $\Psi$. We go from $(H(t)^{F},H^F)$ to $(G(s)^F,G^F)$ in three steps through three restrictions {\it via} 
$$G^F\to H^F,\quad G(s)^{\circ\, F}\to H(t)^{F},\quad G(s)^{\circ\, F}\to G(s)^F$$
Thus we rely $\Psi_{G,b}$  to  $\Psi_{H,B}$ by the condition

\noindent(3.4.2.6) $\quad$ If $\xi\in\II B$ and $\mu\in\II b$ cover $\la\in\II{G(s)^{\circ\, F}}$,  then $\Psi_{H,B}(\xi)$ covers $\Psi_{G,b}(\mu)$.

The three  morphisms above have  common properties : the kernel is $\{1\}$, the image is an invariant subgroup, the cokernel is abelian and there is no multiplicity in restrictions of \irr\ representations (Propositions~1.3.4 and 3.1.1, (B)). In such a morphism $X\to Y$ there is a bijection between sets of orbits $\II Y/(Y/X)^\wedge \leftrightarrow \II X/(Y/X)$ induced by $\eta\in\II {Y\mid \chi}$, where $\chi\in\II X$. Then $\chi$ is covered by $|Y_\chi/X|$ elements of $\II Y$ and $\eta$ covers $|Y/Y_\chi|$ elements of $\II X$.  We say that the ``multiplicative factor" on the number of \irr\ in corresponding orbits from $\II Y$ to $\II X$ is $|Y/X|/|Y_\chi/X|^2$.

The three maps restrict to series.

Consider first the stabilizer in $G(s)^F$ of $\ser{G(s)^{\circ \,F}}{s_0}$, or of the $\cento{G^*}s^F$-conjugacy class of $s_0$ in $\cent{G^*}s^F$. It has image $\cent{\cent{G^*}s}{s_0}^F.\cento{G^*}s/\cento{G^*}s$ in ${\rm A}_{G^*}(s)^F$. We have $\cent{G^*}{ss_0}=\cent{G^*}s\cap \cent{G^*}{s_0}$ as well as $\cento{G^*}{ss_0}=\cento{G^*}s\cap \cento{G^*}{s_0}$ (Proposition~1.2.4). By isomorphism theorem ${\rm A}_{G^*}(s)^F_{(s_0)}$ is isomorphic to $\cent{G^*}s^F\cap\cent{G^*}{s_0}^F/\cento{G^*}s^F\cap \cent{G^*}{s_0}^F$ and   that quotient is  the component ${\rm A}_{\cento{G^*}{s_0}}(s)^F$ of ${\rm A}_{G^*}(ss_0)^F$ we obtained in Proposition 1.2.4. 

A block $b$ of $G(s)^F$ covers a block $b_0$ of $G(s)^{\circ\, F}$ if and only if some element of $\II b$ covers an element of $\II{b_0}$. When a block $b$ of $G(s)^F$ covers $b_0$ and a block $B$ of $H(t)^F$  covers $b_0$, then ${\cal B}_{H,t}(B)$ covers ${\cal B}_{G,s}(b)$ and  there are similar partitions  through series. The three restrictions above  send series in series, so we consider separately sets of \irr\ components 
$$\eta\in \II{{\cal B}_{H,t}(B)}\cap\ser{H^F}{tt_0}\mapsto  \{\chi\in \II{{\cal B}_{G,s}(b)}\cap \ser{G^F}{ss_0}\mid \eta\in\II{H^F\mid \chi}\}$$
$$\xi\in\II{B}\cap\ser{H(t)^{F}}{t_0}\mapsto \{\la\in   \II{b_0}\cap \ser{G(s)^{\circ\, F}}{s_0}\mid \xi\in\II{H(t)^F\mid \la}$$
$$\mu\in\II{b}\cap\ser{G(s)^F}{(s_0)}\mapsto\{\la\in \II{b_0}\cap \ser{G(s)^{\circ\, F}}{s_0}\mid \mu\in\II{G(s)^F\mid\la}\}$$
Thanks to  the definition of $\Psi_{H,B}$, to any  $\eta\in \II{{\cal B}_{H,t}(B)}\cap\ser{H^F}{tt_0}$, or  
$\xi\in \II{B}\cap\ser{H(t)^{F}}{t_0}$,  we may associate some $\beta\in\ser{\cent{H^*}{tt_0}^F}1$ such that $\xi=\Psi_{H(t),t_0}(\beta)$  and $\Psi_{H,B}(\xi)=\Psi_{H,tt_0}(\beta)$. Then, by Propositions~1.3.6,  1.3.7, any component of $ \Res{H^F}{G^F}\eta$ or of $\Res{H(t)^{F}}{G(s)^{\circ F}}\xi$ is associated to an orbit of $\beta$ under  ${\rm A}_{G^*}(ss_0)^F$. To compute the effect of the three restrictions on the number of elements  in  the sets of \irr, and their degrees, we may choice the orbit of $\beta$ and apply 1.3.6, 1.3.7.

(i) From $H^F$ to $G^F$ :

If $\eta=\Psi_{H,tt_0}(\beta)$, $\chi\in\II{G^F}$ and $\eta\in \II{H^F\mid \chi}$, $|{\rm A}_{G^*}(ss_0)_\beta ^F|$ $H^F$-conjugate of $\chi$ are covered by $|{\rm A}_{G^*}(ss_0)^F/{\rm A}_{G^*}(ss_0)_\beta^F| $ elements of $\ser{H^F}{tt_0}$. Therefore $$\chi(1).|{\rm A}_{G^*}(ss_0)^F_\beta |=\eta(1)\leqno{(3.4.2.7)}$$ and  the ``multiplicative factor of $\Res{H^F}{G^F}$ above $\beta$" is $$m_\beta (G^F\to H^F  )={|{\rm A}_{G^*}(ss_0)^F_\beta |^2/|{\rm A}_{G^*}(ss_0)^F|}\leqno{(3.4.2.8)}$$

(ii) From $H(t)^{F}$ to $G(s)^{\circ\,  F}$ :

By Proposition~1.2.4 $\cF(\cent{H^*}t)$ is prime to $\ell$ so that $\cent{\cent{H^*}t}{t_0}$ is connected. By Proposition~1.3.7, if $\xi=\Psi_{H(t),t_0}(\beta)$ as in (A) and $\xi$ covers $\la\in\ser{G(s)^{\circ F}}{s_0}$, then $\la$ has $|{\rm A}_{\cento{G^*}s}(s_0)^F_\beta |$ $H(t)^{F}$-conjugates, so that 
$$\la(1).|{\rm A}_{\cento{G^*}s}(s_0)^F_\beta |=\xi(1)\leqno{(3.4.2.9)}$$
and $\la$ is covered by $|{\rm A}_{\cento{G^*}s}(s_0)^F|/|{\rm A}_{\cento{G^*}s}(s_0)^F_\beta |$ elements of $\ser{H(t)^{F}}{t_0}$. The multiplicative factor of $\Res{H(t)^{F}}{G(s)^{\circ F}}$ above $\beta$ is $$m_\beta(G(s)^{\circ\, F}\to H(t)^F)=|{\rm A}_{\cento{G^*}s}(s_0)_\beta ^F|^2/|{\rm A}_{\cento{G^*}s}(s_0)^F |\leqno{(3.4.2.10)}$$

(iii) From $G(s)^F$ to $G(s)^{\circ \,F}$ :

We need to compute the stabilizer in $G(s)^F$ of $\la\in\ser{G(s)^{\circ F}}{s_0}$. 
To study the action of ${\rm A}_{\cento{G^*}{s_0}}(s)^F$ on $\ser{G(s)^{\circ\,  F}}{s_0}$ we consider a regular embedding 
$$\sigma_s\colon (G(s)^\circ,F) \to ( {\tilde H}(\tilde t),F)= (T\times_{\z {G(s)^\circ} }G(s)^\circ,F)$$ where $T$ is the maximal $F$-stable torus of $G(s)^\circ$ giving rise to the dual root data $\FDD(G^\circ(s),T)$, $\FDD(\cento{G^*}s,T^*)$ that define $G(s)^\circ$ as a dual of $\cento{G^*}s$. The group ${\rm A}_{G^*}(s)$ acts on $T$, $G(s)^\circ$, $\z{G(s)^\circ}$. Hence ${\rm A}_{G^*}(s)$ acts  on $\FDD(\tilde H,\tilde T)$, where $\tilde T=T\times _{\z{G(s)^\circ}}T$, by transport via $\FDD(\sigma_s)$ of its action on $\FDD(G(s)^\circ,T)$. Let $\sigma^*_s\colon \tilde H^*\to \cento{G^*}s$ a  morphism dual of $\sigma_s$. There exist  $\tilde t_0\in (\tilde H^*)^F_\ell$ such that $\sigma^*_s(\tilde t_0)=s_0$ and ${\rm A}_{\cent{G^*}{s_0}}(s)^F$ fixes $\tilde t_0$ : if $\theta=\Psi_{T,s_0}(1)$, define $\tilde t_0$ by $\Psi_{\tilde H,\tilde t_0}(1)=\theta^{-1}\otimes_{\z {G(s)^\circ}^F}\theta$.

 Acting on $\FDD(\tilde H,\tilde T)$, and on $\FDD(\tilde H^*,\tilde T^*)$ by transposition, ${\rm A}_{G^*}(s)$ has an image in the groups of outer automorphisms of $\tilde H$ and $\tilde H^*$. By construction for any $a\in{\rm A}_{G^*}(s)^F$ there exist a  pair $(\tau(a), \tau^*(a))$ of dual automorphisms  of $\tilde H$ and $\tilde H^*$ with  image $(a,\lexp {tr}a)$ such that $\tau(a)$ and $\tau^*(a)$ commute with $F$, $\tau(a)$ (resp. $\tau^*(a)$) stabilizes $G(s)^\circ$ (resp $\cento{G^*}s$). The stabilizer of $\ser{\tilde H(t)^F}{\tilde t_0}$ in ${\rm A}_{G^*}(s)^F$ is ${\rm A}_{\cento{G^*}{s_0}}(s)^F$.
 The isotypic morphism $\sigma^*_{s,\tilde t_0}\colon\cent{\tilde H^*}{\tilde t_0}\to\cento{\cento{G^*}s}{s_0}$ allows us to identify unipotent series, so that $\beta\in \ser{\cent{\tilde H^*}{\tilde t_0}^F}1$. To $\beta$ there correspond  $|{\rm A}_{\cento{G^*}s}(s_0)_\beta^F|$ elements of $\ser{G(s)^{\circ \,F}}{s_0}$, the components of $\Res{\tilde H^F}{G(s)^{\circ F}} (\Psi_{\tilde H,\tilde t_0}(\beta))$ (Proposition 1.3.6) , where is $\la$. 
 
 By (iv) in Proposition~1.3.2, if $a\in {\rm A}_{\cento{G^*}s}(s_0)^F$ and $\lexp a \beta$ is the restriction of $\beta$ to $\cento{\tilde H^*}{t_0}^F$ {\it via} $\tau^*(a)$ then $\Psi_{\tilde H,\tilde t_0}(\beta)$ is the restriction of $\Psi_{\tilde H,\tilde t_0}(\lexp a \beta)$ {\it via} $\tau(a)$. As $(\tau(a),\tau^*(a))$ restricts to  $(G(s)^{\circ\,  F},\cento{G^*}s^F)$, $\Res{\tilde H^F}{G(s)^{\circ F}}(\Psi_{\tilde H,\tilde t_0}(\beta))= \Res{\tilde H^F}{G(s)^{\circ F}}(\Psi_{\tilde H,\tilde t_0}(\lexp a\beta))\circ \tau(a)$. If $\lexp a\beta\neq \beta$, then the orbits of $\beta$ and of $\lexp a\beta$ under ${\rm A}_{\cento{G^*} {s_0}}(s)^F$ are disjoint, thanks to the decomposition in direct product of ${\rm A}_{G^*}(ss_0)$ (Proposition~1.2.4). That implies  ${\rm A}_{G^*}(ss_0)^F_\beta={\rm A}_{\cento{G^*}s}(s_0)^F_\beta\times {\rm A}_{\cento{G^*}{s_0}}(s)_\beta^F$. Hence $G(s)^F_\la/G(s)^{\circ\, F}\subseteq {\rm A}_{\cento{G^*}{s_0}}(s)^F_\beta$ whereas the group ${\rm A}_{\cento{G^*}s}(s_0)^F_\beta$ acts on the set of \irr\ components of $\Res{\tilde H^F}{G(s)^{\circ F}}(\Psi_{\tilde H,\tilde t_0}(\beta))$, that  is a regular orbit under ${\rm A}_{\cento{G^*}s}(s_0)^F_\beta$. As  $[{\rm A}_{\cento{G^*}s}(s_0), {\rm A}_{\cento{G^*}{s_0}}(s)]=\{1\}$, any  element of that orbit is fixed by ${\rm A}_{\cento{G^*}{s_0}}(s)^F_\beta$.
 
 We have find $G(s)^F_\la/G(s)^{\circ \, F}={\rm A}_{\cento{G^*}{s_0}}(s)^F_\beta$. If $\mu\in\II{G(s)^F\mid \la}$,
$$\la(1).|{\rm A}_{G^*}(s)^F| = \mu(1).|{\rm A}_{\cento{G^*}{s_0}}(s)_\beta^F|\leqno{(3.4.2.11)}$$ 
and $
|{\rm A}_{\cento{G^*}{s_0}}(s)^F_\beta|$ elements of $\ser{G(s)^F}{(s_0)}$ are covering $|{\rm A}_{G^*}(s)^F/{\rm A}_{\cento{G^*}{s_0}}(s)_\beta^F|$ elements of $\II{G(s)^{\circ F}}$. If $\gamma$ is in the orbit of $\beta$ under ${\rm A}_{G^*}(ss_0)^F$ we have  $|{\rm A}_{\cento{G^*}{s_0}}(s)^F_\beta|=|{\rm A}_{\cento{G^*}{s_0}}(s)^F_\gamma|$.  We have seen that $|{\rm A}_{G^*}(s)^F/{\rm A}_{\cento{G^*}{s_0}}(s)^F|$ series of $G(s)^{\circ\, F}$ are fused under $G(s)^F$. Hence  the multiplicative factor we are looking for is 
$$m_\beta(G(s)^{\circ\, F} \to G(s)^F)= |{\rm A}_{\cento{G^*}{s_0}}(s)^F|/|{\rm A}_{\cento{G^*}{s_0}}(s)^F_\beta|^2$$
Our claim on the existence of $\Psi_{G,b}$ follows from part (A)  of the proof, (3.4.2.6), and the equality 
$$m_\beta(G^F\to H^F)=m_\beta(G(s)^{\circ\, F}\to H(t)^F)/m_\beta(G(s)^{\circ\, F} \to G(s)^F)$$
easy to verify thanks to (3.4.2.8) and (3.4.2.10) and the preceding equality.

As for degrees, recall (3.4.2.3) for $(H,\xi)$ proved in (A) : $\eta(1).|\cent{H^*}t^F|_{p'}=\xi(1).|H{^* F}|_{p'}$. In the morphism $\sigma^*\colon H^*\to G^*$, $H^{*\, F}$ maps onto $G^{* \,F}$, $\cent{H^*}t$ onto $\cento{G^*}s$ and $\cent{H^*}t^F$ onto $\cento{G^*}s^F$, so that $|H^{*\, F}|/|\cent{H^*}t^F|=|G^{* \,F}|/|\cento{G^*}s^F|= |G^F|/|G(s)^{\circ\,  F}|$. Thus  
$$\eta(1).|G(s)^F|_{p'}=\xi(1).|G^F|_{p'}.|{\rm A}_{G^*}(s)^F|$$
With (3.4.2.7),
(3.4.2.9), (3.4.2.11) and the isomorphism $\cento{G^*}{ss_0}^F_\beta\to \cento{\cento{G^*}s}{s_0}^F_\beta\times \cento{\cento{G^*}{s_0}}{s}^F_\beta$ we obtain $\chi(1).|G(s)^F|_{p'}=\mu(1).|G^F|_{p'}$ that is (3.4.2.3) for $(G,\mu)$.
\bull

\vfill\eject

\noindent{\bf 4.
 Brauer categories}

\medskip\noindent
{\bf 4.1. On defect groups } 

The construction of a maximal subpair for a block
$b_{G^F}(L,\la)$ (see 2.1.7) is given in [15] sections 4 and 5.1. First we describe defect groups.

\medskip\noindent{\bf 4.1.1. Proposition. }{\sl Assumption 2.2.3 on $(G,F,\ell,d$). Let $s$ be a semi-simple element in $(G^*)^F_{\ell'}$. Let $(L^*_s,\al)$ be a unipotent $d$-cuspidal datum  of $(\cento{G^*}s,F)$. Let $(L,\la)$ be one of the
$d$-cuspidal data in  $(G,F)$  in series $(s)$, associated to $(L^*_s,\al)$ by Proposition 2.1.4.  Let $T^*$ be an $F$-stable maximal torus in $L^*_s$ such that 
$T^*\cap\b{G^*}$ is maximally split in
$L^*_s\cap\b{G^*}$.
Let $(T,\theta)$ ($\theta\in \II {T^F}$) be in duality with  $(T^*,s)$ (where $s\in  (T^*)^F$). 

Let $D$ be
an $\ell$-subgroup of
$\nor GT^F_\theta$, maximal for the property 

``for any  $r\in \Phi_{L^*_s}(T^*)$, the one parameter subgroup $Y_r$ is contained in $\cent GD$".

Then $D$ is a defect group of $b_{G^F}(L,\la)$.} 

\preuve  
 If $s=1$, then $\theta=1$ and $L$ is in the dual class of
$L^*_1=L^*$ and the condition is a caracterization of an $\ell$-Sylow subgroup of
$C:=\cento G{[L,L)}^F
$,  that is  an $\ell$-Sylow subgroup of $\nor C{T\cap C}^F$ [12] Theorem~4.4. 
 
 Following 1.1.5 one has decompositions in central products
$G=\a G\b G$,
$G^*=\a{G^*}\b {G^*}$.  
The  isogenies $\b {G^*}\to (\b G)^*$,  $\b G\to (\b{G^*})^*$ (1.1.5.4)  restrict in isomorphisms between finite $\ell$-subgroups of $F$-fixed points. One has a central product
 $$\cento{G^*}s=\cento{\a{G^*}}s\cento{\b
{G^*}}s.\leqno{(4.1.1.1)}$$
If $G=\a G$, then there is only one conjugacy class of unipotent $d$-cuspidal data, say  $(T,1_{T^F})$, in $(G,F)$, where $T$ is a diagonal torus. Hence $L^*_s\cap\cento{\a {G^*}}s$ is a diagonal torus of $\cento{\a{G^*}}s$ (see Proposition 2.1.5). 

In [15] \S~4 a defect-group $D$ of
$b_{G^F}(L,\la)$ is described as follows : 

Let $M$ and $M^*$ be dual $E$-split Levi subgroups  of $G$, $G^*$ respectively such that 
$$M\cap \a G=T\cap \a G,\quad M^*\cap\b {G^*}=\cento{\b{G^*}}{\z{L^*(s)\cap\b
{G^*}}^F_\ell} \leqno{(4.1.1.2)}$$ ($T$ as above, duality around $(T,T^*)$). Then $D$ admits a unique maximal normal abelian subgroup $Z$ such that
$$M=\cento GZ,\quad
Z\cap \a G=(T\cap \a G)^F_\ell,\quad Z\cap \b G=\z {M\cap \b G}^F_\ell\leqno{(4.1.1.3)}$$

Let $Q$ be the subgroup of
${\rm W}(\cento{G^*}s,T^*)$ generated by reflections relatives to the roots $$r\in \Phi_{\cento{G^*}s}(T^*)\cap \Phi_{L^*_s}(T^*)^\perp$$
 let $V$ be an $\ell$-Sylow subgroup of $Q^F$. With these notations 
$\nor{G^F}Z\cap\nor{G} T$ contains a defect group $D$ of
$b_{G^F}(L,\la)$ such that $D\cap T=Z$, $DT/T\subseteq W(G,T)$ is anti-isomorphic to $V$  and the extension $Z\to D\to D/Z$ is split [15] Lemma~4.16.

When $G=\a G$, $L^*_s$ is a diagonal torus of  $\cento{G^*}s$, the caracterization of $Q$
reduces to ``$D$ is an $\ell$-Sylow subgroup of $\nor GT^F\cap \cento G\theta$", (where $\cento G \theta$ is in the dual $G^F$-conjugacy class of  the Levi subgroup $\cento{G^*}s$)
that is an extension of $Z=T^F_\ell$ by an $\ell$-Sylow subgroup of
${\rm W}(\cento G\theta,T)^F\cong {\rm W}(\cento{G^*}s,T^*)^F$.

Assume now $G=\b G$. We claim that 

\centerline
 {\sl $M^*$ is the smallest of
$E$-split Levi subgroups of $G^*$ that  contain $L^*_s$.}

 By (4.1.1.2) one has $L^*_s=(T^*\cap\a{G^*}).(L^*_s\cap \b {G^*})\subseteq M^*$, hence $\z{M^*}^F_\ell\subseteq \z{L^*_s}^F_\ell$.  It follows, using (4.1.1.2) again,  that $\z{L^*_s}^F_\ell=(T^*\cap\a{G^*})^F_\ell.\z{L^*_s\cap\b {G^*}}^F_\ell
\subseteq \z{M^*}^F_\ell$. Thus $\z{M^*}^F_\ell=\z{L^*_s}^F_\ell$. But $\z {M^*}^F_\ell=(\zo{M^*}_{\phi_E})^F_\ell$, $\z{L^*_s}^F_{\ell}=(\zo{L^*_s}_{\phi_E})^F_{\ell}$ by Proposition~1.2.5, and $\zo{M^*}\subseteq T^*\subseteq L^*_s\cap M^*$ so that
$\zo {M^*}_{\phi_E}\subseteq
\zo{L^*_s}_{\phi_E}$. The definition of $E$ (1.1.5.1) implies  $\zo {M^*}_{\phi_E}=
\zo{L^*_s}_{\phi_E}$. 

Now let $Z_1:=T^F_\ell\cap (\cap_{\al^\vee\in\Phi_{L^*_s}(T^*)}\Ker \al)$. We claim that $Z_1=Z$. As $\Phi_{L^*_s}(T^*)\subseteq\Phi_{M^*}(T^*)$, using (4.1.1.3), one has $Z\subseteq \z M^F_\ell\subseteq Z_1$. Let $M_1:=\cento G{Z_1}$. Thanks to our hypotheses ``$\ell$ is good and $G=\b G$" $M_1$ is an $E$-split Levi subgroup of  $G$. By definition of $Z_1$ and $M_1$, $Z_1\subseteq \z {M_1}^F_\ell$ and, if $r^\vee\in\Phi_{L^*_s}(T^*)$, then $r\in\Phi_{M_1}(T)$.
Let $M_1^*$ be the  Levi subgroup of $G^*$ that contains $T^*$ and in duality with $M_1$, then  $M^*_1$ is $E$-split and $\Phi_{L^*_s}(T^*)\subseteq \Phi_{M_1^*}(T^*)$. Thus $L^*_s\subseteq M^*_1$ so that $M^*\subseteq M_1^*$,
$M\subseteq M_1$ and $Z_1\subseteq Z$. 

Now let us consider Weyl's groups : as  $L^*_s$ is a Levi subgroup of $\cento{G^*}s$, $Q$ is the subgroup of elements of ${\rm W}(\cento{G^*}s,T^*)$
that fix any $r\in \Phi_{L^*_s}(T^*)$. Then Steinberg 's relations
show that $Q$ is the image of $N\subseteq \nor G T$ where $N$ centralizes any one parameter subgroup $Y_r$ of $G$ associated to 
some $r\in\Phi_{L^*_s}(T^*)$ and $N\cap T$ is a finite $2$-group.  As $|V|$ is odd there is a subgroup $\dot V$ in $\nor GT^F$, with $\dot V.T/T=V$ and  such that $[\dot V,Y_r]=\{1\}$ for any $\al\in\Phi_{L^*_s}(T^*)$. The semi-direct product $Z.\dot V$ is a 
defect group of $b_{G^F}(L,\la)$.

Going from $\a G$ and $\b G$ to
the central  product
$G=\a G.\b G$ is quite easy. When $s$, $L^*_s$ et $T^*$ are fixed $(T,\theta)$ is defined mod $G^F$-conjugacy, thus we have defined a $G^F$-conjugacy class of $\ell$-subgroups of $G^F$,  it is the set of defect groups of $b_{G^F}(L,\la)$.
\bull

\smallskip 
As $G(s)^F/G(s)^{\circ \,F}$ is prime to $\ell$, if a block $b$ of $G(s)^F$ covers a block $b_0$ of $G(s)^{\circ \,F}$, a defect group of $b_0$ is a defect group of $b$. 

\medskip\noindent{\bf 4.1.2.
Proposition. }{\sl Assumption 2.2.3 on $(G,F,\ell,d)$. Let $s$ be a semi-simple element in $(G^*)^F_{\ell'}$. Let $(L^*_s,\al)$ be a unipotent $d$-cuspidal datum  in $(\cento{G^*}s,F)$. Let $(L,\la)$ (resp. $(L(s),\al)$) be one of the
$d$-cuspidal data in  $(G,F)$  in series $(s)$ (resp. in $(G(s)^\circ, F)$ in unipotent series) associated to $(L^*_s,\al)$ by Proposition 2.1.4. The  defect  groups of   $b_{G^F}(L,\la)$, $b_{G(s)^{\circ F}}(L(s),\al)$ are isomorphic. 
}

\preuve
 The central product in  (4.1.1.1) gives $$G(s)^\circ=\a G(s)^\circ.\b G (s)^\circ$$ (but $\a{G}(s)^\circ$ is not $\a{(G(s)^\circ)}$ and $\b {G}(s)^\circ$ is not exactly in duality with $\cento{\b{G^*}}s$ ...)
 
Let $T$ et $T^*$ be dual tori as in  4.1.1 and let  $T_1$ be a dual of $T^*$ in
$G(s)^\circ$, so that the duality between $G(s)^\circ$ and $\cento{G^*}s$ is defines around $(T_1,T^*)$.
Then
$L^*_s$ is in duality with a Levi subgroup 
$L(s)$ of
$G(s)^\circ$ such that $T_1$ is maximally split in $L(s)$. There exists an $F$-compatible isomorphism 
between the couples
$(T_1,T^*)$ and $(T,T^*)$,
restricting as identity on $T^*$ and ${\rm X}(T^*)$. Let $$\rho\colon T_1\to T$$ be its restriction to $T_1$.

We have
$M^*\cap\cento{\b{G^*}}s=L^*_s\cap\b{G^*}$ because $L^*_s$ is $E$-split in
$\cento{G^*}s$ as is $M^*$ in $G^*$.

 We have seen in the proof of Proposition  4.1.1 how $Z$ is defined from $L^*_s$ : 
on $\a G$ side, $Z\cap \a G=(T\cap\a G)^F_\ell$. On $\b G$ side $\zo{M^*\cap\b{G^*}}_{\phi_E}=\zo{L^*_s\cap\b
{G^*}}_{\phi_E}$ so that $\zo{M\cap\b{G}}_{\phi_E}$ is the biggest $\Phi_{E}$-subgroup of $T\cap \b G$ in the kernel of any $r$ such that $r^\vee\in\Phi_{L^*_s\cap \b{G^*}}(T^*\cap \b{G^*})$.

 Let $Z_1\subseteq(T_1)^F_\ell$  be obtained by the same process from $G(s)^{\circ\, F}$
and the $d$-cuspidal unipotent datum $(L^*_s,\al)$, with $(G(s)^\circ,\cento{G^*}s,1,T_1)$ instead of $(G,G^*,s,T)$. To compare
$Z$ with
$Z_1$, we have  to consider the rational component $\a{\cento{\b{G^*}}s}$ of
$\cento{\b{G^*}}s$. It is a central product $\a{\cento{G^*}s}=\cento{\a{G^*}}s.\a{[\cento{\b{G^*}}s]}$ and in the dual side (dually up to some isogeny) is a similar \dec\ 
  $\a{G(s)^\circ}=\a G(s)^\circ.\a{[\b G(s)^\circ]}$, as well as
$\b{G(s)^\circ}=\b{[\b G(s)^\circ]}$. Clearly 
$T_1\cap
\a {G(s)^\circ}=(T_1\cap \a G(s)^\circ).(T_1\cap \a{[\b G(s)^\circ]})$ and
$$Z_1\cap \a{G(s)^\circ}=(T_1\cap \a
G(s)^\circ)^F_\ell.(T_1\cap\a{[\b G(s)^\circ]})^F_\ell$$  (by Proposition 1.1.3 (b), $\z {\a{[\b G(s)^\circ]}}^F/\zo{\a{[\b G(s)^\circ]}}^F$ is prime to $\ell$). As $(L^*_s\cap
\a{[\cento{\b{G^*}}s]})^F$ admits a cuspidal unipotent,
$L^*_s\cap
\a{[\cento{\b{G^*}}s]}$ is a diagonal torus of $\a{[\cento{\b{G^*}}s]}$.

Therefore
$L^*_s\cap
\a{[\cento{\b{G^*}}s]}=T^*\cap \a{[\cento{\b{G^*}}s]}$ hence
$(T^*\cap
\a{[\cento{\b{G^*}}s]})_{\Phi_{E}}= T^*_{\phi_E}\cap
\a{[\cento{\b{G^*}}s]}$ and $$\z{L^*_s\cap\b
{G^*}}^F_\ell=(T^*\cap\a{[\cento{\b{G^*}}s]})^F_\ell.\z{L^*_s\cap
\b{[\cento{\b{G^*}}s]}}^F_\ell$$ 
The central products 
$$T^*\cap \b{G^*}=(T^*\cap
\a{[\cento{\b{G^*}}s]}).(T^*\cap
\b{[\cento{\b{G^*}}s]}))$$ 
$$ \zo{M^*\cap \b{G^*}}_{\phi_E}= (T^*\cap
\a{[\cento{\b{G^*}}s]})_{\phi_E}.\zo{L^*_s\cap  \b{[\cento{\b{G^*}}s]})}_{\phi_E}$$    
 define on dual side  central products
$$T\cap\b G=S.R,\quad \zo{M\cap\b G}_{\phi_E}=S_{\phi_E}.H$$ where
$H$ is defined from
$L^*_s\cap\b{[\cento{\b{G^*}}s]}$ as $\zo{M\cap\b G}_{\phi_E}$ is defined from
$L^*_s\cap{\b{G^*}}$. As $\rho$ is the identity on $X(T^*)\cong Y(T)\cong Y(T_1)$, one sees that
$\rho(T_1\cap
\a{G(s)^\circ})=(T\cap \a G)S$ and $\rho(T_1\cap \b{G(s)^\circ})=R$, so that $\rho(Z_1\cap
\a{G(s)^\circ})=(Z\cap
\a G)(S_{\phi_E})_\ell^F$ and
$\rho(Z_1\cap\b{G(s)^\circ})=H^F_\ell$. We have obtained  $$\rho(Z_1)=Z$$

We note that  $Z_1=\z{L(s)}^F_\ell=\zo{L(s)}^F_\ell$
(see [13] \S~4.3 and Proposition~3.3). 
Now we know by [13] again that  $\nor{G(s)^{\circ F}}{Z_1}\cap\nor{G(s)^\circ} {T_1}$ contains a
defect
 group $D_1$ of
$b_{G(s)^{\circ F}}(L(s),\al)$ such that $D_1\cap T_1=Z_1$, $D_1T_1/T_1$ is an
$\ell$-Sylow subgroup of  $Q_1^F$,  where $Q_1={\rm W}(T_1.\cento{G(s)^\circ}{[L(s),L(s)]},T_1)$ and the extension$$Z_1\to D_1\to D_1/Z_1$$ is split. 

The anti-isomorphism between $W(G(s)^\circ,T_1)$ and $W(\cento{G^*}s,T^*)$ restricts to an anti-isomorphism between $W(T_1.\cento{G(s)^\circ}{[L(s),L(s)]},T_1) $ and $W(T^*.\cento{\cento{G^*}s}{[L^*_s,L^*_s]},T^*)$. By Proposition 4.1.1 and its proof our group $Q_1$ is anti-isomorphic to  the group we denoted $Q$. Thus the 
 isomorphism from  $Z_1$ onto $Z$ is the restriction of an isomorphism from the defect group of 
$b_{G(s)^{\circ F}}(L(s),\al)$ onto a defect group of
$b_{G^F}(L,\la)$. 

We see also that the defect groups of $b_{G(s)^{\circ F}}(L(s),\al)$ and $b_{\cent{G^*}s^F}(L^*_s,\al)$ are anti-isiomorphic.

\bull

\medskip\noindent{\bf 4.2. Subpairs and Brauer's categories }

Using notation~2.3.1, we enforce Proposition~2.3.2, where $H$ was an $E$-split Levi subgroup of $G$. The proofs are quite similar.

\medskip\noindent{\bf 4.2.1. Proposition. }{\sl Let $Y$ be an $\ell$-subgroup of
$G^F$ and let
$H=\cento GY$. The relation $(L_H,\al_H){\simex {G^F}}(L,\al)$  defines a  bijection
between the set $\simex {G^F}$-classes of unipotent $d$-cuspidal data
$(L,\al)$ in
$(G,F)$ such that $[L,L]\subseteq H^g$ for some $g\in G^F$ and  $\simex {H^F}$-classes
of unipotent $d$-cuspidal data
 $(L_H,\al_H)$ in 
$(H,F)$.}

\preuve Let $G=\a G.\b G$ be the \dec\ defined in 1.1.5. A $d$-cuspidal unipotent datum in $(G,F)$ is provided by the central product of unipotent $d$-cuspidal data $(\a L,\a \al)$, $(\b L,\b \al)$ in $(\a G,F)$ and $(\b G,F)$ respectively and $\a L$ is a diagonal torus in $\a G$, $\a\al=1_{{\a L}^F}$ (Proposition~2.1.5). The set of $\simex {G^F}$-classes of unipotent $d$-cuspidal data in
$(G,F)$ is in bijection with the set of $\simex {\b G^F}$-classes of unipotent $d$-cuspidal data in 
$(\b G,F)$. 

If $\b G\subseteq H$, then $\a H=H\cap \a G$ and $\b H=\b G$, hence our claim.

We assume $\b G\not\subseteq H$ and use induction on the semi-simple rank of $G$. There exists $\ell$-subgroups $\a
Y\subseteq {\a G}^F$ and
$\b Y\subseteq
{\b G}^F$ such that $\cent GY=\cent {\a G}{\a Y}.\cento {\b G}{\b Y}$ (see (D) at the end of the proof of  Proposition~1.2.6).  By assumption on $H$, $\b Y\neq 1$, hence $\z{\b Y}\neq 1$. Let $z\in\z{\b Y}$, $z\neq 1$, so that $z\not\in\z{{\b G}}$, and $\cent Gz$ is a proper $E$-split Levi subgroup of $G$ by Proposition~1.2.5. Clearly we have $Y.\cent GY
\subseteq \a G.
\cent {\b G}z=\cento Gz\neq  G$. There exists a proper $d$-cuspidal Levi subgroup $M$  of $G$ such that  $\cent Gz\subseteq M$ (1.1.5.3). By inductive hypothesis, given a unipotent $d$-cuspidal datum  $(L_H,\al_H)$ in $(H,F)$ there exists $(L_M,\al_M)$ with $(L_H,\al_H)\simex {M^F}(L_M,\al_M)$. But $(L_M,\al_M)$ is a $d$-cuspidal unipotent datum in $(G,F)$, hence $(L_H,\al_H)\simex {G^F}(L_M,\al_M)$. 
\bull

In the following Proposition on unipotent blocks we describe all
Brauer's subpairs, a generalization of  [16] Lemma~23.10, relying inclusion to the relation we have introduced in Proposition~4.2.1.

\medskip\noindent{\bf 4.2.2. Proposition. }{\sl  Assumption 2.2.3 on $(G,F,\ell,d)$. Let $(L,\al)$ be a $d$-cuspidal unipotent datum in $(G,F)$. Let $(Y,b_Y)$ be a Brauer  $\ell$-subpair  of
$G^F$ and let  $(L_Y,\al_Y)$ be  a unipotent $d$-cuspidal datum in  $(\cento GY,F)$. Assume that $b_Y$ covers $b_{\cento GY^F}(L_Y,\al_Y)$.

One has $(\{1\},b_{G^F}(L,\al))\subset (Y,b_Y)$ if and only if
$(L_Y,\al_Y)\simex {G^F}(L,\al)$.}

\preuve  Induction on the semi-simple rank of $G$.

We know that $b_{G^F}(L,\al)$ is defined by $b_Y$ in the inclusion $(1,b_{G^F}(L,\al))\subset (Y,b_Y)$   and that $b_Y$ is defined by $b_{\cento GY^F}(L_Y,\al_Y)$. By Proposition~2.1.6, the $G^F$-conjugacy class of $(L,\al)$ is then defined by the $\cento GY^F$-conjugacy class of $(L_Y,\al_Y)$. Thanks to Proposition~4.2.1 we have only to show that the relation $\simex {G^F}$ between the two unipotent
$d$-cuspidal data implies the inclusion of the so-defined  $\ell$-subpairs.

Case (A) : $\b G$ is not contained in  $\cento GY$.

Let $M$ be a proper $d$-split Levi subgroup of $G$ such that $\a G.\cent GY\subseteq M$ (see the proof of Proposition~4.2.1). 
One has
$M=\cento { G}{\z M^F_\ell}$ and
$M^F=\cent{G^F}{\z M^F_\ell}$  by Proposition  1.2.3. If
$(L_1,\al_1)$  is a $d$-cuspidal unipotent datum in $(M,F)$, induction applies : the relation $(L_Y,\al_Y)\simex {M^F} (L_1,\al_1)$ is equivalent to  the inclusion $(1,b_{M^F}(L_1,\zeta_1))\subset (Y,b_Y)$ in $M^F$. 

That
inclusion is equivalent to  $(\z M^F_\ell,b_{M^F}(L_1,\al_1))\subset (\z M^F_\ell
.Y,b_Y)$ in
$G^F$. 
As $M$ est $d$-split,
$(L_1,\al_1)$ is a $d$-cuspidal datum in $(G,F)$ and the $M^F$-conjugacy class of $(L_1,\al_1)$ defines the $G^F$-conjugacy class of $(L_1,\al_1)$ : $(L_Y,\al_Y)\simex {M^F} (L_1,\al_1)$ is equivalent to $(L_Y,\al_Y)\simex { G^F} (L_1,\al_1)$ (Proposition~4.2.1).  By Propositions~2.1.6, 2.1.7
$\Lu MG(b_{M^F}(L_1,\al_1))=b_{G^F}(L_1,\al_1)$
and $(\{1\},b_{G^F}(L_1,\al_1))\subset (\z
M^F_\ell,b_{M^F}(L_1,\al_1))
$ in $G^F$. Transitivity of inclusion gives $(1,b_{G^F}(L_1,\al_1))\subset (Y,b_Y)
$ in $G^F$.

Case (B) :  $\b G\subseteq \cento GY $.

The map $\ser{G^F}1\to\ser{{\a G}^F}1 \times \ser{{\b G}^F}1$,   defined by
$$\Res{G^F}{{\a G}^F.{\b G}^F} =\a \chi\otimes \b \chi,\quad \chi\mapsto (\a \chi,\b \chi)$$ is one-to-one.

Thus any unipotent block  of $\a G^F.\b
G^F$ is covered by exactly one unipotent block of $G^F$, other covering blocks are in $\ell'$-s\'eries $\ser{G^F}s$ where $s\in\z{G^*}^F$ is in duality with some
$\theta\in(G^F/\a G^F.\b G^F)^\wedge$. 

This applies  to the \dec\  $\cento GY=\cento{\a G}Y.\b G$ :  $L_Y=T_Y.( L_Y\cap \b G)$,  where
$T_Y$ is a  diagonal torus of
${\cento { \a G} Y}$, $L_Y\cap \b G$ is a
  $d$-split Levi subgroup of $\b G$ and $\al_Y\in\ser{(T_Y.( L_Y\cap \b G))^F}1$ covers $1_{T_Y^F}\otimes
\b\al$, $\b \al\in\ser{L_Y\cap {\b G}^F}1$, $\b\al$ is $d$-cuspidal (Proposition~2.1.5, (c)). Let  $B_0(Y)$ be the principal block of $\cento{\a G}Y^F$, $B_0(1)$ be the principal block of ${\a G}^F$ and $b_1=b_{ {\b G}^F}(L_Y\cap\b G,\b \al)$. There is exactly one unipotent block $b_Y^\circ$   of $\cento GY^F$ that covers $B_0\otimes b_1$ and there is exactly one unipotent block  $b_Y$ of $\cent{G^F}Y$ which covers  $b_Y^\circ$.

The inclusion in  $\a G^F.\b G^F$ of subpairs $$(1,B_0(1)
\otimes b_{\b G^F}( L_Y\cap\b G,\b\al))\subset (Y,B_{0}(Y)\otimes b_{\b G^F}(L_Y\cap \b G,\b\al))$$ is clear. By unicity of covering unipotent blocks we have an
inclusion in $G^F$ of unipotent subpairs $$(1,b_{G^F}(\a T.
(L\cap\b G),\al)\subset (Y,b_Y)$$ where $\a T$ is a diagonal torus of $\a G$ and $\al\in\ser{(\a T.( L_Y\cap\b G))^F}1$ covers $1_{\a T^F}\otimes \b \al$.
\bull

Now we describe all Brauer subpairs and their inclusion in any series, under our general hypothesis.

\medskip\noindent{\bf 4.2.3. Proposition. }{\sl Assumption 2.2.3 on $(G,F,\ell,d)$. Let $Y$ be an  $\ell$-subgroup of
$G^F$ and $T$,  $T_Y\subseteq T$, $Y'$, $T_{ Y'}\cong {T_Y}^*\subseteq T^*$ as in
Proposition~1.2.6, so that $\cento GY$ and $\cento{G^*}{Y'}$ are in duality. Let
$(L_Y,\la_Y)$ be a 
$d$-cuspidal datum in
$(\cento GY,F)$  with $T_Y\subseteq L_Y$. Assume
$\la_Y\in\ser{{L_Y}^F}{s}$ with $s\in ({T_Y}^*)^F_{\ell'}$. Let $b_Y$ be the block of
$\cent{G^F}Y$ that covers  
 $b_{\cento GY^F}(L_Y,\la_Y)$. 

Let $(L^*_{Y,s},\al_{Y})$ be a unipotent $d$-cuspidal datum in
$(\cento{\cento GY^*}s,F)$ associated to  $(L_Y,\la_Y)$ by Proposition 2.1.4. Let $(L^*_s,\al)$ be a unipotent $d$-cuspidal datum in  $(\cento
{G^*}s,F)$ such that $$(L^*_{Y,s},\al_{Y})\simex{\cento{G^*}s^F}(L^*_s,\al)\,.$$

 There exists a $d$-cuspidal datum 
$(L,\la)$ in
$(G,F)$ such that $T\subset L$,
 $\la\in\ser{L^F}s$, $(L,\la)$ is associated to $(L^*_s,\al)$ and
$$(\{1\},b_{G^F}(L,\la))\subset (Y,b_Y)\,.\leqno{(4.2.3.1)}$$  

}

\preuve 
Clearly Proposition~4.2.3 goes through direct products. 

\noindent (A)  Assume the center of $G$  connected.

By Proposition~1.2.6, $\cent{G^*}s$, $\cent {G^*}{Y'}$ and  $\z{\cento GY}$ are connected. Then the $\cent{\cent G{Y^*}}s^F$-conjugacy class of $(L^*_{Y,s},\al_{Y})$ is well defined by $(L_Y,\la_Y)$ (Proposition~2.1.10).  
 The $\cent{G^*}s^F$-conjugacy class of 
$(L^*_s,\al)$ is well defined by $(L^*_{Y,s},\al_{Y})$ as in Proposition~4.2.2. We have to prove that the relation between the two unipotent $d$-cuspidal data inside $\cent{G^*}s$ imply (4.2.3.1) in $G^F$. We use induction on the dimension of $G$.

(A.1) Assume $G=\a G$. 

Let $G(s)$ be a Levi subgroup of $G$,  in the dual $G^F$-conjugacy class of $\cent{G^*}s$. There is only one block to consider in series $(s)$, denote it 
$b_{G^F}(s)$ and $b_Y$ covers $b_{\cento GY^F}(s)$.  An $\ell$-Sylow subgroup of $G(s)^F$ is a defect group of $b_{G^F}(s)$. By our hypotheses, $s\in\cent{G^*}{Y'}$ and, up to $G^F$-conjugacy we may assume $Y\subseteq G(s)^F$. The condition $(L^*_{Y,s},\al_{Y})\simex{\cento{G^*}s^F}(L^*_s,\al)$ is always satisfied by $d$-cuspidal data.  It is sufficient  to prove
the inclusion $(1,b_{G^F}(s))\subset (Y,b_Y)$ and this is independant of choices of $Y'$, and of tori defining dualities.

 The $\ell$-subpairs of GL$_n(q)$ are given in [7], where $G$ is identified with $G^*$, the $G^F$-conjugacy class of $(T,\theta)$ ($T$ a maximal $F$-stable torus, $\theta\in (T^F)^\wedge$) corresponds to the $G^{*\,F}$-conjugacy class of $(T^*,s)$ ($T^*$ a maximal $F$-stable torus, $s\in T^{*\,F}$)).  We know  that $\cent GY$ is connected and a direct product of linear groups. From [7] we retain : if $B$ is a block of ${\rm GL}_n(q)$ in series $(s)$ and $Y\subseteq G(s)^F$ there is an inclusion of subpairs
$(1,B)\subset(Y,b_Y)$ in ${\rm GL}_n(q)$ where $b_Y$ is in series $(s)$ in $\cent GY^F$. In our hypothesis $G=\a G$ that gives $(1,b_{G^F}(s))\subset (Y,b_{\cent GY^F}(s))$ in $G^F$, (4.2.3.1) in that case.

Any group $G$ of type $\AA$ (resp. $(G,F)$ such that $G=\a G$) with connected center may be reached from general linear groups GL$_n$ (resp. $G={\rm GL}_n$ such that $G=\a G$) by a sequence of morphisms of three types : 

1. direct product,

2.  $(H,F)\to (G,F)$ is a regular covering between groups with connected centers (and $H=\a H$)

3. $(G,F)\to (H,F)$ is an isotypic embedding, where $\z G$ is connected, $[G,G]$ is simply connected and there is an isotypic embedding $G_0\to H$, with $G_0\cong GL_n$.

We verify that our claim goes from $(H,F)$ to $(G,F)$ in 2. and 3.

In case 2,  $G^F$
is a quotient of $H^F$, the dual map  $G^*\to H^*$ is an embedding, so  let $s\in H^*$. There is some $\ell$-group  $Y_H$ in $H^F$  with image $Y$ in $G^F$ and $\cent {G^F}Y $ is a quotient of
$\cent{H^F}{Y_H}$. There is a regular embedding $\cent GY^*\to \cent H{Y_H}^*$, we assume $s\in \cent
H{Y_H}^*=\zo{H^*}.(\cent GY)^*\subseteq H^*$. The two inclusions between $\ell$-subpairs
$(1,b_{H^F}(s))\subset (Y_H,b_{\cent H{Y_H}^F}(s))$ in $H^F$ and
$(1,b_{G^F}(s))\subset (Y,b_{\cent GY^F}(s))$ in $G^F$ are true or not simultaneously.

In case 3 assuming $Y\subseteq G^F$, let $Y_0=Y.\zo H^F_\ell\cap G_0$. We have
$H^F=\zo H^F.G^F=\zo H^F.G_0^F$ and $\cent HY=\zo H^F.\cent GY=\zo
H^F.\cent{G_0}{Y_0}$. Let
$s_H\in(H^*)^F_{\ell'}$ with images
$s\in G^*$, 
$s_0\in G_0^*$. Then $\cent{G^*}s$ and $\cent{G^*_0}{s_0}$ are  quotients of
$\cent{H^*}{s_H}$ by central torii. Clearly $b_{H^F}(s_H)$ covers $b_{G^F}(s)$ and $b_{G_0^F}(s_0)$.
The morphism $ \cent GY^*\to G^*$ given by Proposition 1.2.6 may be deduced by quotient from a morphism $\cent HY^*\to H^*$ (with duality around $T_Y.\zo H$...), giving $\cent{G_0}{Y_0}^*\to G_0^*$ (duality
around $T_Y.\zo H\cap G_0$...). Then $b_{\cent H{Y}^F}(s_H)$ covers $b_{\cent G{Y}^F}(s)$ and $b_{\cent {G_0}{Y_0}^F}(s_0)$. 
The inclusion $(1,b_{G_0^F}(s_0))\subset (Y_0,b_{\cent{G_0}{Y_0}^F}(s_0))$ in $G_0^F$ implies $(1,b_{H^F}(s_H))\subset (Y,b_{\cent{H}{Y}^F}(s_H))$ which implies (4.2.3.1) in $G^F$.

 (A.2)  Assume $\b G\subseteq \cent GY$.
 
That case reduces to the preceding one by a standard description, with given dualities :

$s\mapsto (\a s,\b s)$ by $G^*\to (\a G)^*(\b G)^*$,

$\cent GY=\cent{\a G}{\a Y}\b G$ where $\a Y\subseteq \a G^F$,

$\cent{G^*}{Y'}$ maps on $\cent{(\a G)^*}{\a Y'}.(\b G)^*$ for some $\a Y'\subset (\a G)^{*\,F}$,

$L_Y=L_{\a Y}M$, where $L_{\a Y}\subseteq \cent{\a G}{\a Y}$, $M\subseteq\b G$, $\Res{L_Y^F}{L_{\a Y}^FM^F}\la_Y=\la_{\a Y}\otimes \mu$ where $\la_{\a Y}$ (resp. $\mu$) is $d$-cuspidal in series $(\a s)$ (resp. $(\b s)$) (see Proposition~2.1.5),

if $(L_{\a Y},\la_{\a Y})$ (resp. $(M,\mu)$)
is associated to the $d$-cuspidal unipotent datum $(L^*_{\a Y,\a s},\al_{\a Y})$ (resp. $(M^*_{\b s},\beta)$) then we may assume that $L^*_{Y,s}$ has image $(L^*_{\a Y,\a s}\times
M^*_{\b s})$ by the map $G^*\to(\a G)^*\times (\b G)^*$ and that $\Res{L^{*\,F}_{\a Y,\a s}\times M^{*\, F}}{L_Y^{*\,F}}( \al_{\a Y}\otimes \beta)=\al_Y$.

the condition $(L^*_{Y,s},\al_{Y})\simex{\cento{G^*}s^F}(L^*_s,\al)$ is then equivalent to :

\noindent `` $L^*_s$ has image $L^*_{\a s}\times M^*_{\b s}$ in $(\a G)^*\times (\b G)^*$,  $\Res{L^{*\,F}_{\a s}\times M^{*\, F}_{\b s}}{L^{*\,F}_s}( \a \al\otimes \beta)=\al$ where $\a \al\in\ser{L^{*\,F}_{\a s}}1$ and  $(L^*_{\a Y,\a s},\al_{\a Y})\simex{\cento{(\a G)^*}{\a s}^F}(L^*_{\a s},\a \al)$ "

If $( L_a,\a \la)$ is a $d$-cuspidal datum in series $(\a s)$ in $(\a G,F)$ associated to $(L^*_{\a s},\a \al) $, then we may assume that $L= L_a.M$ and $\Res{L^F}{L_a^FM^F}\la=\a\la\otimes \mu$.

Then (4.2.3.1) is equivalent to $(1,b_{(\a G)^F}( L_a,\a \la))\subset (\a Y,b_{\a Y})$. That last inclusion follows from (A.1).

(A.3) Assume now  $\b G\not\subseteq \cent GY$.

As in the proof of Proposition~4.2.1, let  $z\neq 1$ such that 
$Y.
\cent GY\subseteq 
\cent Gz$ where $\cent Gz$ is a proper $E$-split Levi  subgroup of
$G$. By Proposition~1.2.6 there exist an $\ell$-subgroup $Y'$ in $G^{*\, F}$ and inclusions of dual groups 
$
\cent GY^*=\cent{G^*}{Y'}\subseteq
\cent Gz^*\subseteq G^*$ with $s\in\cent{G^*}{Y'}$. By Proposition~4.2.2 we have unipotent $d$-cuspidal data $(L^*_{Y,s},\al_{Y})$ in $(\cent{\cent {G^*}{Y'}}s,F)$, $(L^*_{z,s},\al_{z})$ in $(\cent{\cent {G}{z}^*}s,F)$,  $(L^*_s,\al)$ in $(\cent {G^*}s,F)$ such that 
$$(L^*_{Y,s},\al_{Y})\simex{\cent{\cent {G}z^*}s^F}(L^*_{z,s},\al_{z})\simex {\cent
{G^*}s^F}(L_s,\al)$$
These unipotent $d$-cuspidal data define $d$-cuspidal data in series $(s)$ :   
$(L_{Y},\la_{Y})$ in $({\cent GY},F)$, $(L_{z},\la_{z})$ in $({\cent
Gz},F)$,   and $(L,\la)$ in $(G,F)$. By Propositions~2.1.6, 2.1.7 we have an inclusion of subpairs  : 
$$(1,b_{G^F}(L,\la))\subset  (\gen z,b_{\cent Gz^F}(L_z,\la_z))\quad{\rm in}\; G^F$$
By induction we have an inclusion of subpairs:
$$(1,b_{\cent Gz^F}(L_z,\la_z))\subset (Y,b_{\cent GY^F}(L_Y,\la_Y))\quad {\rm  in}\; \cent Gz^F  $$
 equivalent to  
$$(\gen z,b_{\cent Gz^F}(L_z,\la_z)) \subset
(Y,b_{\cent GY^F}(L_Y,\la_Y))\quad {\rm in}\; G^F $$
Transitivity of inclusion of subpairs in $G^F$ gives (4.2.3.1).

\smallskip\noindent (B) $\z G$ connected or not.

 The inclusion of subpairs in $G^F$ may be described by inclusion of so-called ``connected subpairs" in [15], 2.1, where $(1,b_{G^F}(L,\la))\subset (Y,b_Y)$ is  equivalent to its ``restriction" $(1,b_{G^F}(L,\la))\subset b_{\cento GY^F}(L_Y,\la_Y))$.
 
 In a regular embedding $G\to  H$, let $t\in H^*$ of image $s$ by a dual map. Let $(M,\mu)$ be a $d$-cuspidal datum in series $(t) $ in $(H,F)$ such that $L=M\cap G$ and $\mu $ covers $\la$ (see Proposition~2.1.5). The associated $d$-cuspidal data $(L^*_s,\al)$ and $(M^*_ t,\beta)$ are related : $\beta$ is the restriction of $\al$ through the isotypic morphism $M^*_t\to L_s^*$ induced by the restriction $M^*\to L^*$ of $H^*\to G^*$. By Proposition~2.4.2, Proposition~2.4.4 and its proof, a quotient $H^F/K$ of
$H^F/G^F\cong M^F/L^F$ acts regularly on the set of blocks of $G^F$ that are covered by
$b_{H^F}(M,\mu)$, as well as on the set of unipotent $d$-cuspidal elements of $\ser{(L^{*\,F}_s}1$ that are covered by $\mu$. 

If $D$ is a defect group of $b_{G^F}(M,\mu)$, the defect groups of the various $b_{G^F}(L,\la)$, $\la$ covered by $\mu$,  form a $H^F$-conjugacy class and are the $D\cap G^F$, where $D$ is a defect group of $b_{H^F}(M,\mu)$ (see Propositions~5.1.3, 5.1.4 in Appendix). The action of $H^F$ on $G^F$ transforms inclusion of Brauer subpairs in $G^F$ in inclusion of Brauer subpairs.
Let $Y$ be some subgroup of $D\cap G^F$. We have a natural isomorphism $H^F/G^F\cong \cento HY^F/\cento GY^F$ giving an  action of $H^F/G^F$ on the blocks $b_Y$ of $\cent{G^F}Y$ that are covered by a given block $B_Y$ of $\cent{H^F}Y$. By this way $H^F/K$ acts on the set of all Brauer subpairs in $G^F$ containing $(1, b_{G^F}(L,\la))$ such that $\la$ is covered by $\mu$. One sees that (4.2.3.1) follows from the inclusion $(1,b_{H^F}(M,\mu))\subset (Y,B_Y)$.
\bull

\smallskip The Proposition~4.2.4 is the decisive step to show isomorphism between Brauer categories of  a unipotent block $b_s$ of $G(s)^F$  and its Jordan correspondant $b_{G^F}(L,\la)$, a block of $G^F$ in series $(s)$, as described in the proof of Proposition~3.4.1. The link between  $b_s$ and  $b_{G^F}(L,\la)$ is the $d$-cuspidal unipotent $\al\in\ser{\cento{G^*}s^F}1$, so  we may assume, with notations of 3.4.1, that  ${\cal B}_{G,s}(b_s)=b_{G^F}(L,\la)$. Given a subgroup $Y$ of a common defect group, it appears  a similar link between Brauer subpairs, a $d$-cuspidal $\al_Y\in\ser{\cento{G(s)^\circ}Y^F}1$. The inclusion of subpairs is a consequence of the relation we introduced in 2.3.1, by Proposition~4.2.2.

\medskip\noindent{\bf 4.2.4. Proposition. }{\sl Let $(G,F,\ell,d)$, $Y$, $Y'$, $(L_Y,\la_Y)$, $s$, $b_Y$, $(L^*_s,\al)$, $(L^*_{Y,s},\al_Y)$ as in Proposition 4.2.3, hence (4.2.3.1) : 
$$(\{1\},b_{G^F}(L,\la))\subset (Y,b_Y)\,.$$
  Let 
$G^\circ(s)$ and
$G(s)$ be defined as in section 3.1. We may assume $Y\subset G(s)^{\circ\,F}$ by Proposition~4.1.2.
The group $\cento{\cento{G^*}{Y'}}s=\cento{\cento{G^*}s}{Y'}$ is in duality with
$\cento{G(s)^\circ}Y$. Let
$(L_s,\al)$ and 
$(L_{Y,s},\al_Y)$ be unipotent $d$-cuspidal  data of $(G(s)^\circ,F)$ and
$(\cento{G(s)^\circ}Y,F)$  respectively, with dual data  
$(L^*_{s},\al)$ and
$(L^*_{Y,s},\al_Y)$ in $(\cento{G^*}s,F)$ and $(\cento{\cento{G^*}s}{Y'},F)$ respectively (see (1.3.1, 1.3.5)).

Let $b_{Y,s}$ be a block of
$\cent{G(s)^F}Y$ that covers $b_{\cento {G(s)^\circ}Y^F}(L_{Y,s},\al_Y)$  and $b_s$ be the block of $G(s)^F$ such that $$ (\{1\},b_s)\subset (Y,b_{Y,s})\,.$$

The
groups ${\rm E}_{G(s)^F}(Y,b_{Y,s})$ et ${\rm E}_{G^F}(Y,b_Y)$ are isomorphic. }

 \preuve $G(s)$ is defined around a couple of dual maximal tori $(T\subseteq G\cap G(s)^\circ,T^*\subseteq \cento{G^*}s)$.
 The $d$-cuspidal datum $(L,\la)$ in $(G,F)$ is associated to a set of unipotent $d$-cuspidal datum $(L^*_s,\al)$ in $(\cento{G^*}s,F)$ (Propositions~2.1.4, 2.1.7). As a consequence of 
$(L^*_{Y,s},\al_Y)\simex{\cento{G^*}s^F} (L^*_s,\al)$, we assumed in Proposition~4.2.3, one has 
$(L_{Y,s},\al_Y)\simex{G(s)^{\circ  F}} (L_s,\al)$. A subgroup of  ${\rm A}_{G^*}(s)^F$ acts transitively on the  set of $d$-cuspidal $\al$ corresponding to $\la$ by Propositions~1.3.6 and 2.1.4. The groups  ${\rm E}_{G^F}(Y,b_Y)$ defined by different $\al$ in the ${\rm A}_{L^*}(s)^F$-orbit are isomorphic. So we fix $\al$.

 Transport the duality between $G^\circ (s)$ and $\cento{G^*}s$ around dual $F$-stable maximal tori $S\subseteq G^\circ(s)$ and
$S^*\subseteq L^*_{s}$, $L_s$ and $L^*_s$ being $d$-split Levi subgroups in dual conjugacy classes such that, with notations of Proposition~3.4.1, ${\cal B}_{G,s}(b_s)=b_{G^F}(L,\la)$, where $b_s$ is one of the blocks of $G(s)^F$ that cover $b_{G(s)^{\circ\,F}}(L_s,\al)$.  Isomorphic or anti-isomorphic defect groups of $b_{\cento{G^*}s^F}(L^*_s,\al)$, $b_{G(s)^{\circ\,F}}(L_s,\al)$, $b_s$ and $b_{G^F}(L,\la)$ are known by Proposition~4.2.1. We may identify them and assume that, by a suitable choice of $S$,   $Y\subseteq \nor{G^F}S$, $Y\subseteq\nor{G(s)^\circ}S $ and $Y$ centralizes $[L_s,L_s]$. On dual side $Y'\subseteq \nor{(G^*)^F}{S^*}\cap \cento{G^*}s$ and $Y'$ centralizes $[L^*_s,L^*_s]$.

The blocks of $\cento GY^F$ that are covered by $b_Y$ are
$\cent {G^F}Y$-conjugate and $(L_Y,\la_Y)$ is defined up to $\cento GY^F$-conjugacy hence   we have an isomorphism
 $${\rm E}_{G^F}(Y,b_Y)\cong\nor{G}{Y,L_Y,\la_Y}^F/\nor{\cent GY}{L_Y,\la_Y}^F\leqno{(4.2.4.1)}$$

 (A) Assume $\z G$  connected, hence $G(s)$ connected. 
 
By Proposition~2.2.6 , if $(L_{Y,s},\al_Y)$ is a unipotent $d$-cuspidal
datum in  $(\cento{G(s)^\circ}Y,F)$ corresponding to $(L^*_{Y,s},\al_Y)$ in $(\cento{G^*}s,F)$, the groups 
${\rm E}_{\cento{G^*}s^F}( Y',b_{\cento{G^*}s^F}(L^*_Y,\al_{Y}))$ and
${\rm E}_{G(s)^{\circ\,F}}(Y,b_{G(s)^{\circ\,F}}(L_{Y,s},\al_{Y}))$ are anti-isomorphic.

 We may define duality between  $\cento GY$ and $\cento {G^*}{Y'}$ around tori $(T_Y\subseteq \cento GY)$, $(T^*_Y\subseteq \cento {G^*}{Y'})$ as in Proposition~1.2.6.
 We have $\nor{G^F}{Y,b_Y}=\nor{G^F}{Y,L_Y,\la_Y}\subseteq
\nor{G^F}{Y,L_Y,\ser{{L_Y}^F}s}$. 

 By Proposition~1.2.6, (c), the groups $\nor{G}{Y,T_Y}/\nor{\cent {G}Y}{T_Y}$ and $\nor{G^*}{Y',T^*_Y}/\nor{\cent {G^*}{Y'}}{T^*_Y} $  are
anti-isomorphic. So are   
$\nor G{Y,L_Y,T_Y}/\nor {\cent GY} {Y,L_Y}$ and $\nor
{G^*}{ Y',L_{ Y}^*,T^*_Y}/\nor {\cent {G^*}{ Y'}} { Y',{L_{ Y}}^*}$. Here $L^*_Y$ is in the dual $\cent{G^*}{Y'}^F$-conjugacy class of $L_Y$. Then the groups
$L_{ Y}^*$ are  $L^*_{Y,s}$ may be mutually defined, one by the other, as $L_Y$ and $L_{Y,s}$, by the relations
$L^*_{Y,s}=L_Y^*\cap\cent{G^*}s$ and
$L_Y^*=\cent{\cento{G^*}{Y')}}{\zo{L^*_{Y,s}}_{\phi_d}}$ (see Propositions~2.1.4 and 3.4.1). We obtain an anti-isomorphism between $\nor{G}{Y,L_Y,\ser{{L_Y}^F}s}^F/{L_Y}^F$ and
$\nor{\cent{G^*}s}{ Y',L^*_{Y,s}}^F/(L^*_{Y,s})^F$. 

Finally, by
Propositions  1.3.2, 2.1.4,  2.1.7 and 2.2.6 and isomorphism (4.2.4.1),  
${\rm E}_{G^F}(Y,b_Y)$ is anti-isomorphic to 
 $\nor{\cent{G^*}s}{ Y',
L^*_{Y,s}}^F_{\al_Y}/\nor{\cent{\cent{G^*}s}{Y'}}{L^*_{Y,s}}^F_{\al_Y}$,  that is exactly
${\rm E}_{\cent{G^*}s^F}( Y',b_{\cent{G^*}s^F}(L^*_{Y,s},\al_{Y}))$ by (4.2.4.1). 

 As ${\rm E}_{\cent{G^*}s^F}( Y',b_{\cent{G^*}s^F}(L^*_{Y,s},\al_{Y}))$ is anti-isomorphic to ${\rm E}_{G(s)^F}(Y,b_{G(s)^{ F}}(L_{Y,s},\al_Y))$, that gives the claimed isomorphism.

(B) To obtain the claimed isomorphism in case $\z G$ is not connected, consider as usual a regular embedding
 $G\to  H$, defining by restriction regular embeddings
$\cento GY\to
\cento HY$, $L_Y\to M_Y:=\z H.L_Y$ and dual maps $H^*\to G^*$ (where $t\in(H^*)^F_{\ell'}\mapsto s$), $M^*_Y\mapsto L^*_Y$.

The
group 
$\nor{G^F}{Y,L_Y,\ser{{L_Y}^F}s}/{L_Y}^F$ is a split extension of
$ \nor{H^F}{Y, M_Y,\ser{{ M_Y}^F} t}/{ M_Y}^F$ by the stabilizer ${\rm A}_{G^*}(s)^F_{ Y',L^*_{Y,s}}$  of the $G^F$-conjugacy class of $( Y',L^*_{Y,s})$ in ${\rm A}_{G^*}(s)^F$.  But ${\rm A}_{G^*}(s)^F_{ Y',L^*_{Y,s}}$  appears as a subgroup of  $\nor{\cent{G^*}s^F}{ Y'}/\nor{\cento{G^*}s^F}{Y'}$. By part (A) of the proof we have an exact sequence of anti-morphisms
$$\nor{\cento{G^*}s^F}{ Y',L^*_{Y,s}}/L^{*\, F}_{Y,s}\to
\nor{G^F}{Y,L_Y,\ser{L_Y^F}s}/L_Y^F\to {\rm A}_{G^*}(s)^F_{ Y',L^*_{Y,s}}\leqno{(4.2.4.2)}$$
By 
construction of $G(s)$ and duality,
${\rm A}_{G^*}(s)^F_{ Y',L^*_{Y,s}}$ is the image in  $G(s)^F/G(s)^{\circ\,F}$ of
$\nor{G(s)^F}{Y,L_{Y,s}}$.  The group ${\rm A}_{L^*_Y}(s)^F$
acts on $\ser{L^{*\, F}_{Y,s}}1$ and on $\ser{{L_Y}^F}s$ through the Jordan decomposition, as said in section 1.3. The  natural isomorphism between  $\nor{G^F}{Y,L_Y}/{L_Y}^F$
and $\nor{H^F}{Y, M_Y}/{{{M}_Y}}^F$ exchange Jordan \dec\
$\Psi_{M_Y, t}$ and the action of ${\rm A}_{L^*_Y}(s)^F$ (Proposition 1.3.1, (iv)). Furthermore any element of $\nor{\cent{G^*}s^F}{ Y',L^*_{Y,s}}$
or of $\nor{G^F}{Y,L_Y,\ser{L_Y^F}s}$ that fixes $\la_Y$ fixes any of its 
$\nor{H^F}{Y,M_Y}$-conjugate (Proposition~2.4.2). It follows that $\nor{G^F}{Y,L_Y,\la_Y}/{L_Y}^F$ has image in   ${\rm A}_{G^*}(s)^F_{ Y',L^*_{Y,s}}$ the stabilizer  of $ {\al_Y}$. The  sequence (4.2.4.2) restricts in 
$$\nor{\cento{G^*}s}{ Y',L^*_{Y,s}, {\al_Y}}^F/L^{*\, F}_{Y,s}\to
\nor{G}{Y,L_Y,\la_Y}^F/L_Y^F \to {\rm A}_{G^*}(s)^F_{ Y',L^*_Y, {\al_Y}}\to 1\leqno{(4.2.4.3)}$$  
By (4.2.4.1), ${\rm E}_{G^F}(Y,b_Y)$ is a quotient of that extension by $\nor{\cent  GY}{L_Y,\la_Y}^F/{L_Y}^F$. 

Considering $\cento GY$ instead of $G$, we have a split exact subsequence of (4.2.4.3)
$$\nor{\cento{\cento{G^*}{ Y'}}s}{L^*_{Y,s}, {\al_Y}}^F/L^{*\,F}_{Y,s}\to
\nor{\cento GY}{L_Y,\la_Y}^F/{L_Y}^F \to {\rm A}_{\cento{G^*}{Y'}}(s)^F_{L^*_{Y,s}, {\al_Y}}\to 1\leqno{(4.2.4.4)}$$
 Here $\cento{\cento{G^*}s}{Y'}=\cento{\cento{G^*}{Y'}}s$. The quotient group $\nor{\cento{G^*}s}{Y', L^*_{Y,s},{\al_Y}}^F/\nor{\cent{\cento{G^*}s}{Y'}}{ L^*_{Y,s},{\al_Y}}^F$ is anti-isomorphic to ${\rm E}_{G(s)^{\circ\,F}}(Y,b_{\cent{G(s)^\circ}Y^F}(L_{Y,s}, {\al_Y}))$. Furthermore ${\rm A}_{G^*}(s)$ is an
$\ell'$-group and $\cent{G}Y/\cento GY$ and $\cent{G^{* }}{Y'}/\cento {G^*}{Y'}$ are  $\ell$-groups. 
 
From (4.2.4.3) we obtain as quotient a split exact sequence of morphisms
$$1\to{\rm E}_{G(s)^{\circ\,F}}(Y,b_{\cent{G(s)^\circ}Y^F}(L_{Y,s}, {\al_Y}))\to {\rm E}_{G^F}(Y,b_Y)\to  {\rm A}_{G^*}(s)^F_{ Y',L^*_{Y,s}, {\al_Y}}/{\rm A}_{\cento{G^*}{ Y'}}(s)^F_{L^*_{Y,s}, {\al_Y}}\to 1$$  

A similar argument applies in $G(s)^F$, where $G(s)^{\circ\,F}.\nor{G(s)^F}{Y,b_{Y,s}}/G(s)^{\circ\,F}=  {\rm A}_{G^*}(s)^F_{ Y',L^*_{Y,s}, {\al_Y}}$. It gives easily he split extension
$$1\to {\rm E}_{G(s)^{\circ\,F}}(Y,b_{\cent{G(s)^\circ}Y^F}(L_{Y,s}, {\al_Y})) \to {\rm
E}_{G(s)^F}(Y,b_Y) \to {\rm A}_{G^*}(s)^F_{ Y',L^*_Y, {\al_Y}}/{\rm A}_{\cento{G^*}{
Y'}}(s)^F_{L^*_Y, {\al_Y}}\to 1$$ 
hence the claim isomorphism.

\bull

\vfill\eject

\noindent{\bf 5. Appendix }

\bigskip\noindent{\bf 5.1. Self-centralizing Brauer pairs, blocks and normal subgroups}

In that section we collect folklore results on blocks, pairs and Clifford theory with abelian quotient. 

Let $X$ be a finite group.
The abelian group $(X/[X,X])^\wedge$ acts on the set $\II X$ by tensor product, hence on the set of $\ell$-blocks of $X$  and on the set of  Brauer $\ell$-subpairs in a coherent  way as follows.  

Let $B$ be a block of $X$, $Q$ be an $\ell$-subgroup of
$X$, $\theta\in (X/[X,X])^\wedge$ and define $\theta_Q:=\Res X {\cent XQ}\theta$.  

The blocks $B$ and $\theta\otimes B$ are isomorphic and one has $\II{\theta\otimes B}=\{\theta\otimes
\chi\mid \chi\in\II B\}$.

An inclusion of subpairs
$(P,b_P)\subset (Q,b_Q)$ implies $(P,\theta_P\otimes b_P)\subset
(Q,\theta_Q\otimes b_Q)$.
The block  $B$ being defined by a maximal subpair
$(D,B_D)$, $\theta\otimes B=B$ is equivalent to $\theta_D\otimes b_D=b_D$ (note that $\theta_D$ is stable under $\nor XD$). 

The $\ell$-Sylow subgroup of  $(X/[X,X])^\wedge$ fixes any
$\ell$-block of $X$ ([28] Chapter 7, Corollary~5.6). Consider now  a subgroup $Y$ of $X$ with $[X,X]\subseteq Y$.
  Define
${\rm I}^X_Y(B)$ and ${\rm I}^X_Y(\chi)$ for $\chi\in\II X$ by $$(X/{\rm I}^X_Y(B))^\wedge=((X/Y)^\wedge)_B\,,\quad
(X/{\rm I}^X_Y(\chi))^\wedge=((X/Y)^\wedge)_\chi$$ 
Assume that $B$ covers a block $b$ of $Y$, let $X_b$ be the stabilizer of $b$ in $X$. If we consider  blocks as primitive central idempotents in algebras we have 
$\sum_{\theta\in\II{{\rm
I}^X_Y(B)/Y}}\theta\otimes B=\sum_{g\in X/X_b}b^g$.
There exists a unique block $b'$ of $X_b$ such that  $B={\rm Tr}_{X_b}^X(b')$. Then  $\theta\otimes B={\rm Tr}_{X_b}^X((\Res{X}{X_b}\theta)\otimes b')$. Thus ${\rm I}^X_Y(B)={\rm I}^{X_b}_Y(b')$. By a theorem of
Fong-Reynolds, ([28] Chapter~5, Theorem~5.10)
$$\II{{\rm Tr}_{X_b}^X(b')}=\{\Ind{X_b}X(\chi')\mid \chi'\in \II{b'}\}.\leqno{(5.1.1)}$$
  
Hence for any $\theta\in(X/X_b)^\wedge$ and $\chi=\Ind{X_b}X(\chi')\in\II
B$ we have $\theta\otimes \chi=\chi$. That proves
${\rm I}^X_Y(\chi)\subseteq X_b$.   Now if $\zeta\in \II Y$ one has clearly $X_\zeta\subseteq
X_{B_Y(\zeta)}$. The formula (5.1.1) has a twin for \irr\ representations that imply 
  ${\rm I}^X_Y(\chi)\subseteq X_\zeta$ when $\chi\in\II{X\mid\zeta}$. Finally $${\rm I}^X_Y(B_X(\chi))\subseteq {\rm I}^X_Y(\chi)\subseteq X_\zeta\subseteq
X_{B_Y(\zeta)}\quad{\rm when}\quad \chi\in\II{X\mid\zeta}.\leqno{(5.1.2)}$$

If there is no multiplicity in restrictions from $X$ to $Y$, with notations used in (5.1.1), (5.1.2), $\zeta $ extends in some $\zeta'$ to $X_\zeta$, different extensions has the form 
 $\theta\otimes \zeta'$ hence ${\rm I}^X_Y(\chi)=X_\zeta$. 

On defect groups recall a consequence of a theorem of Fong

\medskip\noindent {\bf 5.1.3. }([28] Chapter~5, Theorem~5.16) {\sl Let $Y$ be an invariant subgroup of a finite group $X$ and $b$ be an 
$\ell$-block of $Y$. There exist an $\ell$-block $B$ of $X$ that
cover
$b$ and a defect group $D$ of $B$ contained in $X_b$ such that
$D/D\cap Y$ is an $\ell$-Sylow subgroup of $X_b/Y$.}

\medskip\noindent{\bf 5.1.4. Proposition. }{\sl Let  $\rho\colon X\to X/Y$  be a morphism  of finite groups such that $X/Y$ is an abelian $\ell'$-group, let $b$ be an $\ell$-block of $Y$, $B$ an $\ell$-block of $X$.
 
(a)  (i) Let $(U,b_U)\subset (V,b_V)$ be an inclusion of $\ell$-subpairs in $Y$. If $B_U$ is an $\ell$-block of $\cent XU$ that covers $b_U$, there is an $\ell$-block of $\cent XV$ that covers $b_V$ and such that $(U,B_U)\subset (V,B_V)$ is an inclusion in $X$. All covering blocks of $b_V$ are so obtained.

(ii) Let be
$(D,b_D)$ a maximal subpair of  $Y$,  $(D,B_D)$ a maximal subpair of  $X$ such that $(1,b)\subset (D,b_D)$ and $(1,B)\subset (D,B_D)$. Then $B$ covers $b$ if and only if  some $\nor XD$-conjugate of $B_D$ covers $b_D$.

(b) Let 
$(D,b_D)$ be a maximal subpair in  
$Y$ with canonical caracter $\la\in\II{\cent YD}$ and such that $(1,b)\subset (D,b_D)$. One has $X_b=Y.\nor XD_{\la}$.
 
 (c) Assume $B$ covers $b$ and let $\la$ as in (b). Assume non multiplicity in $\Ind{\cent XD}{\cent Y D} \la$ and let $\mu$ a component of $\Ind{\cent XD}{\cent Y D} \la$.

There is a morphism $f\colon \nor
XD_{\la}\to  \II{\cent XD_{\la}/\cent YD)}$ whose kernel is $\nor X
D_{\la}\cap\nor XD_\mu$. Define  $J$ by the equality
$f( \nor XD_{\la})=\II {\cent
XD_{\la}/J}$ . One has 
$\rho(J)=\rho({\rm I}_Y^X(B))$. The set of blocks of $X$ that cover $b$ is a reguler orbit under $\II{\rho(J)}$. }

\preuve
As $|X/Y|$ is prime to $\ell$, if $B$ covers $b$, then $B$ and $b$ have a common defect group by 5.1.3.

(a) (i) Consider  a ``normal inclusion of subpairs" in $Y$, $(U,b_U)\triangleleft (V,b_V)$ [16], Definition~5.2 :  $U$ and $V$ are $\ell$-subgroups of $Y$ such that $U$ is normal in $V$,  $b_U$, $b_V$,  are $\ell$-blocks of respectively 
$\cent Y U$, $\cent YV$ and $V$ fixes $b_U$. The inclusion of subpairs  in $Y$
 writes with Brauer's morphism  Br$^Y_V$ on the images in  caracteristic $\ell$ :
${\rm Br}^Y_V(\bar {b_U}).\bar {b_V}=\bar {b_V}$ and $b_U$ is  unique when  $b_V$ is given. 
We have
$$\sum_{x\in
\cent XU/\cent XU_{b_U}}{b_U}^x=\sum_jB_{U,j},\quad  \sum_{y\in \cent
XV/\cent XV_{b_V}}{b_V}^y=\sum_kB_{V,k}$$ where $\{B_{U,j}\}_j$ is 
the set of blocks of $\cent XU$ that cover
$b_U$, an orbit under the action of $(\cent XU/\cent YU)^\wedge$, or of $(X/Y)^\wedge$ {\it via} $(\theta\mapsto \theta_ U)$, and 
$\{B_{V,k}\}_k$ is the set of blocks of $\cent XV$ that cover
$b_V$, an orbit under the action of $(\cent XV/\cent YV)^\wedge$, hence of $(X/Y)^\wedge$.

  One has  $\sum_{x\in
\cent XU/\cent XU_{b_U}}{\rm Br}^Y_V(\bar{b_U}^x)=\sum_j{\rm Br}^X_V(\bar{B_{U,j}})$, an equality between  sums of idempotents. As $\bar{b_V}$ appears on the left as primitive in the center of $k\cent YV$ and  the sum is stable under the actions of $\cent XV/\cent YV$ and of $(X/Y)^\wedge$, any $B_{V,k}$ appears on the right side. That proves a normal inclusion $(U,B_{U,j})\triangleleft (V,B_{V,k})$ in $X$, in fact a map $(k\mapsto j)$ that commutes with the action of $(X/Y)^\wedge$.
 
 As inclusion of subpairs is defined by transitivity from normal inclusions, we have (i) in (a).
 
(ii)  Given $B$, $B_D$ is defined  modulo $\nor XD$-conjugacy. (ii) follows from (i).

(b) 
$\la$ has central defect and defines $b_D$. 
By Frattini argument, as all maximal pairs in $Y$ containing $(1,b)$ are $Y$-conjugate we have 
$X_b=Y.\nor X D_{\la}$. 

(c) 
We use (a) and consider the blocks of
$\cent XD$ that cover
$b_D$. They are given by their canonical character
$\mu_j=\Ind{\cent XD_{\la}} {\cent XD}\la'_j$,
  $\la'_j\in \II{\cent XD_{\la}\mid\la}$ ($\la'_j$ and $\mu_j$ have central defect). As each $\la'_j$ extends $\la$ (no-multiplicity hypothesis) we have
$m:=|\cent XD_{\la}/\cent YD|$ distincts blocks $B'_j$ with central defect of $\cent XD$, a regular orbit under
$ (\cent
XD_{\la}/\cent YD)^\wedge$, that cover $|\cent XD/\cent
XD_{\la}|$ blocks of $\cent YD$. As $X/Y$ is abelian we have $[\nor
XD,\cent XD]\subseteq \cent YD$, hence $\nor X D_{\mu_j}=\nor XD_{B'_j}$ does not depend on $j$, denote it  $\nor XD_{\mu}$. If
$x\in\nor XD$ and  ${\mu_j}^x=\mu_k$, then, by restriction to $\cent XD$,  ${B_D}^x=B_D$.
By conjugacy of maximal subpairs containing a block, a block of $X$ that covers  $b$ corresponds to an orbit under $\cent XD.\nor XD_{b_D}$ in the set of blocks of $\cent XD$ that  cover $b_D$. 
Recall that  $\nor XD_{b_D}=\nor XD_{\la}$. We have a regular orbit under $\cent XD.\nor XD_{\la}/\cent XD.\nor
XD_{\la}\cap \nor XD_\mu$, a quotient  isomorphic to $\nor XD_{\la}/\nor
XD_{\la}\cap \nor XD_\mu$.
  
We have obtained an exact sequence
$$1\to \nor XD_{\la}\cap\nor XD_\mu\to \nor XD_{\la}\to \big( {\rm I}^{\cent XD}_{\cent
YD}(B_D)/\cent YD\big)^\wedge\to\big({\rm I}^X_Y(B)/Y\big)^\wedge\to 1$$
\bull

If, in  Proposition 5.1.4,  $\la$ extends to $\mu\in\II{\cent XD}$ and $\mu$ is stable under $\nor XD_\la$ --- and that is the case if $B$ is the principal block --- then
$\rho(\cent XD)^\wedge$ acts regularly on the set of blocks of $X$ that cover $b$.
 
 \medskip\noindent{\bf 5.1.5. Proposition. }{\sl Let $Y\subseteq X$ be finite groups, $Y$ invariant in $X$,  and $b$ a block of $Y$, $B$ a block of $X$ such that $B$ covers $b$. Assume that $B$ has central defect.

(a)  If there is  $\chi\in\II B$ of height zero and whose restriction to $Y$ has multiplicities prime to $\ell$, then $b$ has central defect group and $\z Y_\ell=Y\cap
\z X_\ell$. 

(b) If $X/Y$ is $\ell$-solvable, then (a) applies and $b$ has central defect group.

(c) If $b$ has central defect, there is $\chi\in/II B$ as in (a).

}
\preuve  There is some
 $\chi\in\II B$ whose kernel contains the defect group,  of height 0  : $\chi(1)_\ell=|X/\z X|_\ell$. Let $\xi\in\II b$ such that $ \chi\in\II{X\mid \xi}$ [28], Chapter~5, Lemma~5.7. 
 
 (a)  If $\scal {\Res XY\chi}\xi Y$ is prime to $\ell$, and that is clearly true if $X/Y$ is an $\ell'$-group or if $X/Y$ is cyclic,  then we have $\chi(1)_\ell=\xi(1)_\ell |X : X_\xi|_\ell$, hence $\xi(1)_\ell=|X_\xi/\z
X|_\ell$. We have $\z X\subseteq X_\xi\subseteq X_b$.  By 5.1.3 on defect groups of covering blocks we have
$|X_b/Y\z X|_\ell=1$. Thus $|X_b|_\ell=|X_\xi|_\ell=|Y.\z X|_\ell$,
hence
$\xi(1)_\ell=|Y/Y\cap\z X|_\ell$. But $Y\cap\z X\subseteq \z Y$ and
$\xi(1)$ divides
$|Y/\z Y|$. Hence  $\z Y_\ell=Y\cap \z X_\ell$ and the defect group of  $b$
is $\z Y_\ell$. 

(b) The condition $\z Y_\ell=\z X_\ell\cap Y$
is equivalent to $\z Y_\ell\subseteq\z X$ and to  $\z {Y.\z X_\ell}_\ell=\z X_\ell$. If true that condition implies similar ones for any group $R$
between $Y$ and $X$: we have then  $\z Y_\ell\subseteq \z
R$ and $\z R_\ell\subseteq \z X$. If $X/Y$ is $\ell$-solvable there is a sequence $1=K_0\triangleleft  K_1\triangleleft\dots\triangleleft K_k \triangleleft K_{k+1}=X$ such that each quotient $K_{j+1}/K_j$ is a cyclic $\ell$-group or an $\ell'$-group, hence (a) applies.

(c) If $b$ has central defect as $B$, then $a:=\scal {\Res XY\chi}\xi Y$ is prime to $\ell$ for some $\chi$ :

Let  $\xi$ be the canonical character of $b$ and assume $\chi\in\II{X\mid\xi}$. We have $\xi(1)_\ell=|Y/\z Y|_\ell$ and $\chi(1)=a\xi(1) |X|/|X_\xi|$,
$|X_\xi|_\ell=|X_b|_\ell=|Y.\z X|_\ell$. Hence $\chi(1)_\ell.|\z Y/\z
X\cap Y|_\ell=a_\ell |X/\z X|_\ell$. But  $\chi(1)$ divides $|X/\z X|$ and
$\z Y_\ell=Y\cap
\z X_\ell$, so that  $a_\ell=1$.
\bull

\medskip An alternative proof when $Y/X$ is an $\ell$-group :

As $\z X_\ell$ is in the kernel of the canonical character $\chi$ of $B$,  
$\z X_\ell\cap Y$ is contained in the kernel of any component  $\xi$ of $\II b\cap \Res
XY\chi$, we may assume $\z X_\ell=1$. Then by ???? $X_b=Y$, hence $X_\xi=Y$. It follows that
$\xi$ is projective as is $\chi$, so that $b$ has null defect and $\z Y_\ell=1$. Coming back  to the general case $\z H_\ell\subseteq \z G_\ell\cap H$.

\bigskip

\noindent{\bf 5.2.  On unipotent series}

\medskip 
The following reminds us a result of Geck [25]. If $d=1$, by the theory of Hecke algebras, in an Harish-Chandra series $\ser{G^F}{(L,\la)}$, to the sign representation of the associated Hecke algebra there corresponds a unique \irr\ character whoose degree has maximum $p$-valuation and with multiplicity $1$ in $\Lu LG\la$. In Proposition~5.2, (a) we find in classical types an  \irr\ character whoose degree has minimum $p$-valuation and with multiplicity $1$ in $\Lu LG\la$... Assertion (b) of Proposition~5.2 is true for $d=1$ by Geck's argument. We use notations and definitions of section~2.2. Our reference is [26], Appendix.

\medskip\noindent{\bf 5.2. Proposition. } {\sl Let $(G,F)$ be a connected   reductive  group defined on $\F q$. Let
$d$ be a positive integer and 
$(L,\al)$ a $d$-cuspidal unipotent datum  of $(G,F)$,  defining the Generalized $d$-Harish-Chandra series $\ser{G^F}{(L,\al)}\subseteq \ser{G^F}1$.  Assume $d>1$. 

(a) Assume all components of $G$ have classical type. Let  $\chi_0\in\ser{G^F}{(L,\al)}$ such that $\chi(1)_p$ is minimal on $\chi_0$ for $\chi$ running in the G.$d$-HC series $\chi_0\in\ser{G^F}{(L,\al)}$. Then $\chi_0(1)_p\leq \al(1)_p$, with equality if $p>2$, and $\chi_0(1)\neq \chi(1)$ for any other $\chi\in \ser{G^F}{(L,\al)}$. 
Furthermore $|\scal{\Lu
LG\al}{\chi_0}{G^F}|=1$.

 (b) In any type there exist  $\chi\in\ser{G^F}{(L,\al)}$ such that $\scal{\Lu
{ L} {G} \chi} \al{ {G}^F}\in\{-1,1\}$ and $\xi(1)\neq \chi(1)$ for any $\xi\in\ser{G^F}{(L,\al)}$ different from $\chi$. }

\preuve By inspection. 
As the unipotent Lusztig series is invariant in an isotypic morphism (see 1.3.1), we may assume that $G$ is rationally \irr, finally is \irr. 

(a.1) Type 
$\AA_{n}$.

Assume first that $G$ is split. One has a parametrization of $\ser{G^F}1$ by partitions of $n+1$. A partition $\Lambda=\{0\leq\la_1\leq \dots \leq \la_m\}$  may be defined by a so-called ``$\beta$-set" $B(\Lambda)=\{0<b_1<\dots b_m\}$ where $b_j=\la_j+j$ ($1\leq j\leq m$). One may assume $m$ as large as we want, adding $0$ at the beginning of $\{\la_j\}_j$. Assume $\Lambda$ define $\chi\in\ser{G^F}1$. By the degree formula  given in [26] Appendix, we have $\chi(1)_p=q^{v(B(\Lambda))}$,
where
$$v(B(\Lambda))=\sum_{1\leq j< m}(m-j)b_{j}-\sum_{1\leq j<m-1}{m-j\choose
2}$$ 

If $\Lambda_1$ is a partition of $(n+d+1)$ with the same number of components than $\Lambda$, one says that $\Lambda$ may be deduced from $\Lambda_1$ by deleting a hook of length $d$ if and only if 
 $B(\Lambda_1)=B(\Lambda)\dot + \{b_j,b_j+d\}$ for some $j\in[1,m]$ (here $\dot +$
is boolean addition and $b_j\in B(\Lambda)$, $b_j+d\in B(\Lambda_1)$).  Clearly if $j=m$ then $v(B(\Lambda_1))=v(B(\Lambda))$. Let  $k$ such that $b_k<b_j+d<b_{k+1}$ (with $k=m$ if $b_m<b_j+d$). By the above formula we have 

\centerline{$v(B(\Lambda_1))-v(B(\Lambda))=\sum_{j<s\leq k}b_s+(m-k)d-(k-j)b_j=\sum_{j<s\leq k}(b_s-b_j)+(m-k)d$.}

\noindent The two terms on the right are non negative and they are null only if $j=k$ and $m=k$, that is $j=m$.

A $d$-cuspidal element in $\ser{G^F}1$ is defined by a $d$-core $\kappa$, or partition without any hook of length $d$, equivalently : if $b_j\in B(\kappa)$ and $b_j\geq d$, then $b_j-d\in B(\kappa)$.
 A $d$-cuspidal unipotent datum in $(G,F)$ is defined by a $d$-split Levi subgroup of type $\AA_{n-td}$ ($t\in \N$) and $\la\in\ser{L^F}1$, $\la$ being defined by a $d$-core. The elements of $\ser{G^F}{(L,\la)}$ corresponds to the set of partitions of $(n+1)$ from which one can go down to $\kappa$ by deleting a sequence of hooks of length $d$. By what we saw in the case $t=1$ there is a unique element $\Lambda_0$ in that set with minimal parameter $v$, precisely such that $v(B(\Lambda_0))=v(B(\kappa))$. Furthermore there is only one way to go from $B(\Lambda_0))=\{0<b_1<\dots b_{m-1}<b_m+td\}$ to $B(\Lambda)$ by deleting a sequence of hooks of length $d$, and that implies 
 $|\scal{\Lu
LG\al}{\chi}{G^F}|=1$.

In connection with Ennola's conjecture, there is a bijection from $\ser{{\rm GL}_{n+1}(q)}1$ to $\ser{{\rm GU}_{n+1}(q)}1$, $(\chi\mapsto \chi_1)$ that preserves the $p$-valuation of degrees. There is also a one-to-one map $\N\to\N$, $(d\mapsto d_1)$ and a correspondance $(L\mapsto L_1) $ between ${\rm GL}_{n+1}(q)$-conjugacy classes of $d$-split Levi subgroups of  ${\rm GL}_{n+1}$ and ${\rm GU}_{n+1}(q)$-conjugacy classes of  $d_1$-split Levi subgroup of ${\rm GU}_{n+1}$ that commutes with Lusztig induction so that $\chi$ is $d$-cuspidal if and only if $\chi_1$ is $d_1$-cuspidal, and $\ser{{\rm GL}_{n+1}(q)}{(L,\la)}$ is sent onto $\ser{{\rm GU}_{n+1}(q)}{(L_1,\la_1)}$.

Hence the property we claim goes from type $\AA$ to type $\lexp 2\AA$.

 (a.2) Types $\BB_n$, $\CC_n$. 

We use the following elementary lemma :

\medskip\noindent{\bf Lemma. }{\sl Let $\{b_j\}_{1\leq j\leq N}$ be a sequence of natural integers such that $b_j>b_1$ if $j\neq 1$ and $N$ is even. Let $\epsilon \colon [1,N]\to\{-1,1\}$, $q\in\N$, $q\geq 2$. Then
$$\prod_{1\leq j\leq N}(q^{b_j}-\epsilon(j))\neq   \prod_{1\leq j\leq N}(q^{b_j}+\epsilon(j))$$
}

\noindent{\it Proof of the lemma. } Assume equality. 
The general term in the developpement of each product is defined by a subset $Y$ of $[1,N]$ and writes on the right side $t(Y)=(\Pi_{y\notin Y}\epsilon(y))q^{\Sigma_{y\in Y}b_y}$, on the leftt side $(-1)^{N-|Y|}t(Y)$. If $|Y|$ is even these terms are equal, hence simplify. If $|Y|$ is odd, they differ by a minus sign.

If $q\neq 2$, reducing modulo $q^{b_1+1}$ we obtain $2q^{b_1}\equiv  0 \pmod {q^{b_1+1}}$, a contradiction.

If $q=2$, reducing modulo $2^{b_1+2}$ we obtain $2^{b_1+1}\equiv  0 \pmod {q^{b_1+2}}$, a contradiction.  \bull

An element 
$\chi$ in
$\ser{G^F}1$ is defined by a symbol $(B,C)$ of rank $n$,
where $B\subseteq \N$, $C\subseteq \N$,

\noindent $B=\{b_0=0<b_1<\dots <b_j<\dots b_s\}$, $C=\{c_0=0<c_1<\dots <c_j<\dots c_t\}$ and $(s-t)\in\N$ is odd. By [26] the $p$-valuation of $\chi(1)$ may be written  $(2^{|B\cap
C|})_{p}.q^{v(B,C)+f(s+t)}$ with $$v(B,C)=\sum_{1\leq j\leq s}(s-j+|C\cap
]b_j,\infty [|)b_j+\sum_{1\leq j\leq t}(t-j+|B\cap
]c_j,\infty [|)c_j+\sum_{b\in B\cap C}b $$ 

We note that if $b_s>c_t$ (resp. $c_t>b_s$), then $v(B,C)$ is
independant of $b_s$ (resp. $c_t$), but, as we'll see, is growing up with the other values $b_j$, $c_j$. 

(i) Assume $d$ odd. Then $\al$ is defined by
a symbol $(B_0,C_0)$ of rank $k$, where $(n-k)\in d\N$, and such that the partitions associated to
$B_0$ and $C_0$ respectively are $d$-cores. The symbol $(B,C)$ defines an element of $\ser{G^F}{(L,\al)}$
if and only if one can go from $B$ to $B_0$, and from $C$ to $C_0$ deleting  a sequence of hooks of length $d$. Let $(B,C)$ as above and consider $(B_1,C)$ with $B_1=B\dot+\{b_j,b_j+d\}$ for some $j$ such that $b_j+d\notin B$. Let $r\in[j,s]$ maximum such that $b_r<b_j+d$. One has 

\noindent $v(B_1,C)-v(B,C)=\sum_{j<i\leq r}(b_i-b_j)+(s-r)d+|C\cap]b_j+d,\infty[|d+\sum_{\{i\mid b_j<c_i\leq b_j+d\}}(c_i-b_j)$. 

\noindent  It is a sum of four non negative integers that are all null only under the conditions : $b_j$ is maximum in $B$ (that is $j=s$) and $C\subseteq[0,b_j]$, hence $b_s\geq c_t$. For $p\neq 2$ these are the conditions under which  the $p$-valuation of the degree of corresponding \irr\ unipotent character does not increase when going from $(B,C)$ to $(B_1,C)$.

For $p=2$ we have to take account of the factor  $(2^{|B\cap
C|})_{p}$: to obtain the variation of the $p$-contribution  to the degree, multiply $q^{v(B_1,C)-v(B,C)}$  by $2$ if $b_j+d\in C$ but $b_j\notin C$, and by $1/2$ if $b_j\in C$ but $b_j+d\notin C$. The valuation in $2$ may decrease strictly if $c_s=b_t$ and $j=t$.

Assume ${\rm sup}(C_0)\neq {\rm sup}(B_0)$. By our description in type $\AA$ and the value of $v(B,C)$ above, there exists one and only one symbol
$(B,C)$ of rank $n$ in the  $d$-series above $(B_0,C_0)$ such that
$v(B,C)\leq v(B_0,C_0)$ and then $v(B,C)=v(B_0,C_0)$ : if
$b_s>c_t$, it is 
$(B_1,C_1)=(B_0\dot +\{b_s,b_s+n-k\},C_0)$. If $\chi_1\in\ser{G^F}1$ has parameter $(B_1,C_1)$ $\chi_1(1)_p$ is minimal in the $d$-series $\ser{G^F}{(L,\al)}$ (sketch specially the case $p=2$, assuming $d>1$). Then  $|\scal{\Lu LG\al}{\chi_1}{G^F}|=1$, because there is only one way to go from $(B_1,C_1)$ to $(B_0,C_0)$ by deleting successively $(n-k)/d$ hooks of length $d$.

Assume now $a:={\rm sup}(C_0)= {\rm sup}(B_0)$. There are exactly two different symbols 
$$(B_1,C_1)=(B_0,C_0\dot
+\{a,a+n-k\}),\quad (B_2,C_2)=(B_0\dot +\{a,a+n-k\},C_0)$$ of rank $n$ in the $d$-series of $(B_0,C_0)$ such that  $v(B_i,C_i)\leq v(B_0,C_0)$, ($i=1,2$),
and then $v(B_i,C_i)=v(B_0,C_0)$ (short of the case $p=2$, see above).  We have to compare the degrees of the $\chi_i\in\ser{G^F}1$  corresponding to $(B_i,C_i)$ ($i=1,2$). The degree formula given in [26] shows that if $(B,C)$ is the parameter of $\chi\in\ser{G^F}1$, 
$\chi(1)_{p'}$ is, up to a constant depending only of $(|G|,s,t$) equal to  
$${(2^{|B\cap
C|})_{p'}\over \prod_  \gamma(q^{l(\gamma)}-1)
\prod_\delta(q^{l(\gamma_1)}-1)
\prod_{\gamma'}(q^{l(\delta)}+1)
}$$ where $\gamma$ (resp. $\gamma_1$, $\delta$) runs on the set of hooks of $B$ (resp. hooks of $C$, cohooks of $(B,C)$) and the function $l$ design the length of the hook or cohook. We have 
$$\bigg(2_{p'}.(q^{n-k}-1).\prod_{1\leq u\leq n-k-1}(q^{2u}-1)\bigg)\bigg({\chi_1(1)\over \al(1)}\bigg)_{p'}=\prod_{c<a,c\notin C_0}{q^{a-c}-1\over q^{a+n-k-c}-1}\prod_{b<a,b\notin B_0}{q^{a-b}+1\over q^{a+n-k-b}+1}$$
$$\bigg(2_{p'}.(q^{n-k}-1).\prod_{1\leq u\leq n-k-1}(q^{2u}-1)\bigg)\bigg({\chi_2(1)\over \al(1)}\bigg)_{p'}=\prod_{c<a,c\notin C_0}{q^{a-c}+1\over q^{a+n-k-c}+1}\prod_{b<a,b\notin B_0}{q^{a-b}-1\over q^{a+n-k-b}-1}$$
hence 
$$\bigg({\chi_1(1)\over\chi_2(1)}\bigg)_{p'}=\prod_{c<a,c\notin C_0}{(q^{a-c}-1)(q^{a+n-k-c}+1)\over (q^{a+n-k-c}-1)(q^{a-c}+1)}\prod_{b<a,b\notin B_0}{(q^{a-b}+1)(q^{a+n-k-b}-1)\over (q^{a+n-k-b}+1)(q^{a-b}-1)}$$
Let $\bar B:=\{a-b\mid b<a, b\notin B_0\}$, $\bar C=\{a-c\mid c<a,c\notin C_0\}$, $D:=n-k+\bar B$, $E=n-k+\bar C$. Assume $\chi_1(1)_{p'}=\chi_2(1)_{p'}$. We have 
$$\prod_{x\in \bar C}(q^x-1)\prod_{x\in\bar B}(q^x+1)\prod_{x\in D}(q^x-1)\prod_{x\in E}(q^x+1)=\prod_{x\in \bar C}(q^x+1)\prod_{x\in\bar B}(q^x-1)\prod_{x\in D}(q^x+1)\prod_{x\in E}(q^x-1)$$
in which we may exchange $(\bar B,E)$ and $(\bar C,D)$. 
If $x\in\bar B\cap\bar C$, the factor $(q^{2x}-1)$ appears in each side of that equality. So we may assume $\bar B\cap\bar C=\emptyset$, as well as $D\cap E=\emptyset$. Similarly if $x\in\bar C\cap E$ or $x\in \bar B\cap D$ there is a simplification by $(q^{2x}-1)$. Finally we have an equality as above  with $\bar C\cap\bar B=\bar C\cap E=D\cap E=\bar B\cap D=\emptyset$, $\bar B\neq \emptyset$, $|\bar C|=|E|$, $|\bar B|=|D|$. Then let $y$ be the smallest element of $\bar B\cup\bar C\cup D\cup E$. If $y\in\bar B$ (resp. $y\in\bar C$), then $y\notin D$, (resp. $y\notin E$) because ${\rm Inf} (D)\geq n-k+{\rm Inf} (\bar B)$
(resp. ${\rm Inf} (E)\geq n-k+{\rm Inf} (\bar C)$). The Lemma applies : $\chi_1(1)_{p'}$ and $\chi_2(1)_{p'}$ are different.
 
 The equalities $|\scal{\Lu LG\al}{\chi_1}{G^F}|=1=|\scal{\Lu LG\al}{\chi_2}{G^F}|=1$ are true as in the first case.
 
 (ii) Assume $d$ even.

In that case $\al$ is defined by a symbol $(B_0,C_0)$ with no cohook of length $d/2$. If an element $\chi\in\ser{G^F}1$ defined by a symbol $(B,C)$, it belongs to the $d$-series $\ser{G^F}{(L,\al)}$ if one can go down from $(B,C)$ to $(B_0,C_0)$ by deleting  a sequence of cohooks of length $d/2$. 

As in (i) the degree formula imply that if ${\rm Sup}(B_0)>{\rm Sup}(C_0)$ there is only one symbol $(B_1,C_1)$ of rank $n$ in the $d$-series of $(B_0,C_0)$ such that $v(B_1,C_1)=v(B_0,C_0)$, that is $(B_0\dot +\{b_s\},C_0\dot+\{b_s+n-k\})$ if $2(n-k)/d$ is odd, and $(B_0+\{b_s,b_s+n-k\},C_0)$ if $2(n-k)/d$ is  even. The sequence of cohooks to go down from $(B_1,C_1)$ to $(B_0,C_0)$ is unique so that, if $\chi_1\in\ser{G^F}1$ has parameter $(B_1,C_1)$ then  $\scal{\Lu LG\al}{\chi_1}{G^F}=1$.

Assume ${\rm Sup}(B_0)={\rm Sup}(C_0)=a$. We find two symbols of rank $n$ in the $d$-series of $(B_0,C_0)$ : $(B_0\dot +\{a\},C_0\dot+\{a+n-k\})$ and $(B_0\dot +\{a+n-k\},C_0\dot+\{a\})$. Listing variations on the sets of hooks and cohooks between $(B_0,B_0)$ and $(B_i ,C_i)$ ($i=1,2$)
one finds the same quotient $[\chi_1(1)/\chi_2(1)]_{p'}$ as above, hence $\chi_1(1)_{p'}\neq\chi_2(1)_{p'}$.

(a.3) Types $\DD_n$, $\lexp 2{\DD_n}$. 

An element
$\chi$ in
$\ser{G^F}1$ is defined by a symbol $(B,C)$ of rank $n$,
where  $||B|-|C||\in 4\N$ if $(G,F)$ is split and
$|B|-|C|\equiv 2\pmod 4$ if not. Furthermore, if $B\neq C$, $(B,C)$
and
$(C,B)$ define the same \irr\ character, but
$(B,B)$ defines two $\chi\in \ser{G^F}1$, said to be ``twins". 
A $d$-cuspidal $\la$ is defined par a symbol $(B_0,C_0)$ without any hook of length $d$ (in case $d$ is odd) or without any cohook of length $d/2$ (in case $d$ is even). The $p$-valuation and the $p'$-component of the degrees are given by the same formulas as in type $\BB$, up to a constant factor. If $(B_0,C_0)$ is not symmetric the symbols $(B_i,C_i)$ given in (i) (if $d$ is odd) or (ii) (if $d$ is even) are not symmetrics. The proofs given in (a.2) go on.

Thus we consider a symmetric symbol $(B_0,C_0)$.  The twins element of $\ser{L^F}1$ defined by $(B_0,C_0)$ are $\nor{G^F}L$-conjugate if $L\neq G$.
Then the symbols $(B_i,C_i)$, as defined in (a.2), are not symmetric and verify $(B_1,C_1)=(C_2,B_2)$, hence we obtain a unique $\chi\in\ser{G^F}1$ with claimed properties.

(b) Exceptional types.

One may conjecture that (a) generalizes to exceptional types as well as the following  variation : 

{\sl There exist  $\chi\in\ser{G^F}{(L,\al)}$ such that 
  
 \centerline{   for all $\xi\in\ser{G^F}{(L,\al)}$, $\xi(1)_p$ divides $\chi(1)_p$; 
if $\xi(1)_p=\chi(1)_p$ and $ \xi\neq \chi$ then $\xi(1)<\chi(1)$ .}

Furthermore $|\scal{\Lu
LG\al}{\chi}{G^F}|=1$.
   }

But the verification is quite boresome because in some G.$d$-HC unipotent series one may find different elements with equal generic polynomial degree. In all cases there is a unique element whose generic degree has minimal $q$-valuation, or a unique element whose generic degree has maximal $q$-valuation, and that gives (b).

If $L$ is a maximal torus $T$, so that $\al=1_{T^F}$, $1_{G^F}$ has the property, and the Steinberg character the opposite one, thanks to [20], 12.7, 12.8. 
 
 The $d$-unipotent series $\ser{G^F}{(L,\al)}$ when $L$ is not a torus are given in [10], Table 1, Table 2, with multiplicities $\scal\chi{\Lu LG\al}{G^F}$. The generic degrees, polynomials in $q$, of element of $\ser{G^F}1$ are as in [26, Appendix]. The verification is immediate.

In series that are known only by the sum of two algebraic conjugate, as $15+16$, $40+41$, $42+43$ in notations of [10], two algebraic conjugate \irr\ components (which have equal degree) are always in different series. 
\bull

\bigskip\noindent{\bf 5.3.  On a commutation formula in classical type}

\medskip Here we prove the commutation formula (J3) we considered in section 1.3.3, for groups of classical type. 

\medskip\noindent{\bf 5.3. Proposition. }{\sl Assume $G$ has classical type and a connected center. Let $s\in G^{*\,F}$ be semi-simple, let $L^*$ and $M^*$ be $F$-stable Levi subgroups in $G^*$ such that $s\in L^*\subseteq M^*$, and $L\subseteq M$ be $F$-stable Levi subgroups of $G$ in dual $G^F$-conjugacy classes. Assuming Mackey decomposition formula between $F$-stable Levi subgroups of $G$, one has
$$\Lu LM\circ\Psi_{L,s}=\Psi_{M,s}\circ \Lu{\cent{L^*}s}{\cent {M^*}s}\leqno{(5.3.1)}$$}

\preuve By our general conventions we assume that the duality between $L$ and $L^*$, $M$ and $M^*$, $G$ and $G^*$ may be defined around the same pair of dual maximal tori $(T\subseteq L, T^*\subseteq L^*)$ and so on for any ``dual" sets $\{L_j\}_j$, $\{L^*_j\}_j$ of Levi subgroups of $G$, $G^*$ with a common maximal torus $T\subseteq \cap_jL_j$, $T^*\subseteq\cap_jL^*_j$. Similar restrictions apply to dualities over $\F q$ between $F$-stable groups. It follows that the choice of duality between $L$ and $L^*$, hence the eventual different choices of $\Psi_{L,s}$, if $s\in L^*$, do not affect $\Lu LM\circ\Psi_{L,s}$.  

As a consequence of Mackey formula $\Lu{L\subset P}M\circ\Psi_{L,s}$ is independant of the choice of $P$. By duality (5.3.1) implies 
$$\Psi_{L,s}\circ \slu{\cent{L^*}s}{\cent {M^*}s}=\slu LM\circ\Psi_{M,s} \leqno{(5.3.2)}$$

We prove (5.3.1) by induction on the semi-simple rank of $M$. If $M$ is a torus, $L=M$, there is nothing to prove. If $L$ is a torus the formula follows from (ii) in Proposition~1.3.1.

(a) A " trivial case" is the following : it may happen that $ \cent{M^*}s= \cent{L^*}s$ when $\zo {L^*}\subseteq \zo { \cent{M^*}s}$ (note that $\zo{L^*}\subseteq\cent{G^*}s$ and $\cent{L^*}s=\cent{\cent{G^*}s}{\zo {L^*}}$). In that case $\Lu {L}M$ restricts to and is defined by a bijection $\ser{L^F}s\to \ser{M^F}s$ that commutes with Jordan \dec\ by  (iii) of Proposition~1.3.2. So we have a special case of  (5.3.1) : $\Lu LM\circ \Psi_{L,s}=\Psi_{M,s}$.

 (b) By transitivity of Lusztig induction, if there exists some $F$-stable Levi subgroup $K$ of $G$ with $L\subseteq K\subseteq M$, $L^*\subseteq K^*\subseteq M^*$, such that $\Psi_{M,s}\circ \Lu{\cent{K^*}s}{\cent{M^*}s}=\Lu KM\circ \Psi_{K,s}$ and $\Psi_{K,s}\circ \Lu{\cent{L^*}s}{\cent{K^*}s}=\Lu LK\circ \Psi_{L,s}$, we are done. So we write $\Lu {L}M$ as a product of minimal "steps". 
 
 (c) To a Levi subgroup $L$ with $s\in  L^*$  we associate an eventually smaller Levi subgroup $L_s$ as follows : let $L_s\subseteq L$ such that  $L_s$ is in the dual $G^F$-conjugacy class of $$L^*_s:=\cent{G^*}{\zo {L^*\cap\cent{G^*}s}}$$  As $\zo {L^*}\subseteq \zo{L^*\cap\cent{G^*}s}$, one has $L^*_s\subseteq L^*$, which allows to assume $L_s\subseteq L$. By definition $\cent{L_s^*}s$ is $\cent{\cent{G^*}s}{\zo {L^*\cap\cent{G^*}s}}=\cent{L^*}s$. The "trivial case" we described in (a) shows that 
$\Lu {L_s}L\circ \Psi_{L_s,s}=\Psi_{L,s}$.
Given $L$ and $M$ with $L^*\subseteq M^*$, clearly $L^*_s\subseteq M^*_s$ so we may assume that $L_s\subseteq M_s$. Now it is sufficient to prove our  formula between $L_s$ and $M_s$. Indeed,  assuming that $\Lu {L_s}{M_s}\circ \Psi_{L_s,s}=\Psi_{M_s,s}\circ \Lu{\cent{L^*_s}s}{\cent {M^*_s}s}$, we have by (a) and (b) $$\eqalign{\Lu {L}M\circ\Psi_{L,s} &=
\Lu  {L}M\circ\Lu{L_s}L\circ \Psi_{L_s,s}=\Lu{M_s}M\circ \Lu{L_s}{M_s}\circ \Psi_{L_s,s}\cr&=\Lu{M_s}M\circ \Psi_{M_s,s}\circ \Lu{\cent{L^*_s}s}{\cent {M^*_s}s}=\Psi_{M,s}\circ \Lu{\cent{L^*}s}{\cent {M^*}s}\cr}.$$

(d) So now we assume that $L=L_s\subseteq M=M_s$ with notations of (c) and $L\neq M$. Thanks to (b) we assume that $\cent{L^*}s$ is a maximal $F$-stable Levi subgroup in $\cent{M^*}s$. 

 We claim that $$||\Lu {L}M (\psi_{L,s}(\la))||=||\Lu {\cent{L^*}s}{\cent{M^*}s}\la||.\leqno{(5.3.3)}$$

Denote $\zeta=\Psi_{L,s}(\la)$. The square norm of $\Lu {L}M \zeta$ is $\scal{\slu {L}M( \Lu {L}M  \zeta)}{\zeta}{L^F}$ and may be computed by Mackey formula. It is a sum, indexed on a set of double classes mod $L^F$ in $M^F$, of terms  $\scal{\Lu{L\cap \lexp g L}{L}(\slu{\lexp g L}{L\cap\lexp g L}(\lexp g\zeta))}\zeta{L^F}$. Here $g$ is such that $L$ and $ \lexp g L$ have  a common  maximal torus so that $L\cap \lexp g L$ is a Levi subgroup of $\lexp g L$ and of $L$. When $g_1$ runs in $L^FgL^F$, $L\cap\lexp{g_1}L$ runs in a complete $L^F$-conjugacy class of Levi subgroups of $L$. The scalar product above is 
$$n(g):=\scal{\chi(g)}\zeta{L^F} \quad{\rm with}\; \chi(g):=\Lu{K(g)}{L}(\slu{K( g )}{\lexp g L}(\lexp g\zeta))\;{\rm and}\;  K(g)=L\cap \lexp g L\leqno{(5.3.4)}$$
 As $\zeta\in\ser{L^F}s$, $n(g)$ is non zero only if  for any  couple $(T,\theta)$ with $T\subseteq K(g )$, $\theta\in\II{T^F}$ and 
 
 \noindent $\scal {\slu {K(g)\cap \lexp g L}{\lexp g L}  (\lexp g \zeta)}{\Lu T{K(g)}\theta}{K(g)^F}\neq 0$ then the $L^F$-conjugacy class of $(T,\theta)$ is  associated to the $(L^*)^F$-conjugacy class of some $(S^*,s)$. That condition on $\lexp g\zeta$ is equivalent to $\scal {  \slu{T}{\lexp g L}(\lexp g\zeta)}{\theta}{T^F}\neq 0$. 
As $g\in M^F$ we have $\scal {  \slu{T^g}{L}\zeta}{\theta^g}{(T^g)^F}=\scal {  \slu{T}{\lexp g L}(\lexp g\zeta)}{\theta}{T^F}$. Hence $(T^g,\theta^g)$ and $(T,\theta)$ are $L^F$-conjugate. But $g$ is defined modulo $L^F$, so we may assume  
\smallskip 

\noindent (A) {\sl $\quad g\in \nor{M^F}{T,\theta} $ for some $(T,\theta)$ in the $L^F$-conjugacy class corresponding to the $(L^*)^F$-conjugacy class of $(T^*,s)$. }

and then $$||\Lu {L}M(\Psi_{L,s}(\la)||^2=\sum_{L^FgL^F\mid (\rm A)} n(g)\leqno{(5.3.5)}$$

Let $w=gT\in W(M,T)^F$ and $w^*$ the image of $w$ in the anti-isomorphism $W(M,T)\to W(M^*,T^*)$, $w^*\in W(\cent{M^*}s,T^*)^F$ and $w^*=g^*T^*$ for some $g^*\in\nor{\cent{M^*}s^F}{T^*}$.
Let $\Phi(L,T)\subset X(T)$ be the set of roots of $L$ with respect to $T$. With our choice of $g$ we have $\Phi(K(g),T)=\Phi(L,T)\cap \lexp w\Phi(L,T)$. We define a dual Levi subgroup $K^*(g)$ in $L^*$ by $\Phi({K^*(g)},T^*)=\Phi(L^*,T^*)\cap \lexp {w^*}\Phi(L^*,T^*)$. As $w^*\in W(\cent{M^*}s,T^*)$, we have $\Phi(\cent{K^*(g)}s,T^*)=\Phi(\cent{L^*}s,T^*)\cap \lexp {w^*}\Phi(\cent{L^*}s,T^*)$, that is 
$$K^*(g)=L^*\cap \lexp{g^*}L^*,\quad \cent{K^*(g)}s=\cent{L^*}s\cap\lexp{g^*}{\cent{L^*}s}\leqno{(5.3.6)}$$
 To $K(g)$, $g$ as in (A), we have associate a Levi subgroup $\cent{L^*}s\cap \lexp {g^*}{\cent{L^*}s}$ of $\cent{M^*}s$, hence one of the terms of the Mackey formula for $\slu{\cent{L^*}s}{\cent{M^*}s}\circ \Lu{\cent {L^*}s}{\cent{M^*}s}$. The contribution of that term to the square product $||\Lu {\cent{L^*}s}{\cent{M^*}s}\la||^2$ is 
 $$n(g^*):=\scal{\mu(g^*)}\la{\cent{L^*}s^F}\quad  {\rm for}\;\mu(g^*):=\Lu{\cent{K^*(g)}s}{\cent{L^*}s}(\slu{\cent{K^*( g )}s}{\lexp {g^*} {\cent{L^*}s}}(\lexp {g^*}\la))\leqno{(5.3.7)}$$
with $K^*(g)$ given by (5.3.6).
We have $$||\Lu{\cent{L^*}s}{\cent{M^*}s}\la||^2=\sum_{\cent{L^*}s^Fg^*\cent{M^*}s^F\subseteq M^{*\,F}}n(g^*)\leqno{(5.3.8)}$$

 Assume $h$ satisfy condition (A) with a torus $S$ as $g$ with $T$ and define $(K^*(h),h^*)$ from $(h,S)$ as  $(K^*(g),g^*)$ from $(g,T)$. If $\cent{K^*(h)}s$ is $\cent{M^*}s^F$-conjugate to $\cent{K^*(g)}s$, then 
 
 \noindent 
$\cent{L^*}s^Fh^*\cent{L^*}s^F=\cent{L^*}s^Fg^*\cent{L^*}s^F$ therefore $(L^*)^Fh^*(L^*)^F=(L^*)^Fg^*(L^*)^F$ so that $K^*(g)$ and $K^*(h)$ are $L^{*\,F}$-conjugate. Then $K(h)$ and $K(g)$ are $L^F$-conjugate and $L^FhL^F=L^FgL^F$. 

To any of the double class $L^FgL^F$ $(g\in M^F$) we have to consider in (5.3.5) we  associate a double class $\cent{L^*}s g^*\cent{L^*}s$ ($g^*\in \cent{M^*}s^F$) and that application is injective. Given a double class $\cent{L^*}s g^*\cent{L^*}s$ where $g^*\in\cent{M^*}s^F$ and $\cent{L^*}s\cap \lexp{g^*}{\cent{L^*}s}$ contains a maximal torus $T^*$ we may choice $g^*\in\nor{M^*}{T^*}^F$, so that $g^*T^*=w^*\in W(\cent{L^*}s,T^*)^F\subseteq W(L^*,T^*)^F$. For $w^*\mapsto w$ and $gT=w$ ($g\in M^F$), $w$ in the isomorphic group $W(L,T)$, $\cent{L^*}s g^*\cent{L^*}s$ is associated to $L^FgL^F$.

So we have defined a bijection $L^FgL^F\mapsto \cent{L^*}sg^*\cent{L^*}s$ between terms of the two decomposition formulas (5.3.5) and (5.3.8). Now (5.3.3) will follow from the equality $n(g)=n(g^*)$ for any such couple $(g,g^*)$.

  On $M$-side one has $n(g)=\scal{\slu{K(g)}{\lexp g L}(\lexp g\zeta)}{\slu{K(g)}{L}\zeta}{K(g)^F}$. As $g$ induces an isomorphism $(L,F)\to(\lexp g L,F)$, $g^*$ induces a dual isomorphism  $(L^*,F)\to\lexp {g^*}{L^*}$ that fixes $s$ so that we may assume  $\lexp g\zeta=\Psi_{\lexp g L,s}(\lexp {g^*}\la)$. 

If $g\in\nor{M}L^F$, then $K(g)=L$ and $K^*(g)=L^*$. In that case $n_g\neq 0$ if and only if $\lexp g\zeta=\zeta$, if ond only if $\lexp {g^*}\la=\la$ and then $n(g)=1=\scal{\lexp {g^*}\la}\la{\cent{L^*}s^F}=1=n(g^*)$.

If $L\cap\lexp gL\neq L$, induction hypothesis applies : we have, by (5.3.2), $$\slu{K(g)\cap\lexp gL}{\lexp g L}(\lexp g\zeta)=\Psi_{K(g),s}(\slu{\cent{K^*(g)}s}{\cent{\lexp {g^*}{L^*}}s}(\lexp {g^*}\la)),\quad \slu{K(g)}{L}\zeta= \Psi_{K(g),s}(\slu{\cent{K^*(g)}s}{\cent{L^*}s}(\la))$$
hence $n(g)=\scal {\slu{\cent{K^*(g)}s}{\cent{\lexp {g^*}{L^*}}s}(\lexp {g^*}\la)}{\slu{\cent{K^*(g)}s)}{\cent{L^*}s}(\la)}{\cent{K^*(g)}s^F}=n(g^*)$. 

(5.3.3) is proved.

(e) To conclude we need a

\medskip\noindent{\bf 5.3.9. Lemma. }{\sl Let $(G,F)$ be any algebraic reductive group defined on $\F q$, all components of which have classical type. Let $L$ be a maximal proper $F$-stable Levi subgroup of $G$ and $\la\in\ser{L^F}1$. Then $ \Lu {L\subset P}G \la$ is uniquely defined as the element of  minimal norm among the set of
$\mu\in\Z\ser{G^F}1$ such that $\piu^G(\mu)=\Lu LG(\piu^L(\la))$, so is independant of $P$.}

\noindent{\it (Proof of the Lemma after the end of the proof of the formula (5.3.1))}
 
 We have commutation formulas (1.3.1.4) and (1.3.2.4) : 
 $$\Lu LM\circ\piu^L=\piu^M\circ\Lu LM,\quad \Psi_{M,s}\circ \piu^{\cent{M^*}s}\circ\Lu{\cent{L^*}s}{\cent{M^*}s}=\piu^M\circ\Lu LM\circ \Psi_{L,s}$$ 

 Apply the Lemma between $\cent{L^*}s$ and $\cent{M^*}s$, and isometry $\Psi_{M,s}$, knowing that Jordan decomposition commute with projection on the space of uniform function ((J.1) in Proposition~1.3.2) : $\Psi_{M,s}(\Lu {\cent{L^*}s}{\cent{M^*}s}\la)$ is uniquely defined as the element of minimal norm among the set of $\Lambda\in\Z\ser{M^F}s$ such that $ \piu^M(\Lambda)=\piu^M(\Psi_{M,s}(\Lu {\cent{L^*}s}{\cent{M^*}s}\la))$, its norm being precisely $ ||\Lu {\cent{L^*}s}{\cent{M^*}s}\la||$. 
 
 We have $\piu^M(\Lu LM(\Psi_{L,s}(\la))=\Psi_{M,s}( \piu^{\cent{M^*}s}(\Lu{\cent{L^*}s}{\cent{M^*}s}\la))=\piu^M(\Psi_{M,s}(\Lu {\cent{L^*}s}{\cent{M^*}s}\la))$. By (5.3.3), (5.3.1) is proved.  \bull
 
 \medskip\noindent{\it Proof of the Lemma 5.3.9. }  We use the fact that in classical type the application $(\la\mapsto \piu^G(\la))$ from $\ser{G^F}1$ to $ K\II {G^F}$,  is one-to-one.
 
  Let $[G,G]=G_1.G_2\dots G_k$ be a \dec\ in a central product of rationally \irr\ components.  The Levi subgroup $L$  writes $L=\zo{G}.(L\cap G_1)\dots (L\cap G_k)$ where $L\cap G_i$ is a Levi subgroup of $G_i$. As $L$ is maximal there is only one $i$, say $i=1$, such that $L\cap G_i\neq G_i$. As the set of unipotent \irr\ characters is indifferent to central quotients, we may write $\la=\la_1\otimes \la_2\otimes \dots \otimes \la_k$, where $\la_1\in\ser{(L\cap G_1)^F}1$, $\la_i\in\ser{{G_i}^F}1$ for $i\neq 1$ and $\mu\in \Z[\ser{{G_1}^F}1\otimes \dots\otimes\ser{{G_k}^F}1]$, so that $\Lu LG\la=\Lu{L\cap G_1}{G_1}\la_1\otimes\la_2\otimes  \dots \otimes \la_k$. Hence $||\Lu LG\la||^2=||\Lu{L\cap G_1}{G_1}\la_1||^2$ and $\piu^G(\Lu LG\la)=\piu^{G_1}(\Lu {L\cap G_1}{G_1} \la_1)\otimes \piu^{G_2}(\la_2)\otimes\dots \otimes \piu^{G_k}(\la_k)$. 

(a) If $G_1$ has type $\AA$, every central unipotent function on ${G_1}^F$ is a uniform function. The condition $\piu^{G}(\mu)=\Lu LG(\piu^L(\la))$ gives $\piu^G(\mu)=\Lu {L\cap G_1}{G_1} \la_1\otimes\dots \otimes \piu^{G_k}(\la_k)$ ($L\cap G_1$ maximal  in $G_1$ or not). From the condition $||\mu||^2\leq ||\Lu LG\la||^2$, it follows that $\mu=\Lu {L\cap G_1}{G_1} \la_1\otimes\dots \otimes \mu_k$ with $||\mu_i||^2=1$ for $i>1$ and $\piu^{G_i}(\mu_i)=\piu^{G_i}(\la_i)$, hence $\mu_i=\la_i$ for $i>1$. 

(b) If $G_1$ has classical type $\XX\in\{\BB,\CC,\DD\}$,  then there is some integer $d$ such that $L\cap G_{1}$ has type $\XX_{k-d}$,  (the rational types of $G_1$ and $(L\cap G_1)$ may differ in types $\DD$, $\lexp 2\DD$), $\zo{L\cap G_{1}}$ has polynomial degree $(x^d-1)$ (if $d$ is odd) or $(x^{d/2}+1)$ (if $d$ is even).  Then $\Lu{L\cap G_{1}}{G_1}\la_{1}$ is given by one of Asai's formulas, a $d$-hook (or $(d/2)$-cohook) formula, see  [9], (3,5), (3.9). These formulas show that $\Lu{L\cap G_{1}}{G_{1}}\la_{1}$ is a  sum  $\sum_{1\leq j\leq r} \pm \chi_j$, all $\chi_j$ in distinct families, families defined from $(W(G_{1}),F)$. But $\chi_j$ is uniquely defined in $\ser{{G_1}^F}1$ by $\piu^{G_1}(\chi_j)$. Thus $\Lu{L\cap G_{1}}{G_{1}}\la_{1}$ is unique in $\Z\ser{{G_{1}}^F}1$ with uniform projection $\sum_j \pm\piu^{G_{1}}(\chi_j)$ and minimal square norm $r=||\Lu {L\cap G_1}{G_1}\la_1||^2=||\Lu LG\la||^2$. We conclude as in (a).
\bull

\bigskip\noindent{\bf 5.4. More on Generalized $d$-Harish-Chandra theory }

We use here results of Bonnaf\'e  [3] on type $\AA$ with a Frobenius endomorphism, completed by Cabanes [12] in twisted type. 

\medskip\noindent{\bf 5.4.1. Facts on wreath products }

We need properties of the Weyl groups of centralizers of semi-simple elements in groups of type $\AA$.

In that section $\Sigma$ is a finite set, a finite group $B$ acts by permutations on $\Sigma$, $W^0$ is a direct product  indexed on $\Sigma$ : $W^0=\times_{\beta\in \Sigma}W_\beta$. If $\beta$ belongs to the orbit $B\al\in \Sigma/B$, then $W_\beta$ is isomorphic to $W_\alpha$. We identify the groups $W_\beta$ ($\beta\in B\al$) to $W_\al$ so that $B$ acts on $\times_{\beta\in B\al}W_\beta$ as on the set of functions from $B\al$ to $W_\al$, with its natural group structure. When $\Omega$ is contained in an orbit under $B$ on $\Sigma$ we denote $W^\Omega=\times_{\beta\in \Omega}W(\beta)$, a subgroup of $W^0$, with a projection morphism $W^0\to W^\Omega$, $(w\mapsto w^\Omega)$.

The fundamental property of the action of $B$ on $W^0$ is that $b\in B$ and  $b\beta =\beta$ imply $\beta w=w$ for any $w\in W_\beta$. We say that {\it $W:=W^0\rtimes B$ is a wreath product}. Note that in that sense a direct product $W^0\times B$ is a wreath product! With standard notations, if $B$ acts faithfully and is transitive on  $\Sigma$,  then $W^0\rtimes B\cong W_x \wr B$. The proofs of all properties we recall here may be reduced to the transitive case. To easy references we use notations similar to notations  in [12] or [3].

As $B$ acts on $\II{W^0}\cong \times_{\beta\in \Sigma}\II{W_\beta}$, we identify $\II{W^0}$ with $\times_{\beta\in \Sigma}\II{W_\beta}$ and we use notation $\chi^\Omega\in\II{W^\Omega}$ if $\Omega$ is as above. For any subgroup $C$ of $B$ we have a natural one-to-one map
$$\II{W^{0}}^C\to \II{W^{0\,C}},\quad \chi=\chi_1\mapsto \chi_C\leqno{(5.4.1.1)}$$
and simplify $\chi_{\gen b}$ in $\chi_b$.
Let  $\chi\in\II{W^0}$. Once the complement $B$ is fixed, there is a special  {\it canonical extension} of $\chi$ to $W^0\rtimes B_\chi$, as well as  to any subgroup $W^0\rtimes C$ if $C\subseteq B_\chi$, design it $\chi\rtimes C$, as in [12]. It is caracterized as the unique \irr\ extension of $\chi$ to $W^0\rtimes C$ with natural integer value on $c$ for any $c\in C$. The restriction of $\chi\rtimes C$ to $W^0c$, we denoted $\chi\rtimes c$, may be defined as follows :  for  any $\omega=\gen c\al\in \Sigma/\gen c$, an orbit under $\gen c$ in $\Sigma$, any $w^\omega=(w(\beta))_{\beta\in \omega}\in W^\omega$, put
$$(\chi^\omega\rtimes c)(w^\omega c)=\chi^\omega_c(w(c^{|\omega|-1}\al)\dots w(c^{|\omega|-j}\al)\dots  w(\al))\leqno{(5.4.1.2)}$$
If $(\beta\mapsto w(\beta))$ is constant on $\omega$, then $(\chi^\omega\rtimes c)(w^\omega c)=\chi^\omega_c(w(\al)^{|\omega|})$.
Then $\chi\rtimes C$ is defined by
$$ (\chi\rtimes C)(wc)=(\chi\rtimes \gen c)(wc)=(\chi\rtimes  c)(wc)=[\otimes_{\omega\in \Sigma/\gen c}(\chi^\omega\rtimes c)](w c)\leqno{(5.4.1.3)}$$  
Thus by (5.4.1.2) $(\chi\rtimes C)(c)=(\chi\rtimes c)(c)=\chi_c(1)$. 

Let $\hat \chi$ be in $\II{W\mid \chi}$, $\chi$ is defined by $\hat\chi$ mod $B$-action, $\hat \chi$ is one of the distinct 
$$\Gamma^W(\chi*\theta):=\Ind{W^0\rtimes B_\chi}{W}(\theta\otimes (\chi\rtimes B_\chi)) \quad \chi\in\II{W^0}/B,\;\theta\in\II{B_\chi}\leqno{(5.4.1.4)}$$

\medskip\noindent{\bf 5.4.1.5. }{\sl Let be $W=W^0\rtimes B$, $\chi\in\II{W^0}$ as above.

(i) If $C\subseteq B_\chi$, then $\Res{W^0\rtimes B_\xi}{W^0\rtimes C}(\chi\rtimes B_\chi)=\chi\rtimes C$. For any $b\in B$, $(\chi\rtimes B_\chi)^b=\chi^b\rtimes {B_\chi}^b$. 

(ii) If $B_\chi\subseteq C\subseteq D\subseteq B$ and $\theta\in\II {B_\chi}$, then 

\noindent$\Ind{W^0\rtimes C}{W^0\rtimes D}(\Gamma^{W^0\rtimes C}(\chi*\theta))=\Gamma^{W^0\rtimes D}(\chi*\theta)$. 
 
(iii) Assume $B$ abelian. Let $C$ be a subgroup of $B$, $\psi\in B^\wedge$, $\theta\in(B_\chi)^\wedge$, $\rho=\Res {B_\chi}{C_\chi}\theta$. One has

\noindent  $\psi\otimes\Gamma^W(\chi*\theta)=\Gamma^W(\chi*(\Res B{B_\chi}\psi\otimes \theta))$ and $\Res W {W^0\rtimes C}(\Gamma^W(\chi*\theta))=\sum_{b\in B/C}\Gamma^{W^0\rtimes C}(\chi^b*\rho)$.

(iv) Assume $B$ abelian. Let $V^0\rtimes B_V$ be a sub wreath product of $W^0\rtimes B$, hence $V\cap W^0=V^0$, $B_V\subseteq B$. One has, for any $\theta\in(B_\chi)^\wedge$, $\zeta\in\II{V_0}$, $\psi\in( B_{V,\zeta})^\wedge$,

\noindent $\scal{\Res{W^0\rtimes B}{V^0\rtimes B_V}(\Gamma^{W^0\rtimes B}(\chi*\theta))}{\Gamma^{V^0\rtimes B_V}(\zeta*\psi)}{V^0\rtimes B_V}$

\hfill$=\sum_{b\in B/B_\chi.B_{V,\zeta}}\scal{\Res{W^0}{V^0}{\chi^b}}\zeta {V^0}.\scal{\Res {B_\chi}{B_\chi\cap B_{V,\zeta}}\theta}{\Res{B_{V,\zeta}}{B_\xi\cap B_{V,\zeta}}\psi}{B_\chi\cap B_{V,\zeta}}$}

\noindent{\it On proof. } (i) and (ii) are clear from definition and caracterization of canonical extensions, (iii) is a special case of (iv) with $V^0=W^0$ and is easily deduced from definition (5.4.1.4) and Mackey's formula. Using (iii), (ii) and transitivity of restriction, one sees that (iv) is true if it is true when  $\theta=1=\psi$.  We may extend $\chi\mapsto\chi\rtimes b$ by $\Z$-linearity on $\Z[\II{W^0}^b]$ so that (5.4.1.2) gives $\Res{W^0b}{V^0b}(\chi\rtimes b)= ( \Res{W^0}{V^0}\chi)\rtimes b$. The equality in (iv) follows by Mackey's formula, $B/B_\chi.B_{V,\zeta}$ represents the set of  doubles classes for $W^0\rtimes B_\chi$ and $V^0\rtimes (B_{V,\zeta})$ in $W^0 \rtimes B$.  \bull

\medskip\noindent{\bf  5.4.1.6. } {\sl Assume $B$ abelian. By Mellin transform one defines for any $b\in B$ and $\chi\in \II {W^0}^b$ 
$$\hat \Gamma^W(\chi*b)=\sum_{\theta\in (B_\chi)^\wedge}\theta(b^{-1})\Gamma^W(\chi*\theta)$$
Then $\hat\Gamma(\chi*b)$ extends $|B_\chi|.\Res{W}{W^0b}(\Gamma(\chi*1))$ by $0$ outside $W^0b$ and $\scal{\hat\Gamma(\chi*b)}{\hat\Gamma(\chi*b)}{W^0b}=|B_\chi|.|B|$. The family $\{\chi\rtimes b\}_{\chi\in{\II {W^0}^b}}$ is an orthonormal basis of ${\rm Cent}(W^0b)$.}

\medskip\noindent{\bf 5.4.2. G.$d$-HC in type $\AA$ }

Bonnaf\'e introduced the following conjecture, assuming $p$ is good for $G$,  and $\Gamma_u^M$, $\Gamma_{u_L}^L$,  are Gelfand-Graev characters of $M$ and $L$ respectively associated to regular unipotent elements $u_M$ (of $M$) and $u_L$ (of $L$) [3]  14.E.

\medskip\noindent{\bf 5.4.2.0. Conjecture $\frak G$. }{\sl If $L$ is an $F$-stable Levi complement of a parabolic subgroup $P$ of an $F$-stable  Levi subgroup $M$ of $G$,  with a coherent choice of regular unipotent elements $u_M$, $u_L$  one has 
$\slu {L\subset P}M\Gamma^M_{u_M}=\Gamma^L_{u_L}$.}

\noindent (recall that we have included the sign $\epsilon_L\epsilon_M$ in our definition of $\Lu LM$). 

Conjecture ${\frak G}$ is verified if $p$ is good and [ $\z G$ is connected or $F$ is a Frobenius $\F q$-endomorphism and  $q$ large enough ].

\medskip\noindent{\bf 5.4.2. Proposition. }{\sl Assume all components of $G$ has type $\AA$ and Conjecture ${\frak G}$ holds in $(G,F)$.  Then Generalized $d$-Harish-Chandra theory holds in any series $\ser{G^F}s$, any $d$.}

\preuve In part (A) we describe G.$d$-HC series in type $\AA$ when the center of $G$ is connected. In part (B), thanks to the descripion of $\ser{G^F}s$ in SL$_n(q)$ and SU$_n(q)$ due to Bonnaf\'e and Cabanes, we describe G.$d$-HC series in the general case by their image in the set of \irr\ characters of the Weyl group of $\cent{G^*}s$, see (B.3.2).

  (A) Connected center.

(A.1) The origin of G.$d$-HC for $G:={\rm GL}_n$, $G^F={\rm GL}_n(\epsilon q)$ is Murnaghan-Nakayam formula in symmetric groups. We recall briefly these facts. The Weyl group $W$ with respect to a maximal diagonal torus in $G$ is the symmetric group ${\frak S}_n$ on which $F$ acts trivially. One has a one-to-one map, defined using Deligne-Lusztig induction, [20] 15.8
$$\II W\to \ser{G^F}1,\quad \mu\mapsto R_\mu\leqno{(5.4.2.1)}$$
 and the map $\mu\mapsto R_\mu$ transform induction (on $W$-side) in Lusztig  induction (on $(G,F)$-side).
 
  Let $\delta$ be the order of $\epsilon q$ modulo $\ell$, so that $\delta$ differs from $d$ in twisted type : if $\epsilon=-1$ and $\ell\neq 2$, then $d$ is $\delta, \delta/2, 2\delta$ resp. when $\delta\equiv 0, 2, (1 $ or $3$)$ \pmod 4$ resp.). If $\delta=1$, then there is only one G.$d$-HC series in $\ser{G^F}1$. Assume $\delta\neq 1$. 
  
  Let $S$ be a subgroup of ${\frak S}_n$ generated by m cycles of length $\delta$ with disjoints supports. The product $v$ of these $m$ cycles is a parameter with respect to a diagonal maximal torus of $G$ for a maximal torus $T_v$ in $G$, the centralizer of $v$ in ${\frak S}_n$ has the form ${\frak S}_{n-m\delta}\times (S\rtimes{\frak S}_m)$ and ${\frak S}_{n-m\delta}$ is the Weyl  group  of the $d$-split Levi subgroup $L_v:=\cent G{(T_v)_{\Phi_d}}$  with respect to $T_v$.  Any Levi subgroup of $G$ occurring in a unipotent $d$-cuspidal datum in $(G,F)$ is so obtained. A $d$-split Levi subgroup of $G$ that contains $L_v$ up to $G^F$-conjugacy is defined by a partition of $m$, say $m=\sum_jm_j$, and has Weyl group ${\frak S}_{n-m\delta}\times(\times_j[{\frak S}_{m_j}]^\delta)$ where $F$ acts regularly on each direct product $[{\frak S}_{m_j}]^\delta$. Any $d$-split Levi subgroup of $G$ contains some $L_v$ (with $m$ large enough).
  
  The Murnaghan-Nakayama formula allows to compute $\Res {W}{ W(L_v) v}$. When we consider that restriction to the class $W(L_v)v$ {\it  as a central function on $W(L_v)$, we write it $\Res W {W(L_v)*v}$}.  With these notations all unipotent central functions being uniform functions, $\Res W {W(L_v)*v}\mu$ gives $\slu {L_v}G R_\mu$ (see [20], 15.7, 15.8). As well we may define $\Ind  {W(L_v)*v}W$ so that $\Ind  {W(L_v)*v}W\nu$ gives $\Lu {L_v}G R_\nu$ when $\nu\in\II{W(L_v)}$.
    
  If $\mu=[\la]\in\II {{\frak S}_n}$ is defined by a partition $\la$ of $n$ with no hook of length $\delta$ --- we say that $\la$ is a $\delta$-core --- then, by Murnaghan-Nakayama formula,  $\mu$ has null value on any element whose cycle decomposition has a cycle whose length belongs to $\delta\N$. The  type of any maximal torus in a proper $d$-split Levi subgroup of G contains such a cycle, that's why the corresponding $R_\mu$ is $d$-cuspidal in our first sense (2.4). The set of partitions of $n$ is a disjoint union of ``$\delta$-HC series" defined for each $\delta$-core $\kappa$ of size $(n-m\delta)$ for some $m$, hence $\II {{\frak S}_n}$ is a disjoint union of ``$\delta$-HC series" defined for each ``$\delta$-cuspidal" $\mu=[\kappa]\in\II{{\frak S}_{n-m\delta}}$ ($\kappa$ a $\delta$-core) for some $m$. With $v$ of cycle type $\delta^m$ as above is defined a series in $\II{ {\frak S}_n}$ :
$${\cal S}^{{\frak S}_n}(m,\delta,\mu)=\{\la\in\II{{\frak S}_n}\mid \Res{{\frak S}_n}{{\frak S}_{n-m\delta}*v}\la = c\mu, c\neq 0\}\leqno{(5.4.2.2)}$$
these series  form a partition of $\II{ {\frak S}_n}$ when $\mu=[\kappa]$ and $\kappa$ runs on the set of $\delta$-cores such that $n-||\kappa||\in\delta\N$. In (5.4.2.2) $c\in\Z$.

On $(G,F)$-side, $(L_v,R_{\mu})$ is a unipotent $d$-cuspidal datum. The relation $\scal{\la}{\Ind{{\frak S}_{n-m}\delta*v}{{\frak S}_n}\mu}{{\frak S}_{n}}\neq 0$ is equivalent to $\scal{R_\la}{\Lu {L_v}{G}R_\mu}{G^F}\neq 0$, furthermore the two scalar products are equal up to a sign $\epsilon_\mu\epsilon_\la$ and this sign is 1 in non twisted type. Thus one has a partition of $\ser{G^F}1$ in G.$d$-HC series $\ser{G^F}{(L_v,R_\mu)}$ , (5.4.2.2) gives  $$\ser{G^F}{(L_v,R_\mu)}=\{R_\la\mid \la\in{\cal S}^{{\frak S}_n}(m,\delta,\mu)\}\leqno{(5.4.2.3)}$$

(A.2) Unipotent G.$d$-HC series.

The generalization of G.$d$/$(\epsilon q)$-HC to any direct product of symmetric groups is immediate. 

From the ``invariance" of the set of unipotent \irr\ characters  in isotypic morphisms  and properties of direct products, we may assume that $G$ is rationally \irr. It may happen that the $ \F q$-endomorphism acts non trivially on the Dynkin diagram of $G$ or/and of $L_v$. If it is the case we have to consider $F$-conjugacy classes in $W$, central functions on $W F$, to recover Lusztig induction. 

So we consider the  scalar descent :  the Weyl group $W$ with respect to a maximal diagonal torus is a wreath product $W={\frak S}_n\wr \gen F=W^0\rtimes \gen \phi$.  Here $W^0=\times_{\beta\in \Sigma}W_\beta$, $F$ acts on $\Sigma$ as a transitive permutation, and $\phi$ is defined by $F$ (see 5.4.1 with $B=\gen \phi$). 

One has a one-to-one map  $\II{W^0}^\phi\to \ser{G^F}1$, $(\chi\mapsto \chi\rtimes \phi\mapsto R_\chi)$ that generalizes (5.4.2.1).

 Let $n$ such that $W_\beta\cong {\frak S}_n\cong W^{0\,\phi}$, let $r$ be the order of $\phi$, $r=|\Sigma|$. Let  $m$ be  a positive integer such that  $m\delta<n$, the set of elements $(v_\beta)_{\beta\in \Sigma}$ in $W^0$ such that $v_\beta$ has cycle type $\delta^m$ for any $\beta\in \Sigma$ is an $F$-stable $W^0$-conjugacy class whose intersection with $W^{0\, F}$ is a conjugacy class of $W^0$. Take $v$ inside with $F(v)=v$. We have $\cent{W^0}v=V^0.S^0$ where $S^0$ is generated by the cycles of $v$, $V^0\cong\times_\beta{\frak S}_{n-m\delta}$,  $S^0$ and $V^0$ are $F$-stable subgroups of $W^0$, $S^0\rtimes\gen\phi$ and $V^0\rtimes\gen \phi$ are sub-wreath products of $W$ and $V^0$ is the Weyl group of the $d$-split Levi subgroup $L_v$ of $G$, $\cent G{(T_v)_{\phi_d}}$ with respect to $T_v$. The action of $F$ on $W(L_v,T_v)$ is given by $v\phi$, acting as $\phi$. 
 
 We want to compute $\Res {W^0\phi} {V^0*v\phi}$  (the transpose map $\Ind{V^0v\phi}{W^0\phi}$ is called  {\it induction tordue} in [1], 3.1). With notations introduced in 5.4.1, (5.4.1.4), if $\chi\in \II{W^0}^\phi$, $\theta\in \gen \phi^\wedge$, $h\in V^0$,  then $$\Gamma^W(\chi*\theta)(hv\phi)=\theta(\phi)(\chi\rtimes \phi)(hv\phi)$$  But, as $\phi v=v\phi$ and $[v,V^0]=1$, $(\chi\rtimes\phi)(hv\phi)=\chi_\phi(h'v^r)$
 where  $h'\in V^{0\,\phi}$ is given by (5.4.1.2), with $(\phi,h)$ instead of $(c,w^\omega)$. Let $\delta_1=\delta/(\delta,r)$, where $(\delta,r)$ denote the GCD of $\delta$ and $r$, so that $v^r$ has cycle type $\delta_1^{m(\delta,r)}$.   
  Assume that $\chi_\phi$ belongs to ${\cal S}^{{\frak S}_n}(m_1,\delta_1,\mu)$, with notations of (5.4.2.2), $\mu$ being defined by a  $\delta_1$-core. 
  Let be  $\Res{W^{0}}{V^{0}*v^r}\chi=\sum_{\zeta\in\II{V^0}}t(\zeta)\zeta$.
 Then $t(\zeta)\neq 0$ implies $\zeta\in\II{V^0}^{v\phi}$ and  $$\Res{W^{0\,\phi}}{V^{0\,\phi}*(v^\phi) ^r}\chi_\phi=\sum_{\zeta\in\II{V^0}^{ v\phi}}t(\zeta)\zeta_\phi$$ 
 
 By (5.4.1.2)  $(\chi\rtimes \phi)(hv\phi)=\sum_{\zeta\in\II{V^0}^\phi}c(\zeta)\zeta_\phi(h')$ and $\zeta_\phi(h')=(\zeta\rtimes\phi)(h\phi)=(\zeta\rtimes v\phi)(hv\phi)$. Extending linearly the isometry ($\chi\mapsto\chi\rtimes \phi$) to spaces of central functions Cent$(V^{0\, \phi})\to {\rm Cent}(V^0\phi)$, we obtain  
 $$(\chi\rtimes\phi)(hv\phi)= ((\Res {W^0}{V^0*v^r}\chi)\rtimes v\phi)(hv\phi),\quad \Res{W}{V^0v\phi}(\Gamma^W(\chi*\theta))=\theta(\phi)((\Res{W^0}{V^0*v^r}\chi)\rtimes v\phi)$$
  $$\Res{W^0\phi}{V^0v\phi}(\chi\rtimes \phi)=(\Res{W^0}{V^0*v^r}\chi)\rtimes v\phi \leqno{(5.4.2.4)}$$
The partition in series like (5.4.2.2) applies to the basis $\{\chi\rtimes\phi\}_{\chi\in\II{W^0}^\phi}$ of the space of central functions on $W\phi$ (or $F$-central functions on $W$) (see 5.4.1.6). Assume $\zeta\in \II{V^0}^{v\phi}\cong \II{V^{0,v\phi}}$ is defined by a $\delta_1$-core $\kappa$, partition of $n-m_1\delta_1$, so that $\zeta_\phi=[\kappa]$. Then $\Res W {Vv\phi}(\Gamma^W(\chi*\theta))=\theta(\phi)\Res{W^0\phi}{V^0v\phi}(\chi\rtimes\phi)=c\theta(\phi)\zeta_\phi$ with $c\in\Z$ and $c\neq 0$ if and only if $\chi\in\II{W^0}^\phi$ and $\chi_\phi \in{\cal S}^{{\frak S}_n}(m_1,\delta_1,\zeta_\phi)$ by (5.4.2.4). So are  defined G.$d$-HC unipotent series in $G^F$ from G.$\delta$-HC series in $W^F$ : $$\ser{G^F}{(L_v,R_\zeta)}=\{R_\chi\mid \chi_\phi\in {\cal S}^{{\frak S}_n}(m_1,\delta_1,[\kappa])\},\quad \zeta_\phi=[\kappa]\leqno{(5.4.2.5)}$$
where $\zeta\in{\II {V^0}}^\phi$ and $\zeta_\phi$ is $\delta_1$-cuspidal in ${\frak S}_{n-m_1\delta_1}$.

(A.3) In any Lusztig series.

Assume $s$ is a semi-simple element of $G^{* \,F}$ and $\cent{G^*}s$ is connected.  As an example, think to a central product $G$ of ${\rm GL}_{n_j}$ for some $(n_j)_j$. Then $\cent{G^*}s$ is a central product of connected groups $\prod_{\beta\in \Sigma}C_\beta$ of type $\AA$, the product is indexed on the set $\Sigma=\cup_j\Sigma_j$ of eigenvalues of $s$ in the standard representation. One has $W(s)=\times_{\beta\in \Sigma}W_\beta$. The preceding description applies to $\cent{G^*}s$ by direct product on $j\in J$ and on each $\gen F$-orbit on $\Sigma$. Going down from $\times_{\omega\in \Sigma/\gen F}(\times _{\beta\in\omega}(W_\beta))\rtimes\gen \phi)$ to $W\rtimes \gen\phi$ is easy. That's why  G.$d$-HC theory holds in $\ser{\cent{G^*}s^F}1$. 

Assuming $\z G$  connected we  have commutation formula (J3) : $\Lu LG\circ \Psi_{L,s} =\Psi_{G,s}\circ \Lu{\cent{L^*}s}{\cent{G^*}s}$ (when $s\in L^*$). It gives G.$d$-HC theory in $\ser{G^F}s$ : 

If $L$ is $d$-split in $(G,F)$ then $L^*$ is $d$-split in $(G^*,F)$ hence $\cent{L^*}s$ is $d$-split in $(\cent{G^*}s,F)$. Thanks to (J3), $\Psi_{L,s}(\la)$ is $d$-cuspidal in $\ser{L^F}s$ if and only if $\la$ is $d$-cuspidal in $\ser{\cent{L^*}s^F}1$ and $G^*$ is the only $d$-split Levi subgroup of $G^*$ that contains $\cent{G^*}s$. Any $d$-cuspidal unipotent datum $(L^*_s,\al)$ in $(\cent{G^*}s,F)$ defines a $G^F$-conjugacy class of $d$-cuspidal data $(L,\la)$ in series $s$ of $(G,F)$ as follows :  $L$ is in the dual class of the $d$-split Levi subgroup $\cent{G^*}{\z{L^*_s}_{\phi_d}}$ and $\la=\Psi_{L,s}(\al)$. By (J3) again one obtains a partition of $\ser{G^F}s$ in G.$d$-HC series $$\ser{G^F}{(L,\Psi_{L,s}(\la))}=\Psi_{G,s}[\ser{\cent{G^*}s^F}{(\cent{L^*}s,\la)}], \quad L^*=\cent{G^*}{\z{L^*_s}_{\phi_d}}\leqno{(5.4.2.6)}$$

(B) General case, $\z G$ connected or not.

From section 2 we may consider  a regular embedding $(G,F)\subseteq (H,F)$.  Let $t$ be a semi-simple element of $H^*$ that maps on $s$ by a dual map. We need the following elementary fact :

\smallskip\noindent (B.1) {\sl Let $M^*_t=\times_{\omega\in\Sigma/\gen \phi} M^*_\omega$ be a $d$-split Levi subgroup of $\cent{H^*}t=\times_\omega C_\omega$ ($\Sigma$ as in (A.3)). Define $M^*$ by $M^*=\cent{H^*}{\z {M^*_t}_{\phi_d}}$ and let $L^*$ be the image of $M^*$ in $G^*$. Let $\omega\in\Sigma/\gen\phi$ and assume  $M^*_\omega\neq C_\omega$. Then $\omega$ is stable under ${\rm A}_{L^*}(s)^F$.}
\smallskip\noindent{\it Proof of (B.1). } 

(a) Reduction.

$(G,F)$ is a central product of rationally \irr\ components $G_i$. Let $\nu\colon (G^*,F)\to(\times_iG^*_i,F)$ be a dual map of $(G_0,F)=(\times_iG_i,F)\to (G,F)$. From Proposition~1.1.4, (d)  there exist regular embeddings $G_i\to H_i$, $G\to H$, an  isotypic morphism $\times_iH_i\to H$ and by duality a commutative  diagram of isotypic morphisms
$$\diag { H^*}{\rho}{G^*}{\nu}{G^*_0:=\times_iG^*_i}{\rho_0}{H^*_0:=\times_iH^*_i}{\nu_1}{}$$
So are defined $(t_i)_i=\nu_1(t)$, $(s_i)_i=\rho_0((t_i)_i)=\nu(s)$. Using the four maps above we may identify the Weyl groups $W(t)$, $W((t_i)_i)$, $W^\circ(s)$, $W^\circ((s_i)_i)$, as well as the sets of eigenvalues $\Sigma(t)$, $\cup_i\Sigma(t_i)$ with $\phi$-action (see (A.3)). Furthermore $\nu$ induces morphisms $W(s)\to W((s_i)_i)$, ${\rm A}_{G^*}(s)\to {\rm A}_{G^*_0}((s_i)_i)$ commuting with actions on $\Sigma(t)$ (see Proposition~1.2.3). Given $M^*_t$, we may consider $M^*_{(t_i)_i}=\nu_1(M^*_t)$, and then  $M^*_0=\cent{H^*_0}{\zo{M^*_{(t_i)_i}}_{\phi_d}}=\nu_1(M^*)=\times_iM^*_{0,i}$ so that $\rho_0(M^*_0)=\nu(L^*)=\times_iL^*_{0,i}$, where $\rho(M^*_{0,i})=L^*_{0,i}$. The restriction of $\nu$ induces a $F$-morphism ${\rm A}_{L^*}(s)\to {\rm A}_{\nu(L^*)}((s_i)_i)$, the last one is a subgroup of $\times_i{\rm A}_{L^*_{0,i}}(s_i)$. Now it is clear that if (B.1)  is verified for $(G_i\to H_i, t_i, s_i,M^*_{t_i})$, any $i$, it is true for $(G\to H,t,s,M^*_t)$. So we may assume $G$ rationally \irr.

In a scalar descent (see (A.2)) we may replace $G$ by one of its component and $F$ by some convenient power. Then $H={\rm GL}_n$ for some $n$, so that $G$ contains ${\rm SL}_n$. The above argument shows that if (B.1) is true for ${\rm SL}_n$ it is true for $G$. So we assume $G={\rm SL}_n$. 

(b) $G={\rm SL}_n\subset {\rm GL}_n$. 

By definition $M^*$ and $L^*$ are $d$-split. Thanks to (b) in Proposition~1.2.4, we may assume that $M^*_\omega$ is a maximal proper $d$-split Levi subgroup of $C_\omega$ and that $M^*_{\omega'}=C_{\omega'}$ if $\omega'\neq\omega$. 

Then we have $M^*_\omega=S_\omega\times K_\omega$, a direct product defined on $\F q$, where  $S_\beta\cong [{\rm GL}_m]^{\delta_\omega}$ for some $m$, any $\beta\in\omega$, and $(S_\omega)^F\cong {\rm GL}_m((\epsilon q)^{|\omega|\delta_\omega})$, $K_\beta={\rm GL}_b$, $(K_\omega)^F\cong {\rm GL}_b((\epsilon q)^{|\omega|} )$. We have  $M^*=M^*_1\times M^*_2$, where $M^*_1\cong[{\rm GL}_{m_1}]^\delta$ acts on the space $V_1\oplus V_2\dots \oplus V_\delta$ of fixed points of $K_\omega\times(\times_{\omega'\neq \omega}C_{\omega'})$ in $\bar F{}^n$ (so that $m_1=m\delta_\omega|\omega|/\delta)$) and $(M^*_1)^F\cong ({\rm GL}_{m_1}((\epsilon q)^\delta)$, $M^*_2\cong {\rm GL}_{m_2}$ ($m_2=n-m_1\delta$) and $(M^*_2)^F\cong{\rm GL}_{m_2}(\epsilon q)$.

A dual Levi subgroup $M$ of $M^*$ in $(H,F)$ is isomorphic to $M^*$. 
Note that the center of $L=M\cap G$ may be disconnected : the center of $M$ has rank $(\delta +1)$ and the function determinant on it is $\bar F{}^\times\to \bar F{}$, $(\la_1,\dots,\la_\delta,\la_0)\mapsto (\la_1\dots \la_\delta)^{m_1}{\la_0}^{m_2}$. 

An element $a$ in $A$ acts on $\Sigma$ as an element of $\bar \F {}$ by multiplication and so semi-regularly on the eigenspaces $V_\beta$ of $t$, or analogous in $H$. One sees that $a$ may be found in $M^*$ only if $a\omega=\omega$. Any $\beta\in\omega$, $C_\beta=(C_\beta\cap S_\omega)\times (C_\beta\cap K_\omega)$ and $\phi$ acts regularly on these products, hence such a $a=(a_1,a_2)\in M^*_1\times M^*_2$ may exists, if the spectrum of $s$ in each $V_j$ is $a_1$-stable. The order of $a_1$ has to divide $|\omega|/(|\omega|,\delta)$.
\bull

(B.2) Bonnaf\'e-Cabanes parametrization.

Recall the decomposition $\cent{H^*}t=\times_{\beta\in\Sigma}C_\beta$, $\Sigma$ the set of eigenvalues of $t$ in the standard representation of $H^*$ on $\bar \F{}^n$. We have seen that $\cent{G^*}s$ is a semi-direct product of $\cento{G^*}s$ by the abelian group ${\rm A}_{G^*}(s)$. 

Let $A:= {\rm A}_{G^*}(s)^F$, $W^\circ(s)$, $W(s)$ be the Weyl groups of $\cento{G^*}s$, $\cent{G^*}s$  with respect to a diagonal torus $T^*_s$ of $\cento{G^*}s$, that is $W(s)=\nor{\cent{G^*}s}{T^*_s}/T^*_s$. The map $H^*\to G^*$ restricts to $\cent{H^*}t\to\cento{G^*}s$,  $T^*_t\to T^*_s$ and an isomorphism $W(t)\to W^\circ(s)$. Design by $\phi$ the action of $F$ on  $W^\circ(s)$ and $W(s)$, hence on $\Sigma$, so that $B:=\gen {A,\phi}$ is considered as a subgroup of ${\frak S}(\Sigma)$. Then $W(s)$,  $W^\circ(s)\rtimes \gen {A,F} $ and $W^\circ (s)^\phi\rtimes A$ appear as  wreath products with base group $W^\circ(s):=W(\cento{G^*}s)=\times_{\beta\in \Sigma}W_\beta$ or $W^\circ(s)^\phi=\times_{\omega\in\Sigma/\gen\phi}W_\omega^\phi$. We may have $W_\beta=1$ for some $\beta\in \Sigma$, then the following  combinatorial description of $\ser{G^F}s$ keep sense.

For each orbit $\omega\in\Sigma/\gen \phi$ is given an integer $\delta_\omega$  so that elements of cycle type $(\delta_\omega)^m$ in $W_\beta$ ($\beta\in\omega$) provide parameters with respect to the diagonal torus of $C_\omega$ of  $\phi_d$-split Levi subgroups of $C_\omega$.

Bonnaf\'e [3] Th\'eor\`eme 23.9 (non twisted type) and Cabanes  [12] Theorem~3.6 (twisted type) have shown that, under conjecture ${\frak G}$ (5.4.2.0),
 there is a one-to-one map 
$$R^G[s]=R[s]\colon \pm\II{W(s)^\phi}/\{-1,1\}\to \pm \ser{G^F}s/\{-1,1\}\leqno{(5.4.2.7)}$$
In other words there exists a sign function $\upsilon\colon \II{W(s)^\phi}\to\{-1,1\}$ such that $\upsilon R[s]$ is one-to-one from $\II{W(s)^\phi}$ onto $ \ser{G^F}s$.

Simplify $W^\circ(s)\rtimes B$ (resp. $W^\circ(s)$) in $W$ (resp. $W^0$). If $\chi=(\chi_\beta)_{\beta\in \Omega}\in \times_{\beta\in \Omega}\II{W_\beta}$, we denote $(\chi_\beta)_{\beta\in \Omega}$, considered as an element of $\II{W_\Omega}$, by $\chi^\Omega$  so that $\chi^\Sigma=\chi\in\II{W^0}$.

We first describe the map $R[s]$, following [3], [12].
That map  is defined by restriction of a linear map between spaces of central functions. 

Let ${\rm Cent}(G^F,s)$ be the space of central functions on $G^F$ with basis $\ser{G^F}s$. Using notations of 1.3.5 we knows that in  the regular embedding $G\subseteq H$, $H^F/\tau_{H,s}(A)$ acts on $\ser{G^F}s$. By (1.3.5.1) $(H^F/\sigma_{H,s}(A))^\wedge$ is isomorphic to $A$. So ${\rm Cent}(G^F,s)$ decomposes as a direct sum of isotypic $(H^F/\tau_{H,s}(A))$-spaces ${\rm Cent}(G^F,s,a)$. One obtains an orthogonal decomposition which is independant of the regular embedding [3] (11.15) corresponding on $W$-side to the decomposition $W^{\phi}=\cup_{a\in A}W^{0\,\phi}a$ :
$${\rm Cent}(G^F,s)=\perp_{a\in A}{\rm Cent}(G^F,s,a),\quad {\rm Cent}(W^\phi)=\perp_{a\in A}{\rm Cent}(W^{0\,\phi}a)^A\leqno{(5.4.2.8)}$$
 Isometries are defined in [3] 23.C
$$ R^G[s,a]=R[s,a]\colon {\rm Cent}(W^{0\,\phi}a)^A \to {\rm Cent}(G^F,s,a)\quad R[s]=\oplus_{a\in A}R[s,a] \leqno{(5.4.2.9)}$$
 The scalar product on ${\rm Cent}(W^{0\,\phi}a)^A $ has to be defined as $1/|A|$ times the usual scalar product on $a$-central functions on $W^{0\,\phi}$.
 
When $\cent{G^*}s$ is connected $R^G[s]=R^G[s,1]$ may be obtained as in part (A) thanks to Jordan decomposition and the unipotent case $\II{W(\cent{G^*}s)^F}\to \ser{\cent{G^*}s^F}1$, or simply because $\cent{G^*}s$ is a Levi subgroup of $G^*$ (see 1.3.2 and 1.3.3 (b)). Up to a sign the map $(\mu\mapsto R_\mu)$ in (5.4.2.1) is just $\mu\mapsto R^G[1](\mu)$.

By properties of  wreath products (5.4.1.6), ${\rm Cent}(W^{0\,\phi}a)$ has an orthonormal basis $\{\chi\rtimes a\}_\chi$ where $\chi$ runs in $\II{W^{0\,\phi}}^a$. Using the natural one-to-one maps $\II{W^{0\,\phi }}^a\to \II{W^{0\,\gen{\phi,a}}}\cong\II{W^0}^{\gen{a,\phi}}\to \II{W^{0\,a}}^\phi$ is defined an isometry ${\rm Cent}(W^{0\,\phi }a)\to {\rm Cent}(W^{0\,a}\phi)$, sending $\chi_\phi\rtimes a$ on  $\chi_a\rtimes\phi$ for any $\chi\in\II{W^0}^{\gen {a,\phi}}$, and that isometry commutes with $A$-action by (5.4.1.5) (i). So is defined an isometry
$$\sigma^W_{s,a}=\sigma_{s,a}\colon {\rm Cent}(W^{0\,\phi }a)^A\to {\rm Cent}(W^{0\,a}\phi)^A\leqno{(5.4.2.10)}$$

Let ([3] 23.C) 
$${\cal R}^G[s,a]:=R^G[s,a]\circ (\sigma_{s,a}^W)^{-1} $$  
A fact is that  ${\cal R}^G[s,a]$ may be defined for any type from ${\rm Cent}(W^{0\, a}\phi)$ to ${\rm Cent}(G^F,s,a)$,  and then commute with $A$-action [3] (17.17), Proposition~17.18. So we compose it with the orthogonal projection  $$\pi^a_A\colon {\rm Cent}(W^{0\, a}\phi)\to {\rm Cent}(W^{0\, a}\phi)^A$$
 Furthermore ${\cal R}^G[s,a]$ commute in some sense with Lusztig induction. Let $L$ be an $F$-stable Levi subgroup of $G$, with dual $L^*$ in $G^*$ such that $s\in L^*$, assume that some diagonal maximal torus of $\cento{L^*}s$ has parameter $w_L\in W^0$ with respect to the maximal diagonal torus of $\cento{G^*}s$, let $A(L)={\rm A}_{L^*}(s)^F$, a subgroup of $A$. One has a commutation formula  for any $a\in A(L)$ [3] Proposition~17.24 :
$$ |A(L)|.\Lu LG\circ{\cal R}^L[s,a]=|A|. {\cal R}^G[s,a]\circ \pi^a_A\circ \Ind{W_L^{0\,a}w_L\phi}{W^{0\,a}\phi}\leqno{(5.4.2.11)}$$
Therefore the inverse image  of Lusztig induction $\Lu LG$ by $R[s]$ in (5.4.2.7) is given by    
$${\cal I}_{W_L}^W=|A/A(L)|\oplus_{a\in A(L)}[(\sigma^W_{s,a})^{-1}\circ \pi^a_A\circ \Ind{W_L^{0\, a}w_L\phi}{W^{0\,a}\phi}\circ\sigma_{s,a}^{W_L}]\leqno{(5.4.2.12)}$$
(we don't write the orthogonal projections given by (5.4.2.8)).
We would like to compute its value on $\eta_{w_L\phi}$ when $(L,R^L[s]({\eta_{w_L\phi}}))$ would define a $d$-cuspidal datum in series $[s]$ in $(G,F)$ for some $\eta\in\II{W_L^{0}}^{w_L\phi}$.  
Our claim is that  $R^L[s]({\eta_{w_L\phi}})$ is a component of $\Res{M^F}{L^F}\xi$ where $M=L.\z H$ and $\xi$ is $d$-cuspidal. 

By [3], (23.15)  for any $\zeta\in\II{W^0}^\phi$ and with notations of (5.4.1.4) one has
$$\Res{H^F}{G^F}(R^H[t](\zeta_\phi))=\sum_{\theta\in (A_\zeta)^\wedge}R^G[s](\Gamma^{W^\phi}(\zeta_\phi*\theta))\leqno{(5.4.2.13)}$$

(B.3) $d$-cuspidal data and G.$d$-HC series in type $\AA$.

From our description in part  (A.3)  of that proof and  (5.4.2.6), a $d$-cuspidal datum in series $(t)$ in $(H,F)$ is defined on $W^0$-side as follows : for each orbit $\omega\in\Sigma/\gen \phi$ are given $W_\omega\cong \times_{\beta\in\omega}{\frak S}_{n_\omega}$, an integer $\delta_\omega$, and a ``$\delta_\omega$-cuspidal" element $\zeta^\omega$ in $\II{{\frak S}_{m_\omega}}$ where $\delta_\omega$ divides $(n_\omega-m_\omega)$, $v_\omega\in W_\omega^\phi$, $V_{v_\omega}\subseteq \cent{W_\omega}{v_\omega}$. We obtain a ``$(\delta_\omega)_{\omega\in\Sigma/\gen \phi}$-cuspidal datum  in $(W^0,\phi)$" :  $(\times_\omega V_\omega,(\zeta^\omega_{v_\omega\phi})_\omega)$.  For each $\omega\in\Sigma/{\gen \phi}$ is defined a $d$-split Levi subgroup $M_{v_\omega}$  of $C_\omega:=\Pi_{\beta\in\omega}C_\beta\subseteq \cent{H^*}t$. Then (5.4.2.4) apply in each $W_\omega$, giving a similar formula in $W^0$. So is defined a $d$-cuspidal datum $(M,\la)$ where $M=\times_\omega M_{v_\omega}$, $w_M=(v_\omega)_\omega$, $W(M)=\times_\omega V_{v_\omega}$, $\la=\otimes_\omega R^{M}[t](\zeta^\omega_{v_\omega\phi})$. The G.$d$-HC series in series $(t)$ is a set product on $\Sigma/\gen \phi$ of G.$d$-HC series, each one defined from a series in $W_\omega^\phi$, as (5.4.2.2) gives (5.4.2.5). 

Furthermore, if $\zeta=(\zeta^\omega)_\omega\in\II {W^0}^\phi$, $R^H[t](\zeta_\phi)$ is $d$-cuspidal if and only if $\z {\cent {H^*} t}_{\phi_d}\subseteq \z {H^*}$ and $\zeta_\phi=(\zeta^\omega)_\omega$ is $(\delta_\omega)_\omega$-cuspidal. 
A first step is  to show that $R^G[s](\Gamma^{W^\phi}(\zeta_\phi*\theta))$ is $d$-cuspidal if and only if $R^H[t](\zeta_\phi)$ is $d$-cuspidal, to obtain a generalization of 2.1.5, where the order of $s$ is prime to $\ell$ :

\medskip\noindent (B.3.1) {\sl Let $(G,F)\subseteq (H,F)$ be a regular embedding, $\sigma\colon H^*\to G^*$ a dual morphism, $t$ a semi-simple element of $(H^*)^F$, $s=\sigma^*(t)$, $\mu\in\ser{H^F}t$, $\la$ an \irr\ component of $\Res{H^F}{G^F}\mu$. Then $\la$ is $d$-cuspidal if and only if $\mu$ is $d$-cuspidal.}

\medskip\noindent{\it Proof of (B.3.1). }

(a) If $\mu$ is $d$-cuspidal, $\la$ is $d$-cuspidal :

Let  $\mu=R^H[t](\zeta_\phi)$, assumed to be $d$-cuspidal. Let $L$ be a maximal proper $d$-split Levi subgroup of $(G,F)$, a dual $L^*$ in $(G^*,F)$ such that $s\in L^*$, $M^*=L^*\z {H^*}$, $t\in M^*$. As $\z {\cent {H^*} t}_{\phi_d}\subseteq \z {H^*}$, $\cent{M^*}t\neq \cent{H^*}t$ and $\cent{M^*}t$ (resp. $\cent{L^*}s$) is a proper maximal $d$-split Levi subgroup of $(\cent{M^*}t,F)$ (resp. $(\cent{G^*}s,F)$). (B.1.1) and the description we made in its proof  apply. We have $A(L)(\omega)=\omega$.

Let $\xi\in \II{W_M}^{w_L\phi}$ and $\psi\in (A(L)_\xi)^\wedge$ defining  $\eta=\Gamma^{W_L^{w_L\phi}}(\xi_{w_L\phi}*\psi)$ in $\II{W_L^{w_L\phi}}$. To compute the scalar product of ${\cal I}^W_{W_L}(\eta)$ with $\Gamma^{W^\phi}(\zeta_\phi*\theta)$ using the orthogonal decompositions in (5.4.2.8) we may restrict $a$ in $A(L)_\xi\cap A_\zeta$.  The projection of $\eta$ on ${\rm Cent}(W_L^{w_L\phi}a)$ is $\psi(a)\sum_{b\in A(L)/A(L)_\xi}(\xi_{w_L\phi}^b\rtimes a)$ and one has 
$$\sigma^{W_L}_{s,a}(\eta)=\psi(a)\sum_{b\in A(L)/A(L)_\xi}({\xi_a}^b \rtimes w_L\phi)$$
To compute $\Ind{W_L^{0\, a}w_L\phi}{W^{0\,a}\phi}(\sigma^{W_L}_{s,a}(\eta))$ in (5.4.2.12) we use (5.4.2.4) applied in $W^{0\, a}$.  Then $w_L^\omega:=(v(\beta))_{\beta\in \omega}$ and $\zeta^\omega:=(\zeta(\beta))_{\beta\in\omega}$ are constant on $\omega$, so are $v_a$ and $\zeta_a$ on $\omega_a:=\omega/\gen a$. The cardinal $r_a$ of $\omega_a$ is a divisor of $r=|\omega|$. In the computation we made to obtain (5.4.2.4), if we consider $(W^{0\,a},W_L^{0\,a},\zeta_a)$ instead of $(W^{0},V^{0},\chi)$, we have to replace $r$ by $r_a$. Then $\delta_\omega=\delta/(\delta,r)$ divides $\delta_a:=\delta/(\delta,r_a)$. The description of partitions and hooks by so-called $\beta$-sets shows that any $\delta_\omega$-core is a $\delta_a$-core, because deleting a hook of length $\delta_a$ may always be obtained by deleting successively $\delta_a/\delta_\omega$ hooks of length $\delta_\omega$. As $\zeta^\omega_\phi$ is defined by a $(\delta_\omega)$-cuspidal element, then $\Res{W^{0\,a}\phi}{W_L^{0\, a}w_L\phi}({\zeta_a}\rtimes \phi)=0$. That implies $\slu LG[R^G[s](\Gamma^{W^\phi}(\chi*\theta))]=0$ hence  $R^G[s](\Gamma^{W^\phi}(\chi*\theta))$ is $d$-cuspidal for any $\theta\in (A_\chi)^\wedge$, a part of (B.3.1).

(b) If $\mu$ is not cuspidal, $\la$ is not $d$-cuspidal :

By (A.3), if $\z {H^*}_{\phi_d}\neq \z{\cent{H^*}t}_{\phi_d}$, then $R^H[t](\zeta_\phi)$ is not $d$-cuspidal whenever $\zeta\in\II{W^0}^\phi$.

In that case let $M^*=\cent{H^*}{\z{\cent{H^*}t}_{\phi_d}}$, $L^*$ the image of $M^*$ in $G^*$, $M$ a $d$-split Levi subgroup of $(H,F)$ in duality with $M^*$, $L=M\cap G$ : we have $W(\cent{M^*}t)=W(\cent{H^*}t)$, $W(\cent{G^*}s)=W(\cent{L^*}s)$ so that  $R^H[t](\zeta_\phi)=\Lu {M}H(R^M[t](\zeta_\phi))$ and ${\rm A}_{L^*}(s)\cong {\rm A}_{G^*}(s)$. By (1.3.1.1) and (5.4.2.13) 

\noindent $\Res{H^F}{G^F}(R^H[t](\zeta_\phi))=\Lu LG(\Res{M^F}{L^F}(R^M[t](\zeta_\phi))=\Lu LG(\sum_{\theta\in A^\wedge} R^L[s](\Gamma^{W^\phi}(\zeta_\phi\otimes\theta)))$. 

By Clifford theory the isomorphic quotient groups $H^F/G^F$ and $ M^F/L^F$ act regularly and the same way on the sets $\{R^G[s](\Gamma^{W^\phi}(\zeta_\phi \otimes\theta)\}_\theta$ and $\{R^L[s](\Gamma^{W^\phi}(\zeta_\phi*\theta)\}_\theta$ [3], (10.4). Thus $\Lu LG$ induces a one-to-one map between these two sets (one may show that $\Lu LG(R^L[s](\Gamma^{W^\phi}(\zeta_\phi*1)))=R^G[s](\Gamma^{W^\phi}(\zeta_\phi*1))$
hence $\Lu LG(R^L[s](\Gamma^{W^\phi}(\zeta_\phi*\theta)))=R^G[s](\Gamma^{W^\phi}(\zeta_\phi*\theta))$). The components of $R^H[t](\zeta_\phi)$ are not $d$-cuspidal. If $\zeta$ is $(\delta_\omega)_\omega$-cuspidal, $(M,R^M[t](\zeta_\phi))$ is a $d$-cuspidal data in series $(t)$ in $(H,F)$. 

More generally (b) is implied by  the following descripion of G.$d$-HC series in series $[s]$ in $G$ :

\medskip\noindent{(B.3.2) {\sl Let $(M,\mu)$ be a $d$-cuspidal datum in $(H,F)$, where $\mu=R^M[t](\xi_{w_M\phi})$, $\xi\in \II{W_M}^{w_M\phi}$, and  let $L=M\cap G$, $\la=R^L[s](\Gamma^{W_L^{w_L\phi}}(\xi*\psi))$ ($W_M=W_L^0=W(\cent{H^*}t)$, $\psi\in (A(L)_\xi)^\wedge$). Then $R^G[s](\Gamma^{W^\phi}(\zeta*\theta))$ is a component of $\Lu LG\la$ if and only if $R^H[t](\zeta_\phi)\in\ser{H^F}{(M,\mu)}$ and $\Res{A_\theta}{A(L)_\xi}\theta=\psi$.}

\medskip\noindent{\it Proof of (B.3.2).}   Let $\chi:=R^H[t](\zeta_\phi)\in\ser{H^F}{(M,R^M[t](\xi_{w_M\phi})}$ ($\zeta\in\II{W^0}^\phi$). We may assume $M\neq H$, thanks to (a) and (5.4.2.13).

The last condition in (B.3.2) assume $A(L)_\xi\subseteq A_\zeta$, a consequence of (B.1) :

 We have $\cent{G^*}t=\times_\omega C_\omega$, $\cent {M^*}t=\times_\beta M_\beta=\times_\omega M_\omega$. For any $\omega\in\Sigma/\gen \phi$, if $M_\omega=C_\omega$, $\zeta^\omega=\xi^\omega$ and $A(L)_\xi$ fixes $\zeta^\omega$. Assume $M_\omega\neq C_\omega$, then $M_\omega$ is a maximal proper $d$-split Levi subgroup of $C_\omega$ because $M_\omega=S_\omega\times K_\omega$  where $S_\omega$ is a torus (as in the proof of (B.1.1) with $m=1$). By (B.1), $A(L)\omega=\omega$.  As $\zeta\in\II{W^0}^\phi$, $\zeta(\beta)$ is constant on $\omega$ hence $\zeta^\omega$ is fixed by $A(L)$. Thus we have $A(L)_\xi\subset A_\zeta$.

 We have by (5.4.2.13)
$$\Res{M^F}{L^F}(R^M[t](\xi_{w_M\phi}))=\sum_{\psi\in(A(L)_\xi)^\wedge}R^L[s](\Gamma^{W_L^{w_L\phi}}(\xi_{w_L\phi}*\psi))$$
and know by (a) that  $R^L[s](\Gamma^{W_L^{w_L\phi}}(\xi_{w_L\phi}*\psi))$ is $d$-cuspidal for any $\psi\in(A(L)_\xi)^\wedge$. Note that for any $nM^F$ in $ \nor{H^F}{M^F,[t]}/M^F$, the isomorphisms $\nor{H^F}{M^F,[t]}/M^F\cong \nor{W^{0\,\phi}}{W_Mw_M}/W_M^{w_M\phi}\cong \nor{G^F}{L^F,[s]}/L^F$ imply $A(L)_\xi=A(L)_{\xi^n}$ and allow to write

$R^M[t](\xi_{w_M\phi})^n=R^M[t](\xi^n_{w_M\phi})$, 
$\Res{M^F}{L^F}(R^M[t](\xi^n_{w_M\phi}))=\sum_{\psi\in(A(L)_\xi)^\wedge}R^L[s](\Gamma^{W_L^{w_L\phi}}(\xi^n_{w_L\phi}*\psi))$.

Thanks to the functorial behaviour of formula (5.4.2.4)  with respect to the direct product on $\Sigma/\gen {\phi,A}$, 
$\zeta_a\rtimes\phi$ is a component of  $\Ind{{W_L}^{0\,a}*w_L\phi}{W^{0\,a}\phi}(\xi_a\rtimes w_L\phi)$ only if, for any $\Omega\in\Sigma/{\gen {a,\phi}}$, $\zeta_a^{\omega_a}=(\zeta^\Omega)_a$ appears as a component of some $\Ind{V^{\omega_a}*v_{\omega_a}^{r_a}}{W^{\omega_a}}(\xi^b_a)^{\omega_a}$, where $b\in A$ and  $\omega_a\in(\Sigma/\gen a)/\gen \phi\cong \Sigma/{\gen {a,\phi}}$. 

 Coming back to the $d$-cuspidal datum $(L,R^L[s](\eta_{w_L\phi}))$ we know that if $\zeta_\phi\notin {\cal S}^{W^{0\,\phi}}(W_L^{0\,w_L\phi},\xi^b_{v\phi})$ for some $b\in A$, then $R^G[s](\Gamma^{W^\phi}(\chi*\theta))$ is not a component of $\Lu LG(R^L[s](\eta_{w_L\phi}))$. 
 
  So assume $\zeta_\phi\in {\cal S}^{W^{0\,\phi}}(W_L^{0\,\phi},{\xi_{w_L\phi}}^a)$, $a\in A(L)_{\xi}$, so that $R^H[t](\zeta_\phi)\in\ser{H^F}{(M,\mu)}$. 
  On $W^{0\,\phi}$ we have a decomposition 
  $\Ind{W_L^{0\,\phi}*w_L}{W^{0\,\phi}}\xi_{w_L\phi}=\sum_\chi d(\xi,\chi)\chi_\phi$   where  $d(\xi,\chi)=\scal{\Res{W^{0\,\phi}}{W_L^{0\,\phi}*w_L}\chi_\phi}{\xi_{w_L\phi}}{W_L^{0\,w_L\phi}}$.
   
 As $(W_L,w_L,\xi)=(W_L^a,w_L^a,\xi^a)$, $d(\xi,\chi)=d(\zeta,\chi^a)$. On $W^{0\, a}$  we have $\Ind{W_L^{0\, a}*w_L}{W^{0\, a}}\xi_a=\sum_{\chi_a} d(\xi,\chi)\chi_a$. Furthermore $d(\xi^b,\chi)=d(\xi,\lexp b \chi)$ for any $b\in A$ and $\pi^a_A(\chi_a\rtimes \phi)=1/|A|\sum_{b\in A}((\chi^b)_a\rtimes\phi)$. Finally we obtain
$$|A/A(L)|(\oplus_{a\in A(L)}\sigma^W_{s,a}\circ {\cal I}_{W_L}^W)(\eta_{w_L\phi})=\oplus_{a\in A(L)}\psi(a)\sum_{\chi_a} \sum_{b\in A/A(L)_\xi}d(\xi,\chi^b)((\chi^b)_a\rtimes \phi)$$

Let $\theta\in(A_\zeta)^\wedge$. If $a\notin A_\zeta$, the projection of $\Gamma^{W^\phi}(\zeta*\theta)$ on ${\rm Cent}(W^{0\,\phi}a)$ is null. If  $a\in A_\zeta$, the projection of $\Gamma^{W^\phi}(\zeta*\theta)$ on ${\rm Cent}(W^{0\,\phi}a)$ is $\theta(a)\sum_{b\in A/A_\zeta}({\zeta_\phi}^b\rtimes a)$, whose image by $\sigma^W(s,a)$ is $\theta(a)\sum_{b\in A/A(L)_\xi}({\zeta_a}^b\rtimes \phi)$. As $\sigma^W_{s,a}$ is an isometry, the scalar product of $(\Gamma^{W^\phi}(\zeta*\theta))$ with ${\cal I}_{W_L}^W(\eta_{w_L\phi})$ is  equal to

\noindent $|A/A_V|.|A_\zeta/A(L)_\xi|\sum_{a\in A_\zeta}\theta(a)\psi(a)^{-1}d(\xi,\zeta)$. One sees that it is non zero and only if $\Res{A(L)_{\zeta_{v\phi}}}{A_\zeta}\psi=\theta$. \bull

\vfill\eject

\centerline{REFERENCES.}
\bigskip
\item{[0]} M. Enguehard, Vers une d\'ecomposition de Jordan des blocs des groupes r\'eductifs finis, {\it J. Algebra} {\bf  319} (2008), 1035--1115.

\item{[1]} C. Bonnaf\'e, Produits en couronne de groupes lin\'eaires, {\it J. Algebra} {\bf 211} (1999), 57--98.

\item{[2]} C. Bonnaf\'e,   Quasi-isolated elements in finite reductive groups,
{\it Comm. Algebra} {\bf 33} (2005),  2315--2337.

\item{[3]} C. Bonnaf\'e,   Sur les caract\`eres des groupes r\'eductifs finis \`a centre non connexe : applications aux groupes sp\'eciaux lin\'eaires et unitaires, {\it Ast\'erisque} {\bf 306}, 2006.

\item{[4]} C. Bonnaf\'e and J. Michel, Computational proof of the Mackey formula for $q>2$, {\it J. Algebra} {\bf 327} (2011), 506--526.

\item{[5]} C. Bonnaf\'e et R. Rouquier, Vari\'et\'es de Deligne-Lusztig et
cat\'egories d\'eriv\'ees, {\it Publ. Math. I.H.E.S.} {\bf 97} (2003),  1--59.

\item{[6]} C. Bonnaf\'e and R. Rouquier, Coxeter orbits and modular representations, {\it Nagoya Math. J.}  {\bf 183} (2006), 1--34.

\item{[7]} M. Brou\'e, Les $\ell$-blocs des groups GL$(n,q)$ et
U$(n,q^2)$ et leurs structures locales, {\it S\'eminaire Bourbaki}, {\bf 640} (1984),
28 pages.

\item{[8]} M. Brou\'e, Blocs, isom\'etries parfaites, cat\'egories d\'eriv\'ees,
{\it C.R. Acad.  Sci. Paris} {\bf 307} (1988),  13--18. 

\item{[9]} M. Brou\'e et J. Michel, Blocs et s\'eries de Lusztig dans un groupe
reductif fini, {\it J. reine angew. Math.} {\bf 395} (1989),  56--67.

\item{[10]} M. Brou\'e, G. Malle  and J. Michel, Generic blocks of finite
reductive groups, {\it Ast\'erisque} {\bf 212} (1993),  7--92.

\item{[11]} M. Cabanes, Unicit\'e du sous-groupe ab\'elien distingu\'e maximal dans certains sous-groupes de Sylow, {\it Cr. Acad. Sci. Paris, S\'er. I}  {\bf 318} (1994), 889--894.

\item{[12]} M. Cabanes, On Jordan decomposition of characters for 
${\rm SU}(n,q)$, {\it arXiv} :1112.6035v1, 2011, {\it Jl of Algebra}, {\bf 374} (2013), 216--230.

\item{[13]} M. Cabanes and M. Enguehard, On unipotent blocks and their ordinary
characters, {\it Invent. Math.}  {\bf 117} (1994),  149--164.

\item{[14]} M. Cabanes and M. Enguehard, Local methods for the  blocks of 
reductive groups over finite fields, {\it Progr. Math.} {\bf 141} (1996), 141--164.

\item{[15]} M. Cabanes and M. Enguehard, On Blocks of Finite Reductive Groups and Twisted Induction,  {\it Advances in Math.}  {\bf 145} (1999), 189--229.

\item{[16]} M. Cabanes and M. Enguehard,
{\it Representation Theory  of Finite Reductive Groups}, new math.
monographs {\bf 1}, Cambridge University Press, 2004.

\item{[17]} R. W. Carter,  {\it Finite Groups of Lie Type: Conjugacy Classes and
Complex Characters}, Wiley, 1985.

\item{[18]} P. Deligne and G. Lusztig, Representations of reductive groups
over finite fields, {\it Ann. of Math.}  {\bf  103} (1976), 103--161.

\item{[19]} F. Digne et J. Michel, On Lusztig's parametrization of characters of
finite groups of Lie type, {\it Ast\'erisque}  {\bf 181-182} (1990),  113--156.

\item{[20]} F. Digne and J. Michel, {\it Representations of Finite Groups
of Lie Type}, London Math. Soc. Student Texts {\bf 21}, Cambridge
University Press, 1991.

\item{[21]} F. Digne et J. Michel, Groupes r\'eductifs non connexes, {\it
Ann. scient. Ec. Norm. Sup., $4^e$ s\'erie} {\bf 27} 1994,  345--406.

\item{[22]} F. Digne et J. Michel, Points fixes des automorphismes
quasi-semi-simples, C.R.A.S., Ser. I  {\bf 334} ( 2002),  1055--1060. 

\item{[23]} M. Enguehard, Sur les $\ell$-blocs unipotents des groupes
r\'eductifs finis quand $\ell$ est mauvais, {\it J.  Algebra}, {\bf 230}
(2000),  334--377. 

\item{[24]} P. Fong and B. Srinivasan, The blocks of finite classical groups, {\it
J. reine  angew. Math.}  {\bf 396} (1989), 122--91. 

\item{[25]} M.  Geck, A note on Harish-Chandra induction. {\it
Manuscripta math.}  {\bf 80} (1993),  393--401.

\item{[26]} G. Lusztig, {\it Characters  of reductive groups over a finite field},
Ann. Math. Studies {\bf 107}, Princeton University Press, 1984.

\item{[27]} G. Lusztig, On the representations of reductive groups with disconnected
center, {\it Ast\'erisque}  {\bf 168} (1988),  157--166.

\item{[28]} H. Nagao and Y. Tsushima,  {\it Representations of finite Groups},
Academic Press, 1989.

\item{[29]} T. Shoji, Some generalizations of Asai's results for classical groups, in {\it Algebraic Groups and Related Topics, Advances Studies in Pure Mathematics} {\bf 6} (1985).q

\item{[30]} T. Shoji,  Shintani descent for special linear groups, {\it J. Algebra} 
{\bf 199} (1998),  175--228.

\item {[31]} R. Steinberg, {\it Endomorphisms of linear algebraic groups},
Memoirs of the American Mathematical Society  {\bf 80}, 1968.

\item{[32]} J. Th\'evenaz, {\it $G$-algebras and Modular Representation Theory}, Oxford
University Press, 1995.

\end